\documentclass[12pt]{report}
\usepackage{graphicx}					
\usepackage{mathrsfs}
\usepackage{verbatimbox}
\usepackage{listofitems}
\usepackage{amsmath, amssymb}
\usepackage{wasysym} 
\usepackage{arydshln}
\usepackage{amsfonts, ifthen, amsthm, vmargin}
\usepackage[usenames]{color}
\usepackage{ulem}
\usepackage[colorlinks=true, linkcolor=blue3, 
urlcolor=gray2, citecolor=blue3, pagebackref=true]{hyperref} 
\newcommand{\sref}[2]{\hyperref[#2]{#1 \ref*{#2}}}
\usepackage[T1]{fontenc}
\usepackage{lmodern}
\newcommand{\dotcup}{\mathbin{\mathaccent\cdot\cup}}
\newcommand{\widesim}[2][4.2]{\mathrel{\overset{\scalebox{#1}[1]{$\sim$}}{#2}}} 
\newcommand{\keywords}[1]{{\textbf{Key words and phrases.}} #1}
\newcommand{\msc}[1]{{\textbf{MSC2020 subject classifications.}} #1}
\newcommand{\dci}[1]{{\textbf{Declaration of competing interest.}} #1}

\renewcommand\labelenumi{\textup{(\roman{enumi})}}
\renewcommand\theenumi\labelenumi
\DeclareFontFamily{OT1}{rsfs}{} \DeclareFontShape{OT1}{rsfs}{m}{n}{
<-7> rsfs5 <7-10> rsfs7 <10-> rsfs10}{}
\DeclareMathAlphabet\mathcurl{OT1}{rsfs}{m}{n}

\newtheorem{theorem}{Theorem}[section]
\newtheorem{lemma}[theorem]{Lemma}
\newtheorem{corollary}[theorem]{Corollary}
\newtheorem{proposition}[theorem]{Proposition}

\theoremstyle{definition}
\newtheorem{example}[theorem]{Example}

\theoremstyle{definition}
\newtheorem{definition}[theorem]{Definition}

\theoremstyle{definition}
\newtheorem{remark}[theorem]{Remark}

\theoremstyle{plain}
\newtheorem*{assump*}{Assumption}

\theoremstyle{plain}
\newtheorem*{lemma*}{Lemma}

\theoremstyle{plain}
\newtheorem*{cor*}{Corollary}

\theoremstyle{plain}
\newtheorem*{gt*}{The Grothendieck Inequality (Lindenstrauss-Pe{\l}czy\'{n}ski style)}

\theoremstyle{plain}
\newtheorem*{thmGI*}{The Grothendieck equality}

\theoremstyle{plain}
\newtheorem*{thmHI*}{The Haagerup equality}

\theoremstyle{plain}
\newtheorem*{obsv*}{Observation}

\theoremstyle{plain}
\newtheorem*{fact*}{Fact}

\theoremstyle{definition}
\newtheorem*{df*}{Definition}

\theoremstyle{definition}
\newtheorem*{ex*}{Example}

\theoremstyle{definition}
\newtheorem*{ack*}{Acknowledgments}

\theoremstyle{plain}
\newtheorem*{exkrv*}{Example (Krivine's constant)}

\theoremstyle{plain}
\newtheorem*{gi*}{Matrix Version of the Grothendieck equality}

\theoremstyle{plain}
\newtheorem*{guidex*}{Guiding Example}

\theoremstyle{plain}
\newtheorem*{prp*}{Proposition}

\theoremstyle{plain}
\newtheorem*{cprp*}{Proposition (correlation matrix version of GT and little GT)}

\theoremstyle{plain}
\newtheorem*{thm*}{Theorem}

\theoremstyle{remark}
\newtheorem*{thmTS*}{\textup{\textbf{\normalsize T\scriptsize HEOREM  \normalsize 
(T\scriptsize SIREL'SON, \normalsize 1980)}}}

\theoremstyle{plain}
\newtheorem*{thmgn*}{Theorem (Grothendieck, 1953; Niemi, 1983)}

\theoremstyle{plain}
\newtheorem*{c*}{Conjecture}

\theoremstyle{plain}
\newtheorem*{p*}{Problem}

\theoremstyle{definition}
\newtheorem*{rem*}{Observation}

\theoremstyle{remark}
\newtheorem*{rem2*}{Remark}

\theoremstyle{plain}
\newtheorem*{arcsinex*}{Example (Stieltjes, 1889)}

\theoremstyle{plain}

\theoremstyle{plain}

\theoremstyle{plain}

\newcount\colveccount
\newcommand*\colvec[1]{
        \global\colveccount#1
        \begin{pmatrix}
        \colvecnext
}
\def\colvecnext#1{
        #1
        \global\advance\colveccount-1
        \ifnum\colveccount>0
                \\
                \expandafter\colvecnext
        \else
                \end{pmatrix}
        \fi
}

\definecolor{amaranth}{rgb}{0.9, 0.17, 0.31}
\definecolor{blue3}{rgb}{0.32, 0.09, 0.98}
\definecolor{blue4}{rgb}{0.13, 0.67, 0.8}
\definecolor{tangerine}{rgb}{1.0, 0.8, 0.0}
\definecolor{taupegray}{rgb}{0.55, 0.52, 0.54}
\definecolor{coolblack}{rgb}{0.0, 0.18, 0.39}
\definecolor{gray2}{rgb}{0.5, 0.5, 0.5}
\definecolor{cadetblue}{rgb}{0.37, 0.62, 0.63}
\definecolor{green2}{rgb}{0.0, 0.5, 0.0}
\definecolor{myred}{rgb}{0.8, 0.0, 0.0}

\definecolor{blue2}{rgb}{0.08, 0.38, 0.74}
\definecolor{mypink}{rgb}{0.858, 0.188, 0.478}
\definecolor{cite}{rgb}{1.0, 0.13, 0.32}

\newcommand{\floor}[1]{\left\lfloor #1 \right\rfloor}
\newcommand{\ceil}[1]{\left\lceil #1 \right\rceil}

\newcommand\ignore[1]{}

\newcommand\C{\mathbb C}
\newcommand\D{\mathbb D}
\newcommand\R{\mathbb R}
\newcommand\T{\mathbb T}
\newcommand\M{\mathbb M}
\newcommand\N{\mathbb N}
\newcommand\Z{\mathbb Z}
\newcommand\F{\mathbb F}
\newcommand\vc[2]{{\text{vec}}(#1,#2)}
\newcommand\vcfour[4]{{\text{vec}}(#1,#2,#3,#4)}
\newcommand\E[2][\P]{\mathbb E_{#1}\left[ #2\right]}
\renewcommand\E{\mathbb E}
\renewcommand\P{\mathbb P}
\renewcommand\H{\mathbb H}
\newcommand\Q{\mathbb Q}

\renewcommand{\S}{\mathbb{S}}

\renewcommand*{\Re}{\operatorname{Re}}
\renewcommand*{\Im}{\operatorname{Im}}

\newcommand\ind{1\hspace{-2.5mm}1}
\newcommand\indsm{1\hspace{-1mm}1}

\newcommand\ldju[2]{\bigcup_{#1}^{#2}\hspace{-8.75mm\cdot}\hspace{8mm}}

\makeatletter
\renewcommand{\@seccntformat}[1]{\csname the#1\endcsname.\hspace{1em}}%

\makeatother

\makeatletter
  \renewcommand*\env@matrix[1][*\c@MaxMatrixCols c]{%
    \hskip -\arraycolsep
    \let\@ifnextchar\new@ifnextchar
  \array{#1}}
\makeatother

\makeatletter
\def\widebreve{\mathpalette\wide@breve}
\def\wide@breve#1#2{\sbox\z@{$#1#2$}%
     \mathop{\vbox{\m@th\ialign{##\crcr
\kern0.08em\brevefill#1{0.8\wd\z@}\crcr\noalign{\nointerlineskip}%
                    $\hss#1#2\hss$\crcr}}}\limits}
\def\brevefill#1#2{$\m@th\sbox\tw@{$#1($}%
  \hss\resizebox{#2}{\wd\tw@}{\rotatebox[origin=c]{90}{\upshape(}}\hss$}
\makeatletter

\cleardoublepage 
\begin{document}
\title{Upper bounds for Grothendieck constants, quantum correlation matrices and
CCP functions 
}
\author{Frank Oertel  \\
	Philosophy, Logic \& Scientific Method\\
	Centre for Philosophy of Natural and Social Sciences (CPNSS)\\
	London School of Economics and Political Science\\
	Houghton Street, London WC2A 2AE, UK\\
	}

\date{} 

\maketitle

{
\hypersetup{linkcolor=coolblack} 
\tableofcontents
\cleardoublepage 
\pagenumbering{arabic} 
}


\begin{abstract}
\noindent Within the framework of the search for the still unknown {best 
possible} 
value of the real and complex Grothendieck constant $K_G^\F$ in the famous 
Grothendieck inequality (unsolved since 1953), where $\F$ denotes either the 
real or the complex field, we concentrate our search on their smallest upper 
bound. To this end, we establish a basic framework, built on functions which map 
correlation matrices to correlation matrices entrywise by means of the Hadamard 
product, such as the Krivine function in the real case or the Haagerup function 
in the complex case. By making use of multivariate real and complex Gaussian 
analysis, higher transcendental functions, integration over spheres and 
combinatorics of the inversion of Maclaurin series, we provide an approach by 
which we also recover all famous upper bounds of Grothendieck himself ($K_G^\R 
\leq \sinh(\pi/2) \approx 2.301$ - \cite{G1953}), Krivine ($K_G^\R \leq \frac{\pi}
{2 \ln(1 + \sqrt{2})} \approx 1.782$ - \cite{Kr1977_2}) and Haagerup ($K_G^\C  
\leq 1.405$, numerically approximated - \cite{H1987}); each of them as a 
special case. In doing so, we aim to unify the real and complex cases as much 
as possible and apply our results to several concrete examples, including the 
Walsh-Hadamard transform (``quantum gate'') and the multivariate Gaussian 
copula - with foundations of quantum theory and quantum information theory 
in mind. 
{Futhermore we give a {slightly simplified} proof of the best non 
computer-aided approximation known up to now, i.e., $K_G^\R < \frac{\pi}
{2 \ln(1 + \sqrt{2})}$ ({\cite{BMMN2013, Kr2023}}).}
We summarise our key results in form of an algorithmic scheme and shed light 
on related open problems for future research works.
\end{abstract}
\noindent\keywords{Grothendieck inequality, Grothendieck real constant, Grothendieck 
complex constant, Gram matrix, quantum correlation, Bell inequality, Hadamard 
product, Kronecker product, Pearson correlation coefficient, completely correlation 
preserving function, Schoenberg's theorem, Hermite polynomial, Ornstein-Uhlenbeck 
semigroup, Gamma function, Stein's lemma, spherical integration, Gaussian 
hypergeometric function, real and complex Gaussian random vector, Gaussian 
copula, noise stability, operator ideal, Taylor series inversion, Bell 
polynomial}\\[0.5em]
\noindent\msc{Primary 05A10, 15A60, 33C05, 33C45, 41A58, 62H05, 62H20; 
secondary 46A20, 47B10, 81P45.}   
\chapter{Introduction and motivation: {the outstanding story of 
Grothendieck's theorem}}
{\section{Historical perspective and theoretical framework}}
{
In order to better categorise the primary subject of our monograph, we briefly 
examine the history of an essential part of the classical theory of Banach spaces 
in functional analysis, which emerged from A. Grothendieck's seminal article 
``R\'{e}sum\'{e} de la Th\'{e}orie M\'{e}trique des Produits Tensoriels 
Topologiques'' (\cite{G1953}). In this regard, we would like to highlight the 
article \cite{P2015}, which provides a detailed overview of this remarkable 
development. 

In the late thirties, tensor products entered the area of functional analysis, 
due to the works of F. J. Murray, J. von Neumann and R. Schatten. However, it 
was Grothendieck who revealed the structural richness of tensor products of 
Banach spaces and who used their various norms to construct related classes of 
bounded linear operators. In this context, he actually established the origin 
of the ``local'' theory of Banach spaces, i.e., the study of the structure of 
Banach spaces in terms of their finite-dimensional subspaces. Here, it should be 
noted that Grothendieck also \textit{introduced} and characterised the famous 
and highly consequential approximation property of Banach spaces (cf. 
\cite[Section 1.5]{G1955}).

Despite its emergence more than six decades ago, the techniques and results of 
the pioneering work of Grothendieck in \cite{G1953} are still not widely known 
or appreciated, including his main result, mostly known as ``the Grothendieck 
inequality''. Grothendieck himself called it ``the fundamental theorem of the 
metric theory of tensor products''. It is likely that his work was rejected 
because he provided almost no proofs and relied on the (duality) theory of the 
rather abstract (yet very powerful) notion of tensor products of Banach spaces. 
(cf. \cite{AK2016, DF1993, DFS2008, G1953}). \cite{DFS2008} gives a very 
readable and comprehensive account of the tensor product theory, developed in 
\cite{G1953} while maintaining the symbolic language of Grothendieck. The 
culmination in \cite{DFS2008} is Chapter 4, where the Grothendieck inequality 
and its consequences are considered in detail. Actually, \cite{G1953} appeared 
in 1956. 

It was not until the end of the sixties when the scientific community gave 
\cite{G1953} more recognition. The interest in Grothendieck's work namely 
revived when J. Lindenstrauss and A. Pe{\l}czy\'{n}ski recast its main results 
in the more traditional language of operators and matrices, including the 
Grothendieck inequality (\sref{Theorem}{thm:GT_matrix_form}), on which our 
monograph is based (cf. \cite[Theorem A.3.1]{DFS2008} and \cite{LP1968}). They 
presented important applications to the theory of absolutely $p$-summing 
operators and translated results, which were written in terms of tensor 
products by Grothendieck, into properties of linear operators and operator 
ideals. 

Almost at the same time, a general theory of operator ideals on the class of 
Banach spaces was developed by A. Pietsch and his academic school in Jena, yet 
without the use and the abstract language of Grothendieck's tensor norms. Due to 
Pietsch's seminal book ``Operator Ideals'' (\cite{P1980}), that theory became 
a central theme in Banach space theory. Particularly during this time, theory 
and applications of operator ideals had a greater prevalence, as opposed to 
the tensor norm theory of Grothendieck. A comprehensive overview of Pietsch's 
theory and application of operator ideals - including a corresponding 
reformulation of Grothendieck's seminal inequality - is given in 
\cite{DJP2001}.

In 1993, A. Defant and K. Floret published their pathbreaking and comprehensive 
monograph ``Tensor Norms and Operator Ideals'' (\cite{DF1993}). Here, 
deep interconnections between operator ideals, the ``local'' theory of Banach 
spaces and tensor norms are revealed with a high level of attention to detail. 
They made very clear that tensor products and operator ideals are closely 
connected and showed in detail how to transform tensor products to operator 
ideals and conversely, revealing that normed tensor products of Banach spaces 
(in the sense of Grothendieck) and Banach operator ideals (in the sense of Pietsch) 
are ``two sides of the same coin''! Nowadays, many researchers follow the 
approach of Defant and Floret and make use of both languages simultaneously, 
just like we do (cf. \cite{Oe1992, Oe1998, Oe2002, Oe2003}). 
The monographs \cite{DF1993, DFS2008, J1987, R2002} are very valuable sources 
which strongly help to make Grothendieck's approach accessible to a wider 
community.

In conclusion, the} Grothendieck inequality had a profound influence on the 
geometry of Banach spaces and operator theory; particularly between 1970 and 
1990. We highly recommend the readers who have a solid knowledge of functional 
analysis to study Chapter 8 of the superb monograph \cite{AK2016}. Here is 
worked out in great clarity, step-by-step (even without the use of tensor 
products of Banach spaces, and without the use of operator ideal 
theory), how the Grothendieck inequality can be equivalently characterised, 
including Grothendieck's key result, that the inequality 
is equivalent to the deep fact that \textit{any} bounded linear operator $T \in 
{\mathfrak{L}}(L^1(\mu), l_2)$ (where the measure $\mu$ lives on a 
$\sigma$-finite measure space) already is absolutely $1$-summing and 
satisfies the norm inequality $\Vert T \Vert_{{\mathfrak{P}}_1} \leq 
K_G^\F\,\Vert T \Vert$ (cf. \cite[Remark 8.3.2 (b)]{AK2016}, 
\cite[Theorem 23.10]{DF1993}, {\cite{DJ1991, DJP2001}}, 
\cite[Theorem 10.7]{J1987} and \sref{Remark}{rem:GT_vs_abs_summing_operators}). 
An exceptional proof of the latter result (which is built on a factorisation of 
$T \in {\mathfrak{L}}(l_1, l_2)$ over the disc algebra $A(\D)$) is given in 
\cite[Theorem III.F.7]{W1991} (cf. also 
\sref{Remark}{rem:disc_algebra_and_positive_Wiener_algebra} below).

Meanwhile, in addition to this impact, the Grothendieck inequality 
exhibits deep  applications in different fields ({such as algorithmic 
complexity in theoretical computer science}, analysis of Boolean functions, 
random graphs (including the mathematics of the systemic risk in financial networks, 
analysis of nearest-neighbour interactions in a crystal structure (Ising model), 
correlation clustering and image segmentation in the field of computer vision), 
NP-hard combinatorial optimisation, non-convex optimisation and semidefinite 
programming (cf. \cite{GW1995}), foundations and philosophy of quantum mechanics, 
quantum information theory, quantum correlations (cf. 
\sref{Section}{sec:Gram_matrices_and_quantum_correlation}), quantum 
cryptography, communication complexity protocols and even high-dimensional 
private data analysis (cf. \cite{DNT2015})! Also in these fields there exist 
many challenging related open problems.

\begin{theorem}[\textbf{Grothendieck inequality in matrix form}]
\label{thm:GT_matrix_form}
Let $\F \in \{\R, \C\}$. There is an absolute constant $K>0$ such that for 
any $m, n \in \N$, for any $A = (a_{ij}) \in \M_{m,n}(\F)$, any 
$\F$-Hilbert space $H$, and any $(u_1, \ldots, u_m) \in B_H^m$, 
$(v_1, \ldots, v_n) \in B_H^n$, the following inequality is satisfied:
\begin{align*}
\big\vert \sum_{i = 1}^{m}\sum_{j = 1}^{n} a_{ij}\langle u_i, v_j\rangle_H 
\big\vert \leq K\,\sup{\big\{ \big\vert \sum_{i = 1}^{m}\sum_{j = 1}^{n} a_{ij} 
p_i q_j\big\vert : \vert p_i \vert \leq 1, \vert q_j \vert \leq 1 \, 
\forall i,j\big\}}.
\end{align*}
\end{theorem}
\noindent The smallest possible value of the corresponding absolute constant $K$ 
is called the \textit{Grothendieck constant $K_G^\F$} (cf. also 
\autoref{thm:GT_rewritten_form}). The superscripts $\R$ and $\C$ are used 
to indicate the different values in the real and complex cases. Regarding functional 
analytic key reformulations of the Grothendieck inequality, involving the 
infinite-dimensional Banach spaces of type $C(K)$, $C(L)^\prime$ and $L^1(\mu)$, 
we highly recommend the readers to study \cite[Section 2]{P2012}, including the 
detailed and very helpful proof of the equivalence of \cite[Theorem 2.3]{P2012} 
and the Grothendieck inequality in matrix form (on which our paper is based). 
Observe that in the case $m = n = 1$, already $K=1$ satisfies the Grothendieck 
inequality. However, it is well-known that $K_G^\R > K_G^\C > 1$ (cf. also  
\sref{Corollary}{cor:Krivine_KGR_of_2}). In his seminal paper \cite{G1953}, 
Grothendieck proved that $K_G^\R \leq \sinh(\frac{\pi}{2}) \approx 2.301$ 
(within our framework recovered as special case in 
\sref{Example}{ex:Krivine_recovered}). In 1974, Grothendieck's result could be improved 
by R. E. Rietz, who showed that $K_G^\R < 2.261$ (cf. \cite{R1974}). Until 
present (rounded to three digits) the following encapsulation of $K_G^\R$ 
holds; rounded to 3 digits (cf. {\cite{BMMN2013, Kr2023}},  
\sref{Example}{ex:Krivine_recovered} and \sref{Example}{ex:Naor_et_al_s_result}):
\[
1.676 < K_G^\R \stackrel{(!)}{<} \frac{1}{\tfrac{2}{\pi}\,
\sinh^{-1}(1)} = \frac{\pi}{2 \ln(1 + \sqrt{2})} \approx 1.782\,.
\]
The complex constant is strictly smaller than the real one. Namely, if we merge 
the values of the upper bounds of $K_G^\C$ achieved to date (cf. 
\cite{H1987, Kr1977_2, P1978, P2012}, \autoref{thm:GT_Niemi_PSD_case} and our 
approximative calculation of the number $\frac{1}{c^\ast} \approx 1.40449$ at 
the end of \sref{Example}{ex:Haagerup}), we obtain (rounded to three digits): 
\begin{align}\label{eq:known_upper_bounds_of_the_complex_GT_constant}
1 < \frac{4}{\pi} < 1.338 \leq K_G^\C \leq 1.405 < \sqrt{2} < 
e^{1-\gamma} < \frac{\pi}{2} < K_G^\R \leq \sqrt{2}\,K_G^\C\,,
\end{align}
where $\gamma : = \sum_{n=2}^\infty (-1)^n \frac{\zeta(n)}{n} = -\Gamma'(1) 
\approx 0.577$ denotes the Euler-Mascheroni constant. Until present, the 
best-known lowest upper bound of $K_G^\C$ is given by $K_G^\C \leq 1.40491$, 
carried out by U. Haagerup in \cite{H1987} (approximatively achieved again in 
\sref{Example}{ex:Haagerup}). 

Regarding apparently surprising equivalent formulations of 
\autoref{thm:GT_matrix_form} (including their detailed verifications), revealing the 
depth of the structure beneath the ``surface of the inequality'', we refer to 
\cite[Equivalent formulations, p. 109 ff]{J1987}.
 
Computing the {best possible} numerical value of the constants $K_G^\R$ and 
$K_G^\C$ is still an open problem (unsolved since 1953). This is where our 
own research continues. We look for a general framework (primarily build on 
methods originating from (block) matrix analysis (cf. \cite{HJ2013}), 
multivariate statistics with real and complex Gaussian random vectors, theory 
of special functions, modelling of statistical dependence with 
copulas and combinatorics, whose complexity increases rapidly in dimension, 
though) which allows either to give the value of $K_G^\R$, respectively 
$K_G^\C$ explicitly or to approximate these values from above and from below 
at least. However, our approach - which in particular allows a short 
proof of the real and complex Grothendieck inequality, even with J.-L. Krivine's 
upper bound of $K_G^\R$ - confronts us strongly with the question whether the 
seemingly non-avoidable combinatorial complexity actually allows us to determine 
the values of $K_G^\R$, respectively $K_G^\F$ explicitly, or not. A detailed 
description of this research problem can be studied in 
\sref{Section}{sec:RP_1} of our monograph.

If either the size $m \times n$ of the arbitrarily chosen matrix $A \in 
\M_{m,n}(\F)$ or the dimension $d$ of the finite-dimensional Hilbert space $\F_2^d$ 
is predefined, we obtain the corresponding two weakened forms of Theorem 
\ref{thm:GT_matrix_form}: 
\begin{proposition}\label{prop:GT_induced_by_fixed_matrix_size_or_fixed_HS_dimension}
Let $\F \in \{\R, \C\}$. 
\begin{enumerate}
\item For any $d \in \N$ there is a constant $K^\F(d) > 1$ such 
that
\begin{align}\label{eq:GT_for_given_l_2_d}
\big\vert \sum_{i = 1}^{m}\sum_{j = 1}^{n} a_{ij}
\langle u_i, v_j\rangle_{\F_2^d}\big\vert \leq K^\F(d)\,
\sup{\big\{ \big\vert \sum_{i = 1}^{m}\sum_{j = 1}^{n} a_{ij} p_i q_j\big\vert :
\vert p_i \vert \leq 1, \vert q_j \vert \leq 1 \, \forall i,j\big\}}
\end{align} 
for any $m, n \in \N$, for any $A \in \M_{m,n}(\F)$, for any 
$(u_1, \ldots, u_m) \in B_d^m$, and for any $(v_1, \ldots, v_n) \in B_d^n$.
\item For any $(m,n) \in \N^2$ there is a constant $K^\F(m,n) > 1$ 
such that
\begin{align}\label{eq:GT_for_given_m_and_n}
\big\vert \sum_{i = 1}^{m}\sum_{j = 1}^{n} a_{ij}
\langle u_i, v_j\rangle_{H}\big\vert \leq K^\F(m,n)\,
\sup{\big\{ \big\vert \sum_{i = 1}^{m}\sum_{j = 1}^{n} a_{ij} p_i q_j\big\vert :
\vert p_i \vert \leq 1, \vert q_j \vert \leq 1 \, \forall i,j\big\}}
\end{align} 
for any Hilbert space $H$, for any $A \in \M_{m,n}(\F)$, for any 
$(u_1, \ldots, u_m) \in B_H^m$, and for any $(v_1, \ldots, v_n) \in B_H^n$.
\end{enumerate}
\end{proposition}
\noindent Let $K_G^\F(d)$ denote the smallest possible value of the 
corresponding constant $K^\F(d)$, introduced by Krivine (cf. 
\cite[Proposition 20.17]{DF1993}), and let $K_G^\F(m,n)$ be the smallest possible 
value of the constant $K^\F(m,n)$, introduced by B. S. Tsirel'son for $\F = \R$ 
(cf. \cite{T1987} and the detailed elaboration in \cite{L2018, L2020}). 
Consequently, $K_G^\F(d) \leq K_G^\F$ for all $d \in \N$, whence $\sup\limits_
{d \in \N}K_G^\F(d) \leq K_G^\F$. Similarly, it follows that $\sup\limits_
{(m,n) \in \N^2}K_G^\F(m,n) \leq K_G^\F$. It seems to us that the numbers 
$K_G^\F(m,n)$ and $K_G^\F(d)$ in general do not stand in relation to each other. 
Hence, to avoid any risk of confusion, it is important to understand whether 
authors refer to $K_G^\F(m,n)$ or to $K_G^\F(d)$ (or even to $K_G^\F$) in their 
work, when they talk about ``the Grothendieck constant'' (such as it is the 
case in \cite{BDGIKP2023, BBT2011, FHLLZZ2015, FR1994, K2017}). For any $d \in 
\N_3$, explicit \textit{lower} bounds of $K_G^\R(d)$ in closed, analytic form 
are provided in \cite[Theorem 1]{BBT2011} and \cite[Theorem 2.2]{FHLLZZ2015}. 
Very recently, the lower bound of $K_G^\R(3)$ (which is precisely the threshold 
value for the nonlocality of the two-qubit Werner state for projective 
measurements in quantum information theory (cf. \cite{AGT2006}, 
\cite[Section 3]{FHLLZZ2015} and \sref{Example}{ex:Werner_state})) could be 
improved. \cite{BDGIKP2023} namely reveals that $1.4367 \leq K_G^\R(3) \leq 
1.4546$. To achieve this result, however, a high computing power was required. 
In \cite{K2017}, an application of duality in semidefinite programming 
(implemented via the so-called ``convex hull algorithm'' in MATLAB) lead to the 
following values of $K_G^\R(m, n)$: $K_G^\R(5,5) = K_G^\R(4,n) = \sqrt{2}$, 
where $n \in \{4,5,6,7\}$. 
  
Note that $K_G^\F(1) = 1$. Since the sequence $(K_G^\F(d) )_{d \in \N}$ is 
non-decreasing it even follows that $K_G^\F = \lim\limits_{d \to \infty}K_G^\F(d) 
= \sup\{K_G^\F(d): d \in \N\}$ (see  
\sref{Proposition}{prop:K_G_R_resp_K_G_C_as_sup_of_nuclear_norms}). Moreover, 
we may add (see \sref{Corollary}{cor:K_GR_2d_vs_K_GC}): 
\[
K_G^\R(2d) \leq K_G^\R(2)\,K_G^\C(d) = \sqrt{2}\,K_G^\C(d) \text{ for all } d \in \N\,.
\]
In particular, by taking the limit $d \to \infty$, we reobtain $K_G^\R \leq \sqrt{2}
\,K_G^\C$. 

Another important special case of the Grothendieck inequality (known as 
\textit{the little Grothendieck inequality}) appears if just positive semidefinite 
matrices $A$ are considered. Let $k_G^\F$ denote the Grothendieck constant, derived 
from the Grothendieck inequality restricted to the set of all positive semidefinite 
$n \times n$ matrices, with entries in $\F$. Then (cf. 
{\cite[Theorem 3.5.9]{DFS2008}, \cite[Th\'{e}or\`{e}me 4, p. 41]{G1953}, 
\cite[Theorem II and Remark, p. 179]{N1983}, \cite{P1985}} and 
\sref{Remark}{rem:GT_vs_abs_summing_operators}):
\begin{theorem}
\label{thm:GT_Niemi_PSD_case}
{
\[
k_G^\R = \frac{\pi}{2} 
\,\text{ and }\, k_G^\C = \frac{4}{\pi}.
\]
} 
\end{theorem}
\noindent An approximation of the largest lower bounds of both Grothendieck constants 
(which is not the subject of our current work) can be found in 
\cite{D1984}. The real case is studied in \cite{R1992} as well.
\section{Preliminaries, terminology and notation} 
This section serves to provide the foundation upon which our whole work is 
built. To this end, we list the basic notation and symbolic abbreviations used 
throughout our paper. More specific terminology, including terms introduced for 
the first time and related symbolic shortcuts will be introduced on the spot. 
The few remaining symbolic shortcuts which are not explicitly described, are 
either self-explanatory or can be found in any well-established and relevant 
undergraduate textbook in mathematics.

\textit{Numbers and sets} -- As is usual, we denote the set of complex 
numbers by $\C$ and the set of real numbers by $\R$. $\Z$ represents the set of 
all integers and $\N$ stands for the subset of positive integers. We will use the 
symbol $\F$ to denote either the real field $\R$ or the complex field $\C$. 
If we wish to state a definition or a result that is satisfied for either real 
or complex numbers (i.e., if $\F \in \{\R, \C\}$), we simply will make use of 
the letter symbol $\F$. Where there is no risk of confusion, we suppress the 
symbol $\F$. In order to save unnecessary case distinctions, we constantly view 
the set $\R$ as a subset of $\C$, so that $\R = \{z \in \F : z = \overline{z}\} = 
\{z \in \F : \Im(z) = 0\}$. $\T : = \{z \in \C : \vert z \vert = 1\}$ denotes 
the unit circle (``one-dimensional torus''), $\D : = \{z \in \C : \vert z \vert 
< 1\}$ the open unit disk and $\overline{\D} : = \{z \in \C : \vert z \vert 
\leq 1\}$ the closed unit disk. If $F$ is an arbitrary normed space, then 
$S_F : = \{x \in F : \Vert x \Vert_F = 1\}$ denotes its unit sphere and 
$B_F : = \{x \in F : \Vert x \Vert_F \leq 1\}$ its closed unit ball. Thus, $S_\F =
\F \cap \T$. In particular, $S_\R = \{-1,1\}$ and $S_\C = \T$. 
$\N_0 : = \{0\} \cup \N$ denotes the set of all non-negative integers 
(often also somewhat unhappily denoted as $\Z_+$). If $m \in \N$, we put 
$[m] : = \N \cap [1, m] = \{1, 2, \ldots, m\}$ and $\N_m : = \N \cap [m, \infty) = 
\{m, m+1, m+2, \ldots\}$.
%
Fix $n \in \N_0$. In addition to the factorial $n! : = \prod\limits_
{i=1}^n (n-i) \in \N$, the double factorial $n!! \in \N$ will play a dominant 
role. The latter is defined as 
\[
n!! : = \prod\limits_{i=0}^{\floor{\frac{n-1}{2}}} (n-2i),
\] 
where $\R \ni x \mapsto \floor{x} : = \max\{\nu \in \Z : \nu \leq x\}$ denotes 
the floor function (and $\R \ni x \mapsto \ceil{x} : = \min\{\nu \in \Z : x \leq 
\nu\}$ the ceiling function). We adopt the usual approach to include $(-1)!! 
: = 1$ as well. A straightforward proof by induction on $n \in \N_0$, including 
the well-known fact that $\Gamma(\frac{1}{2}) = \sqrt{\pi}$ shows that 
\[
n!! = 
\begin{cases} \hspace{0.7cm}2^{\frac{n}{2}}\,\Gamma(\frac{n}{2}+1) &\text{if } n \text{ is 
even}\\ 
\sqrt{\frac{2}{\pi}}\,2^{\frac{n}{2}}\,\Gamma(\frac{n}{2}+1) &\text{if } n 
\text{ is odd}
\end{cases}\,,
\]
where $\{z \in \F : \Re(z) > 0\} \ni z \mapsto \Gamma(z) : = 
\int_{0}^\infty e^{-t}\,t^{z-1}\,\textup{d}t$ denotes the Gamma function which will play 
an important role in our paper (cf. \cite[Chapter 6.1]{SS2003} and 
\sref{Lemma}{lem:double_factorial_fact}).

\textit{Vectors, matrices, norms and linear operators in general} -- 
Fix $m, n \in \N$. The set of all $m \times n$-matrices with entries in a given 
non-empty subset $S \subseteq \F$ is denoted by $\M_{m,n}(S)$. The matrix ring 
$\M_{n,n}(\F)$ is abbreviated as $\M_n(\F)$. As usual, $e_i \in \F^n$ denotes 
the column vector having a $1$ in the $i$th place and zeros elsewhere. If we wish 
to emphasize the dependence on the dimension $n$ of the vector space $\F^n$, then 
we speak of the set $\{e_1^{(n)}, e_2^{(n)},\dots, e_n^{(n)}\} \subseteq \F^n$ 
(cf., e.g., \eqref{eq:tensor_prod_rep_II}). $I_n : = (e_1 \,\brokenvert\, e_2 \,
\brokenvert\, \cdots \,\brokenvert\, e_n) \in \M_n(\F)$ describes the identity 
matrix. Initially, if not indicated otherwise, any vector (deterministic or 
random) $x \in \F^n$ is set as column vector, so that the allocated row vector 
is decribed by transposition ($x \mapsto x^\top$). Translated into Dirac's bra-ket 
language, which is also used in quantum information theory, it holds that 
$e_i = \vert i-1 \rangle$ and $e_i^\top = \langle i-1 \vert$ ($i \in [n]$). In 
particular, $\vert 0 \rangle = e_1$, $\vert 1 \rangle = e_2$ and $\vert n-1\rangle
\langle 1 \vert = e_n e_2^\top \in \M_{n,2}(\F)$ (cf. \cite{NC2000, SH2018} and 
\eqref{eq:bra_ket_version_of_qubit_tensor_prod}). If $A \in 
\M_{m,n}(S)$ is given, it is sometimes very fruitful to represent the entries of 
$A$ as $A_{ij} : = e_i^\top A e_j = e_j^\top A^\top e_i = (A^\top)_{ji}$, so 
that $A = (a_{ij})$, where $a_{ij} : = A_{ij}$. $\overline{A} \in \M_{m,n}(\F)$ 
is defined as $\overline{A}_{ij} : = \overline{A_{ij}}$, implying that 
$A^\ast : = \overline{A}^\top = \overline{A^\top}$ and $x^\ast : = 
\overline{x}^\top = \overline{x^\top}$. Recall that the Euclidean norm is given 
by $\Vert x \Vert_2 : = \sqrt{x^\ast\,x} = \sqrt{\sum\limits_{i=1}^n 
\vert x_i \vert^2}$ for any $x = (x_1, \ldots, x_n)^\top \in \F^n$. If we equip 
the $n$-dimensional vector space $\F^n$ with the Euclidean inner product, we 
obtain the $n$-dimensional Hilbert space  $\F_2^n : = (\F^n, \langle \cdot, 
\cdot\rangle_2)$, where the inner product on $\F^n$ is given by $\langle x, 
y\rangle_2 : = \langle x, y\rangle_{\F_2^n} : = y^\ast x = \sum\limits_{i=1}^n x_i 
\overline{y_i}$. In particular, $\langle z, w\rangle_{\F_2^1} = z\cdot
\overline{w}$ for all $z, w \in \F$. As usual in mathematics, we adopt the 
convention that any inner product $\langle \cdot, \cdot\rangle_H$, defined on 
an arbitrary Hilbert space $H$, is linear in the first argument and conjugate 
linear in the second one, implying that $\langle x, y \rangle_H^{[P]} : = 
\overline{\langle x, y \rangle_H} = \langle y, x \rangle_H$ ($x,y \in H$) is 
conjugate linear in the first argument and linear in the second one; a rather 
common approach in (quantum) physics. An orthonormal basis in $\F_2^n$ is given 
by the set of vectors $\{e_1, e_2,\dots, e_n\}$; i.e., by the standard basis of 
$\F^n$. Occasionally, if $1 \leq p \leq \infty$, we put $\F_p^n : = 
(\F^n, \Vert\cdot\Vert_p)$, where
\[
\Vert x \Vert_p : = \begin{cases}
(\sum_{i=1}^n \vert x_i\vert^p)^{1/p}& \text{ if } 1 \leq p < \infty\\
\max\{\vert x_i \vert: i \in [n]\}& \text{ if } p=\infty
\end{cases}
\]
denotes the $p$-norm of $x = (x_1, \ldots, x_m)^\top \in \F^n$. If there is no 
risk of confusion regarding $\F$, we simply speak of the space $l_p^n$ (as usual).
As usual, if $n \in \N_2$, then $\S^{n-1}$ denotes the unit sphere in $\R_2^n$. 
Throughout the paper, we also identify any linear operator $T : \F^n  \longrightarrow 
\F^m$ with its representing matrix with respect to the respective standard bases: 
$T \equiv (T_{ij})_{(i,j) \in [m] \times [n]}$. In particular, we have
\[
T_{ij} \equiv e_i^\top Te_j = \langle Te_j, e_i\rangle_{\F_2^m} = 
\langle e_j, T^\ast e_i\rangle_{\F_2^n} \equiv (T^\ast)_{ji} \text{ for all } 
i,j \in [m] \times [n]\,,
\]
where $T^\ast : \F^m_2 \longrightarrow \F^n_2$ is the adjoint operator. 
Furthermore, in the case $n=1$, $\F$ is considered throughout as the 
one-dimensional Hilbert space $(\F_2^1, \langle \cdot, \cdot\rangle_2)$, where 
$\langle z, w\rangle_2 : = z\,\overline{w}$ for all $z,w \in \F$. As usual, $O(n)$ 
denotes the orthogonal group, consisting of all invertible matrices $A \in 
\M_n(\R)$ such that $A^{-1} = A^\top$. $U(n)$ describes the unitary group, 
consisting of all invertible matrices $A \in \M_n(\C)$ such that $A^{-1} = A^\ast$. 
$SO(n) : = \{A \in O(n) : \text{det}(A) = 1\}$ is the special orthogonal group, 
and $SU(n) : = \{A \in U(n) : \text{det}(A) = 1\}$ describes the special unitary 
group.

An important inner product on the $\F$-vector space $\M_{m,n}(\F)$ of 
all $m \times n$-matrices with entries in $\F$, which turns $\M_{m,n}(\F)$ into 
an $mn$-dimensional Hilbert space, is the Frobenius inner product, which is 
defined as follows: if $A = (a_{ij}) \in \M_{m,n}(\F)$ and $B = (b_{ij}) \in 
\M_{m,n}(\F)$, then
\[
\langle A, B\rangle_F : = \text{tr}(A B^\ast) = \text{tr}(B^\ast A) = 
\sum\limits_{i=1}^m \sum\limits_{j=1}^n a_{ij}\,\overline{b_{ij}} = 
\overline{\langle B, A\rangle_F}\,,
\]
where 
\[
\text{tr}(C) : = \sum\limits_{i=1}^n \langle C e_i, e_i\rangle_2 = 
\sum\limits_{i=1}^n c_{ii} = \text{tr}(C^\top) = \overline{\text{tr}(C^\ast)} = 
\overline{\text{tr}(\overline{C})}
\] 
denotes the trace of a given (quadratic) matrix $C = (c_{ij}) \in \M_{n,n}(\F)$. 
One can easily verify the well-known fact that the set of all elementary matrices 
$\{e_i\,e_j^\top : (i,j) \in [m] \times [n]\}$ is an orthonormal basis in the 
$mn$-dimensional $\F$-Hilbert space $(\M_{m,n}(\F), \Vert \cdot \Vert_F)$ 
(since $(e_i\,e_j^\top)_{\alpha\beta} = \delta_{i\alpha}\,\delta_{j\beta}$ for 
all $(\alpha,\beta) \in [m] \times [n]$ and $\text{tr}(xy^\top) = y^\top x$ for 
all $x, y \in \F^n$). We adopt the symbolic notation of the ``$\mathfrak{L}$-community'' 
to represent the set of all bounded linear operators between two normed spaces 
$(E, \Vert\cdot\Vert_E)$ and $(F, \Vert\cdot\Vert_F)$ by $\mathfrak{L}(E,F)$.
As usual, the Banach space $E^\prime : = \mathfrak{L}(E,\F)$ denotes the dual 
space of $E$. (The ``$\mathfrak{B}$-community'', often encountered among 
researchers in the field of $C^\ast$-algebras, uses $\mathfrak{B}(E,F)$ instead, 
so that for example, $\mathfrak{B}(H) = \mathfrak{L}(H,H)$, if $H$ is a given 
Hilbert space (cf., e.g., \cite[Chapter 1.3]{H2006})). Remember that every 
linear operator $T : E_0 \longrightarrow F$ from a finite-dimensional normed 
space $E_0$ to an arbitrary normed space $F$ already is bounded. It should also 
be noted that actually $\langle A, B\rangle_F = \langle A, B\rangle_{HS}$ 
coincides with the Hilbert-Schmidt inner product, defined on the Hilbert space 
$\mathfrak{L}(\F_2^n, \F_2^m)$. To this end, recall that if $H$ and $K$ are 
arbitrarily given $\F$-Hilbert spaces, $T \in 
\mathfrak{L}(H,K)$ is a Hilbert-Schmidt operator if and only if $(\Vert Te_\iota\Vert_K)_
{\iota \in J} \in l^2(J)$ for some orthonormal basis $(e_\iota)_{\iota \in J}$ 
in $H$. Here, $J$ denotes an arbitrary index set which must be neither finite nor 
at most countable (cf., e.g., \cite[Proposition 20.2.7]{J1981}). Hence, if 
$S, T \in \mathfrak{S}_2(H,K)$ are two Hilbert-Schmidt operators, the 
Cauchy-Schwarz inequality implies that
\[
\langle S, T \rangle_{HS} : = \sum\limits_{\iota \in J} \langle Te_\iota, 
Se_\iota\rangle_K = \sum\limits_{\iota \in J} \langle e_\iota, 
T^\ast Se_\iota\rangle_K = \text{tr}(T^\ast S)\,,
\]
is a well-defined inner product on the $\F$-vector space $\mathfrak{S}_2(H,K)$ of 
all Hilbert-Schmidt operators. In fact, it turns $\mathfrak{S}_2(H,K)$ into a 
Hilbert space itself (cf. \cite[Exercises IX.2.19, IX.2.20]{C1990} and 
\cite[Proposition 20.2.7]{J1981}). Let us also note the easy-to-prove fact that
\[
\Pi : (\M_{m,n}(\F), \Vert \cdot \Vert_F)  \stackrel{\cong}{\longrightarrow} 
(\M_{m,n}(\F), \Vert \cdot \Vert_F)^\prime, A \mapsto (B \mapsto \text{tr}(B A^\top))
\]
is an isometric isomorphism, whose inverse is given by $\Pi^{-1} = \Psi$, where
\[
\M_{m,n}(\F) \ni \Psi(t) : = (t(e_i e_j^\top))_{i,j} \text{ for all } t \in 
(\M_{m,n}(\F), \Vert \cdot \Vert_F)^\prime\,.
\]
Obviously, also the canonically defined mapping
\[
\Theta : (\M_{m,n}(\F), \Vert \cdot \Vert_F)^\prime \stackrel{\cong}{\longrightarrow}
(\M_{n,m}(\F), \Vert \cdot \Vert_F)^\prime, t \mapsto (M \mapsto \langle M^\top, 
t \rangle) 
\]
is an isometric isomorphism, implying that the composition of these two isometric 
isomorphisms lead to the finite-dimensional version of \textit{trace duality} 
with respect to the norm $\Vert \cdot \Vert_F = \Vert \cdot \Vert_{HS}$ (cf. 
\cite[Theorem 6.4]{DJT1995}):
\[
\Pi \circ \Theta : (\M_{m,n}(\F), \Vert \cdot \Vert_F)  \stackrel{\cong}{\longrightarrow} 
(\M_{n,m}(\F), \Vert \cdot \Vert_F)^\prime, A \mapsto (B \mapsto \text{tr}(B A)).
\]
\noindent Although it is our intention that the main ideas developed in our 
paper can be captured {without having any knowledge of} advanced functional 
analysis and related operator theory, we will add a few text passages which 
should show how also our approach extends into the area of functional analysis 
and operator ideal theory. Related references will be listed, of course. 
In particular - despite its elegance and power - we intentionally avoid the 
explicit use of the language of abstract tensor products of Banach spaces and 
related tensor norms (originally coined by A. Grothendieck in his seminal 
paper \cite{G1953}) as far as possible. Of course, any attentive reader will 
recognise that tensor products occasionally also are lurking in our framework 
(primarily in the form of concrete Kronecker products of matrices). Remarks in 
this regard could be skipped at the first reading. However, for particularly 
stubborn readers and authors, we strongly refer to \cite{C1990, DF1993, H2006, 
J1987, DJT1995, P1980, W2018}. 

Since (symmetrically) partitioned random vectors and block matrices play 
a key role in our analysis, it is sometimes very useful to transform matrices 
into column vectors by making use of a technique known as matrix vectorisation 
(cf. \cite[Chapter 10]{AM2005}). If  $A = (a_1 \,\brokenvert\, a_2 \,\brokenvert\, 
\cdots \,\brokenvert\, a_n) \in \M_{m,n}(\F)$, with columns $a_j \in \F^m$ 
$(j \in [n])$, then 
\[
\text{vec}(A) : = \text{vec}(a_1, \ldots, a_n) : = (a_1^\top\,\brokenvert\, 
\ldots \,\brokenvert\, a_n^\top)^\top \in \F^{mn}
\]
denotes the \textit{column} vector constructed by stacking the columns of $A$ 
on top of each other. A concise \textit{entrywise} implementable construction 
(built on Euclidean division with remainder) of $\text{vec}(A)$ will be 
studied at the beginning of \sref{Section}{sec:Hadamard_gate_and_Krivine} 
(cf. \eqref{eq:vec_op}). Obviously, 
\[
\text{vec} : (\M_{m,n}(\F), \langle\cdot, \cdot \rangle_F) 
\stackrel{\cong}{\longrightarrow} \F_2^{mn}, A \mapsto \text{vec}(A)
\]
is an isometric isomorphism (between finite-dimensional Hilbert spaces). In 
particular, 
\begin{align}\label{eq:vec_as_isometric_isomorphism}
\langle \text{vec}(A), \text{vec}(B)\rangle_2 = \text{tr}(B^\ast A) \text{ and }
\Vert\text{vec}(A)\Vert_2 =\sqrt{\text{tr}(A^\ast A)}  
\end{align}
for all $A, B \in \M_{m,n}(\F)$. We also need $\text{vec}$'s cousin, the Kronecker 
product of matrices (cf. \cite[Chapter 10]{AM2005}), on which the construction 
of a matrix is based which delivers $\sqrt{2}$ as a lower bound of the real 
Grothendieck constant $K_G^\R$ \textit{and} plays a key role in the foundations 
of quantum mechanics, quantum information and even in evolutionary biology: the 
Walsh-Hadamard transform (cf. \sref{Example}{ex:Hadamard_gate}, 
\sref{Remark}{rem:CHSH_inequalities} and 
\sref{Remark}{rem:appl_in_evolutionary_biology})! The Kronecker product is 
constructed as follows: if $A \in \M_{m,n}(\F)$ and 
$B \in \M_{p,q}(\F)$, then 
\[
\M_{mp, nq}(\F) \ni A \otimes B = \begin{pmatrix}
M_{11} & M_{12} & \cdots & M_{1n}\\
M_{21} & M_{22} & \cdots & M_{2n}\\
\vdots & \vdots & \vdots & \vdots\\
M_{m1} & M_{m2} & \cdots & M_{mn} 
\end{pmatrix},
\]
where $M_{ij} : = a_{ij}B \in \M_{p,q}(\F)$. If $C \in \M_{n,r}(\F)$ and $D \in 
\M_{q,s}(\F)$ are two further matrices, then elementary block matrix 
multiplication instantly results into the well-known fact that
\begin{align}\label{eq:multipl_of_Kronecker_products}
(A \otimes B)(C \otimes D) = AC \otimes BD\,. 
\end{align}
We will also provide a rigorous \textit{entrywise} implementable 
construction of the Kronecker product at the beginning of 
\sref{Section}{sec:Hadamard_gate_and_Krivine} (again built on Euclidean 
division with remainder). In the context of partitioned random vectors, we 
will apply the $\text{vec}$ operator in the following sense: if $x = 
(x_1, \ldots, x_n, x_{n+1}, \ldots, x_{2n})^\top \in \F^{2n}$, then $x = 
\text{vec}(a_1, a_2)$, where $a_1 : = (x_1, \ldots, x_n)^\top \in \F^n$ and 
$a_2 : = (x_{n+1}, \ldots, x_{2n})^\top \in \F^n$.

 Fix $v = (v_1, \ldots, v_l)^\top \in \F^l$ and $A = (a_{ij}) \in 
\M_{m,n}(\F)$. Let $\Psi : [m] \times [n] \longrightarrow  [l]$ and $\Lambda : 
[l] \longrightarrow [m] \times [n]$ be given. In the context of our 
analysis, $\Psi$ should be viewed as a mapping which maps an index $(i,j) 
\in [m] \times [n]$ of the matrix element $a_{ij} \in \F$ to the index $\Psi(i,j) 
\in [l]$ of an allocated vector element. Conversely, $\Lambda$ should be regarded 
as a mapping which maps the index $\alpha \in [l]$ of the vector element 
$v_\alpha \in \F$ to an index $\Lambda(\alpha) = (\Lambda_1(\alpha), 
\Lambda_2(\alpha))\in [m] \times [n]$ of an allocated matrix element. More 
precisely formulated, if we consider the (linear) composition operator 
$C_\psi : \F^l \longrightarrow \M_{m,n}(\F)$, we map the given vector $v \in 
\F^l$ to a matrix $C_\Psi(v) : = v_\Psi \in \M_{m,n}(\F)$ as follows:
\[
(v_\Psi)_{i, j} : = v_{\Psi(i,j)} \text{ for all } (i,j) \in [m] \times [n].
\]
Analogously, we map the matrix $A = (a_{ij}) \in \M_{m,n}(\F)$ to a vector 
$A_\Lambda \in \F^l$, according to the rule
\[
(C_\Lambda(A))_\alpha : = (A_\Lambda)_\alpha : = a_{\Lambda(\alpha)} = 
a_{(\Lambda_1(\alpha), \Lambda_2(\alpha))} \text{ for all } \alpha \in [l],
\]
where $C_\Lambda : \M_{m,n}(\F) \longrightarrow \F^l$ denotes the related 
linear composition operator. Similarly, the matrix $A \in \M_{m,n}(\F)$ can be 
mapped to the matrix $A_\sigma \in \M_{r,s}(\F)$, where $\sigma$ is now a given
mapping of type $\sigma : [r] \times [s] \longrightarrow [m] \times [n]$. 
Observe that $A_\sigma$ consists of entries of the originally given matrix $A$, 
so that we could view $A_\sigma$ as a subordinated matrix of $A$ (cf. 
\sref{Remark}{rem:Hadamard_subord_to_Kronecker}). For example, $A^\top =
A_\tau = C_\tau(A)$, where the mapping $\tau: [n] \times [m] \longrightarrow 
[m] \times [n]$ is defined as transposition: $\tau(\nu,\mu) : = (\mu,\nu)$ 
(cf. \eqref{eq:vec_op}). A very recent application of vectorisation within the 
framework of single quantum systems (where entangled states do not play a 
role), including a related application of the Grothendieck inequality can be 
found in \cite{V2022, V2023}.

\textit{Partitioning, $\C^{2n}$ versus $\R^{4n}$, matrices with special 
form and positive semidefiniteness} -- With regard to the study of a class of 
crucially important partitioned multivariate complex Gaussian random vectors for 
our analysis (cf. \sref{Chapter}{chp:multivar_complex_Gaussian_law}), we firstly 
have to shed some light on the structure of the following two important mappings,
which we will encounter many times in this paper. Similar constructions and 
particular cases are listed in \cite[Chapter 1.2]{AHSE1995} and 
\cite[Problem 1.3.P20 and Problem 1.3.P21]{HJ2013}. To this end, fix $m,n \in \N$ 
and $w \in \C^n$. Put
\begin{equation}\label{eq:canonical_1}
\C^n \ni w = \Re(w) + i \Im(w) \mapsto J_2(w) \equiv J_2^{[n]}(w) := 
\vc{\Re(w)}{\Im(w)} \in \R^{2n} 
\end{equation}
and 
\begin{equation}\label{eq:canonical_2}
M_{m,n}(\C) \ni A \mapsto R_2(A) : = 
\begin{pmatrix}
\Re(A) & -\Im(A)\\
\Im(A) & \Re(A)
\end{pmatrix}
\in \M_{2m,2n}(\R) 
\end{equation}
Observe, that if $n \geq 2$, $J_2(w) = {\text{vec}}(\Re(w), \Im(w)) =
{\text{vec}}(\Re(w_1), \ldots, \Re(w_n), \Im(w_1), \ldots, \Im(w_n))$ in general 
does not coincide with ${\text{vec}}(J_2(w_1), J_2(w_2), \ldots, J_2(w_n))$. 
However, given arbitrary $z, w \in \C^n$, we obtain an important equality which 
will be applied several times in our paper; namely:   
\begin{align}\label{eq:impact_of_G}
\begin{split}
J_2^{[2n]}\vc{z}{w} &= {\text{vec}}(\Re(z), \Re(w), \Im(z), \Im(w))\\ 
&= G\,{\text{vec}}(\Re(z), \Im(z), \Re(w), \Im(w)) = G\,\vc{J_2(z)}{J_2(w)}\,,
\end{split}
\end{align}
 where
\begin{align}\label{eq:orthog_matrix}
G \equiv G_n : =
\begin{pmatrix}
I_n & 0 & 0 & 0\\
0 & 0 & I_n & 0\\
0  & I_n & 0 & 0\\
0 & 0 & 0 & I_n
\end{pmatrix} = G^\top = G^{-1} \in O(4n) \text{ is orthogonal}\,.
\end{align}
Observe that the matrix $G_1 \in O(4)$ precisely coincides with the ``swap gate'' 
(also known as ``flip operator''), used in quantum information theory (cf. 
\cite[Problem 28]{SH2018} and \sref{Example}{ex:Werner_state}). In general, 
if ${\text{vec}}(x_1, x_2, x_3, x_4) \in \R^{4n}$ is given, then $G \equiv G_n$ 
swaps the column vectors $x_2 \in \R^{n}$ and $x_3 \in \R^{n}$ and maps 
${\text{vec}}(x_1, x_2, x_3, x_4)$ to ${\text{vec}}(x_1, x_3, x_2, x_4)$. Also 
the matrix $G$ will be needed repeatedly; for example, in the proofs of  
\sref{Lemma}{lem:correl_Gaussians_I} and 
\sref{Corollary}{cor:inheritage_prop_of_Sigma_2n_rho}.

Not only in this context, we occasionally need the following construction. Let 
$f, g : \F^k \longrightarrow \F$ be two arbitrary functions ($k \in \N$). Then 
the function $f \otimes g : \F^{2k} \longrightarrow \F$ is (purely symbolic) 
defined as 
\[
(f \otimes g)({\text{vec}}(x,y)) : = f(x)g(y)\,\text{ for all }\, x,y \in \F^k\,.
\]
Clearly the mapping $J_2 : \C^n \longrightarrow \R^{2n}$ is bijective, 
and $\Vert J_2(a)\Vert_{\R_2^{2n}} = \Vert a \Vert_{\C_2^n}$ for any $a \in \C^n$. 
Moreover,
\begin{align}\label{eq:J_2_of_conj_complex}
J_2(\overline{a}) = \begin{pmatrix}
I_n & 0\\
0 & -I_n
\end{pmatrix}J_2(a) \text{ for all } a \in \C^n\,.
\end{align}
Recall that $\S^{\nu-1}$ denotes the unit sphere in $\R_2^\nu$ ($\nu \in \N_2$). Thus, 
$S_{\C_2^n} = J_2^{-1}(\S^{2n-1})$ describes the unit sphere in $\C_2^n$. Let 
$z, w \in \C^n$. Then, 
\begin{align}\label{eq:Re_part}
\Re(\langle w, z \rangle_{\C_2^n}) = \Re(z^\ast\,w ) = J_2(z)^\top J_2(w) =
\langle J_2(w), J_2(z) \rangle_{\R_2^{2n}}
\end{align}
induces an inner product which turns $\C^n$ into a \textit{real} 
finite-dimensional Hilbert space and $J_2$ into an isometric isomorphism 
between the \textit{real} Hilbert spaces $(\C^n, \Re(\langle \cdot, \cdot 
\rangle_{\C_2^n}))$ and $\R_2^{2n}$. Clearly, $J_2$ cannot be extended to a linear 
mapping between $\C^n$ and $\C^{2n} \supseteq \R^{2n}$. However, by construction $J_2^{-1} : 
\R^{2n} \longrightarrow \C^n$ clearly satisfies
\begin{align}\label{eq:complex_matrix_extension}
x_1 + i\,x_2 = J_2^{-1}x = (I_n \, \brokenvert\,i I_n)x \text{ for all } x = 
\vc{x_1}{x_2} \in \R^{2n} \cong \R^n \times \R^n\,,
\end{align}
implying that the linear and non-injective mapping between the complex vector spaces 
$\C^{2n}$ and $\C^n$, induced by the matrix $(I_n \, \brokenvert\,i I_n) \in 
\M_{n, 2n}(\C)$ actually is a linear extension of $J_2^{-1}$. Moreover, it 
follows that
\begin{align}\label{eq:Im_part}
\Im(\langle w, z \rangle_{\C_2^n}) = \Im(z^\ast\,w ) = 
\Re(z^\ast(-i w) ) = J_2(z)^\top J_2(-i w) =
\langle J_2(-i w), J_2(z) \rangle_{\R_2^{2n}}
\end{align}
and
$J_2(Aw) = R_2(A)\,J_2(w)$ and $R_2(rA) = r R_2(A)$ for all $r \in \R$, $z, w \in \C^n$ and 
$A \in \M_{m,n}(\C)$. In particular, the following diagram commutes
\begin{figure}[h]
\centering
\setlength{\unitlength}{0.4cm}
\begin{picture}(20,5)\thicklines 
\put(6.8,4){\vector(1,0){5.6}} 
\put(5.5,3.8){$\C^n$}
\put(12.7,3.8){$\C^m$}
\put(6,3.5){\vector(0,-1){3.5}} 
\put(5.4,-0.9){$\R^{2n}$}
\put(13,3.5){\vector(0,-1){3.5}} 
\put(12.8,-0.9){$\R^{2m}$\hspace{1mm}.}
\put(7.1,-0.6){\vector(1,0){5.4}} 
\put(8.9,4.3){$A$}
\put(13.2,1.7){$J_2$}
\put(8.3,-0.1){$R_2(A)$}
\put(4.7,1.7){$J_2$}
\end{picture}
\caption{\small{$R_2(A) = J_2 \circ A \circ J_2^{-1}$}}
\label{fig:CD}
\end{figure}

\noindent Note also that $R_2(I_n) = I_{2n}$ and
\[
R_2(GH) = R_2(G)R_2(H) \text{ for all } (G, H) \in \M_{k,m}(\C) \times \M_{m,n}(\C)\,.
\]
Moreover, 
\begin{align}\label{eq:from_Hermitian_to_symmetric}
R_2(A^\ast) = R_2(A)^\top \text{ for any } A \in \M_{m,n}(\C)\,.
\end{align}
In particular,
\[
R_2(A^\top)^\top = R_2(\overline{A}) = \begin{pmatrix}
I_n & 0\\
0 & -I_n
\end{pmatrix}R_2(A)\begin{pmatrix}
I_n & 0\\
0 & -I_n
\end{pmatrix} \text{ for any } A \in \M_{m,n}(\C)\,.
\]
Thus, from the algebraic viewpoint, $R_2 : M_n(\C) \longrightarrow \M_{2n}(\R)$ 
is an injective unitary ${}^\ast$-ring homomorphism. In particular, $A \in M_n(\C)$ 
is invertible if and only if $R_2(A) \in \M_{2n}(\R)$ is invertible (since 
$R_2(A)R_2(A^{-1}) = I_{2n} = R_2(A^{-1})R_2(A)$ for any $A \in GL(n; \C)$). In 
the case $n=1$, we reobtain the well-known Abelian group isomorphism $R_2 : \T 
\stackrel{\cong}{\longrightarrow} SO(2)$. Consequently, it follows that
\begin{align}\label{eq:complex_form}
\begin{split}
\Re(\langle A\,w,z\rangle_{\C_2^m}) &= \Re(z^\ast A\,w) = 
J_2(z)^\top R_2(A)J_2(w) = \langle R_2(A)J_2(w), J_2(z)\rangle_{\R_2^{2m}}\text{ and }\\
\Im(\langle A\,w,z\rangle_{\C_2^m}) &= \Im(z^\ast A\,w) =
-J_2(z)^\top R_2(i\,A)J_2(w) = -\langle R_2(i\,A)J_2(w), J_2(z)\rangle_{\R_2^{2m}}
\end{split} 
\end{align}
for all $(z,w) \in \C^m \times \C^n$ and $A \in \M_{m,n}(\C)$. 
\begin{lemma}\label{lem:skew_sym}
Let $n \in \N$ and $C \in M_n(\R)$. Then the following statements are equivalent
\begin{enumerate}
\item $C$ is skew symmetric.
\item $x^\top Cy = -y^\top Cx$ for all $x,y \in \R^n$.
\item $z^\top Cz = 0$ for all $z \in \C^n$.
\item $x^\top Cx = 0$ for all $x \in \R^n$.
\end{enumerate}
\end{lemma}
\begin{comment}
Since $C$ is skew symmetric, it follows that $x^\top Cy = (x^\top Cy)^\top = 
y^\top (-C) x$ for all $x,y \in \R^n$, wherefrom the implication $(i) 
\Rightarrow (ii)$ follows. $(ii)$ obviously implies condition $(iii)$: we just 
have to develop $(x^\top + i y^\top)C(x + iy)$ for arbitrary $x,y \in \R^n$ and 
apply (ii) to the four consecutive factors. $(iii) \Rightarrow (iv)$ is trivial. 
Assume that the hypothesis $(iv)$ holds. Let $x, y \in \R^n$ be given. Then 
$x^\top Cx = 0$, $y^\top Cy = 0$ and $(x^\top + y^\top)C(x+y) = 0$. Consequently,
\[
0 = (x^\top + y^\top)C(x+y) = x^\top Cy + y^\top Cx = x^\top (C+C^\top)y\,,
\]
and (i) follows.
\end{comment}
\noindent Combining \sref{Lemma}{lem:skew_sym} with 
\eqref{eq:from_Hermitian_to_symmetric} and \eqref{eq:complex_form}, we 
immediately obtain another neat result, including a full characterisation of 
Hermitian matrices $A=A^\ast \in M_n(\C)$ by their symmetric real representation 
$R_2(A) = R_2(A)^\top \in \M_{2n}(\R)$.
\begin{proposition}\label{prp:Hermitian_vs_symmetric}
Let $n \in \N$ and $\Gamma = \Re(\Gamma) + i\,\Im(\Gamma) \in M_n(\C)$. Then the following statements are 
equivalent:
\begin{enumerate}
\item $\Gamma$ is Hermitian.
\item $i\Gamma$ is skew Hermitian.
\item $\Re(\Gamma) \in M_n(\R)$ is symmetric and $\Im(\Gamma) \in M_n(\R)$ is skew symmetric.
\item $R_2(\Gamma)$ is symmetric.
\item $z^\ast \Gamma\,z \in \R$ for all $z \in \C^n$.
\item $R_2(i\,\Gamma)$ is skew symmetric.
\end{enumerate}
 
In particular, if $\Sigma = 
\begin{pmatrix}
A & B\\
C & D
\end{pmatrix} \in \M_{2}(M_n(\R))$, then the following applies: 
\begin{align}\label{eq:char_of_R_2_Gamma}
\begin{split}
{}&\Sigma = R_2(\Gamma) \text{ for some Hermitian matrix } \Gamma \in 
\M_n(\C) \text{ if and only if } A = D\\ 
{}&\text{and } B + C = 0 \text{ and } \Sigma \text{ is symmetric}.
\end{split}
\end{align}
Thereby, the uniquely defined Hermitian matrix is given by $\Gamma = A + i C$.
\end{proposition}
\noindent Another important implication refers to the role of the matrix $R_2(A)$
regarding a full clarification of the reason for the difference between the 
structure of positive semidefinite matrices in $\M_n(\C)$ and the structure of 
positive semidefinite matrices in $\M_n(\R)$ (cf. e.g. \cite[Theorem 4.1.10]{HJ2013} 
and the field-independent definition in the form of 
\sref{Lemma}{lem:uniform_psd_char} below): 
\begin{corollary}\label{cor:char_of_complex_PSD}
Let $n \in \N$ and $A = \Re(A) + i\,\Im(A) \in M_n(\C)$. Then the following statements are 
equivalent:
\begin{enumerate}
\item $z^\ast A\,z \geq 0$ for all $z \in \C^n$.
\item $z^\ast A\,z \geq 0$ for all $z \in \C^n$ and $A$ is Hermitian.
\item $x^\top R_2(A)\,x \geq 0$ for all $x \in \R^{2n}$ and $R_2(A)$ is symmetric.
\end{enumerate}
If in addition $\Im(A) = 0$, then (i) is equivalent to
\begin{enumerate}
\item[\textup{(i${}^\prime$)}] $x^\top A\,x \geq 0$ for all $x \in \R^n$ and 
$A$ is symmetric.
\end{enumerate}
\end{corollary}
\noindent \sref{Corollary}{cor:char_of_complex_PSD} reveals the role of symmetry 
in the established definition of a positive semidefinite real matrix. For example, 
if we consider the non-symmetric real matrix $A : = \begin{pmatrix}
0 &  1\\
-1 & 0
\end{pmatrix}$ and $z : = (1,i)^\top \in \C^2$, then $x^\top A\,x = 0$ for all 
$x \in \R^2$, but $z^\ast A\,z = 2i$. Thus, throughout the paper, we 
apply the following characterisation of positive semidefinite (respectively 
positive definite) matrices, which does not depend on the choice of the 
field $\F \in \{\R, \C\}$: 
\begin{lemma}\label{lem:uniform_psd_char}
Let $\F \in \{\R, \C\}, n \in \N$ and $A \in M_n(\F)$. $A$ is 
positive semidefinite (respectively positive definite) in $\M_n(\F)$ if the 
following two conditions are satisfied:
\begin{enumerate}
\item $A = A^\ast$.
\item $z^\ast A z \geq 0$ (respectively $z^\ast A z > 0$) for all 
$z \in {\F^n}\setminus\{0\}$.
\end{enumerate}
In particular, $A \in M_n(\R)$ is positive semidefinite in $M_n(\R)$, if and 
only if $A \in M_n(\R) \subseteq M_n(\C)$ is positive semidefinite in $M_n(\C)$.
\end{lemma}
\begin{remark}
If we identify (bounded) linear operators $A \in {\mathfrak{L}}(\F^n, \F^n)$ and 
matrices $A \in M_n(\F)$, then $A \in {\mathfrak{L}}(\F^n, \F^n)$ is positive 
semidefinite if and only if $A$ is a positive self-adjoint operator. Since any 
positive operator $A \in {\mathfrak{L}}(\C^n, \C^n)$ is self-adjoint, positivity 
coincides with positive semidefiniteness on ${\mathfrak{L}}(\C^n, \C^n)$, in 
contrast to positivity on ${\mathfrak{L}}(\R^n, \R^n)$; i.e., there are positive 
non-symmetric operators $B \in {\mathfrak{L}}(\R^n, \R^n)$ (such as $B : = 
\begin{pmatrix}
0 & 1\\
-1 & 0
\end{pmatrix}$), implying that these operators cannot be positive semidefinite 
in $M_n(\C)$.  
\end{remark}
\noindent Given $\emptyset \neq S \subseteq \F$, we put (cf. \cite{GMS2015})
\[
\M_n(S)^+ : = \{A : A \in \M_n(S) \text{ and } A \text{ is psd in } \M_n(\F)\}\,.
\]
Thus, $M_n(\C)^+ = \{A \in M_n(\C) : z^\ast A z \geq 0 \text{ for all } z \in 
\C^n\}$ and $\M_n(\R)^+ = \{A \in M_n(\R) : A = A^\top \text{ and } x^\top A x 
\geq 0 \text{ for all } x \in \R^n\}$. Moreover, $A \in \M_n(\C)^+$ if and only 
if $R_2(A) \in \M_{2n}(\R)^+$. Here, the subclass of all correlation matrices, 
i.e., of all psd matrices with ones on their diagonal (cf. 
\cite[Definition 2.14.]{K2022} and \sref{Lemma}{lem:charact_of_corr_matrices}) 
plays the main role in our work. Only through their structure, including the 
deep impact of correlation-preserving mappings (cf. 
\sref{Definition}{def:nCP_and_CCP}) our main results could be developed. We 
actually work with exactly those correlation matrices that are used in 
statistics. So, our approach could also be interesting for the statistical 
community; especially for those researchers who are working in spatio-temporal 
modelling and functional data analysis (FDA).
 
All basic properties of positive semidefinite (respectively, positive 
definite) matrices including the ``striking if not almost magical'' structure of 
related $2 \times 2$ block matrices, used and listed throughout our paper (without 
giving any proof) can be found in \cite[Chapter 1]{B2007}. A further, very 
detailed analysis of the convex psd cone and its geometry, considered from the 
point of view of convex optimisation is listed in \cite[Chapter 2.9]{D2019}. 
(Note also that in addition to the symbol ``\,$\M_n(\R)^+$\,'', the terms 
``\,${\S_n^+}$\,'' and ``\,${\P_n}$\,'' are often found in the literature.)

\textit{Measurability, probability, random vectors} -- If not 
specified differently, $(\Omega, \mathscr{F})$ always denotes a measurable space 
which is not specified in more detail. However, we have to make use of different 
probability spaces, including $(\F^k, \mathcal{B}(\F^k), \gamma_k^\F)$, where 
$\gamma_k^\F$ denotes the real or complex Gaussian measure, described in detail 
in \sref{Section}{sec:complex_Gaussian_rvcs}. As usual, if $\P$ is a given 
probability measure on some $(\Omega, \mathscr{F})$ and $X : \Omega 
\longrightarrow \F$ a $\P$-integrable $\F$-valued random variable, then 
(cf., e.g., \cite{B1996})
\[
\E_\P[X] : = \E_\P[\Re(X)] + i\,\E_\P[\Im(X)] : = 
\int_{\Omega} \Re(X)\,\textup{d}\P + i\,\int_{\Omega} \Im(X)\,\textup{d}\P\,.
\]
In order not to unnecessarily complicate readability, we use the symbols 
$\textup{d}^n x$ and $\lambda_n$ interchangeably to denote the real 
$n$-dimensional Lebesgue measure (e.g., $\int_{\R^n} f\,
\textup{d}\lambda_n = \int_{\R^n} f(x)\lambda_n(\textup{d}x) = 
\int_{\R^n} f(x)\,\textup{d}^n x$). Unless otherwise stated, random \textit{variables} 
will be denoted by capital letters (such as $X : \Omega \longrightarrow \R$, or 
$Z : \Omega \longrightarrow \C$), whereas random \textit{vectors} will be denoted 
by bold capital letters (such as $\textbf{X} : \Omega \longrightarrow \R^n$ or 
$\vc{\textbf{Z}}{\textbf{W}} : \Omega \longrightarrow \C^{2n}$). $\textbf{X} 
\stackrel{d}{=} \textbf{Y}$ stands for the equality $\P_{\textbf{X}} = 
\P_\textbf{Y}$ of the respective probability laws.

 Within the framework of standard measure theory (including 
classical $L^p$-spaces), we will tacitly assume that we always are working 
with equivalence classes of almost everywhere coinciding $\F$-valued functions,
respectively vector valued measurable mappings on some underlying measure space 
$(\Omega, {\mathscr{F}}, \mu)$. However, since paths of stochastic processes will 
not play a role in this paper, we do not have to pay special attention to the 
structure of null sets. In this regard, a typical example is the real-valued 
signum function:
\[
\text{sign}(x) : = \begin{cases}
1& \text{if } x > 0\\
0& \text{if } x=0\\
-1& \text{if } x < 0
\end{cases}\,.
\]
If we namely view $\text{sign}$ as an element of $L^\infty(\R)$ (with
$\Vert \text{sign} \Vert_\infty = 1$), it follows that
\[
\text{sign} = 2\,\ind_{[0, \infty)} - 1 \text{ in } L^\infty(\R)
\]
(since $\{0\}$ is a Lebesgue null set), where $\ind_A$ denotes the 
indicator function of $A$. Observe that $H : = \ind_{[0, \infty)}$ 
is also well-known as Heaviside step function, which is especially used for 
applications of Fourier analysis in electricity engineering. This perspective 
will become an important part of our approach (cf. 
\sref{Example}{ex:nth_deriv_of_Heaviside_step_fct}).   

Finally, let us remark, that we also make use of the purely symbolic 
notation $x \equiv y$ to indicate that $x$ can be canonically identified with 
the quantity $y$ (such as $\M_{n,1}(\F) \equiv \F^n$) or that it is just a 
shortcut for the previously rigorously defined quantity $y$ (cf., e.g., 
\eqref{eq:complex_Lebesgue_measure}).
\chapter{Complex Gaussian random vectors and the probability law 
$\C N_{2n}(0, \Sigma_{2n}(\zeta))$}\label{chp:multivar_complex_Gaussian_law}
\section{General complex Gaussian random vectors in $\C^n$ and their probability 
distribution}\label{sec:complex_Gaussian_rvcs}
 
Regarding a deeper analysis of the underlying structure of the Grothendieck 
inequality in the complex case, including the Haagerup equality, it is very 
helpful to work with centred random vectors whose probability law follows 
the multivariate \textit{complex} Gaussian distribution, fully characterised 
through certain correlation matrices, whose entries are elements of $\overline{\D}$. 
This approach allows us to generalise the Haagerup equality by \textit{substituting} the complex sign function, chosen by Haagerup (see 
\cite{H1987}), through arbitrary ``circularly odd'' functions $b : \C^k \longrightarrow \T$, 
where $k \in \N$ (see \sref{Corollary}{cor:Haagerup_incl}).

It is far beyond the scope of the present contribution, to recall the 
rich structure of the multivariate complex Gaussian distribution in detail. 
However, for the convenience of the readers we list and describe the whole 
properties of that probability law which are implemented in some of our following 
proofs in relation to the complex version of the Grothendieck inequality and 
beyond in a self-contained way. We highly recommend the readers to study related chapters in 
the references \cite{AHSE1995, Gn1963}, respectively \cite[Chapter 2.1]{HKPV2009} and 
\cite[Appendix E.2]{HvNVW2017}, where that class of random vectors and their distribution 
functions is comprehensively and rigorously introduced including the related 
symbolic (mostly self-explaining) notation. Significant facts about real Gaussian random 
vectors are also listed and discussed thoroughly in \cite[Chapter 5.II.1]{LQ2018} and 
\cite[Chapter 1.10]{S2013}. In \cite[Chapter 30]{B1996}, a real Gaussian random vector  
is viewed and studied as a special case of a measurable mapping between a probability space 
and a measurable space. Recall the powerful general characterisation of the 
Gaussian law of random vectors with values in $\R^n$ (cf. e.g. 
\cite[Theorem 30.2]{B1996}):  
\begin{proposition}\label{prop:n_dim_real_Gaussian_law}
Let $(\mu, \Sigma) \in \R^n \times \M_n(\R)^+$ and $\textbf{X} = 
{\text{vec}}(X_1, X_2, \ldots, X_n)$ be a random vector in $\R^n$. Then the 
following statements are equivalent:
\begin{enumerate}
\item \[
\textbf{X} \sim N_n(\mu, \Sigma);
\]
\item For all $a \in \R^n$, 
\[
a^\top\textbf{X} = \sum_{i=1}^n a_i\,X_i \sim N_1(a^\top\mu, a^\top\Sigma\,a);
\]
\item The characteristic function of $\textbf{X}$ is given by
\[
\R^n \ni a \mapsto \phi_\textbf{X}(a) : = \E[\exp(ia^\top \textbf{X})] =
\exp(i a^\top \mu - \tfrac{1}{2}a^\top\Sigma a).
\] 
\end{enumerate}
\end{proposition}
\noindent Consequently, the following important fact (which we apply in this 
paper frequently) follows at once:  
\begin{align}\label{eq:linear_transf_of_Gaussian_law}
\text{If }\textbf{X} \sim N_n(\mu, \Sigma), \text{ then } A\textbf{X} + b 
\sim N_n(A\mu+b, A\Sigma A^\top) \text{ for all } m \in \N \text{ and } 
(b, A) \in \R^m \times \M_{m,n}(\R). 
\end{align}  
Furthermore, recall that a random vector $\textbf{Z} = 
{\text{vec}}(Z_1, Z_2, \ldots, Z_n)$ which maps into $\C^n$ is a complex random vector 
if for all $\nu \in [n]$ $Z_\nu = X_\nu + i\,Y_\nu$, where $X_\nu = \Re(Z_\nu)$ and 
$Y_\nu = \Im(Z_\nu)$ both are real random variables (each one 
defined on the same probability space). Along the lines of the notation for (deterministic) 
vectors in $\C^n$ one puts 
\[
\textbf{X} \equiv \Re(\textbf{Z}) : = {\text{vec}}(\Re(Z_1), \Re(Z_2), \ldots, 
\Re(Z_n)) = {\text{vec}}(X_1, X_2, \ldots, X_n)
\]  
and
\[
\textbf{Y} \equiv \Im(\textbf{Z}) : = {\text{vec}}(\Im(Z_1), \Im(Z_2), \ldots, 
\Im(Z_n)) = {\text{vec}}(Y_1, Y_2, \ldots, Y_n),
\]
implying that $\textbf{Z} = \textbf{X} + i\,\textbf{Y} = \Re(\textbf{Z}) + 
i\,\Im(\textbf{Z})$. Let 
\[
\lambda_n^{\C} : = (J_2^{-1})_\ast\lambda_{2n}
\] 
be the Lebesgue measure on $\C^n$ (i.\,e., the image measure of the real Lebesgue measure 
$\lambda_{2n}$). Fix $0 < p < \infty$. If $z = x + i y \in \C$, then
\[
\vert z\vert^p = (x^2 + y^2)^{p/2} \leq \max\{2^{(p/2)-1},1\}(\vert x\vert^p + 
\vert y\vert^p) \leq \max\{2^{p/2},2\}\vert z\vert^p
\] 
(see \cite[2.10.E]{J1981}). Consequently, the change-of-variables formula 
(cf. e.g. \cite[Chapter 19]{B2001}), applied to the image measure 
$\lambda_n^{\C}$, implies that 
\[
h \in L^p(\C^n, \lambda_n^{\C}) \text{ if and only if } 
\Re(h)\circ J_2^{-1} \in L^p(\R^{2n}, \lambda_{2n}) \text{ and } 
\Im(h)\circ J_2^{-1} \in L^p(\R^{2n}, \lambda_{2n}).
\]
The construction of $\lambda_n^{\C}$ namely implies that 
\begin{align}\label{eq:complex_Lebesgue_measure}
\int_{\C^n} \vert h(z)\vert^p\,\lambda_n^{\C}(\textup{d}z) \equiv 
\int_{\C^n} \vert h\vert^p\,\textup{d}\lambda_n^{\C} = 
\int_{\R^{2n}} \vert h\vert^p \circ J_2^{-1}\,\textup{d}\lambda_{2n}
\equiv \int_{\R^{2n}} \vert h(x+iy)\vert^p\,\lambda_{2n}(\textup{d}(x,y)). 
\end{align}
Equipped with these basic, well-known facts about the Lebesgue measure on 
$\C^n$, we reintroduce complex Gaussian random vectors in the following, 
seemingly elementary way:
\begin{definition}[\textbf{Complex Gaussian random vector}]
An $n$-dimensional complex random vector $\textbf{Z}$ is a complex Gaussian 
random vector if the real $2n$-dimensional random vector $J_2(\textbf{Z}) =
\vc{\Re(\textbf{Z})}{\Im(\textbf{Z})}$
is a real Gaussian random vector.
\end{definition}
\noindent Although that definition of complex Gaussian random vectors seems to be a quite 
inconspicuous one, it encapsulates a rich underlying structure which strongly differs from 
that one of real Gaussian random vectors. Firstly, without having to know any further details 
about the structure of complex Gaussian random vectors, the change-of-variables formula 
(cf. e.g. \cite[Chapter 19]{B2001}) implies that 
\[
\E[g(\textbf{Z})] \equiv \E[g \circ \textbf{Z}] = 
\int_{\Omega} g\circ \textbf{Z}\,\textup{d}\P = 
\int_{\C^n}g\,\textup{d}\P_{\textbf{Z}} = 
\int_{\C^n} g\,\textup{d}\,(J_2^{-1})_\ast\P_{J_2(\textbf{Z})}
\] 
can be written as 
\[
\E[g(\textbf{Z})] = \E_{\P_{\textbf{X}}}[g \circ J_2^{-1}] =
\E[\Re(g(J_2^{-1}(\textbf{X})))] + 
i\,\E[\Im(g(J_2^{-1}(\textbf{X})))]
\]
for any $\P_{\textbf{Z}}$-integrable function $g = \Re(g) + i\,\Im(g)$ and any 
complex Gaussian random vector $\textbf{Z}$, where $\textbf{X} \stackrel{d}{=} 
J_2(\textbf{Z})$ is a real $2n$-dimensional Gaussian random vector. Consequently, 
the expectation vector $\mu \equiv \E[\textbf{Z}] : = {\text{vec}}
(\E[Z_1], \E[Z_2], \ldots, \E[Z_n])$ as well as the variance matrix $\Gamma \equiv 
\text{var}(\textbf{Z}) : = \E[(\textbf{Z}-\mu)(\textbf{Z}-\mu)^\ast] = 
(\E[(Z_i-\mu_i)(\overline{Z_j} - \overline{\mu_j})])_{1 \leq i,j \leq n}
\in \M_n(\C)^+$ and the cross-covariance matrix $C \equiv \text{cov}(\textbf{Z}, 
\overline{\textbf{Z}}) : = \E[(\textbf{Z}-\mu)(\overline{\textbf{Z}-\mu})^\ast] = 
\E[(\textbf{Z}-\mu)(\textbf{Z}-\mu)^\top] = (\E[(Z_i-\mu_i)(Z_j-\mu_j)])_
{1 \leq i,j \leq n} \in M_n(\C)$ are well-defined. Let $S \in \M_{2n}(\R)^+$ be 
the variance matrix of $J_2(\textbf{Z})$. Then 
\[
J_2(\textbf{Z}) \sim N_{2n}(J_2(\mu), S)\,,
\]
where $S = \E[J_2(\textbf{Z}-\mu)J_2(\textbf{Z}-\mu)^\top]$. 
A straightforward computation of $C + \Gamma = 
2\E[(\textbf{Z}-\mu)(\Re(\textbf{Z}-\mu))^\top]$ and 
$C - \Gamma = 2i\,\E[(\textbf{Z}-\mu)(\Im(\textbf{Z}-\mu))^\top]$ 
implies that
\begin{align}\label{eq:the_matrix_S}
\begin{split}
2\,S 
&= \begin{pmatrix}
\Re(C+\Gamma) & \Im(C-\Gamma)\\
\Im(C+\Gamma) & -\Re(C-\Gamma)
\end{pmatrix}
=
R_2(\Gamma) + 
\begin{pmatrix}
\Re(C) & \Im(C)\\
\Im(C) & -\Re(C)
\end{pmatrix}\\
&=
\Lambda_{2n}^\ast 
\begin{pmatrix}
\Gamma & C\\
\overline{C} & \overline{\Gamma}
\end{pmatrix}
\Lambda_{2n}\,,
\end{split}
\end{align}
where
$
\Lambda_{2n} : = 
\frac{1}{\sqrt{2}}
\begin{pmatrix}
I_n & i\,I_n\\
I_n & -i\,I_n
\end{pmatrix} 
\in U(2n)$ is an unitary matrix (with $\det(\Lambda_{2n}) = (-i)^n$). 
Observe that 
\[
\begin{pmatrix}
\Gamma & C\\
\overline{C} & \overline{\Gamma}
\end{pmatrix} = \E[\textbf{W}\textbf{W}^\ast],
\]
where the complex random vector $\textbf{W} : = 
\vc{\textbf{Z}-\mu}{\overline{\textbf{Z}-\mu}} = 
\vc{\textbf{Z}}{\overline{\textbf{Z}}} - \widetilde{\mu}$, with 
$\widetilde{\mu} := \vc{\mu}{\overline{\mu}}$, maps into $\C^{2n}$. 
Consequently, 
$\begin{pmatrix}
\Gamma & C\\
\overline{C} & \overline{\Gamma}
\end{pmatrix}$ is the variance matrix of the random vector 
$\vc{\textbf{Z}}{\overline{\textbf{Z}}}$. Observe that
\begin{align}\label{eq:not_a_complex_Gaussian_rvc}
\R^{4n} \ni J_2(\vc{\textbf{Z}}{\overline{\textbf{Z}}}) = 
A J_2(\textbf{Z}) \sim 
N_{4n}(A J_2(\mu), A S A^\top),
\end{align}
where $A: = \begin{pmatrix}
I_n & 0\\
I_n & 0\\
0 & I_n\\
0 & -I_n\\
\end{pmatrix}
\in \M_{4n, 2n}(\R)$. Thus, since
$a^\top S a = \frac{1}{2}(\Lambda_{2n}\,a)^\ast
\begin{pmatrix}
\Gamma & C\\
\overline{C} & \overline{\Gamma}
\end{pmatrix}
\Lambda_{2n}\,a =
\frac{1}{2}\,\E[(\Lambda_{2n}\,a)^\ast \textbf{W}\textbf{W}^\ast\,
\Lambda_{2n}\,a] = \E[(\textbf{W}^\ast\,\Lambda_{2n}\,a)^\ast 
\textbf{W}^\ast\,\Lambda_{2n}\,a]$ for all $a \in \R^{2n}$, it follows that 
the real matrix $S$ is always positive semidefinite (cf. 
\sref{Corollary}{cor:char_of_complex_PSD}); i.e., $S \in \M_{2n}(\R)^+$. 
\textit{If in addition} $C \Gamma = \Gamma C$, then (see \cite[Theorem 3]{S2000})
\begin{equation*}
\det(S) = \frac{1}{4^n}\,\det\Big(R_2(\Gamma) + 
\begin{pmatrix}
\Re(C) & \Im(C)\\
\Im(C) & -\Re(C)
\end{pmatrix}\Big) = 
\frac{1}{4^n}\,\det(\Gamma \overline{\Gamma} - C \overline{C})\,.
\end{equation*}
In particular, if $C = 0$, then $\frac{1}{2}R_2(\Gamma) = S$ is positive 
semidefinite, implying that 
\begin{align}\label{eq:Re_and_Im_matrix_parts}
\Re(\Gamma) = \Re(\Gamma)^\top \text{ and } -\Im(\Gamma) = \Im(\Gamma)^\top
\end{align}
(cf. \sref{Proposition}{prp:Hermitian_vs_symmetric}-\ref{(iii)}) and 
\begin{equation}\label{eq:det}
\det(R_2(\Gamma)) = \vert\det(\Gamma)\vert^2\,.
\end{equation}
Hence, $\det(\sqrt{R_2(\Gamma)}) = \vert\det(\Gamma)\vert$.
In particular, $\Re(\textbf{Z}) \sim N_n(\Re(\mu), \tfrac{1}{2}\Re(\Gamma))$, 
$\Im(\textbf{Z}) \sim N_n(\Im(\mu), \tfrac{1}{2}\Re(\Gamma))$ and 
$\E[\Re(Z_i)\Im(Z_i)] = \Im(\Gamma_{ii}) = 0$ for all $i \in [n]$. However, 
because of \sref{Proposition}{prop:Re_Im_independence_condition}, the
random \textit{vectors} $\Re(\textbf{Z})$ and $\Im(\textbf{Z})$ in general are 
not independent!   

 Since the distribution of $J_2(\textbf{Z})$ is fully specified by $\mu$, $\Gamma$
and $C$ (due to \eqref{eq:the_matrix_S}), we write 
$\textbf{Z} \sim {\C}N_n(\mu, \Gamma, C)$ if $\textbf{Z}$ is an $n$-dimensional 
complex Gaussian random vector. Thus, if $\textbf{X}$ and $\textbf{Y}$ are 
$n$-dimensional real random vectors, then
\[
\tfrac{1}{\sqrt{2}}{\textbf{X}} + i\,\tfrac{1}{\sqrt{2}}{\textbf{Y}} \sim 
{\C}N_n(\mu, \Gamma, C) \text{ if and only if } \vc{\textbf{X}}{\textbf{Y}}
\sim N_{2n}\Big(J_2(\mu), R_2(\Gamma) + 
\begin{pmatrix}
\Re(C) & \Im(C)\\
\Im(C) & -\Re(C)
\end{pmatrix}\Big),
\]
or equivalently that 
\[
\textbf{Z} \sim {\C}N_n(\mu, \Gamma, C) \text{ if and only if }
J_2(\sqrt{2}\,\textbf{Z})\sim N_{2n}\Big(J_2(\mu), R_2(\Gamma) + 
\begin{pmatrix}
\Re(C) & \Im(C)\\
\Im(C) & -\Re(C)
\end{pmatrix}\Big).
\]
Regarding the main topic of our monograph, we only need to work with $C=0$, what 
will happen from now on. In this case, we just write $\textbf{Z} \sim 
{\C}N_n(\mu, \Gamma)$. Thus, $\textbf{Z} = \tfrac{1}{\sqrt{2}}{\textbf{X}} + 
i\,\tfrac{1}{\sqrt{2}}{\textbf{Y}}\sim {\C}N_n(\mu, \Gamma)$, if and only
if 
\[
\vc{\textbf{X}}{\textbf{Y}} = \sqrt{2}J_2(\textbf{Z}) = 
\sqrt{2}{\text{vec}}(\Re(Z_1), \ldots, 
\Re(Z_n),\Im(Z_1), \ldots, \Im(Z_n)) \sim N_{2n}(J_2(\mu), 
R_2(\Gamma)),
\]
implying that \eqref{eq:linear_transf_of_Gaussian_law} carries over to the 
complex case:   
\begin{align}\label{eq:linear_transf_of_complex_Gaussian_law}
\text{if }\textbf{Z} \sim {\C}N_n(\mu, \Gamma), \text{ then } A\textbf{Z} + b 
\sim {\C}N_n(A\mu+b, A\Gamma A^\ast) \text{ for all } m \in \N \text{ and } 
(b, A) \in \C^m \times \M_{m,n}(\C). 
\end{align}
\begin{remark}
Let $n \in \N$ and $0 \not= \Sigma \in \M_n(\R)^+$ be given. Fix some $\textbf{X} \sim 
N_n(0, \Sigma)$. A natural question would be, to ask whether $\textbf{Z} : = 
\frac{1}{\sqrt{2}}\textbf{X} + i\,\textbf{0} \sim {\C}N_n(0, \Sigma)$ in 
particular is a complex Gaussian random vector? However, if this were the case, 
it would follow that
\[
\vc{\textbf{X}}{\textbf{0}} = J_2(\sqrt{2}\textbf{Z}) \sim N_{2n}(0, R_2(\Sigma))
= N_{2n}\Big(0, \begin{pmatrix}
    \Sigma & 0 \\
    0 & \Sigma
  \end{pmatrix}\Big),
\]   
implying that $0 \not= \Sigma = \E[\textbf{X} \textbf{X}^\top] = 
\E[\textbf{0} \textbf{0}^\top] = 0$, which is absurd. 
\end{remark}
\begin{remark}
In general, the random vector 
$\vc{\textbf{Z}}{\overline{\textbf{Z}}}$ is not a complex Gaussian one, even if 
$\textbf{Z}$ is. In order 
to recognise this, let e.g. $Z \sim {\C}N_1(0, 1)$ be given. Then 
$J_2(\sqrt{2} Z) \sim N_2(0, I_2)$. Thus, \eqref{eq:not_a_complex_Gaussian_rvc} 
implies that
\[
J_2(\sqrt{2}\,\vc{Z}{\overline{Z}}) \sim N_4(0, A A^\top).
\]
However, since the matrix
\[
A A^\top = \begin{pmatrix}[cc:cc]
1 & 1 & 0 & 0\\ 
1 & 1 & 0 & 0\\ \hdashline
0 & 0 & 1 & -1\\
0 & 0 & -1 & 1
\end{pmatrix}
\] 
does not coincide with a block matrix of type $\begin{pmatrix}[c:c]
A & -B\\ \hdashline
B & A
\end{pmatrix} \in M_2(M_2(\R))$, $\vc{Z}{\overline{Z}}$ is not a complex 
Gaussian random vector. That observation also holds in the multi-dimensional 
case (see \sref{Lemma}{lem:gen_of_CF_and_MGF}-\ref{(iii)}).
\end{remark}
\noindent Occasionally, in view of embedding both, the complex and the real case
into a single statement, we also unambiguously say that $\textbf{Z} \sim 
{\F}N_n(\mu, \Gamma)$ if the random vector $\textbf{Z}$ maps into $\F^n$, where 
$\F \in \{\R, \C\}$. 
Note that for both fields, we explicitly include the case of variance matrices 
$\Gamma \in \M_n(\F)^+$ which are not invertible, so that a probability 
\textit{density} function of $\textbf{Z} \sim {\F}N_n(\mu, \Gamma)$ would not 
have to exist; as opposed to the characteristic function of \textbf{Z} which 
completely determines the probability law $\P_{\textbf{Z}}$. The characteristic 
function of $\textbf{Z} \sim {\C}N_n(\mu, \Gamma)$ can be reduced to the well-known 
characteristic function of the $2n$-dimensional real Gaussian random vector 
$J_2(\textbf{Z}) \sim N_{2n}(J_2(\mu), \frac{1}{2}R_2(\Gamma))$. This 
follows from 
\begin{align*}
\C^n \ni c \mapsto \phi_{\textbf{Z}}(c) :&= 
\E[\exp(i \Re(c^\ast{\textbf{Z}}))] =
\E[\exp(i J_2^\top(c)J_2({\textbf{Z}}))]\\
&= \exp(i J_2(c)^\top J_2(\mu))
\exp(- \tfrac{1}{4}J_2(c)^\top R_2(\Gamma)J_2(c))\\
& = \exp(i\Re(c^\ast \mu))
\exp(-\tfrac{1}{4}\Re(c^\ast\Gamma c))\\ 
& = \exp(i\Re(c^\ast \mu))\exp(-\tfrac{1}{4}c^\ast\Gamma c)
\end{align*}
(cf. \cite[Theorem 2.7]{AHSE1995}, \cite[Theorem 30.2]{B1996}, 
\cite[Definition E.1.13 and Theorem E.1.16]{HvNVW2017} and 
\sref{Lemma}{lem:gen_of_CF_and_MGF}-\ref{(ii)} below). 
\begin{proposition}\label{prop:Re_Im_independence_condition}
Let $n \in \N$, $\mu \in \C^n$, $\Gamma \in \M_n(\C)^+$ and 
$\textbf{Z} \sim {\C}N_n(\mu, \Gamma)$. Then $\Re({\textbf{Z}})$ and 
$\Im({\textbf{Z}})$ are independent if and only if $\Im(\Gamma) = 0$.
\end{proposition}
\begin{comment}
\eqref{eq:Re_and_Im_matrix_parts} and \sref{Lemma}{lem:skew_sym} imply that
\begin{align*}
\phi_{\vc{\Re({\textbf{Z}})}{\Im({\textbf{Z}})}}(\vc{\Re(c)}{\Im(c)})
&= \phi_{\textbf{Z}}(c)\\
&= \phi_{\Re({\textbf{Z}})}(\Re(c))\,\phi_{\Im({\textbf{Z}})}(\Re(c))
\exp(-\tfrac{1}{2}\Im(c)^\top\Im(\Gamma)\Re(c))
\end{align*}
for all $c \in \C^n$. Hence, $\Re({\textbf{Z}})$ and $\Im({\textbf{Z}})$ are 
independent if and only if 
\[
\phi_{\Re({\textbf{Z}})}(\Re(c))\,\phi_{\Im({\textbf{Z}})}(\Re(c))
(1 - \exp(-\tfrac{1}{2}\Im(c)^\top\Im(\Gamma)\Re(c)) = 0
\]
for all $c \in \C^n$. Consequently, if we apply the absolute value of the 
latter equality to any vector $e_l + i e_k \in \C^n$, where $k,l \in [n]$, 
it follows that
\[
\vert 1 - \exp(-\tfrac{1}{2}\Im(\Gamma_{kl}))\vert = 0
\]
for all $k,l \in [n]$, and the claim follows.
\end{comment}
\noindent Similarly, if $\Gamma$ (respectively $R_2(\Gamma)$) is invertible, the 
complex density function under the Lebesgue measure $\lambda_n^{\C}$ on $\C^n$ 
can be constructed as 
\begin{align*}
\C^n \ni a \mapsto \varphi_{\mu, \Gamma}(a) &: = 
\varphi_{J_2(\mu), \frac{1}{2}\,R_2(\Gamma)}(J_2(a))\\ 
&= \frac{1}{\pi^n\sqrt{\det(R_2(\Gamma))}}
\exp(-(J_2(a-\mu))^\ast R_2(\Gamma^{-1})(J_2(a-\mu)))\\
&\stackrel{\eqref{eq:det}}{=} \frac{1}{\pi^n\,\vert\det(\Gamma)\vert}
\exp(-(a-\mu)^\ast\,\Gamma^{-1}\,(a-\mu))\,.
\end{align*}
These facts, including \eqref{eq:linear_transf_of_Gaussian_law},
\eqref{eq:Re_and_Im_matrix_parts}, 
\eqref{eq:linear_transf_of_complex_Gaussian_law} and 
\sref{Proposition}{prop:n_dim_real_Gaussian_law}, immediately imply the 
following comprehensive characterisation of the probability law 
$\C N_n(\mu, \Gamma)$.
\begin{proposition}\label{prop:complex_Gaussian_law}
Let $n \in \N$, $\mu \in \C^n$, $\Gamma \in \M_{n}(\C)^+$ and $\textbf{Z}$ a complex 
$n$-dimensional random vector. 
Then the following statements are equivalent:
\begin{enumerate}
\scalebox{0.88}{
\vbox{
\item For all $c \in \C^n$, $\phi_{\textbf{Z}}(c) = 
\E[\exp(i\,\Re(c^\ast \textbf{Z})] = \exp(i\Re(c^\ast \mu))
\exp(-\tfrac{1}{4}c^\ast\Gamma\,c)$.
\item $\sqrt{2}\,J_2(\textbf{Z}) =
{\text{vec}}(\sqrt{2}\Re(Z_1), \ldots, \sqrt{2}\Re(Z_n),\sqrt{2}\Im(Z_1), \ldots, 
\sqrt{2}\Im(Z_n)) \sim N_{2n}(\sqrt{2}\,J_2(\mu), R_2(\Gamma))$.
\item $\textbf{Z} \sim \C N_n(\mu, \Gamma)$.
\item For all $\alpha \in \T$, $\alpha\textbf{Z} \sim 
\C N_n(\alpha \mu, \Gamma)$.
\item For all $c \in \C^n$, $c^\ast\textbf{Z} \sim 
{\C}N_1(c^\ast \mu, c^\ast\Gamma\,c)$.
\item For all $c \in \C^n$, $\sqrt{2}\,J_2(c^\ast \textbf{Z}) 
\sim N_2(\sqrt{2}\,J_2(c^\ast \mu), R_2(c^\ast\Gamma\,c))$. 
\item For all $c \in \C^n$, $\sqrt{2}\Re(c^\ast \textbf{Z}) \sim 
N_1(\sqrt{2}\Re(c^\ast \mu), c^\ast\Gamma\,c)$.  
}}
\end{enumerate}
In particular, if $\textbf{Z} \sim \C N_n(\mu, \Gamma)$, then 
$\Re(\textbf{Z}) \stackrel{d}{=} \Im(\textbf{Z}) \sim N_n(0, \Re(\Gamma))$, and
$\Re(Z_i)$ and $\Im(Z_i)$ are independent for all $i \in \N$. $\textbf{Z} \sim 
\C N_n(\mu, \Gamma)$ if and only if $\overline{\textbf{Z}} \sim 
\C N_n(\overline{\mu}, \overline{\Gamma})$. Moreover, if $\mu=0$ and $\Gamma = 
I_n$, then $\{\Re(Z_1), \ldots, \Re(Z_n),\Im(Z_1), \ldots, \Im(Z_n)\}$ are 
pairwise independent. 
\end{proposition}
\noindent In relation to a further investigation of the complex Grothendieck 
constant (built on complex Hermite polynomials), we need a further analysis of 
the structure of the random vector $\vc{\textbf{Z}}{\overline{\textbf{Z}}}$ 
(cf. \autoref{thm:correlated_complex_Hermite_polynomials}). That analysis, 
encapsulated in the next lemma, includes a short proof of a generalisation of a 
result of L. J. Halliwell (see \cite[Appendix B]{H2015}), built on a 
``change of mean trick'', that allows us to pull both, the characteristic function 
and the moment generating function of a real Gaussian random vector out of a single 
formula, without having to assume the existence of a density function. In doing so, 
we will recognise again that in general, the complex random vector 
$\vc{\textbf{Z}}{\overline{\textbf{Z}}}$ in $\C^{2n}$ is not Gaussian, even
if the random vector $\textbf{Z}$ in $\C^n$ were a complex Gaussian one. 
\begin{lemma}\label{lem:gen_of_CF_and_MGF}
Let $n \in \N, \Sigma \in \M_n(\R)^+$ and $\Gamma \in \M_n(\C)^+$. Let $\textbf{X} \sim 
N_n(0, \Sigma)$ and $\textbf{Z} \sim {\C}N_n(0, \Gamma)$. Then
\begin{enumerate}
\item $\E[\exp(c^\top\textbf{X})] = 
\exp(\frac{1}{2}\,c^\top \Sigma c )$ for all $c \in \C^n$.
\item $\E[\exp(a^\ast\textbf{Z} + b^\ast\overline{\textbf{Z}})] = 
\exp(a^\ast\Gamma \overline{b})$ for all $a, b \in \C^n$. 
\item For all $a,b \in \C^n$,
\begin{align*}
\phi_{\vc{\textbf{Z}}{\overline{\textbf{Z}}}}(\vc{a}{b}) &=
\exp(-\tfrac{1}{4} a^\ast \Gamma a)
\exp(-\tfrac{1}{4} b^\ast \overline{\Gamma} b) 
\exp(-\tfrac{1}{2}\Re(a^\ast \Gamma \overline{b}))\\
&= \exp\big(-\tfrac{1}{4}{\vc{a}{b}}^\top
\begin{pmatrix}
\Gamma & 0\\
0 & \overline{\Gamma}
\end{pmatrix}
\vc{a}{b}\big)\exp(-\tfrac{1}{2}\Re(a^\ast \Gamma \overline{b})).
\end{align*}
\end{enumerate}
In particular, $\vc{\textbf{Z}}{\overline{\textbf{Z}}} \sim 
{\C}N_{2n}\big(0, \begin{pmatrix}
\Gamma & 0\\
0 & \overline{\Gamma}
\end{pmatrix}\big)$ if and only if $\Gamma = 0$. Moreover,
$\E[\exp(a^\ast\textbf{Z})] = 1$ for all $a \in \C^n$. 
\end{lemma}
\begin{remark}
The main difficulty in the proof of \sref{Lemma}{lem:gen_of_CF_and_MGF}-\ref{(i)} arises from the 
fact that the normal random variables $\alpha^\top\textbf{X}$ and $\beta^\top\textbf{X}$ are 
correlated, so that we cannot simply represent $\E_{\P_{\textbf{X}}}
[\exp(\alpha^\top\textbf{X})\exp(i \beta^\top\textbf{X})]$ as a product of two 
expectations. However, it is possible to construct a completely different proof of 
\sref{Lemma}{lem:gen_of_CF_and_MGF}-\ref{(i)}, which is built on (an application of 
the one-dimensional case of) \autoref{thm:h_fg_real_case}. We strongly encourage 
the readers to work out the details.
\end{remark} 
\noindent Let $p \in \N$. It is well-known that the image measure 
$\P_\textbf{X}$ of a real Gaussian random vector $\textbf{X} \sim N_p(0, I_p)$ 
actually coincides with the Gaussian measure $\gamma_{p}$ on $\R^p$ (cf. e.g. 
\cite[Proposition 1.2.2.]{B1998}), constructed via
\[
{\mathcal{B}}(\R^p) \ni B \mapsto \gamma_{p}(B) : = 
(2\pi)^{-p/2}\,\int_B \exp(-\frac{1}{2}\Vert x \Vert^2_2 )\lambda_p(\textup{d}x)
= \P(\textbf{X} \in B).
\]
Due to \sref{Proposition}{prop:complex_Gaussian_law} 
this fact can be easily transferred to the 
complex case. To this end, let $n \in \N$, $\textbf{Z} \sim \C N_n(0, I_n)$ and $b = \Re(b) + 
i\,\Im(b) : \C^n \longrightarrow \C$. Consider the two mappings
\[
r(b) : = \Re(b)\circ \frac{1}{\sqrt{2}}\,J_2^{-1} : \R^{2n} \longrightarrow 
\R  \text{ and } s(b) : = \Im(b)\circ \frac{1}{\sqrt{2}}\,J_2^{-1} : \R^{2n} 
\longrightarrow \R\,.
\]
By construction, it follows that for any $x,y \in \R^n$,
\begin{align}\label{eq:real_valued_fcts_r_b_and_s_b}
r(b)({\text{vec}}(x,y)) = \Re(b(\tfrac{1}{\sqrt{2}}\,x + 
i\,\tfrac{1}{\sqrt{2}}\,y)) \text{ and } s(b)({\text{vec}}(x,y)) = 
\Im(b(\tfrac{1}{\sqrt{2}}\,x + i\,\tfrac{1}{\sqrt{2}}\,y)).
\end{align}
Obviously, $s(b) = r(-i\,b)$, $b \circ \frac{1}{\sqrt{2}}\,J_2^{-1} = 
r(b) + i s(b)$, $\Re(b) = r(b) \circ \sqrt{2}J_2$, $\Im(b) = s(b) \circ 
\sqrt{2}J_2$ and $r(\alpha b) = \alpha r(b)$ for any $\alpha \in \R$.

 Let $p,q \in [1, \infty)$, such that $\frac{1}{p} + \frac{1}{q} = 1$. 
A direct application of H\"{o}lder's inequality (to the vectors 
$(r(b), s(b))^\top \in \R^2$ and $(1,1)^\top \in \R^2$) implies that
\begin{align*}
\max\{\vert r(b) \vert^p, \vert s(b) \vert^p\} \leq r(\vert b \vert)^p = 
(\vert b\vert^p \circ \frac{1}{\sqrt{2}}\,J_2^{-1}) \leq 
(\vert r(b) \vert\cdot 1 + \vert s(b) \vert\cdot 1)^p \leq 2^{p/q}\,
(\vert r(b) \vert^p + \vert s(b) \vert^p).
\end{align*}
Hence, $\max\{\vert r(b) \vert^p, \vert s(b) \vert^p\} \in 
L^p(\R^{2n},\gamma_{2n})$, if and only if $b \in L^p(\C^n, \P_{\textbf{Z}})$, 
and
\begin{align}\label{eq:complex_Gaussian_integral}
\begin{split}
\E[b(\textbf{Z})] &= \int_{\C^n}b\,\textup{d}\P_{\textbf{Z}} 
= \int_{\C^n} b\,\textup{d}\,(\tfrac{1}{\sqrt{2}} J_2^{-1})_\ast\P_
{\sqrt{2} J_2(\textbf{Z})} = \int_{\R^{2n}}b(\tfrac{1}{\sqrt{2}}(x + i y))
\gamma_{2n}(\textup{d}(x,y))\\
&= \int_{\R^{2n}}r(b)\,\textup{d}\gamma_{2n} + i\,\int_{\R^{2n}}s(b)\,
\textup{d}\gamma_{2n} = \E[r(b)(\textbf{X})] + i\,\E[s(b)(\textbf{X})],
\end{split}
\end{align}
where $\textbf{X} \stackrel{d}{=} \sqrt{2}\,J_2(\textbf{Z}) \sim 
N_{2n}(0, I_{2n})$. In particular,

\scalebox{0.88}{
\vbox{
\begin{align*}
{\mathcal{B}}(\C^n) \ni A \mapsto \gamma_n^{\C}(A) : &= \P(\textbf{Z} \in A) = 
\E[\ind_{A}(\textbf{Z})] \stackrel{\eqref{eq:complex_Gaussian_integral}}{=} 
\gamma_{2n}(J_2(\sqrt{2}A))\\ 
&=
\pi^{-n}\int_{J_2(A)}\exp(-\Vert x \Vert^2_2)\,\lambda_{2n}(\textup{d}x) 
\stackrel{\eqref{eq:complex_Lebesgue_measure}}{=} 
\pi^{-n}\int_{A}\exp(-\Vert J_2(z) \Vert^2_2)\,\lambda^\C_{n}(\textup{d}z)
\end{align*}
}}

\noindent emerges as the Gaussian measure on $\C^n$, implying that
\[
\gamma_n^{\C} = (\tfrac{1}{\sqrt{2}}J_2^{-1})_\ast \gamma_{2n} = 
\P_\textbf{Z}
\] 
is absolutely continuous with respect to $\lambda^\C_{n}$, with 
Radon-Nikod\'{y}m derivative\\ $\frac{\textup{d}\gamma_n^{\C}}{\textup{d}
\lambda^\C_{n}} = \pi^{-n}\exp(-\Vert J_2(z) \Vert^2_2)$. 
\begin{remark}
In \cite[Section 8.7]{DF1993}, the complex Gaussian measure on $\C^n$ is defined 
in such a manner that it coincides with the real Gaussian measure $\gamma_{2n}$ 
on $\R^{2n}$. Given that construction, the important factor $\sqrt{2}$ - which 
actually emerges from the underlying structure of the probability law of a 
complex Gaussian random vector - is ignored. In our view, that approach creates 
a bit of dissonance. For example, \sref{Corollary}{cor:little_GT_and_Gaussian_structure}, 
which shows us that for both fields, $\F = \R$ and $\F = \C$, the little 
Grothendieck constant $k_G^\F$ actually emerges from a common source, can no 
longer be maintained.
\end{remark}
\noindent Hence, $b \in L^p(\C^n,\gamma_n^{\C})$ if and only if 
$\max\{\vert r(b) \vert^p, \vert s(b) \vert^p\} \in L^p(\R^{2n},\gamma_{2n}^{\R})$, 
and 
\begin{align}\label{eq:complex_Gaussian_integral_II}
\begin{split}
\int_{\C^n}\Re(b)\,\textup{d}\gamma_n^{\C} + i\int_{\C^n}\Im(b)\,\textup{d}\gamma_n^{\C} =
\int_{\C^n}b\,\textup{d}\gamma_n^{\C} = 
\E[b(\textbf{Z})] \stackrel{\eqref{eq:complex_Gaussian_integral}}{=} 
\int_{\R^{2n}}r(b)\,\textup{d}\gamma_{2n} + i\,\int_{\R^{2n}}s(b)\,
\textup{d}\gamma_{2n}\,.
\end{split}
\end{align}

\noindent In particular, $b \in L^2(\C^n,\gamma_n^{\C})$ if and only if 
$r(b) \in L^2(\R^{2n},\gamma_{2n})$ and $s(b) \in L^2(\R^{2n},\gamma_{2n})$, 
so that (in either case) 
\begin{equation}\label{eq:complex_Gaussian_L2}
\E[b(\textbf{Z})\overline{b}(\textbf{Z})] 
= \Vert b \Vert_{\gamma_n^{\C}}^2 = \int_{\C^n} \vert b \vert^2 
\textup{d}\gamma_n^{\C} \stackrel{\eqref{eq:complex_Gaussian_integral_II}}{=}  
\int_{\R^{2n}}r(\vert b\vert^2)\,\textup{d}\gamma_{2n} = 
\Vert r(b) \Vert_{\gamma_{2n}}^2 + \Vert s(b) \Vert_{\gamma_{2n}}^2
\end{equation}
(since $r(\vert b\vert^2) = r(b)^2 + s(b)^2$). In particular, for any function
$f : \R^{2n} \longrightarrow \C$, it follows that $f \circ \sqrt{2}J_2 
\in L^2(\C^n, \gamma_n^{\C})$ if and only if $\Re(f) \in L^2(\R^{2n},\gamma_
{2n})$ and $\Im(f) \in L^2(\R^{2n},\gamma_{2n})$, whence 
\begin{equation}\label{eq:real_L_2_as_complex_L2}
\Vert f \circ \sqrt{2}J_2 \Vert_{\gamma_n^{\C}}^2 = 
\Vert \Re(f) \Vert_{\gamma_{2n}}^2 + \Vert \Im(f) \Vert_{\gamma_{2n}}^2\,.
\end{equation}
in either case. Consequently, 
\begin{equation}\label{eq:real_L_2_as_complex_L2_product_case}
\Vert g^\C \Vert_{\gamma_n^{\C}} = \Vert g \Vert_{\gamma_{n}}
\end{equation}
for any real-valued function $g : \R^{n} \longrightarrow \R$, where 
$g^\C : = (g \otimes 1)\circ \sqrt{2}\,J_2$ (since $\gamma_{2n} = \gamma_n 
\otimes \gamma_n$). Moreover, if $b, c \in L^2(\C^n,\gamma_n^{\C})$, then 
$r(b\overline{c}) = r(b) r(c) + s(b) s(c)$ and $s(b\overline{c}) = s(b) r(c) - 
r(b) s(c)$. \eqref{eq:complex_Gaussian_integral_II} therefore implies that 
\begin{align}\label{inner_product_of_complex_Gaussian_L2}
\begin{split}
\langle b, c \rangle_{\gamma_n^{\C}} &= 
\langle r(b), r(c) \rangle_{\gamma_{2n}} + 
\langle s(b), s(c) \rangle_{\gamma_{2n}}\\ 
&+ \,i\,\langle s(b), r(c) \rangle_{\gamma_{2n}} -
i\,\langle r(b), s(c) \rangle_{\gamma_{2n}}.
\end{split}
\end{align} 
\begin{remark}
If in addition the function $b : \C^n \longrightarrow \C$ is holomorphic, then $b$ actually 
is an element of the Segal-Bergmann space (cf. \cite[Chapter 3.10]{H2017})
\[
\mathcal{H}L^2(\C^n) : = \{c : \C^d \longrightarrow \C : c \text{ is holomorphic and }
\Vert c \Vert_{\gamma_n^{\C}} < \infty\}.
\]
\end{remark}
 
\noindent In a similar vein, 
one can now prove easily the more comprehensive 
\begin{corollary}\label{cor:complex_Gaussian_integral}
Let $n \in \N$, $p \in [1, \infty)$, $\mu \in \C^n$, $\Gamma \in M_n(\C)^+$ 
and $\textbf{Z} \sim \C N_n(\mu, \Gamma)$. Let $b = \Re(b) + i\,\Im(b) : 
\C^n \longrightarrow \C$, such that 
\[
\R^{2n} \ni y \mapsto r(b)(y + \sqrt{2}\,J_2(\mu)) \in 
L^p(\R^{2n}, \gamma_{2n})
\]
and
\[
\R^{2n} \ni y \mapsto s(b)(y + \sqrt{2}\,J_2(\mu)) \in 
L^p(\R^{2n}, \gamma_{2n}).
\]
Then $b \in L^p(\C^n, \P_{\textbf{Z}})$, and
\begin{align*}
\E[b(\textbf{Z})] = \vert \det(\Gamma)\vert(
\int_{\R^{2n}}r(b)(y + \sqrt{2}\,J_2(\mu))\gamma_{2n}(\textup{d}y) + 
i\,\int_{\R^{2n}}s(b)(y + \sqrt{2}\,J_2(\mu))\gamma_{2n}(\textup{d}y)).
\end{align*}
\end{corollary}
\begin{comment}
We just have to observe that
\begin{align*}
\E[b(\textbf{Z})] &= \E[r(b)(\textbf{Y})] 
+ i\,\E[s(b)(\textbf{Y})]\\
&= \int_{\R^{2n}}r(b)(\sqrt{R_2(\Gamma)}x + \sqrt{2}\,J_2(\mu))\gamma_{2n}
(\textup{d}x) + i\,\int_{\R^{2n}}s(b)(\sqrt{R_2(\Gamma)}x + \sqrt{2}\,J_2(\mu))
\gamma_{2n}(\textup{d}x)\\
&\stackrel{\eqref{eq:det}}{=} \vert \det(\Gamma)\vert(\int_{\R^{2n}}r(b)
(y + \sqrt{2}\,J_2(\mu))\gamma_{2n}(\textup{d}y) + 
i\,\int_{\R^{2n}}s(b)(y + \sqrt{2}\,J_2(\mu))\gamma_{2n}(\textup{d}y)),
\end{align*}
where $\textbf{Y} \stackrel{d}{=} \sqrt{2}\,J_2(\textbf{Z}) \sim 
N_{2n}(\sqrt{2}\,J_2(\mu), R_2(\Gamma))$.
\end{comment}
\section{Partitioned complex Gaussian random vectors in $\C^{2n}$ and 
the probability law ${\C}N_{2n}(0, \Sigma_{2n}(\zeta))$}
 
For the remainder of the paper we put, without loss of generality, $\mu = 0$, 
so that we are working with centred Gaussian random vectors (with respect to both fields, 
$\R$ and $\C$). Moreover, we make use of a specific class of partitioned correlation 
matrices, which turns out to be of crucial importance regarding the topic of the paper. To 
this end, let $\F \in \{\R, \C\}$, $n \in \N$ and $\zeta \in \F \cap 
\overline{\D}$. Put
\begin{align}\label{eq:GT_corr_matrix}
\begin{split}
\Sigma_{2n}(\zeta) : &= 
\begin{pmatrix}
    I_n & \zeta\,I_n \\
    \overline{\zeta}\,I_n & I_n
  \end{pmatrix}
=	\Sigma_2(\zeta) \otimes I_n =		
\begin{pmatrix}
   \begin{matrix}
    1 & 0 & \ldots & 0  \\
    0 & 1 & \ldots & 0 \\
		\vdots & \vdots & \ddots & \vdots \\
		0 & 0 & \ldots & 1
  \end{matrix} 
	& 
	\begin{matrix}
    \zeta & 0 & \ldots & 0  \\
    0 & \zeta & \ldots & 0 \\
		\vdots & \vdots & \ddots & \vdots \\
		0 & 0 & \ldots & \zeta
  \end{matrix}\\	
   \begin{matrix}
    \overline{\zeta} & 0 & \ldots & 0  \\
    0 & \overline{\zeta} & \ldots & 0 \\
		\vdots & \vdots & \ddots & \vdots \\
		0 & 0 & \ldots & \overline{\zeta} 
  \end{matrix}				
				& 
	\begin{matrix}
    1 & 0 & \ldots & 0  \\
    0 & 1 & \ldots & 0 \\
		\vdots & \vdots & \ddots & \vdots \\
		0 & 0 & \ldots & 1
  \end{matrix}
\end{pmatrix}\\ 
&= \Sigma_{2n}(\Re(\zeta)) + R_2(-i\,\Im(\zeta)I_n).
\end{split}
\end{align}
Since $\vert \zeta \vert \leq 1$, it follows that  
\[
\big\langle 
\begin{pmatrix}
    I_n & \zeta\,I_n \\
    \overline{\zeta}\,I_n & I_n
\end{pmatrix}\!\!\colvec{2}{a}{b},\colvec{2}{a}{b} \big\rangle_{\F^{2n}}
= \Vert a + \zeta\,b\Vert_{\F^n}^2 + (1 - \vert \zeta \vert^2)\,
\Vert b \Vert_{\F^n}^2 \geq 0 
\]
for all $a, b \in \F^n$, implying that $\Sigma_{2n}(\zeta)$ in fact is positive 
semidefinite and hence a correlation matrix (due to 
\sref{Lemma}{lem:charact_of_corr_matrices}). 
\sref{Lemma}{lem:charact_of_corr_matrices} also clearly implies that 
\begin{align}\label{eq:the_2_dim_corr_matrix}
C(2; \F) = \big\{\Sigma_{2}(\zeta) = \begin{pmatrix}
1 & \zeta\\
\overline{\zeta} & 1 
\end{pmatrix} 
: \zeta \in \F \cap \overline{\D}\big\}.
\end{align} 
Moreover, the determinant of the Kronecker product $\Sigma_{2n}(\zeta) = 
\Sigma_2(\zeta) \otimes I_n$ is calculated as (cf. \cite[Problem 4.2.1]{HJ1991}):
\[
\det(\Sigma_{2n}(\zeta)) = \det(I_n)^2\,\det(\Sigma_{2}(\zeta))^n =
(1- \vert \zeta \vert^2)^n\,.
\]
If in addition $\vert \zeta \vert < 1$, then $\Vert a + \zeta\,b\Vert_{\F^n}^2 + 
(1 - \vert \zeta \vert^2)\Vert b \Vert_{\F^n}^2 > 0$ for all $(a,b) \in 
(\F^{n} \times \F^{n})\setminus\{(0,0)\}$, implying that in this case 
$\Sigma_{2n}(\zeta)$ even is positive definite and hence invertible, with inverse 
$\Sigma_{2n}(\zeta)^{-1} = \frac{1}{1-\vert \zeta \vert^2}{\Sigma_{2n}(-\zeta)} = 
\frac{1}{1-\vert \zeta \vert^2}\,
\begin{pmatrix}
  I_n & -\zeta\,I_n \\
  -\overline{\zeta}\,I_n & I_n
\end{pmatrix}$, implying that also $(1-\vert \zeta \vert^2)\,\Sigma_{2n}(\zeta)^{-1} = 
\Sigma_{2n}(-\zeta)$ is a correlation matrix of rank $2n$.
 
In particular, for any $\rho \in (-1,1)$, the density function of the 
$2n$-dimensional random vector $\vc{{\textbf{X}}_1}{{\textbf{X}}_2} \sim 
N_{2n}(0, \Sigma_{2n}(\rho))$, where both, ${\textbf{X}}_1$ and 
${\textbf{X}}_2$ are $n$-dimensional random vectors, exists. It is given by 
\begin{align}\label{eq:Mehler_kernel_rep}
\begin{split}
\varphi_{0, \Sigma_{2n}(\rho)}(x_1,x_2) &= 
\frac{1}{(2\pi)^n\,(1 - \rho^2)^{n/2}}\exp\big(-\frac{1}{2(1-\rho^2)}
\langle \Sigma_{2n}(-\rho)x, x\rangle_{\R_2^{2n}}\big)\\
&= \frac{1}{(2\pi)^n\,(1 - \rho^2)^{n/2}}\exp\big(-\frac{\Vert x_1 \Vert^2 + 
\Vert x_2 \Vert^2 - 2\rho\langle x_1,x_2 \rangle}{2(1-\rho^2)}\big)\\
&= M_\rho(x_1,x_2; n)\,\varphi_{0, I_{2n}}(x_1,x_2)\\ 
&= M_\rho(x_1,x_2; n)\,\varphi_{0, I_n}(x_1)\,\varphi_{0, I_n}(x_2)\\
&= \varphi_{0, \Sigma_{2n}(\rho)}(x_2,x_1)\,,
\end{split}
\end{align}
where $x_1, x_2 \in \R^n$, $x : = \vc{x_1}{x_2}$ and 
\begin{align}\label{def:Mehler_kernel}
M_\rho(x_1,x_2; n) : = \frac{1}{(1-\rho^2)^{n/2}}\exp\big(\frac{2\rho\langle x_1, 
x_2\rangle - \rho^2(\Vert x_1 \Vert^2 + \Vert x_2 \Vert^2)}{2(1-\rho^2)}\big) 
= M_{-\rho}(x_1,-x_2; n)
\end{align} 
denotes the $n$-dimensional Mehler kernel (cf. \cite{H1999}).
Moreover, since $\det(\Sigma_{2n}(\zeta)) = 0$ if and only if $\zeta \in \D$, 
we achieve the following result: 
\begin{proposition}\label{prop:structure_of_Sigma_2n_zeta}
Let $\F \in \{\R, \C\}$, $n \in \N$ and $\zeta \in \overline{\D} \cap \F$. Then
\[
\Sigma_{2n}(\zeta) : = 
\begin{pmatrix}
    I_n & \zeta\,I_n \\
    \overline{\zeta}\,I_n & I_n
  \end{pmatrix} \in C(2n; \F).
\]
Moreover, the following statements are equivalent
\begin{enumerate}
\item $\Sigma_{2n}(\zeta)$ is a correlation matrix of rank $2n$.
\item $\Sigma_{2n}(\zeta)$ is invertible.
\item $\zeta \in \D$.
\end{enumerate}
If one of these equivalent statements is given, then $\Sigma_{2n}(\zeta)^{-1} = 
\frac{1}{1-\vert\zeta\vert^2}\,\Sigma_{2n}(-\zeta)$. In particular, also the 
$2n \times 2n$-matrix $(1-\vert\zeta\vert^2)\Sigma_{2n}(\zeta)^{-1}$ is a correlation matrix 
of rank $2n$.
\end{proposition}
\noindent It is not obvious that the probablity law 
$N_{2n}(0, \Sigma_{2n}(\rho))$ can also be described as follows 
(cf. \cite[Definition 11.6 and Definition 11.10]{O'D2014}):
\begin{proposition}\label{prop:char_of_the_Sigma_2n_prob_law}
Let $n \in \N$ and $\rho \in [-1,1]$. Let $\textbf{X} = (X_1, \ldots, X_n)^\top$ and
$\textbf{Y} = (Y_1, \ldots, Y_n)^\top$ be two $\R^n$-valued random vectors. Then 
the following statements are equivalent:
\begin{enumerate}
\item $\vc{\textbf{X}}{\textbf{Y}} \sim N_{2n}(0, \Sigma_{2n}(\rho))$.
\item $\vcfour{X_i}{Y_i}{X_j}{Y_j} 
\sim N_{2n}\big(0, \begin{pmatrix}
\Sigma_{2}(\rho) & 0\\
0 & \Sigma_{2}(\rho)
\end{pmatrix}\big)$ for all $i \not=j \in [n]$.
\item The random component vector pairs $\vc{X_i}{Y_i}, \ldots, 
\vc{X_n}{Y_n}$ are mutually independent, and 
$\vc{X_i}{Y_i} \sim N_{2}(0, \Sigma_{2}(\rho))$ for all $i \in [n]$.
\end{enumerate} 
\end{proposition}
\begin{comment}
We just sketch the main ideas and leave the elaboration of the proof 
to the readers. Suppose that (i) holds. Since
\[
\vcfour{X_i}{Y_i}{X_j}{Y_j} = A(i,j)\vc{\textbf{X}}{\textbf{Y}}\,,
\]
where $\M_{4, 2n}(\R) \ni A(i,j) : = \begin{pmatrix}
e_i^\top & 0^\top\\
0^\top & e_i^\top\\
e_j^\top & 0^\top\\
0^\top & e_j^\top
\end{pmatrix}$, it follows that
\[
A(i,j)\Sigma_{2n}(\rho)A(i,j)^\top = \begin{pmatrix}
\Sigma_{2}(\rho) & 0\\
0 & \Sigma_{2}(\rho)
\end{pmatrix}
\]
is the correlation matrix of the Gaussian random vector $\vcfour{X_i}{Y_i}{X_j}{Y_j}$
(cf. \eqref{eq:linear_transf_of_Gaussian_law}), which proves (ii). If (ii) is given, the structure of the correlation matrix
$\begin{pmatrix}
\Sigma_{2}(\rho) & 0\\
0 & \Sigma_{2}(\rho)
\end{pmatrix}$ implies that for $i \not=j$ the 
characteristic function of $\vcfour{X_i}{Y_i}{X_j}{Y_j}$ can be written as the 
product of the two characteristic functions of the pairs $\vc{X_i}{Y_i}$ and 
$\vc{X_j}{Y_j}$, respectively. Finally, if (iii) is given, it follows that 
for all $a, b \in \R^n$ 
\[
(\vc{a}{b})^\top\vc{\textbf{X}}{\textbf{Y}} = \sum_{i=1}^n V_i\,,
\]
where $V_i : = a_i\,X_i + b_i\,Y_i \sim N_1(0, a_i^2 + 2\rho a_i b_i + b_i^2)$ 
for all $i \in [n]$ and $V_1, \ldots, V_n$ are mutually independent. 
Consequently, it follows that $(\vc{a}{b})^\top\vc{\textbf{X}}{\textbf{Y}} \sim 
N_1(0, \Vert a \Vert^2 + 2\rho a^\top b + \Vert b \Vert^2)$, and (i) is achieved.
\end{comment}
\noindent Regarding the underlying structure of the Haagerup equality (and its 
generalisation - see \autoref{thm:h_b_c_complex_correl_case}) an analysis of 
the structure of partitioned complex $2n$-dimensional Gaussian random vectors 
whose probability law is induced by the correlation matrix $\Sigma_{2n}(\zeta)$ 
leads to another important
\begin{lemma}\label{lem:correl_Gaussians_I}
Let $n \in \N$, $\zeta = x + iy \in \overline{\D}$ and 
$\vc{\textbf{Z}}{\textbf{W}} \sim {\C}N_{2n}(0, \Sigma_{2n}(\zeta))$, where the 
complex random vectors $\textbf{Z}$ and \textbf{W} both map into $\C^n$. Then the following 
statements hold:
\begin{enumerate}
\item For any $\alpha, \beta \in \T$, $\vc{\alpha\textbf{Z}}{\beta\textbf{W}} \sim 
{\C}N_{2n}(0, \Sigma_{2n}(\alpha\,\overline{\beta}\,\zeta))$.
\item $\vc{\textbf{W}}{\textbf{Z}} \sim 
{\C}N_{2n}(0, \Sigma_{2n}(\overline{\zeta}))$.
\item $\sqrt{2}\,\vc{\Re(\textbf{Z})}{\Re(\textbf{W})} \sim 
N_{2n}(0, \Sigma_{2n}(\Re(\zeta)))$ and $\sqrt{2}\,\vc{\Im(\textbf{Z})}
{\Im(\textbf{W})} \sim N_{2n}(0, \Sigma_{2n}(\Re(\zeta)))$.
\item If $\zeta \in [-1,1]$, then 
$\sqrt{2}\,{\text{vec}}(\Re(\textbf{Z}), \Im(\textbf{Z}),  
\Re(\textbf{W}), \Im(\textbf{W})) =
\sqrt{2}\,\vc{J_2(\textbf{Z})}{J_2(\textbf{W})} \sim N_{4n}(0, \Sigma_{4n}(\zeta))$, 
and $\vc{\Re(\textbf{Z})}{\Re(\textbf{W})}$ and 
$\vc{\Im(\textbf{Z})}{\Im(\textbf{W})}$ are independent.
\item If $\zeta \in \overline{\D}\setminus\{0\}$, then 
$\vc{{\text{sign}}(\overline{\zeta})\textbf{Z}}{\textbf{W}} \stackrel{d}{=}
\vc{{\text{sign}}(\zeta)\textbf{Z}}{\textbf{W}} \sim 
{\C}N_{2n}(0, {}{\Sigma_{2n}(\vert \zeta \vert)})$. Moreover,
$\sqrt{2}\,\vc{J_2({\text{sign}}(\overline{\zeta})\textbf{Z})}{J_2(\textbf{W})} 
\stackrel{d}{=}
\sqrt{2}\,\vc{J_2({\text{sign}}(\zeta)\textbf{Z})}{J_2(\textbf{W})} \sim 
N_{4n}(0, \Sigma_{4n}(\vert \zeta \vert))$.
\end{enumerate}
\end{lemma}
\noindent Next, we will recognise that for any $\rho \in [-1,1]$, the real probability 
law $N_{4n}(0, \Sigma_{4n}(\rho))$ actually originates from the complex 
probability law ${\C}N_{2n}(0, \Sigma_{2n}(\rho))$! Since:
\begin{corollary}\label{cor:inheritage_prop_of_Sigma_2n_rho}
Let $\rho \in [-1,1]$. Let ${\textbf{X}}_1$, ${\textbf{Y}}_1$, ${\textbf{X}}_2$ and
${\textbf{Y}}_2$ be four $\R^n$-valued random vectors. 
Then the following statements are equivalent:
\begin{enumerate}
\item $\vc{{\textbf{X}}_1}{{\textbf{X}}_2} \sim 
N_{2n}(0, \Sigma_{2n}(\rho))$, $\vc{{\textbf{Y}}_1}{{\textbf{Y}}_2} 
\sim N_{2n}(0, \Sigma_{2n}(\rho))$, and 
$\vc{{\textbf{X}}_1}{{\textbf{X}}_2}$ and 
$\vc{{\textbf{Y}}_1}{{\textbf{Y}}_2}$ are independent.
\item $\text{vec}({\textbf{X}}_1, {\textbf{X}}_2, {\textbf{Y}}_1, 
{\textbf{Y}}_2) \sim N_{4n}(0, R_2(\Sigma_{2n}(\rho)))$.
\item $\text{vec}({\textbf{X}}_1, {\textbf{Y}}_1, {\textbf{X}}_2, 
{\textbf{Y}}_2) \sim N_{4n}(0, \Sigma_{4n}(\rho))$.
\item $\text{vec}({\textbf{X}}_1, {\textbf{Y}}_1, {\textbf{X}}_2, 
{\textbf{Y}}_2) = \sqrt{2}\text{vec}(\Re({\textbf{Z}}_1), 
\Im({\textbf{Z}}_1), \Re({\textbf{Z}}_2), 
\Im({\textbf{Z}}_2)) = \vc{J_2({\textbf{Z}}_1)}{J_2({\textbf{Z}}_2)}$, 
where $\vc{{\textbf{Z}}_1}{{\textbf{Z}}_2} \sim 
{\C}N_{2n}(0, \Sigma_{2n}(\rho))$.
\end{enumerate}
\end{corollary}
\noindent Similarly, under the inclusion of 
\sref{Corollary}{cor:the_random_norm_case}, respectively 
\cite[Proposition 8.7]{DF1993} (including a minor adjustment of the factor 
$c_1$ of the complex Gaussian probability measure in their proof, required due 
to the shape of the complex density function induced by the law $\C N_1(0, 1) = 
\R N_2(0, \frac{1}{2}I_2)$), we obtain 
\begin{lemma}\label{lem:correl_Gaussians_II}
Let $n \in \N$, $z \in \overline{\D} \cap \F$, $\textbf{X} = 
{\text{vec}}(X_1, \ldots, X_n) \sim {\F}N_n(0, I_n)$, $\vc{\textbf{Y}}
{\textbf{Z}} \sim {\F}N_{2n}(0, \Sigma_{2n}(z))$ and $u, v \in S_{\F^n}$. Then
\[
u^\top\textbf{X} \stackrel{d}{=} u^\ast\textbf{X} \sim \F N_1(0, 1) \,  
\text{ and } \, \colvec{2}{u^\ast\textbf{Y}}{v^\ast\textbf{Z}} \sim 
\F N_{2}(0, \Sigma_2((u^\ast v)z))\,.
\]
In particular, $\colvec{2}{u^\ast\textbf{X}}{v^\ast\textbf{X}} \sim
\F N_{2}(0, \Sigma_2(u^\ast v))$ and 
\[
\E\big[\big\vert \sum_{k=1}^n a_k\,X_k\big\vert^p\big] = 
\E\big[{\vert a^\top\textbf{X} \vert}^p\big] = 
\Vert a \Vert_{\F_2^n}^p\,\E[\vert X_1 \vert^p] = 
\Vert a \Vert_{\F_2^n}^p\,C_p^\F
\]
for all $a \equiv (a_1, \ldots, a_n)^\top \in \F^n$ and $p \in (-1, \infty)$, where 
$C_p^\R : = \frac{(\sqrt{2})^{p}}{\sqrt{\pi}}\,\Gamma(\frac{p+1}{2})$ and 
$C_p^\C : = \Gamma(1+\frac{p}{2})$.
\end{lemma}
\chapter{A quantum correlation matrix version of the Grothendieck inequality}
\section{Gram matrices, quantum correlation, and beyond}
\label{sec:Gram_matrices_and_quantum_correlation}
In this section, we aim at another equivalent reformulation of the Grothendieck 
inequality (occasionally abbreviated by ``GT'') for both fields, 
built on the inclusion of correlation matrices; i.e., positive semidefinite 
matrices with entries in $\F \in \{\R, \C\}$ whose diagonal is occupied with 
1's only. Since we need equivalent descriptions of a correlation matrix including 
its fundamental representation as a Gram matrix (cf. \cite[Theorem 2.7.10]{HJ2013}), 
we proceed with a fundamental acronym, to indicate a comprehensive 
class of matrices with entries in $\F$ which properly contains the class of 
all Gram matrices (cf. \cite[page 441]{HJ2013}) and reveals a deep connection 
to the foundations and philosophy of quantum mechanics. 
\begin{definition}\label{def:Gamma_H_u_v}
Let $\F \in \{\R, \C\}$, $m,n \in \N$ and $H$ be an $\F$-inner product space 
with inner product $\langle \cdot, \cdot \rangle \equiv \langle \cdot, \cdot 
\rangle_H$. Let $u \equiv (u_1, u_2, \ldots, u_m) \in H^m$ and $v \equiv 
(v_1, v_2, \ldots, v_n) \in H^n$. We put
\[
\Gamma_H(u,v) : = 
(\langle v_j, u_i \rangle)_{(i,j) \in [m] \times [n]} =
\begin{pmatrix}
    \langle v_1, u_1 \rangle & \langle v_2, u_1 \rangle & \ldots & \langle v_n, u_1 \rangle \\
    \langle v_1, u_2 \rangle & \langle v_2, u_2 \rangle & \ldots & \langle v_n, u_2 \rangle \\
		\vdots & \vdots & \ddots & \vdots \\
		\langle v_1, u_m \rangle & \langle v_2, u_m \rangle & \ldots & \langle v_n, u_m \rangle
  \end{pmatrix} \in \M_{m,n}(\F).
\]
\end{definition}
 
\noindent Observe that $\Gamma_H(u,v)^\ast = \Gamma_H(v,u)$ for all $(u, v) \in 
H^m \times H^n$. $\Gamma_H(u,v)$ should be viewed as an element of the 
image of the matrix-valued sesquilinear mapping
\[
\Gamma_H : H^m \times H^n \longrightarrow \M_{m,n}(\F).  
\]
In particular, if $\Gamma_H$ were restricted to the product $S_H^m \times 
S_H^n$, we would obtain a matrix-valued sesquilinear mapping, which is defined on the 
infinite-dimensional $C^\infty$-manifold $S_H^m \times S_H^n$; an interesting 
fact, which actually underlies our chosen notation (cf. 
\cite[Example 1.34 and Definition 1.36]{S2023}. If $m = n$ and $u=v \in H^m$, 
we get again the Gram matrix of the vectors $u_1, \ldots, u_m \in H$:
\[
\Gamma_H(u,u) = 
\begin{pmatrix}
    \Vert u_1 \Vert^2 & \langle u_2, u_1 \rangle & \ldots & \langle u_m, u_1 \rangle \\
    \langle u_1, u_2 \rangle & \Vert u_2 \Vert^2 & \ldots & \langle u_m, u_2 \rangle \\
		\vdots & \vdots & \ddots & \vdots \\
		\langle u_1, u_m \rangle & \langle u_2, u_m \rangle & \ldots & \Vert u_m \Vert^2
  \end{pmatrix}
	\in M_m(\F)^+\,.
\]
If $(r, s) \in \F^m \times \F^n$, then
\[
\Gamma_\F(r,s) \equiv \Gamma_{\F_2^1}(r,s) = \overline{r} s^\top = \overline{r} 
\otimes s^\top =
	\begin{pmatrix}
     \overline{r_1}\,s_1   &  \overline{r_1}\,s_2  & \ldots &  \overline{r_1}\,s_n  \\
     \overline{r_2}\,s_1  &  \overline{r_2}\,s_2  & \ldots &  \overline{r_2}\,s_n  \\
		\vdots & \vdots & \vdots & \vdots \\
		 \overline{r_m}\,s_1  &  \overline{r_1}\,s_2  & \ldots &  \overline{r_m}\,s_n 
  \end{pmatrix}\,.
\]
More generally, if $H = \F_2^d$ for some $d \in \N$, it follows that 
$\langle x, y \rangle_H = y^\ast x$ for all $x, y \in H$. Consequently, we obtain 
an important factorisation:
\begin{align}\label{eq:factorisation_of_Gamma_H_if_H_fin_dim}
\Gamma_{\F_2^d}(u,v) = U^\ast\,V \text{ for all } (u, v) \in (\F_2^d)^m \times (\F_2^d)^n,
\end{align}
where $U : = (u_1 \,\brokenvert\, u_2 \,\brokenvert\, \cdots \,\brokenvert\, u_m) \in 
\M_{d,m}(\F)$ and $V : = (v_1 \,\brokenvert\, v_2 \,\brokenvert\, \cdots \,\brokenvert\, v_n) 
\in \M_{d,n}(\F)$. Thus,
\begin{align}\label{eq:fact_property_of_Q_mn_and_the_left_side_of_GT}
{\text{tr}}(A^\ast\Gamma_{\F_2^d}(u, v)) = 
\overline{{\text{tr}}(\Gamma_{\F_2^d}(u, v)^\ast A)} =
\overline{{\text{tr}}(V^\ast (U A))} = \overline{\langle U A, V\rangle_{F}}
\end{align}
for all $A \in \M_{m,n}(\F)$ (cf. also 
\sref{Proposition}{prop:char_of_Q_m_n_for_both_fields}). 
A straightforward proof shows that any Gram matrix is positive semidefinite (cf. 
\cite[Theorem 2.7.10]{HJ2013}). Moreover, if $(a_{ij}) \equiv A = (A^{1/2})^2 = 
A^{1/2}(A^{1/2})^\top \in M_n(\F)^+$ is positive semidefinite and 
$\textbf{X} \sim {\F}N_n(0, I_n)$, then $\textbf{Z} : = A^{1/2}\textbf{X} \sim 
{\F}N_n(0, A)$, implying that $A = A^{1/2}\E[\textbf{X}\textbf{X}^\ast]A^{1/2} = 
\E[\textbf{Z}\textbf{Z}^\ast]$ (cf., e.g., \cite[Lemma 12.10.]{PR2016}) and 
$a_{ij} = \langle A^{1/2} e_j, A^{1/2} e_i \rangle_{\F_2^n}$ for all $i,j \in 
[n]$. Let us also recall the following characterisation of the set $C(n; \F)$ 
of all $n \times n$ correlation matrices (with entries in $\F$): 
\begin{lemma}\label{lem:charact_of_corr_matrices}
Let $\F \in \{\R, \C\}$, $n \in \N$ and $\Sigma = (\sigma_{ij})  \in 
\M_n(\F)$. Then the following statements are equivalent:
\begin{enumerate}
\item $\Sigma \in C(n; \F)$.
\item $\Sigma \in \M_n(\F)^+$ and $\sigma_{ii} = 1$ for all $i \in [n]$. 
\item There exist vectors $x_1, \ldots, x_n \in S_{\F_2^n}$ such that 
\[
\sigma_{ij} = \langle x_j, x_i\rangle_{\F_2^n} = x_i^\ast\,x_j \text{ for all } i, j \in [n].
\]
\item There exist an $\F$-Hilbert space $L$ and $x = (x_1, x_2, \ldots, x_n) 
\in L^n$ such that $\Vert x_i \Vert_L = 1$ for all $i \in [n]$ and 
\[
\Sigma = \Gamma_L(x,x). 
\]
\item $\Sigma = \E[\textbf{Z}\textbf{Z}^\ast]$ for some $n$-dimensional Gaussian 
random vector $\textbf{Z} \sim {\F}N_n(0, \Sigma)$, and $\sigma_{ii} = 1$ for 
all $i \in [n]$.  
\end{enumerate}
In particular, the set $C(n; \F)$ is convex, and $\vert \sigma_{ij} \vert \leq 1$ 
for all $i, j \in [n]$.
\end{lemma}
\begin{remark}[\textbf{The elliptope ${\mathcal{E}}_n \equiv C(n; \R)$}]
In the real case, the set of all $n \times n$-correlation matrices, which is 
very rich in geometrical and combinatorial structure, is also known as 
the so-called \textit{elliptope} (standing for \textit{ellip}\!\,soid and 
poly\!\,\textit{tope}) ${\mathcal{E}}_n$, studied in detail by M. Laurent and 
S. Poljak (cf. \cite[Example 5.44.]{BPT2013}, \cite[Chapter 5.9.1]{D2019} and 
\cite[Chapter 31.5]{DL1997}). From these sources, we learn, among many other 
deep facts, that for any $n \in \N$ the set ${\mathcal{E}}_n \equiv C(n; \R)$ is 
a convex polytope, which in general is not a polyhedron, so that it cannot be 
described as a finite intersection of weak half spaces (cf. \cite[Chapter 5.10]{AB2006}). 
Since both sets, $C(n; \R)$ and $C(n; \C)$ are compact (with respect to the 
topology of pointwise convergence), and since norms on finite-dimensional 
$\F$-vector spaces are equivalent, it follows that any linear functional on the 
finite-dimensional Hilbert space $(\M_n(\F), \Vert \cdot \Vert_F) \cong 
\F_2^{n^2}$ attains its maximum (and minimum) on the compact set $C(n; \F) 
\subseteq \M_n(\F)$, including the linear functional $\text{tr}(A^\ast\,\cdot) : 
\M_n(\F) \longrightarrow \F$, where $A \in \M_n(\F)$ is given. From the point 
of view of real (convex) semidefinite optimisation, both, the primal SDP 
\[
\sup\limits_{\Sigma \in C(m+n; \R)}\text{tr}(A^\top\Sigma) =
\sup\big\{\text{tr}(A^\top\Sigma) : \text{tr}(e_\nu e_\nu^\top \Sigma) = 1 
\text{ for all } \nu \in [m+n], \Sigma \in \M_{\nu}(\R)^+\big\}
\] 
and its dual SDP have non-empty, compact sets of optimal solutions and hence 
attain their respective optima (cf. \cite[Theorem 2.15, Theorem 2.29 and 
Exercise 2.41]{BPT2013} and \cite[Chapter 5]{BV2004}). We do not know whether this 
``strong duality'' is also valid in the complex case (cf. 
\sref{Lemma}{lem:correl_matrix_version_as_equality_of_sets}).
\end{remark}
\noindent Of particular relevance is the set 
\[
C_1(n; \F) : = \{\Theta : \Theta \in C(k; \F) \text{ and } 
{\text{rk}}(\Theta) = 1 \}
\]
of all $n \times n$ correlation matrices of rank 1. The Gram matrix structure 
implies a neat characterisation of $C_1(n; \F)$. At the same time, we 
recognise again the structure of all pure states on the Hilbert space 
$\F_2^n$. To this end, recall that $S_\R = \{-1,1\}$ and $S_\C = \T = 
\{z \in \C : \vert z \vert = 1\}$.
\begin{proposition}\label{prop:rep_of_C_1_n_F}
Let $m, n \in \N$. Then the following statements hold:
\begin{enumerate}
\item
\begin{align*}
\{A : A \in \M_{m,n}(\F) \text{ and } {\text{rk}}(A) = 1 \} &= 
\{\overline{p} q^\top : (p,q) \in \F^m\setminus\{0\} \times \F^n\setminus\{0\}\}\\ 
&= \{\Gamma_{\F}(p,q) : (p,q) \in \F^m\setminus\{0\} \times \F^n\setminus\{0\}\}.
\end{align*}
\item
\[
\{A : A \in \M_n(\F)^+ \text{ and } {\text{rk}}(A) = 1 \} = 
\{x x^\ast : x \in \F^n\setminus\{0\}\} = \{\Gamma_{\F}(z,z) : z \in 
\F^n\setminus\{0\}\}.
\]
In particular,
\[
\{A : A \in \M_n(\F)^+ \text{ and } {\text{rk}}(A) = 1 \text{ and } \text{tr}(A) 
= 1\} = \{x x^\ast : x \in S_{\F^n}\} = \{\Gamma_{\F}(z,z) : z \in S_{\F^n}\}\,.
\]
\noindent Moreover,
\[
\M_{n}(S_\F)^+ = \{x x^\ast : x \in S_\F^n\} = \{\Gamma_{\F}(z,z) : z \in 
S_\F^n\} = C_1(n; \F).
\] 
\end{enumerate}
\end{proposition}
\noindent Note that for all $m, n, \mu, \nu \in \N$, for all Hilbert spaces 
$H$, for all $(u, v) \in H^{m+n} \equiv H^{m} \times H^{n} $ and for all $(w, z)  
\in H^{\mu+\nu} \equiv H^{\mu} \times H^{\nu}$, the following block matrix 
representation in $\M_{m+n, \mu + \nu}(\F)$ always is satisfied:
\[
\begin{pmatrix}
\Gamma_H(u,w) & \Gamma_H(u,z) \\
\Gamma_H(v,w) & \Gamma_H(v,z)
\end{pmatrix}
= \Gamma_H((u,v),(w,z)).
\]
In particular,
\begin{align}\label{eq:l2d_rep_of_m_plus_n_correl_matrix}
\Gamma_{l_2^d}((u,v),(u,v)) = 
\begin{pmatrix}
U^\ast U & U^\ast V\\
V^\ast U & V^\ast V
\end{pmatrix} =
\begin{pmatrix}
U^\ast & 0\\
V^\ast & 0
\end{pmatrix}
\begin{pmatrix}
U & V\\
0 & 0
\end{pmatrix} = 
{\begin{pmatrix}
U & V\\
0 & 0
\end{pmatrix}}^\ast
\begin{pmatrix}
U & V\\
0 & 0
\end{pmatrix}
\end{align}
for any $d \in \N$, where $U : = (u_1 \,\brokenvert\, u_2 \,\brokenvert\, \cdots 
\,\brokenvert\, u_m) \in \M_{d,m}(\F)$ and $V : = (v_1 \,\brokenvert\, v_2 \,\brokenvert\, 
\cdots \,\brokenvert\, v_n) \in \M_{d,n}(\F)$. (due to 
\eqref{eq:factorisation_of_Gamma_H_if_H_fin_dim}). Regarding the topic of our 
work, the block structure of the elements of $C(m+n; \F)$ is of particular 
interest. To this end, put
\[
{\mathcal{Q}}_{m,n}(\F) : = \{S : S = \Gamma_H(u,v) \text{ for some $\F$-Hilbert space } H 
\text{ and } (u,v) \in S_H^m \times S_H^n\}.
\]
Any Hilbert space $(H, \langle \cdot, \cdot \rangle_H)$ can be isometrically 
embedded into the Hilbert space 
\[
(\widetilde{H}, \langle \cdot, \cdot \rangle_
{\widetilde{H}}) : = (H \oplus \F_2^2, \langle \cdot, \cdot \rangle_H + 
\langle \cdot, \cdot \rangle_{\F_2^2})
\]
(external direct sum of $H$ and $\F_2^2$). Thus, for any $(h,k) \in B_H \times 
B_H$, it follows by construction that
\begin{align}\label{eq:extension_from_B_H_to_S_K}
\widetilde{h} : = (h, 0, \sqrt{1-\Vert h \Vert_H}) \in S_{\widetilde{H}},\,
\widetilde{k} : = (k, \sqrt{1-\Vert k \Vert_H}, 0) \in S_{\widetilde{H}} \text{ and }
\langle \widetilde{k}, \widetilde{h} \rangle_{\widetilde{H}} = 
\langle k, l \rangle_H\,.
\end{align}
Note that $\text{dim}({\widetilde{H}}) \geq 3$. Consequently,
\begin{align}\label{eq:extension_from_B_H_to_S_K_II}
{\mathcal{Q}}_{m,n}(\F) = \{S : S = \Gamma_H(x,y) \text{ for some $\F$-Hilbert 
space } H \text{ and } (x,y) \in B_H^m \times B_H^n\}.
\end{align}
Observe also that for $C(\nu; \F) \subseteq {\mathcal{Q}}_{\nu,\nu}(\F)$ for all 
$\nu \in \N$. By ``inflating'' the set ${\mathcal{Q}}_{m,n}(\F)$ to a 
$(m+n)\times(m+n)$\,-\,correlation matrix, we obtain the non-trivial fact that 
${\mathcal{Q}}_{m,n}(\F)$ is absolutely convex:  
\begin{corollary}\label{cor:block_matrix_structure_of_C_m_plus_n}
Let $m,n \in \N$ and $\Sigma \in \M_{m+n}(\F)$. Then the following statements 
are equivalent:
\begin{enumerate}
\item $\Sigma \in C(m+n; \F)$.
\item
$\Sigma = 
\begin{pmatrix}
\Gamma_H(u,u) & \Gamma_H(u,v)\\
\Gamma_H(u,v)^\ast & \Gamma_H(v,v)
\end{pmatrix}$, for some Hilbert space $H$ over $\F$ and some $(u,v) \in S_H^m 
\times S_H^n$.
\item
$\Sigma = 
\begin{pmatrix}
\Gamma_H(u,u) & \Gamma_H(u,v)\\
\Gamma_H(u,v)^\ast & \Gamma_H(v,v)
\end{pmatrix}$, for some Hilbert space $H$ over $\F$ and some $(u,v) \in B_H^m 
\times B_H^n$.
\item There exist $d \in \N$, $U \in \M_{d,m}(\F)$ and $V \in 
\M_{d,n}(\F)$, such that $u_i : = Ue_i \in S_{\F_2^d}$ for all $i \in [m]$, 
$v_j : = V e_j \in S_{\F_2^d}$ for all $j \in [n]$ and
\[
\Sigma = 
\begin{pmatrix}
U^\ast U & U^\ast V\\
V^\ast U & V^\ast V
\end{pmatrix} =
\begin{pmatrix}
U^\ast & 0\\
V^\ast & 0
\end{pmatrix}\!\!
\begin{pmatrix}
U & V\\
0 & 0
\end{pmatrix}.
\]
\end{enumerate}
In particular, $S \in {\mathcal{Q}}_{m,n}(\F)$ if and only if there exist 
correlation matrices $A \in C(m; \F)$ and $B \in C(n; \F)$, such that
$\begin{pmatrix}
A & S\\
S^\ast & B
\end{pmatrix} \in C(m+n; \F)$ is a correlation matrix. The set 
${\mathcal{Q}}_{m,n}(\F)$ is absolutely convex.
\end{corollary}
\begin{remark}[\textbf{Tsirel'son's characterisation of quantum 
correlation matrices}]\label{rem:quantum_correlation_matrix}
In the real case, i.e., if $\F = \R$, ${\mathcal{Q}}_{m,n} \equiv 
{\mathcal{Q}}_{m,n}(\R)$ \textit{coincides} with the class of so-called 
\textit{quantum correlation matrices}. These matrices are particularly essential 
in the foundations and philosophy of quantum mechanics (cf. 
\cite[Chapter 11]{AS2017}, \cite{L2018, L2020}, \cite[Theorem 1]{T1980} and 
\cite[Section 4]{T1987}). To this end, recall that in quantum mechanics a 
\textit{matrix} $\rho \in \M_{n,n}(\C)$ is called a \textit{(quantum) state} if 
$\rho \in \M_{n,n}(\C)^+$ and $\text{tr}(\rho) = 1$ (cf., e.g., 
\cite[Chapter 0.10]{AS2017}). In fact, we have (cf. \cite[Chapter 11]{AS2017}):
\begin{thmTS*}
\textsl{
Let $m, n \in \N$ and $S \equiv (s_{ij}) \in \M_{m,n}(\R)$. Then the following 
statements are equivalent:
\begin{enumerate}
\item $S \in {\mathcal{Q}}_{m,n}(\R)$.
\item There exists a unital $C^\ast$-algebra $\mathfrak{A}$, self-adjoint 
elements $A_1, \ldots, A_m$, $B_1, \ldots, B_n$ and a state $\tau$ on $\mathfrak{A}$, 
such that
$A_i\,B_j = B_j\,A_i$, $\max\{\Vert A_i\Vert, \Vert B_j\Vert\} \leq 1$ and 
$s_{ij} = \tau(A_i B_j)$ for all $(i,j) \in [m] \times [n]$.
\item There is a state $\rho \in \M_{d_1\cdot d_2}(\C)$ (for some $d_1, d_2 
\in \N$), Hermitian matrix families  $(W_1, \ldots, W_m) \in 
{\M_{d_1}(\overline{\D})}^m$ and $(Z_1, \ldots, Z_n) \in 
{\M_{d_2}(\overline{\D})}^n$, such that
\[
s_{ij} = \text{tr}((W_i \otimes Z_j)\rho)\,\text{ for all }\, (i,j) \in 
[m] \times [n]\,. 
\] 
\end{enumerate}
}
\end{thmTS*}
\noindent Here, according to the construction, the Hermitian matrix 
$W_i \otimes Z_j \in \M_{d_1\cdot d_2}(\overline{\D})$ is given by the Kronecker 
product of $W_i \in \M_{d_1}(\overline{\D})$ and $Z_j \in \M_{d_2}(\overline{\D})$.

In particular, if $k \in \N_3$, then \textit{every} real standard 
$(k \times k)$-correlation matrix -- used in everyday statistical calculations -- 
actually contains a \textit{quantum} correlation matrix block part 
$\in \mathcal{Q}_{m, k-m}$ ($m \in [k-1]$) and its transpose\,! Although, the 
Grothendieck inequality actually ``compares'' the set 
${\mathcal{Q}}_{m,n}(\F)$ of all real (respectively complex\,!) quantum 
correlation matrices with their extreme counterparts of rank 1 (cf. 
\sref{Proposition}{prop:rep_of_C_1_n_F}-\ref{(i)}, \eqref{eq:max_at_the_boundary_1}, 
\eqref{eq:max_at_the_boundary_2}, \autoref{thm:GT_rewritten_form}, 
\sref{Corollary}{cor:GT_equivalent_statements} and 
\autoref{thm:nuclear_unit_ball}), Tsirel'son's groundbreaking result, 
\textit{per se}, won't be discussed in detail in this paper, though. Regarding a 
detailed introduction to this fascinating subject including full and detailed 
proofs of Tsirel'son's results, we particularly refer to 
\cite[Ch. 11.2]{AS2017} and \cite{GLPV2017, L2018, L2020}, and the references 
therein. 
\end{remark}
\section{The Grothendieck inequality, correlation matrices and the matrix 
norm $\Vert \cdot \Vert_{\infty, 1}^\F$}
Fix $\F \in \{\R, \C\}$ and $m, n \in \N$. Let $A \in \M_{m,n}(\F)$ and 
$H$ be an arbitrary Hilbert space. If $(x,y) \in H^m \times H^n$, then
\[
\overline{\sum_{i = 1}^{m}\sum_{j = 1}^{n} a_{ij}
\langle x_i, y_j\rangle_H} = \sum_{i = 1}^{m}\sum_{j = 1}^{n}\overline{a_{ij}}\,
\langle y_j, x_i\rangle_H = \text{tr}(A^\ast\Gamma_H(x, y)) = 
\overline{\text{tr}(A\,\Gamma_H(y, x))}.
\]
Moreover, observe that (cf. \cite[Lemma 2.2.]{FLZ2018}, \cite[Remark 10.1]{J1987}
and \sref{Remark}{rem:attainability_of_max} regarding the existence of the 
respective maxima) 
\begin{align}\label{eq:max_at_the_boundary_1}
\begin{split}
\Vert A \Vert_H^{\text{G}} &:= 
\max_{\Vert u_i \Vert = 1, \Vert v_j \Vert = 1}\big\vert \sum_{i = 1}^{m}\sum_{j = 1}^{n} 
a_{ij}\,\langle u_i, v_j\rangle_H \big\vert 
= \max_{(u, v) \in S_H^m \times S_H^n}\vert\text{tr}(A^\ast\Gamma_H(u, v))\vert\\
&\,= \max_{\Vert u_i \Vert \leq 1, \Vert v_j \Vert \leq 1}\big\vert \sum_{i = 1}^{m}
\sum_{j = 1}^{n} a_{ij}\langle u_i, v_j\rangle_H \big\vert\,.
\end{split}
\end{align}
In particular (if $H = \F$, where $\langle z, w\rangle_H =\overline{w} z$ for 
all $z, w \in \F$), we have:
\begin{align}\label{eq:max_at_the_boundary_2}
\Vert A \Vert_\F^{\text{G}} = \max\limits_{(p,q) \in S_\F^m \times 
S_\F^n}\vert {\text{tr}}(A^\ast\Gamma_\F(p, q))\vert =
\max_{\vert p_i \vert \leq 1, \vert q_j \vert \leq 1}
\big\vert \sum_{i = 1}^{m}\sum_{j = 1}^{n} a_{ij} p_i \overline{q_j}\big\vert\,,
\end{align}
Consequently, if the matrix $A \in \M_{m,n}(\F)$ is viewed as a bounded 
linear operator from $l_\infty^n$ into $l_1^m$, then 
\eqref{eq:max_at_the_boundary_2} and the fact that $l_\infty^m$ 
is isometrically isomorphic to the dual space of $l_1^m$ 
(via the linear map $\chi : l_\infty^m \longrightarrow (l_1^m)^\prime$, 
defined as $l_1^m \ni z \mapsto \langle z, \chi(q) \rangle := q^\top z 
= \sum_{i=1}^m q_i z_i$), implies that
\begin{align}\label{eq:matrix_infty_1_norm}
\begin{split}
\Vert A \Vert_\F^{\text{G}} 
&= \max\limits_{(p,q) \in S_\F^m \times S_\F^n}
\vert {\text{tr}}(A^\ast\Gamma_\F(p, q))\vert 
\stackrel{\eqref{eq:max_at_the_boundary_2}}{=}
\max_{(p,q) \in B_\infty^m \times B_\infty^n}\vert {\text{tr}}(A^\ast\Gamma_
\F(p, q))\vert\\
&= \max_{(p,q) \in B_\infty^m \times B_\infty^n}\vert\text{tr}(A\,\overline{q} 
p^\top)\vert = \max_{(p,q) \in B_\infty^m \times B_\infty^n}\vert\text{tr}(A\,q 
p^\top)\vert = \max_{(p,q) \in B_\infty^m \times B_\infty^n}\vert\langle A q, 
\chi(p)\rangle\vert\\
&= \Vert A \Vert_{{\mathfrak{L}}(l_\infty^n, l_1^m)} =: \Vert A \Vert_{\infty, 1}^\F\,,
\end{split}
\end{align}
where $B_\infty^\nu : = B_{l_\infty^\nu}$, $\nu \in \N$. Consequently, Theorem 
\ref{thm:GT_matrix_form} is equivalent to
\begin{theorem}
\label{thm:GT_rewritten_form}
Let $\F \in \{\R, \C\}$. There is an absolute constant $K>0$ such that for 
any $m, n \in \N$, for any $A \in \M_{m,n}(\F)$ and any $\F$-Hilbert space $H$, 
the following inequality is satisfied:
\begin{align*}
\max\limits_{(u,v) \in S_H^m \times S_H^n}
\vert{\text{tr}}(A^\ast\,\Gamma_H(u,v))\vert = \Vert A \Vert_H^{\text{G}} \leq K\,
\Vert A \Vert_{\infty, 1}^\F\,. 
\end{align*}
$K_G^\F > 1$ is the smallest possible value of the corresponding absolute constant 
$K$.
\end{theorem}
\noindent The operator norm $\Vert \cdot \Vert_{\infty, 1}^\F$ on the right side 
of the Grothendieck inequality is a particular example of a \textit{mixed 
subordinate matrix norm} (cf., e.g., \cite[A.1.5]{BV2004} and \cite{R2000, 
S2005}). 
Here, we have to recall two key results, particularly regarding the 
computational complexity of $\Vert A \Vert_{\infty, 1}^\R$ (cf. 
\cite{HO2010, RS2009, R2000, S2005}): 
\begin{theorem}[\textbf{Rohn, 2000}]
Computing $\Vert A \Vert_{\infty, 1}^\R$ is NP-hard in the class of Maximum Cut 
Matrices. 
\end{theorem}
\noindent Even an approximation of $\Vert A \Vert_{\infty, 1}^\R$ is NP-hard 
(see also \cite[Theorem 6]{R2000}):
\begin{theorem}[\textbf{Hendrickx and Olshevsky, 2010}]
Unless $P = NP$, there is no polynomial time algorithm which, given a real 
matrix $A$ with entries in $\{-1, 0, 1\}$, approximates 
$\Vert A \Vert_{\infty, 1}^\R$ to some fixed error with polynomial running time 
in the dimensions of the matrix. 
\end{theorem}
\noindent These observations immediately result in another important well-known 
fact which will be used later in this paper to show that \textit{for both fields} 
the calculation of $K_G^\F$ can also be elaborated by means of semidefinite 
programming, which is a \textit{convex} optimisation problem (cf. 
\sref{Corollary}{cor:GT_equivalent_statements}, 
\sref{Proposition}{prop:correlation_matrix_version_of_GT} and 
\cite[Chapter 4.6.2]{BV2004}). Namely,   
\begin{lemma}\label{lem:correl_matrix_version_as_equality_of_sets}
Let $\F \in \{\R, \C\}$, $m, n \in \N$ and $A \in \M_{m,n}(\F)$. Let $\texttt{HIL}^\F$ 
denote the class of all $\F$-Hilbert spaces. Put
\begin{align}\label{eq:Delta_transform_of_A}
\Delta(A) \equiv \Delta^\F(A) : = 
\tfrac{1}{2}\begin{pmatrix}
0 & A\\
A^\ast & 0
\end{pmatrix}.
\end{align}
Then $\Delta(A)$ is Hermitian and
\begin{align}\label{eq:quadratic_form_and_Delta_A} 
\text{vec}(x,y)^\ast\Delta(A)\text{vec}(x,y) = \Re(x^\ast Ay) \text{ for all } 
(x, y) \in \F^m \times \F^n\,. 
\end{align}
In particular, $\Delta(A) \in \M_{m+n}(\F)^+$ if and only if $A = 0$. Moreover,
\begin{align}\label{eq:Delta_A_and_equality_of_traces}
{\text{tr}}\big(\Delta(A) 
\begin{pmatrix}
C & S \\
R^\ast & D
\end{pmatrix}\big) = \tfrac{1}{2}(\overline{{\text{tr}}(A^\ast R)} + 
{\text{tr}}(A^\ast S))
\end{align}
for all $(C, D) \in \M_m(\F) \times \M_n(\F)$ and $S, R \in \M_{m,n}(\F)$. $0 \leq 
\max\limits_{\Sigma \in C(m+n; \F)}\text{tr}(\Delta(A)\Sigma) < \infty$,
and
\begin{align}\label{eq:m_n_GT_reformulated}
\begin{split}
\sup\limits_{H \in \texttt{HIL}^\F}\Vert A \Vert_H^{\text{G}} 
&= \max\limits_{S \in {\mathcal{Q}}_{m,n}(\F)}\vert\text{tr}(A^\ast S)\vert = 
\max\limits_{S \in {\mathcal{Q}}_{m,n}(\F)}\Re(\text{tr}(A^\ast S))\\
&= \max\limits_{\Sigma \in C(m+n; \F)}\text{tr}(\Delta(A)\Sigma) \leq K_G^\F(m,n)
\Vert A \Vert_{\infty, 1}
\end{split}
\end{align}
\end{lemma}
\begin{comment}
\eqref{eq:quadratic_form_and_Delta_A} and its impact on the positive 
semidefiniteness of $\Delta(A)$ follows immediately. Since
\[
\Delta(A)\begin{pmatrix}
C & S \\
R^\ast & D
\end{pmatrix} = \tfrac{1}{2}\begin{pmatrix}
A R^\ast & A D \\
A^\ast C & A^\ast S
\end{pmatrix},
\]
it follows that
\[
\text{tr}\big(\Delta(A)\begin{pmatrix}
C & S \\
R^\ast & D
\end{pmatrix}\big) = \tfrac{1}{2}(\text{tr}(A R^\ast) + 
\text{tr}(A^\ast S)) = \tfrac{1}{2}(\overline{\text{tr}(A^\ast R)} + 
\text{tr}(A^\ast S)).
\]
Observe, that \eqref{eq:Delta_A_and_equality_of_traces} does not at all 
depend on the choice of the matrices $C$ and $D$\,! Regarding the proof of 
\eqref{eq:m_n_GT_reformulated}, we firstly verify the second equality. 
To this end, fix $S \in {\mathcal{Q}}_{m,n}(\F)$. Since any $z \in \F$ satisfies 
$\vert z \vert = \Re(\zeta z)$, where $\zeta : = 1$ if $z=0$ and $\zeta : = 
\frac{\overline{z}}{\vert z \vert}$ if $z \not= 0$, it follows 
in particular that $\vert{\text{tr}}(A^\ast S)\vert = 
\Re(\alpha\,{\text{tr}}(A^\ast S)) = \Re({\text{tr}}(A^\ast \alpha S)) \geq 0$ 
for some $\alpha \in S_\F$. Since $\alpha S \in {\mathcal{Q}}_{m,n}(\F)$, the 
completion of the proof of \eqref{eq:m_n_GT_reformulated} now follows by applying 
\eqref{eq:Delta_A_and_equality_of_traces} twice, including the fact that 
$\vert \text{tr}(\Delta(A)\Sigma)\vert \leq \max\limits_
{V \in {\mathcal{Q}}_{m+n,m+n}(\F)}\vert\text{tr}(\Delta(A)V)\vert$ for any 
$\Sigma \in C(m+n; \F) \subseteq {\mathcal{Q}}_{m+n,m+n}(\F)$. 
\end{comment}
\noindent Let $A \in \M_n(\F)$. Put $d_1(A) \equiv d_1^\F(A) : = 
\sup\limits_{\Theta \in C_1(n; \F)}\vert{\text{tr}}(A^\ast \Theta)\vert$. 
Because of \sref{Proposition}{prop:rep_of_C_1_n_F}-\ref{(ii)} it follows that
\[
d_1(A) = \max\limits_{x \in S_\F^n}\vert x^\ast Ax\vert = \max\limits_{x \in S_\F^n}
\vert\text{tr}(A^\ast x x^\ast)\vert \leq \Vert A \Vert_{\infty, 1}^\F\,.
\] 
Observe that $d_1 \not= \Vert \cdot \Vert_{\infty, 1}^\F$ (due to 
\cite[Corollary 2.11]{FL2020}). 
If $B$ is symmetric, respectively Hermitian, then $d_1(B)$ coincides with 
the seminorm $\Vert B \Vert_{\gamma, 1}$ of S. Friedland and L.-H. Lim (cf. 
\cite[Proposition 2.5]{FL2020}). Moreover, in the positive semidefinite case 
\cite[Proposition 2.8]{FL2020} implies that
\begin{align}\label{eq:d_1_psd_case}
d_1(M) = \max\limits_{x \in S_\F^n}\vert x^\ast Mx\vert = 
\max\limits_{x \in [-1,1]^n}\vert x^\ast Mx\vert \text{ for all } M \in 
\M_n(\F)^+\,.
\end{align}
Recall from \sref{Lemma}{lem:correl_matrix_version_as_equality_of_sets} the 
Hermitian matrix $\Delta(A) \equiv \Delta^\F(A) : = 
\frac{1}{2}\begin{pmatrix}
0 & A\\
A^\ast & 0
\end{pmatrix}$, where $m, n \in \N$ and $A \in M(m\times n; \F)$. 
Observe that in the following inequalities, which are an immediate application of 
\cite[Corollary 2.6., (29) and Proposition 2.8., (31)]{FL2020}, seemingly 
no Hilbert space presence is required (cf. also \eqref{eq:m_n_GT_reformulated} 
and \sref{Corollary}{cor:Nesterov_as_special_case}). 
\begin{proposition}\label{prop:correlation_matrix_version_of_GT}
Let $\F \in \{\R, \C\}$ and $m, n \in \N$. 
Then $K_G^\F(m,n)$ is the smallest constant $c > 0$, satisfying
\begin{align}\label{eq:symmetric_block_matrix_version_of_GT}
\vert{\text{tr}}(\Delta(A)\,\Sigma)\vert \leq c\,\Vert A \Vert_{\infty, 1} 
\text{ for all } \Sigma \in C(m+n;\F) \text{ and } A \in \M_{m,n}(\F). 
\end{align}
The little Grothendieck constant $k_G^\F$ is the smallest constant $\gamma > 0$, 
such that
\begin{align}\label{eq:corr_matrix_version_of_GT_in_the_PSD_case}
\vert{\text{tr}}(B\,\Sigma)\vert \leq \gamma\,\Vert B \Vert_{\infty, 1} 
\text{ for all } \Sigma \in C(n;\F) \text{ and } B \in \M_n(\F)^+. 
\end{align}
\end{proposition}
\noindent In fact, if we allow the implementation of a possibly strictly larger 
absolute constant than $K_G^\F$, our approach leads to a further, more general 
inequality, which encompasses the real and the complex Grothendieck inequality 
as a special case (cf. \autoref{thm:GT_generalising_Krivine} (real case), 
respectively \autoref{thm:odd_complex_bdd_CCP_general_version} 
(complex case)). Moreover, it extends the symmetric Grothendieck equality of 
Friedland and Lim in \cite{FL2020} from symmetric $\F$-matrices to arbitrary 
$\F$-matrices (cf. \eqref{eq:abs_cx_hull} and \sref{Remark}{rem:on_Delta_A}): 
\begin{theorem}
Let $\F \in \{\R, \C\}$. Then there exists an absolute constant $K_\ast^\F > 1$ 
such that
\[
\vert{\text{tr}}(B^\ast\,\Sigma)\vert \leq 
K_\ast^\F d_1(B)\,. 
\]
for any $k \in \N$, any $\Sigma \in C(k; \F)$ and any $B \in \M_k(\F)$. Moreover, 
\[
C(k; \F) \subseteq K_\ast^\F \,{\text{acx}}(C_1(k; \F)) \text{ for all } k \in \N\,,
\]
$K_\ast^\R \in [K_G^\R, \sinh(\frac{\pi}{2})]$ and $K_\ast^\C \in 
[K_G^\C, \frac{8}{\pi}-1]$.
\end{theorem}
\section{Characterisation of $K_G^\F$ through operator ideals and violation 
of Bell inequalities (a brief digression)} 
 
Readers who are familiar with operator ideals in the sense of A. Pietsch 
(cf. \cite{DF1993, DJT1995, P1980}) should take notice of \sref{Remark}{rem:GT_norm} 
below regarding the Grothendieck norm \eqref{eq:max_at_the_boundary_1} on the 
left side of the Grothendieck inequality. To round out the picture, we list a 
rather elementary result (\sref{Proposition}{prop:char_of_Q_m_n_for_both_fields}), 
which however unveals the link between matrices in ${\mathcal{Q}}_{m,n}(\F)$ and 
a well-known functional analytic formulation of the Grothendieck inequality; 
namely as an inequality between Banach ideal matrix norms, when matrices are 
viewed as bounded linear operators from $l_\infty^n$ into $l_1^m$ (respectively 
from $l_1^n$ into $l_\infty^m$), being elements of certain $1$-Banach ideals 
(in the sense of A. Pietsch), or equivalently as an inequality between certain 
tensor product norms on tensor products of Banach spaces (cf. 
\sref{Remark}{rem:tensor_norm_rep_of_GT}). The latter approach 
was developed by Grothendieck. Regarding the underlying functional analytic 
details, we refer the readers to \cite{DF1993, P1980}. Primarily, we need an 
intuitive understanding of bounded linear operators between (finite-dimensional) 
Banach spaces, \textit{factoring} through a Hilbert space and the basics of 
nuclear operators (cf. \cite[Chapter 6.3]{P1980}).
 
 
So, let $E$ and $F$ be $\F$-Banach spaces and $S \in {\mathfrak{L}}(E, F)$. By 
definition, $S \in {\mathfrak{L}}_2(E, F)$ if and only if there exist an 
$\F$-Hilbert space $H$ and bounded linear operators $R \in \mathfrak{L}(H, F)$, 
$T \in \mathfrak{L}(E, H)$, such that $S = RT$. $S \in {\mathfrak{L}}_2(E,F)$ is 
said to be \textit{2-factorable} (cf., e.g., \cite[Corollary 18.6.2]{DF1993}).
\setlength{\unitlength}{0.4cm}
\begin{center}
\begin{picture}(20,0)\thicklines
 \put(6.8,0){\vector(1,0){5.6}}
 \put(5.6,-0.2){$E$} 
 \put(12.8,-0.2){$F$}
 \put(6.5,-0.5){\vector(1,-1){2.5}}
 \put(9.0,-3.5){$H$} 
 \put(10,-3){\vector(1,1){2.5}}
 \put(9.3,0.25){$S$}
 \put(11.6,-2.2){$R \in \mathfrak{L}$}
 \put(4.2,-2.2){$\mathfrak{L} \ni T$}
 \end{picture}
\end{center}
\vspace{1.5cm}
It can be shown that $({\mathfrak{L}}_2(E, F), \Vert \cdot 
\Vert_{{\mathfrak{L}}_2(E, F)})$ is an $\F$-Banach space. The norm is defined as
\[
\Vert S : E \longrightarrow F\Vert_{{\mathfrak{L}_2}} \equiv \Vert S \Vert_
{{\mathfrak{L}}_2(E, F)} : = \inf\Vert R\Vert\,\Vert T\Vert,
\]
where the infimum is taken over all factorisations $S = RT$ through any Hilbert 
space. 

The unit ball of the Banach space ${\mathfrak{L}}_2(l_1^n, l_\infty^m)$ completely 
characterises the convex set ${\mathcal{Q}}_{m,n}(\F)$, since:
\begin{proposition}\label{prop:char_of_Q_m_n_for_both_fields}
Let $m, n \in \N$ and $S \in \M_{m,n}(\F)$. Then the following statements are equivalent:
\begin{enumerate}
\item $S \in B_{{\mathfrak{L}}_2(l_1^n, l_\infty^m)}$.
\item There exist $d \in \N$
and $(u, v) \in {(S_{\F_2^d})}^m \times {(S_{\F_2^d})}^n$ such that $S = 
\Gamma_{\F_2^d}(u, v)$.
\item There exist an $\F$-Hilbert space $H$ 
and $(u, v) \in B_H^m \times B_H^n$ such that $S = \Gamma_H(u, v)$.
\item There exist $d \in \N$, $U \in \M_{d,m}(\F)$ and $V \in \M_{d,n}(\F)$ such 
that $\Vert U : l_1^m \longrightarrow l_2^d \Vert = 1$, $\Vert V : l_1^n \longrightarrow 
l_2^d \Vert = 1$ and $S = U^\ast V$.
\end{enumerate}
\end{proposition}
\noindent In other words, if $m, n \in \N$ and $\F \in \{\R, \C\}$, then
\begin{align}\label{eq:char_of_quantum_corr_matrices}
B_{{\mathfrak{L}}_2(l_1^n, l_\infty^m)} = {\mathcal{Q}}_{m,n}(\F) = 
\bigcup_{d=1}^\infty\{\Gamma_{\F_2^d}(u, v) : (u,v) \in (S_{\F_2^d})^m \times 
(S_{\F_2^d})^n\}.
\end{align}
\begin{corollary}\label{cor:real_vs_complex_quantum_correlations}
Let $d, m, n \in \N$, $(z, w) \in (\C_2^d)^m \times (\C_2^d)^n$ and 
$(a, b) \in (\R_2^{2d})^m \times (\R_2^{2d})^n$. Then
\[
\Gamma_{\C_2^d}(z, w) = \Gamma_{\R_2^{2d}}(x, y) + i\,\Gamma_{\R_2^{2d}}(x, y^\prime)\,,
\]
and
\[
\Gamma_{\R_2^{2d}}(a, b) = \Re(\Gamma_{\C_2^d}(\zeta, \xi)),
\]
where $(\zeta, \xi) \in (\C_2^d)^m \times (\C_2^d)^n$, $(x, y) \in (\R_2^{2d})^m 
\times (\R_2^{2d})^n$ and $y^\prime \in (\R_2^{2d})^n$ are given as $x_i : = 
J_2(z_i), y_j : = J_2(w_j)$, $y_j^\prime : = J_2(-i\,w_j) = R_2(-i\,Id_d)y_j$, 
$\zeta_i : = J_2^{-1}(a)$ and $\xi_j : = J_2^{-1}(b)$ $((i,j) \in [m]\times[n])$. 
In particular, $\{\Re(S) : S \in {\mathcal{Q}}_{m,n}(\C)\} \subseteq 
{\mathcal{Q}}_{m,n}(\R)$, $\{\Im(S) : S \in {\mathcal{Q}}_{m,n}(\C)\} \subseteq 
{\mathcal{Q}}_{m,n}(\R)$ and ${\mathcal{Q}}_{m,n}(\C) \subseteq {\mathcal{Q}}_{m,n}(\R) + 
i\,{\mathcal{Q}}_{m,n}(\R)$. Moreover, $\{\Re(\Sigma) : \Sigma \in C(n; \C)\} 
\subseteq C(n; \R)$.
\end{corollary}
\begin{comment}
We just have to apply \eqref{eq:char_of_quantum_corr_matrices}, \eqref{eq:Re_part} 
and \eqref{eq:Im_part}, taking into account the algebraic properties of the 
mappings $J_2$ and $R_2$.
\end{comment}
\noindent Let $(u, v) \in S_H^m \times S_H^n$. Since $\Gamma_H(u,v)^\ast = 
\Gamma_H(v,u)$, it follows that also $\Gamma_H(u,v)^\ast \in {\mathcal{Q}}_
{n,m}(\F)$. 
Consequently, if we recall \autoref{thm:GT_rewritten_form} and
\sref{Lemma}{lem:correl_matrix_version_as_equality_of_sets}, we arrive at the 
following crucial implication of 
\sref{Proposition}{prop:char_of_Q_m_n_for_both_fields} (cf. \cite{GLPV2017}):
\begin{corollary}\label{cor:GT_equivalent_statements}
Let $\F \in \{\R, \C\}$ and $m, n \in \N$. Let $A \in \M_{m,n}(\F)$ such 
that $\Vert A \Vert_{\infty, 1}^\F \leq 1$. Then the following 
statements are equivalent to each other and to \eqref{eq:GT_for_given_m_and_n}:
\begin{enumerate}
\item
\[
\vert{\text{tr}}(A^\ast\,\Gamma_H(u,v))\vert
\leq K_G^\F(m,n) \text{ for all Hilbert spaces }
H \text{ and for all } (u,v) \in S_H^m \times S_H^n\,.
\]
\item
\begin{align*}
\max\limits_{S \in {\mathcal{Q}}_{m,n}(\F)}
\vert{\text{tr}}(A^\ast\,S)\vert = 
\max\limits_{R \in B_{{\mathfrak{L}}_2(l_1^m, l_\infty^n)}}
\vert{\text{tr}}(A\,R)\vert \leq K_G^\F(m,n)\,.
\end{align*}
\item
\begin{align*}
\max\limits_{\Sigma \in C(m+n; \F)}\text{tr}(\Delta(A)\Sigma) \leq K_G^\F(m,n).
\end{align*}
\end{enumerate}
\end{corollary}
\noindent Recall that for any $\nu \in \N$, the canonical isometric isomorphisms 
$(l_1^\nu)^\prime \cong l_\infty^\nu$, and $(l_\infty^\nu)^\prime \cong l_1^\nu$ 
(since $l_\infty^\nu$ is finite-dimensional) explicitly characterise the 
respective dual spaces. Put 
\begin{align}\label{eq:spaces}
X := {\mathfrak{L}}_2(l_1^n, l_\infty^m), Y : = {\mathfrak{L}}(l_\infty^m, l_1^n),
Z : = {\mathfrak{N}}(l_1^n, l_\infty^m) \text{ and } 
W : = {\mathfrak{L}}(l_\infty^n, l_1^m),
\end{align} 
where $({\mathfrak{N}}, \Vert \cdot \Vert_{\mathfrak{N}})$ denotes the 
$1$-Banach ideal of all nuclear operators, which is the smallest $1$-Banach 
ideal (originally created by A. Grothendieck in his famous thesis \cite{G1955}). 
Let us quickly recall that a linear operator $T : E \longrightarrow F$ 
between two $\F$-Banach spaces $E$ and $F$ is said to be \textit{nuclear} if 
there exist sequences $(a_n)_{n \in \N} \subseteq B_{E^\prime}$, 
$(y_n)_{n \in \N} \subseteq B_F$ and $(\lambda_n)_{n \in \N} \in l_1$ such that
\[
T = \sum_{n=1}^\infty \lambda_n\langle \cdot, a_n\rangle y_n\,,
\]
and the series converges in $\mathfrak{L}(E,F)$ (cf. \sref{Remark}{rem:GT_vs_trace_class}, 
\cite{G1955} and \cite[Chapter 6.3 and Theorem 9.2.1.]{P1980}). It is a 
well-known fact that the Banach space $Z$ is isometrically 
isomorphic to the dual Banach space $Y^\prime$. The isometric isomorphism 
$\tau : Z \longrightarrow Y^\prime$ is given by canonical trace duality
\begin{align}\label{eq:trace_duality_fin_dim_case} 
Z \ni S \mapsto \tau_S \equiv \tau(S), \text{ where }
\langle B, \tau_S\rangle : = {\text{tr}}(BS) \text{ for all } B \in Y
\end{align}
(see \cite[Theorem 9.2.1]{P1980}). Readers, who are familiar with tensor norms 
and Banach ideals could verify the above trace duality very quickly, 
(since $Y^{\,\prime} \cong (l_1^m \otimes_{\varepsilon} l_1^n)^\prime 
\cong (l_1^n \otimes_{\varepsilon} l_1^m)^\prime \cong 
\mathfrak{I}(l_1^n, l_\infty^m) = Z$ (cf. 
\cite[Corollary 5.7.1 and Proposition 16.7]{DF1993})). Observe that
\sref{Proposition}{prop:char_of_Q_m_n_for_both_fields} implies that 
$B_Z \subseteq B_X = {\mathcal{Q}}_{m,n}(\F)$. However, since $X$ consists of 
elementary operators only (i.e., linear operators between finite-dimensional 
$\F$-vector spaces), it follows that we may identify the equivalently 
normed finite-dimensional Banach spaces $(Z, \Vert \cdot \Vert_Z)$ and 
$(X,\Vert \cdot \Vert_X)$ topologically (yet not isometrically\,!), implying 
that ${\mathcal{Q}}_{m,n}(\F) = B_X \subseteq X \subseteq Z$. Fix $B \in Y$. 
Then $A : = B^\ast \in W$, and $\Vert B \Vert_Y = \Vert A^\ast \Vert_Y = 
\Vert A \Vert_W = \Vert A \Vert_{\infty, 1}$. 
\eqref{eq:m_n_GT_reformulated} implies that
\[
\vert{\text{tr}}(B S)\vert = \vert{\text{tr}}(A^\ast S)\vert 
\leq K_G^\F(m,n)\,\Vert A \Vert_{\infty, 1} = 
K_G^\F(m,n)\,\Vert B \Vert_Y \text{ for all } (B, S) \in 
Y \times {\mathcal{Q}}_{m,n}(\F)\,. 
\]
By making use of polarisation with respect to the dual pairing $(Y, Z)$ 
(cf., e.g., \cite[Chapter 8.2]{J1981}), the latter inequality is equivalent to
\[
{\mathcal{Q}}_{m,n}(\F) 
\subseteq K_G^\F(m,n)\,B_Y^\circ = K_G^\F(m,n)\,
B_Z.
\]
So, we get again a well-known norm inequality variant of the Grothendieck 
inequality (cf. also \cite[Corollary 14.3]{DF1993} and 
\cite[Section 1.2]{GLPV2017}); namely:
\begin{align}\label{eq:geometric_interpretation_of_GT}
{B_{{\mathfrak{L}}_2(l_1^n, l_\infty^m)} 
\stackrel{\eqref{eq:char_of_quantum_corr_matrices}}{=}}{\mathcal{Q}}_{m,n}(\F) 
\subseteq K_G^\F(m,n)\,B_{{{\mathfrak{N}}(l_1^n, l_\infty^m)}} 
\subseteq K_G^\F\,B_{{{\mathfrak{N}}(l_1^n, l_\infty^m)}}
\end{align}
and $\Vert S \Vert_{{{\mathfrak{N}}(l_1^n, l_\infty^m)}} \leq K_G^\F(m,n) 
\leq K_G^\F$ for all $m, n \in \N$ and $S \in {\mathcal{Q}}_{m,n}(\F)$. 

Thus, we recognise that the isometric isomorphism $\tau : Z \longrightarrow 
Y^\prime$ can be \textit{extended} to the well-defined bounded linear operator 
$\widetilde{\tau} : X \longrightarrow Y^\prime$, where the latter is given as 
$\widetilde{\tau}(R) : = K_G^\F(m,n)\,\tau\big(\frac{1}{K_G^\F(m,n)}\,R\big)$ 
for all $R \in X$; i.e.,
\begin{align}\label{eq:linear_extension_from_Z_to_X}
Y \times X \ni (B, R) \mapsto \langle B, \widetilde{\tau}(R)\rangle : = 
K_G^\F(m,n)\,\big\langle B, \tau\big(\frac{1}{K_G^\F(m,n)}\,R\big)\big\rangle = 
\text{tr}(B R)\,. 
\end{align}
Note that $\Vert\widetilde{\tau}(R)\Vert = K_G^\F(m,n)\,\Vert \frac{1}{K_G^\F(m,n)}\,
R \Vert_Z \leq K_G^\F(m,n)\,\Vert R \Vert_X \leq K_G^\F\,\Vert R \Vert_X$ for 
all $R \in X$ (due to \eqref{eq:geometric_interpretation_of_GT}), whence 
$\Vert \widetilde{\tau} \Vert \leq K_G^\F(m,n) \leq K_G^\F$.

Consequently, since $K_G^\F(m,n)$ is the smallest 
constant which satisfies inequality \eqref{eq:GT_for_given_m_and_n} (or equivalently 
\eqref{eq:m_n_GT_reformulated}), it even follows that 
\begin{align}\label{eq:rep_of_K_G_F_m_n}
K_G^\F(m,n) = \sup\limits_{S \in {\mathcal{Q}}_{m,n}(\F)}\Vert S \Vert_Z
\end{align} 
(since $\vert\text{tr}(A^\ast S)\vert 
\stackrel{\eqref{eq:linear_extension_from_Z_to_X}}{=}
\vert\langle A^\ast, \widetilde{\tau}(S)\rangle\vert \leq
K_G^\F(m,n)\,\Vert \tau({\frac{1}{K_G^\F(m,n)}\,S})\Vert\,
\Vert A^\ast\Vert_Y \stackrel{\eqref{eq:trace_duality_fin_dim_case}}{\leq}
K_G^\F(m,n)\,\Vert {\frac{1}{K_G^\F(m,n)}\,S}\Vert_Z$ for all $S \in 
{\mathcal{Q}}_{m,n}(\F)$ and $A \in B_W$).
\begin{remark}(\textbf{Attainability of maximum in GT})
\label{rem:attainability_of_max}
For any $A \in W$ the linear functional $f_{A} : X \longrightarrow \R, 
R \mapsto \langle A^\ast, \widetilde{\tau}(R)\rangle = \text{tr}(A^\ast R)$ 
satisfies $\vert f_{A}(R) \vert \leq \Vert \widetilde{\tau}(R)\Vert\,
\Vert A \Vert_W \leq K_G^\F \Vert R \Vert_X\,\Vert A \Vert_{\infty, 1}$ 
for all $R \in X$. Hence, $f_{A} : X \longrightarrow \R$ is continuous and 
attains it maximum on the \textit{compact} unit ball $B_X = {\mathcal{Q}}_{m,n}(\F)$ 
(since $X$ is finite-dimensional). Thus, we may indeed replace the supremum by 
the maximum in \autoref{thm:GT_rewritten_form}. Similarly, the maximum is 
attained in \eqref{eq:rep_of_K_G_F_m_n}, whence
\[
K_G^\F(m, n) = \Vert S_0 \Vert_Z > 1
\]
for some $S_0 \in B_X\setminus B_Z = {\mathcal{Q}}_{m,n}(\F)$. Recall that 
$\tau : Z \stackrel{\cong}{\longrightarrow} Y^\prime$ is an isometric 
isomorphism (cf. \eqref{eq:trace_duality_fin_dim_case}) and observe that
\[
\langle S, (\tau^\prime j_Y)(B) \rangle = \langle B, \tau(S) \rangle =
{\text{tr}}(BS) \text{ for all } (S, B) \in Z \times Y\,.  
\]
Thus, $\tau^\prime j_Y : Y \stackrel{\cong}{\longrightarrow} Z^\prime$
again is an isometric isomorphism (since $Y$ is finite-dimensional). It 
therefore follows the existence of some $B_0 \in S_Y$, such that 
\[
K_G^\F(m, n) = \Vert S_0 \Vert_Z = \vert{\text{tr}}(A_0^\ast S_0)\vert\,,
\]
where $A_0 : = B_0^\ast \in S_W$. Moreover, the polar of $B_Z$ satisfies 
\begin{align}\label{eq:polar_reprs}
B_Z^\circ = B_{Z^\prime} = \tau^\prime j_Y(B_Y).
\end{align}
Consequently, for \textit{any} matrix $A \in \M_{m,n}(\F)$, the following 
equivalence holds: $\vert \text{tr}(A^\ast R)\vert \leq 1$ for all $R \in B_Z$ 
\textit{if and only if} $\Vert A \Vert_W = \Vert A \Vert_{\infty, 1} \leq 1$.
In summary, we have:
\begin{align}\label{eq:quantum_viol_of_a_Bell_inequality}
K_G^\F(m, n) = \vert{\text{tr}}(A_0^\ast S_0)\vert > 1\,\text{ and }\,
\vert \text{tr}(A_0^\ast R)\vert \leq 1\,\text{ for all }\, R \in B_Z\,;
\end{align}
a fact which plays a key role in quantum mechanics (cf. 
\sref{Remark}{rem:quantum_violation_of_a_Bell_inequality}). 
Moreover, due to \eqref{eq:char_of_quantum_corr_matrices}, $S_0 = 
\Gamma_{\F_2^{d_0}}(u_0, v_0)$, for some $d_0 \equiv d_0(m,n) \in \N$ and 
$(u_0,v_0) \in (S_{\F_2^{d_0}})^m \times (S_{\F_2^{d_0}})^n$. Hence, 
$K_G^\F(m, n) \leq \sup\limits_{(u,v)\in (S_{\F_2^{d_0}})^m \times 
(S_{\F_2^{d_0}})^n}\vert\text{tr}(A_0^\ast \Gamma_{\F_2^{d_0}}(u, v))\vert \leq 
K_G^\F(d_0)$.
\end{remark}
\begin{remark}[\textbf{Adjoining GT}]\label{rem:GT_norm}
Readers who are familiar with adjoint normed operator ideals and trace duality 
in general (cf. \cite[Chapter 9.1]{P1980}) immediately recognise that 
\sref{Corollary}{cor:GT_equivalent_statements}-\ref{(ii)} implies that the 
Grothendieck inequality actually is equivalent to an inequality between two 
matrix norms, induced by two (adjoint) Banach ideals; namely:
\[
\sup\limits_{H \in \texttt{HIL}^\F}\Vert A \Vert_H^{\text{G}} 
\stackrel{\eqref{eq:m_n_GT_reformulated}}{=} \Vert A \Vert_{{\mathfrak{D}}_2} 
\leq K_G^\F(m,n)\, \Vert A \Vert_{\infty, 1} \leq K_G^\F\,\Vert A \Vert_{\infty, 1}
\]
for all matrices $A \in \M_{m,n}(\F), m, n \in \N$, or, equivalently,
\begin{align}\label{eq:unit_ball_version_of_GT}
B_{{\mathfrak{L}}(l_\infty^n, l_1^m)} \subseteq 
K_G^\F(m,n)\, B_{{\mathfrak{D}}_2(l_\infty^n, l_1^m)} \subseteq
K_G^\F\, B_{{\mathfrak{D}}_2(l_\infty^n, l_1^m)}
\end{align}
for all $m, n \in \N$, where $({\mathfrak{D}}_2, \Vert \cdot \Vert_
{{\mathfrak{D}}_2}) = ({\mathfrak{L}}_2^\ast, \Vert \cdot \Vert_
{{\mathfrak{L}}_2^\ast}) = ({\mathfrak{P}}_2^d \circ {\mathfrak{P}}_2^{}, 
\Vert \cdot \Vert_{{\mathfrak{P}}_2^d \circ {\mathfrak{P}}_2^{}})$ characterises the 
Banach ideal of 2-dominated operators (cf. \cite[Table 17.12, Theorem 17.14 
and Chapter 19]{DF1993} and \sref{Remark}{rem:GT_vs_abs_summing_operators}). 
The latter inclusion also follows directly from ``adjoining'' 
\eqref{eq:geometric_interpretation_of_GT} above. 
\noindent In fact, if $K$ and $L$ are arbitrary compact sets, the following 
deep result of Grothendieck holds:
\[
{\mathfrak{L}(C(K), C(L)^\prime)} \subseteq {\mathfrak{D}_2(C(K), C(L)^\prime)} 
\subseteq {\mathfrak{P}_2(C(K), C(L)^\prime)} \subseteq 
{\mathfrak{L}_2(C(K), C(L)^\prime)}\,,
\]
and
\[
\Vert T \Vert_{{\mathfrak{L}_2}} \leq \Vert T \Vert_{{\mathfrak{P}_2}} 
\leq \Vert T \Vert_{{\mathfrak{D}_2}} \leq 
K_G^\F\,\Vert T \Vert
\]
for all $T \in \mathfrak{L}(C(K), C(L)^\prime)$. To recognise this highly 
noteworthy statement, we just have to note that \cite[Theorem 2.1]{P2012} 
implies that any $T \in \mathfrak{L}(C(K), C(L)^\prime)$ can be represented as 
$T = (J_{\P_L})^\prime U J_{\P_K}$, where for $\Delta \in \{K,L\}$, $\P_\Delta$ 
is a well-defined \textit{probability} measure on $\Delta$, $J_{\P_{\Delta}} : 
C(\Delta) \hookrightarrow L^2(\P_\Delta)$ denotes the canonical (norm 1) 
inclusion and $\Vert U \Vert \leq K_G^\F\,\Vert T \Vert$. Each of the two 
operators $J_{\P_{\Delta}}$ is absolutely $2$-summing (such as their biduals - 
cf. \cite[Corollary 17.8.4]{DF1993}) and satisfies 
$\Vert J_{\P_{\Delta}}\Vert_{\mathfrak{P}_2} = \P_\Delta(\Delta)^{1/2} = 1$ 
(cf. \cite[Subsection 11.2]{DF1993})). In particular, we reobtain 
\cite[Corollary 2.2]{P2012}. It is quite instructive to compare this result 
with \cite[Corollary 14.5.2 and Theorem 17.14]{DF1993}, \cite[Section 5]{H1994} 
and \cite[Theorem G]{K1976}.
\end{remark}
{
\begin{remark}(\textbf{Tensor norm representation of GT})
\label{rem:tensor_norm_rep_of_GT}
Let $m, n \in \N$. Readers who are familiar with both, Banach operator ideals 
and tensor norms, very likely re-recognise the following \textit{equivalent} 
tensor norm representations of the inequality 
\eqref{eq:geometric_interpretation_of_GT} (respectively 
\eqref{eq:unit_ball_version_of_GT}) at once (cf. 
\cite[Theorem 14.4 and Theorem 20.17]{DF1993}):
\begin{enumerate}
\item $\pi(\cdot; l_\infty^n, l_\infty^m) \leq K_G^\F(m,n)\, 
w_2(\cdot; l_\infty^n, l_\infty^m)\leq K_G^\F\, 
w_2(\cdot; l_\infty^n, l_\infty^m)$.
\item $w_2^\ast(\cdot; l_1^n, l_1^m) \leq K_G^\F(m,n)\, 
\varepsilon(\cdot; l_1^n, l_1^m)\leq K_G^\F\,\varepsilon(\cdot; l_1^n, l_1^m)$.
\end{enumerate}
We just have to apply the representation theorem for minimal Banach operator 
ideals to the minimal kernels of the maximal Banach operator ideals ${\mathfrak{I}} 
\sim \pi$ and ${\mathfrak{L}}_2 \sim w_2$ (cf. \cite[Corollary 22.2.1]{DF1993}). 
(i) then follows from the isometric equalities 
\[
l_\infty^n \otimes_\pi l_\infty^m \cong {\mathfrak{I}}^{\mbox\small{min}}
(l_1^n, l_\infty^m) \stackrel{1}{=} {\mathfrak{N}}(l_1^n, l_\infty^m)  
\]
and 
\[
l_\infty^n \otimes_{w_2} l_\infty^m \cong {\mathfrak{L}}_2^{\mbox\small{min}}
(l_1^n, l_\infty^m) \stackrel{1}{=} {\mathfrak{L}_2}(l_1^n, l_\infty^m)\,.
\]
Since $(M \otimes_\alpha N)^\prime \cong M^\prime \otimes_{\alpha^\prime} 
N^\prime$ for all finitely generated tensor norms $\alpha$ and all 
finite-dimensional Banach spaces $M$ and $N$ (prove it!), it follows from trace 
duality that (i) and (ii) in fact are equivalent. Here, it should be noted that 
quite often the tensor norm $w_2$ is also denoted as $\gamma_2$, and the maximal 
Banach ideal ${\mathfrak{L}_2}$ is also known as $\Gamma_2$.
\end{remark}
}
\noindent Given that trace duality view, the role of the ``free parameters'' 
$m, n, d \in \N$, where the pair $(m,n) \in \N^2$ describes the size of 
the matrices and $d$ is the dimension of the underlying
finite-dimensional Hibert space $\F_2^d$ is explicitly described in        
\begin{proposition}\label{prop:K_G_R_resp_K_G_C_as_sup_of_nuclear_norms}
Let $\F \in \{\R, \C\}$ and $m, n, d \in \N$. Put
\[
K_G^\F(m,n; d) : = 
\sup\big\{\Vert \Gamma_{\F_2^d}(u, v) : l_1^n \longrightarrow l_\infty^m 
\Vert_{{\mathfrak{N}}}: (u, v) \in {(S_{\F_2^d})}^m \times 
{(S_{\F_2^d})}^n\big\}.
\]
Then
\begin{align}\label{eq:K_G_F_d}
K_G^\F(d) = \sup\limits_{(m,n) \in \N^2}K_G^\F(m,n; d).
\end{align}
and
\begin{align}\label{eq:K_G_F_m_n}
K_G^\F(m,n) = \sup\limits_{S \in {\mathcal{Q}}_{m,n}(\F)}
\Vert S \Vert_{{\mathfrak{N}}} = \sup\limits_{d \in \N}K_G^\F(m,n; d).
\end{align}
The sequence $(K_G^\F(d))_{d \in \N}$ is non-decreasing, and
\begin{align}\label{eq:K_G_F_rep}
K_G^\F = \sup\big\{\Vert S \Vert_{{\mathfrak{N}}}: m, n \in \N, S \in 
{\mathfrak{Q}}_{m,n}(\F)\big\} = \sup\limits_{(m,n) \in \N^2}K_G^\F(m,n) = 
\sup\limits_{d \in \N}K_G^\F(d) = \lim\limits_{d \to \infty}K_G^\F(d).
\end{align}
In particular, $K_G^\F(1,1) = 1$ and $K_G^\F(m,n;1) = K_G^\F(1) = 1$ for all 
$m,n \in \N$.
\end{proposition}
\begin{remark}\label{rem:GT_vs_trace_class}
Even in the matrix case, the Banach space $({\mathfrak{N}}(l_1^n, l_\infty^m), 
\Vert \cdot \Vert_{\mathfrak{N}})$ should not be confused with the Banach space 
$({\mathfrak{N}}(\F_2^n, \F_2^m), \Vert \cdot \Vert_{\mathfrak{N}}) = 
({\mathfrak{S}}_1(\F_2^n, \F_2^m), \sigma_1)$\,! The latter 
space namely consists of matrices - viewed as operators - 
between (finite-dimensional) \textit{Hilbert spaces}, contained in the so-called 
\textit{Schatten-von Neumann class of index 1}, also known as 
\textit{trace-class operators} (cf., e.g., \cite{HJ1991} and 
\cite[Chapter 20.2]{J1981}). In particular, if we view a given matrix 
$M \in \M_{m,n}(\F)$ as linear operator \textit{from $l_1^n$ to $l_\infty^m$}, 
the norm $\Vert M \Vert_{\mathfrak{N}}$ in general does \textit{not} coincide 
with the so-called \textit{trace norm} of $M$ (also known as 
\textit{nuclear norm}). The latter is given by $\Vert M \Vert_\ast : = 
\text{tr}(\vert M \vert)$, where $\vert M \vert : = (M^\ast\,M)^{1/2}$. Since 
the trace norm of $M$ coincides with the Schatten 1-norm $\sigma_1(M)$, it 
equals the sum of the singular values of the matrix $M$ (cf., e.g., 
\cite[Chapter IV.2]{Bh1997}, \cite[Exercises IX.2.19, IX.2.20 and IX.2.21]{C1990} 
and \cite[Chapter 5.6 and Chapter 7.4.7]{HJ2013}).
\end{remark}
\noindent Let us recall the isometric isomorphism $\tau : Z \stackrel{\cong}
{\longrightarrow} Y^\prime$, induced by the Banach spaces $Y : = 
{\mathfrak{L}}(l_\infty^m, l_1^n)$ and $Z : = {\mathfrak{N}}(l_1^n, l_\infty^m) 
\cong Y^\prime$ (cf. \eqref{eq:trace_duality_fin_dim_case}). Put $\delta : = 
\tau^\prime j_Y : Y \stackrel{\cong}{\longrightarrow} Z^\prime$. As we have seen, 
also $\delta$ is an isometric isomorphism (cf. 
\sref{Remark}{rem:attainability_of_max}). Its inverse is given by $\delta^{-1} = 
j_Y^{-1}(\tau^{-1})^\prime$ (since $Y$ is finite-dimensional). An application of 
the bipolar theorem to the dual pairing $\langle Z, Z^\prime\rangle = \langle Z, 
\delta(Y)\rangle$, induced by the bilinear form $Z \times Z^\prime \longrightarrow 
\F, (R, z^\prime) \mapsto \langle R, z^\prime \rangle : = {\text{tr}}
(R\,\delta^{-1}(z^\prime))$ (cf. \cite[V.1.8]{C1990} and \cite[Chapter 8.2]{J1981}) 
implies the following explicit representation result for the so-called 
``local correlation polytope''. To the best of our knowledge, the outcome for 
the complex case (i.e., if $\F = \C$) is new. A (different) part of our proof 
for the real case can be found in the proof of \cite[Proposition 11.7]{AS2017}. 
In particular, we are going to shed some light on the geometry of the unit ball of 
${\mathfrak{N}}(l_1^n, l_\infty^m)$. To this end, recall that for any subset $S$ 
of an $\F$-vector space 
\[
\text{cx}(S) : = \bigcup_{n \in \N}\big\{\sum_{i = 1}^n \alpha_i x_i : x_1, \ldots, 
x_n \in S, \alpha_1, \ldots, \alpha_n \geq 0 \text{ and } \sum_{i=1}^n 
\alpha_i = 1\big\}
\]
denotes the convex hull of $S$ and that
\begin{align}\label{eq:abs_cx_hull}
\text{acx}(S) : = \bigcup_{n \in \N}\big\{\sum_{i = 1}^n \alpha_i x_i : x_1, \ldots, 
x_n \in S, \alpha_1, \ldots, \alpha_n \in \F \text{ and } \sum_{i=1}^n 
\vert\alpha_i\vert \leq 1\big\}
\end{align}
marks the absolute convex hull of $S$. Here we adopt the notation, introduced 
right below \cite[Proposition 6.1.3]{J1981}. Recall also that $\text{acx}(S) 
= {\text{cx}}(\widebreve{S})$, where $\widebreve{S} : = (\F \cap \overline{\D})
S$ denotes the circled hull of $S$ (cf. \cite[Proposition 6.1.4]{J1981}). 
\begin{theorem}\label{thm:nuclear_unit_ball}
Let $\F \in \{\R, \C\}$ and $m,n \in \N$. Put 
\[
{\mathcal{G}}_{m,n}(\F) := 
\{\overline{p}q^\top : (p, q) \in S_\F^m \times S_\F^n\} = 
\big\{\Gamma_\F(p,q) : (p,q) \in S_\F^m \times S_\F^n \big\}
\]
and 
\[
{\mathcal{H}}_{m,n}(\F) := 
\{\overline{x}y^\top : (x, y) \in (\F \cap \overline{\D})^m \times 
(\F \cap \overline{\D})^n\} = \big\{\Gamma_\F(x,y) : (x,y) \in 
(\F\cap \overline{\D})^m \times (\F\cap \overline{\D})^n \big\}.
\]
Then
\begin{align}\label{eq:the_nuclear_infty_1_ball}
B_{{\mathfrak{N}}(l_1^n, l_\infty^m)} = {\text{acx}}({\mathcal{G}}_{m,n}(\F)) =
{\text{cx}}(\widebreve{{\mathcal{G}}_{m,n}(\F)}) = {\text{cx}}
({\mathcal{H}}_{m,n}(\F)).
\end{align}
\begin{enumerate}
\item If $\F = \R$, then the set ${\mathcal{G}}_{m,n}(\R)$ even 
coincides with the set of all extreme points of $B_{{\mathfrak{N}}(l_1^n, l_\infty^m)}$ 
and
\begin{align*}
B_{{\mathfrak{N}}(l_1^n, l_\infty^m)} &= {\text{cx}}({\mathcal{G}}_{m,n}(\F)) =
\big\{\sum_{i=1}^{mn+1}\lambda_i\,x_i : x_i \in {\mathcal{G}}_{m,n}(\R), 
0 \leq \lambda_i \leq 1, \sum_{i=1}^{mn+1}\lambda_i = 1\big\}\\
&= \big\{\E[\textbf{X}\textbf{Y}^\top] : \max\{\vert X_i\vert, \vert Y_j\vert\} 
\leq 1 \text{ a.s., }\text{for all } (i,j) \in [m] \times [n]\big\}\,.
\end{align*}
\item If $\F = \C$, then
\[
B_{{\mathfrak{N}}(l_1^n, l_\infty^m)} = 
\big\{\sum_{i=1}^{2mn+1}\lambda_i\,x_i : x_i \in \widebreve{{\mathcal{G}}_{m,n}(\C)}, 
0 \leq \lambda_i \leq 1, \sum_{i=1}^{2mn+1}\lambda_i = 1\big\}.
\]
\end{enumerate}
\end{theorem}
\noindent Based on \cite[Definition 11.5 and remark right below]{AS2017}, we 
have shown that at least in the real case (i.e., if $\F = \R$), $B_Z$ precisely 
coincides with the set of all \textit{classical (or local) ``correlation'' 
matrices}.
 
Consequently, \eqref{eq:geometric_interpretation_of_GT} implies 
that for all $m, n \in \N$, for all $\F$-Hilbert spaces 
$H^\F$ and $(u, v) \in S_{H^\F}^m \times S_{H^\F}^n$ there are $x_1^\F, \ldots, 
x_{k^\F}^\F \in (\F\cap\overline{\D})^m$, $y_1^\F, \ldots, y_{k^\F}^\F \in 
(\F\cap\overline{\D})^n$ and $(\lambda_1^\F, \ldots, \lambda_{k^\F}^\F) \in 
[0,1]^{k^\F}$, such that $\sum_{\nu=1}^{k^\F} \lambda_{\nu}^\F = 1$ 
and
\begin{align}\label{eq:QM_is_GT_constant_times_acx_hull_of_rank_1}
\Gamma_{H^\F}(u,v) =
K_G^\F\,\sum_{\nu=1}^{k^\F}\lambda_\nu^\F\,\overline{x_\nu^\F}\,(y_\nu^\F)^\top,
\end{align}
where $k_\R : = m n+1$ and $k_\C : = 2m n+1$. If $\F = \R$, we may assume that 
$\vert x_\nu^\R\vert = 1$ and $\vert y_\nu^\R\vert = 1$ for all $\nu \in [k^\R]$. 
In particular (if $m = n = 1$),
\[
\langle u, v\rangle_{H^\R} = 
K_G^\R\, (\lambda_1^\R\,x_1^\R\,y_1^\R + (1-\lambda_1^\R)\,x_2^\R\,y_2^\R) 
\] 
and
\[
\langle b, a\rangle_{H^\C} = K_G^\C\,\sum_{\nu=1}^{3} \lambda_\nu^\C\,
\overline{x_\nu^\C}\,y_\nu^\C
\]
for all $u, v \in S_{H^\R}$ and $a, b \in S_{H^\C}$. Since $K_G^\R > 1$, it 
follows that $\text{sign}(x_1^\R\,y_1^\R) \not= \text{sign}(x_2^\R\,y_2^\R)$. 
Thus, if $\lambda_1 \not= \frac{1}{2}$, then
\[
\frac{\pi}{2} < K_G^\R \leq \frac{1}{\vert 2\,\lambda_1 - 1 \vert}, 
\text{ respectively } \vert \lambda_1 \vert \leq \frac{1+K_G^\R}{2K_G^\R} < 
\frac{\pi+2}{2\pi} \approx 0.818\,.
\]
\begin{remark}[\textbf{Quantum violation of a Bell inequality}]
\label{rem:quantum_violation_of_a_Bell_inequality}
Fix an \textit{arbitrarily given} $S \in \mathcal{Q}_{m,n}(\F)\setminus B_Z$. 
Recall again the construction of the isometric isomorphism 
$\tau : Z \stackrel{\cong}{\longrightarrow} Y^\prime$ via trace duality (cf. 
\eqref{eq:spaces} and \eqref{eq:trace_duality_fin_dim_case}). Since $B_Z$ 
is a non-empty closed and absolutely convex subset of the Banach space $Z$, we 
may apply \cite[Corollary 7.3.6]{J1981}. The latter is an implication of hyperplane 
separation (which is a geometric version of the Hahn-Banach theorem and thus 
an implication of Zorn's lemma). Hence, due to \eqref{eq:polar_reprs}, it 
follows the existence of a matrix $B_0 \in B_Y$, such that
\[
\text{tr}(A_0^\ast S) = \text{tr}(B_0 S) > 1 \,\text{ and }\, 
\vert \text{tr}(A_0^\ast R)\vert \leq 1 \text{ for all } R \in B_Z\,,
\]
where $A_0 : = B_0^\ast \in B_W$ (compare also with 
\eqref{eq:quantum_viol_of_a_Bell_inequality}). In particular, we have been 
provided with a matrix $A_0 \in \M_{m,n}(\F)$, such that 
$\Vert A_0 \Vert_{\infty, 1} \leq 1$ and
\[
1 < \sup\limits_{H \in \texttt{HIL}^\F}\Vert A_0 \Vert_H^{\text{G}} = \max\limits_
{S \in {\mathcal{Q}}_{m,n}(\F)}\vert\text{tr}(A_0^\ast S)\vert = \max\limits_
{S \in {\mathcal{Q}}_{m,n}(\F)}\Re(\text{tr}(A_0^\ast S)) = 
\max\limits_{\Sigma \in C(m+n; \F)}\text{tr}(\Delta(A_0)\Sigma)
\]
(due to \eqref{eq:m_n_GT_reformulated}). In quantum mechanics (if $\F = \R$), 
the latter inequality is somewhat vaguely referred to as the ``maximal quantum 
violation of a Bell (correlation) inequality''. 
The ``maximal violation of the related Bell inequality'' coincides precisely 
with the inequality
\[
\Vert f_{A_0} \Vert = \max\limits_{S \in {\mathcal{Q}}_{m,n}(\F)}
\vert f_{A_0}(S)\vert = \max\limits_{S \in {\mathcal{Q}}_{m,n}(\F)}
\vert\text{tr}(A_0^\ast S)\vert > 1\,,
\] 
where $f_{A_0}$ is defined as in \sref{Remark}{rem:attainability_of_max} 
(for both fields). 
The term ``Bell inequality'' (now with respect to any given $A \in B_W$, of 
course) is therefore to be understood as the inequality
\[
\max\limits_{R \in B_Z}
\vert f_{A}(R)\vert = \max\limits_{R \in B_Z}
\vert\text{tr}(A^\ast R)\vert\leq 1\,,
\]
which holds for all $A \in B_W$ (due to \eqref{eq:polar_reprs}). That 
Bell inequality is ``violated'', if and only if $\vert f_{A}(S_0)\vert > 1$ 
for some quantum correlation matrix $S_0 \in \mathcal{Q}_{m,n}(\F)\setminus 
B_Z$. (cf. also \sref{Remark}{rem:attainability_of_max}, respectively 
\autoref{thm:nuclear_unit_ball}-\ref{(i)}, together with 
\cite[Lemma 2]{L2020} if $\F = \R$). In this respect, the maximum value 
$\Vert f_{A} \Vert$ is then referred to as ``maximal violation of the 
Bell inequality''. Consequently, due to \eqref{eq:char_of_quantum_corr_matrices} 
and \sref{Remark}{rem:attainability_of_max}, there are $S_\ast \in 
\mathcal{Q}_{m,n}(\F)\setminus B_Z$ and - hence - $d_\ast \in \N$, such that
``the maximal violation of the Bell inequality'' (i.e., $\Vert f_{A} \Vert$) 
is uniformly bounded above by $K_G^\F(d_\ast) \leq K_G^\F$. More precisely, if 
$S_\ast = \Gamma_{\F_2^{d_\ast}}(u_\ast, v_\ast)$, for some $d_\ast \in \N$ and 
$(u_\ast,v_\ast) \in (S_{\F_2^{d_\ast}})^m \times (S_{\F_2^{d_\ast}})^n$, then
\[
1 < \Vert f_{A} \Vert = \vert\text{tr}(A^\ast S_\ast)\vert =
\sup\limits_{(u,v)\in (S_{\F_2^{d_\ast}})^m \times (S_{\F_2^{d_\ast}})^n}
\vert\text{tr}(A^\ast \Gamma_{\F_2^{d_\ast}}(u, v))\vert \leq 
K_G^\F(d_\ast) \leq K_G^\F\,.
\]
\end{remark}
\begin{remark}\label{rem:on_Delta_A}
We don't know whether in \eqref{eq:symmetric_block_matrix_version_of_GT} 
we may substitute the block matrix $\Delta(A) \in \M_{m+n}(\F)$ through an 
\textit{arbitrary} matrix $B \in \M_{m+n}(\F)$. If this were the case, a 
further application of the bipolar theorem (cf. \cite[Theorem 8.2.2]{J1981}) 
shows that the latter would be equivalent to
\[
C(k; \F) \subseteq K_G^\F \, {\text{acx}}(C_1(k; \F)) \text{ for all } k \in \N_2\,. 
\]
\autoref{thm:nuclear_unit_ball} therefore would imply that
\[
C(k; \F) \subseteq K_G^\F \, {}{{\text{acx}}}(\{\overline{p}q^\top : (p, q) \in 
S_\F^k \times S_\F^k\}) \stackrel{\eqref{eq:the_nuclear_infty_1_ball}}
{=} K_G^\F \, B_{{}{\mathfrak{N}(l_1^k, l_\infty^k)}} \text{ for all } k \in \N_2\,,
\]
where $B_{{}{\mathfrak{N}(l_1^k, l_\infty^k)}}$ again denotes the unit ball of 
the Banach space of nuclear operators between $l_1^k$ and $l_\infty^k$, 
equipped with the nuclear norm. 
\end{remark}
\section{$K_G^\R(2)$ and the Walsh-Hadamard transform: Krivine's 
approach revisited}\label{sec:Hadamard_gate_and_Krivine}
 Regarding explicit constructions of elements of
${\mathcal{Q}}_{m,n}(\F)$, a rigorous description of the \textit{entries} of 
the Kronecker product of matrices proves to be a very useful tool (cf. 
\sref{Example}{ex:Hadamard_gate}). To this end, we consider the mapping: 
\[
\Z \times \N \ni (\nu, n) \mapsto r_n(\nu) : = 
\begin{cases} n &\text{if } n \text{ is a divisor of } \nu\\ 
\text{rem}_n(\nu) &\text{if } n \text{ is not a divisor of } \nu\,,
\end{cases}
\]
where $\text{rem}_n(\nu) \in \{0,1, \ldots, n-1\}$ denotes the uniquely 
determined remainder in Euclidean division of $\nu$ by $n$, implying that 
$r_n(\nu) \in [n]$ (by construction). Thus, if $p \in \Z$ and $\nu = p n + 
\text{rem}_n(\nu)$, then $p+1 \geq \frac{\nu+1}{n}$ and
\[
f_n(\nu) : = \frac{\nu - r_n(\nu)}{n}+1 = 
\begin{cases}
p+1 < \frac{\nu}{n}+1& \text{ if } n \text{ is not a divisor of } \nu\\
\frac{\nu}{n}& \text{ if } n \text{ is a divisor of } \nu\,.
\end{cases}
\]
Consequently, if $l \in \N$ and $\nu \in [ln]$, then $f_n(\nu) \in [l]$. 
In particular, $r_n(\nu) = \nu$ if $\nu \in [n]$.

 Especially with regard to \sref{Example}{ex:Hadamard_gate} the ``Boolean'' 
case n=2 is of particular importance to us. Here, we obviously obtain:   
\begin{align}\label{eq:Boolean_case}
f_2(\nu) = \left\lceil{\frac{\nu}{2}}\right\rceil =
\begin{cases}
\frac{\nu}{2}& \text{ if } \nu \text{ is even}\\
\frac{\nu+1}{2}& \text{ if } \nu \text{ is odd}
\end{cases}
\,\text{ and }\,b_1(\nu) : = r_2(\nu)- 1 = \ind_{2\N}(\nu) = 
\begin{cases}
1& \text{ if } \nu \text{ is even}\\
0& \text{ if } \nu \text{ is odd}\,.
\end{cases}
\end{align}
Note that the structure of $f_2$ implies that $f_2([2^i]) = [2^{i-1}]$ for all
$i \in \N$. In particular, for any $m \in \N_2$ and $i \in \{2,3, \ldots, m\}$, 
the well-defined function $b_i : = b_1 \circ \underbrace{f_2 \circ \cdots \circ f_2}_
{(i-1)-\text{times}}$ is the $i$th component of the $\{0,1\}^m$-valued function 
$\pi_m : \Z \longrightarrow \{0,1\}^m$, defined as
\begin{align}\label{eq:Boolean_vector_function}
\pi_m(\nu) : = (b_1(\nu), b_2(\nu), \ldots,b_m(\nu))^\top \text{ for all } \nu \in \Z\,.
\end{align}
 The actual role of the sequence of functions $(r_n)_{n \in \N}$ is 
encoded in 
\begin{lemma}\label{lem:remainder_bijection}
Let $n \in \N$. Then the mapping
\[
\Psi_n : \Z \times [n] \stackrel{\cong}{\longrightarrow} \Z\,,
(i,j) \mapsto (i-1)n + j
\]
is bijective. Its inverse is given by $\Psi_n^{-1} = \Lambda_n$, where
\[
\Lambda_n : \Z \longrightarrow \Z \times [n], \nu \mapsto 
(f_n(\nu),  r_n(\nu)).
\]
Moreover, $\Psi_n([l] \times [n]) = [ln]$ for all $l \in \N$. 
\end{lemma}
\begin{comment}
Fix $(i,j) \in \Z \times [n]$. If $\Psi(i,j) = (i-1)n+j = ln$ for some $l \in \Z$, 
it follows that $0 < j = (l-(i-1))n \leq n$, whence $0 < l-(i-1) \leq 1$. Thus, 
$l=i$, and hence $j = n$. Thus, $n$ is a divisor of $(i-1)n+j$ if and only if 
$j = n$, implying that $r_n((i-1)n+j) = n$ (by construction of the mapping $r$). 
Hence, if $n$ is not a divisor of $(i-1)n+j$, then $j < n$, and it follows that
$r_n((i-1)n+j) = j$ is the uniquely determined remainder of $(i-1)n+j$.
Now, we only have to use a bit of elementary algebra, to show that 
$\Lambda_n \circ \Psi_n = {\text{id}}_{\Z \times [n]}$ and $\Psi_n \circ 
\Lambda_n = {\text{id}}_{\Z}$, where $\Lambda(\nu) : = 
\big(\frac{\nu - r_n(\nu)}{n}+1,  r_n(\nu)\big)$ for all $\nu \in 
\Z$. Consequently, $\Lambda_n = \Psi_n^{-1}$. Finally, if $\Psi_n(i,j) = (i-1)n + j 
\leq ln$, it follows that $n(l-(i-1)) \geq j > 0$, whence $i-1 < l$. Thus, 
$i \in [l]$, and it follows that $[ln] = [ln] \cap \Psi_n(\Z \times [n]) 
\subseteq \Psi_n([l] \times [n])$ (since $\Psi_n$ is onto). Conversely, if 
$i \in [l]$ and $j \in [n]$, then $\Psi_n(i,j) = (i-1)n + j \leq in \leq ln$.
\end{comment}
\noindent Equipped with
the remainder mapping $r$ and the both bijections $\Psi_n : \Z \times [n] 
\longrightarrow \Z$ and $\Psi_m : \Z \times [m] \longrightarrow \Z$, 
we are now able to describe both, the bijective linear operator $\text{vec} : 
\M_{m,n}(\F) \longrightarrow \F^{mn}$ and the Kronecker product explicitly 
\textit{entrywise}. So, fix $m, n, p, q \in \N$. If $A = (a_{ij}) \in \M_{m,n}(\F)$, 
$B = (b_{kl}) \in \M_{p,q}(\F)$, $\alpha \in [mp]$, $\beta \in [nq]$ and 
$\gamma \in [mn]$, put 
\begin{align}\label{eq:Kronecker_product}
(A \otimes B)_{\alpha, \beta} : = 
a_{f_p(\alpha), f_q(\beta)}\cdot b_{r_p(\alpha), r_q(\beta)} = 
((A^\top \otimes B^\top)^\top)_{\alpha, \beta} 
\end{align}
and
\begin{align}\label{eq:vec_op}
\text{vec}(A)_\gamma \equiv {\text{vec}_m}(A)_\gamma : = a_{r_m(\gamma), f_m(\gamma)} = 
(A^\top)_{\Psi_m^{-1}(\gamma)} = 
\big(A_{\tau\circ\Psi_m^{-1}}\big)_\gamma\,,
\end{align}
where $(i,j) \mapsto \tau(i,j) : = (j,i)$ denotes transposition. In other words,
$\text{vec} \equiv {\text{vec}}_m = C_{\tau\circ\Psi_m^{-1}}$. In particular,
\[
\text{vec}\big(e_i^{(m)} (e_j^{(n)})^\top\big) = e_{(j-1)m+i}^{(nm)} 
\text{ for all } (i,j) \in [m] \times [n].
\]
Consequently, \sref{Lemma}{lem:remainder_bijection} implies that
\begin{align}\label{eq:tensor_prod_rep_I}
e_{\nu}^{(nm)} = \text{vec}\big(e_{r_m(\nu)}^{(m)} 
e_{f_m(\nu)}^{(n)\top}\big) \text{ for all } \nu \in [nm]\,.
\end{align} 
Moreover, if $n=1$, it follows that $\text{vec}(A)_i = a_{i, 1}$ for all $i \in [m]$. 
Observe that the vectorisation of the matrix $A$ involves its transpose $A^\top = 
A_\tau \in \M_{n,m}(\F)$. Not too surprisingly, it follows that 
\[
\text{vec}_n \circ C_\tau(A) = \text{vec}_n(A^\top) = \text{vec}_n(A_\tau) = 
A_{\Psi_n^{-1}} = C_{\Psi_n^{-1}}(A)
\]
(since $\tau \circ \tau = \text{id}$). The construction implies at once the 
non-trivial fact that $\text{vec}_n(A^\top) = C_{\Psi_m \circ \tau \circ 
\Psi_n^{-1}}(\text{vec}_m(A))$, where $C_{\Psi_m \circ \tau \circ \Psi_n^{-1}} : 
\F^{mn} \longrightarrow \F^{nm} = \F^{mn}$ is the related linear composition 
operator. Consequently,
\[
\text{vec}_n(A^\top) = K_{m,n}\,\text{vec}_m(A)\,,
\]  
where the matrix $K_{n,m}^\top = K_{m,n}\in O(mn)$ satisfies 
\begin{align*}
(K_{m,n})_{\nu \mu} &\,\,\,= \delta_{(\Psi_m \circ \tau \circ \Psi_n^{-1})(\nu),\,\mu} = 
\delta_{(r_n(\nu)-1)m + f_n(\nu),\,\mu}\\
&\,\,\,\stackrel{(!)}{=} 
\delta_{f_m(\mu), r_n(\nu)}\cdot \delta_{f_n(\nu), r_m(\mu)}\\ 
&\stackrel{\eqref{eq:Kronecker_product}}{=} 
\big(\sum_{i=1}^m\sum_{j=1}^n e_i\,e_j^\top \otimes e_j\,e_i^\top\big)_{\nu \mu}
\end{align*}
for all $\nu, \mu \in [mn]$. The second equality follows from 
\sref{Lemma}{lem:remainder_bijection}: $(r_n(\nu)-1)m + f_n(\nu) = \mu$ if and 
only if $\Psi_m(r_n(\nu), f_n(\nu)) = \mu = \Psi_m(f_m(\mu), r_m(\mu))$. Since 
$f_n(\nu) \in [m]$, it follows that $r_n(\nu) = f_m(\mu)$ and $f_n(\nu) = 
r_m(\mu)$. Therefore, \sref{Lemma}{lem:remainder_bijection} (respectively 
Euclidean division with remainder) allows to extend the results in \cite[Chapter 11, 
including Exercise 11.8]{AM2005} by an explicit entrywise (and hence 
implementable) description of the commutation matrix $K_{m,n} = 
\sum_{i=1}^m\sum_{j=1}^n e_i\,e_j^\top \otimes e_j\,e_i^\top$; namely in form 
of a product of two Kronecker delta symbols.
 
Moreover, $\text{vec}^{-1} = \text{mat}$, where the ``matrixation 
operator'' $\text{mat} : \F^{mn} \longrightarrow \M_{m,n}(\F)$ is given by 
$\text{mat} \equiv \text{mat}_m : = C_{\Psi_m\circ \tau}$; i.e.,
\[
\text{mat}(x) : = x_{\Psi_m\circ \tau} = 
\big(x_{\Psi_m(\tau(i,j))}\big)_{(i,j)} = 
\big(x_{(j-1)m + i}\big)_{(i,j)} =
\begin{pmatrix}
x_1 & x_{m+1} & \ldots & x_{(n-1)m+1}\\
x_2 & x_{m+2} & \ldots & x_{(n-1)m+2}\\
\vdots & \vdots & \ldots & \vdots\\
x_m & x_{2m} & \ldots & x_{mn}
\end{pmatrix} 
\]
for all $x \in \F^{mn}$. Again, if $n=1$, we recognise that $\text{mat}(x)_{i1} = x_i$ 
for all $i \in [m]$. Consequently, we may identify $\text{mat}(\F^m) \equiv 
\text{vec}(\M_{m,1}(\F)) \equiv \F^m$ for all $m \in \N$, so that we may assume without 
loss of generality that $(m, n) \in \N_2 \times \N_2$. That assumption also 
avoids necessarily the review, whether $\text{vec}$ maps $\F^{mn}$ into 
$\M_{m,n}(\F)$, or into $\M_{mn, 1}(\F) \equiv \F^{mn}$, or into $\M_{1, mn}(\F) 
\equiv\{x^\top : x \in \F^{mn}\}$!
\begin{example}
Consider $A : =
\begin{pmatrix}
a_{11} & a_{12} & a_{13}\\
a_{21} & a_{22} & a_{23}
\end{pmatrix} \in \M_{2,3}(\F)$. Then $4 \in [2\cdot 3]= [6]$ and $\Psi_3^{-1}(4) = 
(\frac{4 - r_3(4)}{3} + 1, r_3(4)) = (2,1)$. Thus, $\text{vec}(A^\top)_4 = 
a_{21}$. 
\end{example}
\noindent Therefore, we reobtain the following well-known, easily computable 
statements
\begin{align*}
x y^\top 
\equiv x \otimes y^\top \text{ for all } (x,y) \in \F^m \times \F^n\,.
\end{align*}
and
\begin{align*}
\text{vec}(x y^\top) 
\equiv y \otimes x \text{ for all } (x,y) \in \F^m \times \F^n\,.
\end{align*}
In particular,
\begin{align}\label{eq:tensor_prod_rep_II}
\begin{split}
&e_{\nu}^{(nm)} \stackrel{\eqref{eq:tensor_prod_rep_I}}{=} 
e_{f_m(\nu)}^{(n)} \otimes e_{r_m(\nu)}^{(m)}\\
&\parallel \hspace{3cm} \parallel \\ 
&e_{(j-1)m+i}^{(nm)} \hspace{1cm} e_{j}^{(n)} \otimes e_{i}^{(m)}\\ 
\end{split}
\end{align}
for all $\nu \in [nm] = \Psi_m([n] \times [m]) = \{(j-1)m+i : (j,i) \in [n] \times 
[m]\}$. Equivalently (now translated into Dirac's bra-ket language):
\begin{align}\label{eq:bra_ket_version_of_qubit_tensor_prod}
\vert \alpha\,\beta \rangle \equiv \vert \alpha\rangle\,
\vert \beta\rangle = \vert \alpha m + \beta \rangle \text{ for all } 
(\alpha,\beta)\in \{0, 1, \ldots, n-1\} \times \{0, 1, \ldots, m-1\}.
\end{align}
Consequently, if $C = \sum\limits_{j=1}^n\sum\limits_{k=1}^p c_{jk}\,e_j e_k^\top 
\in \M_{n,p}(\F)$ is a third given matrix, then $ACB = \sum\limits_{j=1}^n
\sum\limits_{k=1}^p c_{jk} (Ae_j)(B^\top e_k)^\top$ and $\text{vec}(C) = 
\sum\limits_{j=1}^n\sum\limits_{k=1}^p c_{jk}\,e_k \otimes e_j$, leading to another 
important, well-known matrix equality:
\[
\text{vec}(ACB) = (B^\top \otimes A)\text{vec}(C).
\]
\begin{example}[\textbf{Werner state}]\label{ex:Werner_state}
Let $p \in [0,1]$. Put 
\[
\R^4 \ni \psi^{-} : = \frac{1}{\sqrt{2}}(0,1, -1, 0)^\top 
\stackrel{\eqref{eq:bra_ket_version_of_qubit_tensor_prod}}{=} 
\frac{1}{\sqrt{2}}(\vert 0\,1 \rangle - \vert 1\,0 \rangle). 
\]
Consider the matrix
\[
\M_4(\R) \ni \rho_p^\text{W} : = p\,\psi^{-}(\psi^{-})^\top + \frac{1-p}{4}\,I_4 =
\frac{1}{4}\begin{pmatrix}
1-p & 0 & 0 & 0\\
0 & 1+p & -2p & 0\\
0 & -2p & 1+p & 0\\
0 & 0 & 0 & 1-p\\
\end{pmatrix}.
\]
Since the rank one matrix $p\,\psi^{-}(\psi^{-})^\top$ is positive semidefinite 
and $\frac{1-p}{4}\,I_4$ is positive definite if $p < 1$, it immediately follows 
that $\rho_p^\text{W}$ is positive definite if $p < 1$ and $\rho_1^\text{W} = 
\psi^{-}(\psi^{-})^\top$ is positive semidefinite (yet not invertible). Moreover, 
by construction, it trivially follows that in general 
$\text{tr}(\rho_p^\text{W}) = 1$. Note that $\rho_p^\text{W}$ can also be written 
as
\[
\rho_p^\text{W} = \frac{3-2\lambda(p)}{6} I_4 - \frac{3-4\lambda(p)}{6} G_1 = 
\frac{3-2\lambda(p)}{6}\big(I_4 - \frac{3-4\lambda(p)}{3-2\lambda(p)}\,G_1\big)\,,
\]
where $\lambda(p) : = \frac{3}{4}(1-p) \in [0,1]$ and $G_1 \in O(4)$ is the 
``flip operator'' (cf. \eqref{eq:orthog_matrix}). In quantum physics, the matrix 
$\rho_p^\text{W} \in \M_4(\R)^+$ is known as the so-called ``two-qubit Werner 
state''. A very detailed discussion of the origin and the meaning of Werner states 
in the foundations and philosophy of quantum mechanics, particularly in 
relation to the topic of entanglement and ``local hidden-variable theories'' 
can be found, for example, in \cite{AS2017, B2021, FHLLZZ2015} and in the 
relevant references therein. 
\end{example}  
\begin{proposition}\label{prop:tensor_product_of_quant_corr_is_quant_corr}
Let $m, n, p, q \in \N$, $S \in \M_{m,n}(\F)$ and $R \in \M_{p,q}(\F)$.
If $S \in {\mathcal{Q}}_{m,n}(\F)$ and $R \in {\mathcal{Q}}_{p,q}(\F)$, then
$S \otimes R \in {\mathcal{Q}}_{mp,nq}(\F)$. 
\end{proposition}
\begin{comment}
We only have to apply \sref{Proposition}{prop:char_of_Q_m_n_for_both_fields} and
the entrywise description \eqref{eq:Kronecker_product} of the Kronecker product
$S \otimes R$. So, choose $\F$-Hilbert spaces $H_1, H_2, (u^{(1)},v^{(1)}) 
\in S_{H_1}^m \times S_{H_1}^n$ and $(u^{(2)},v^{(2)}) \in S_{H_2}^p 
\times S_{H_2}^q$, such that $S = \Gamma_{H_1}(u^{(1)},v^{(1)})$ and 
$R = \Gamma_{H_2}(u^{(2)},v^{(2)})$. Fix $(\alpha, \beta) \in [mp] \times [nq]$.
\eqref{eq:Kronecker_product} consequently implies that
\[
(S \otimes R)_{\alpha, \beta} = 
\big\langle v^{(1)}_{\frac{\beta - r_q(\beta)}{q}+1},
u^{(1)}_{\frac{\alpha - r_p(\alpha)}{p}+1}\big\rangle_{H_1} \cdot 
\big\langle v^{(2)}_{r_q(\beta)}, u^{(2)}_{r_p(\alpha)}\big\rangle_{H_2} =
\langle v_\beta, u_\alpha \rangle_H\,,
\]
where the Hilbert space $H : = H_1 \otimes H_2$ denotes the standard tensor 
product of the Hilbert spaces $H_1$ and $H_2$, $v_\beta : = 
v^{(1)}_{\frac{\beta - r_q(\beta)}{q}+1} \otimes v^{(2)}_{r_q(\beta)}$ and 
$u_\alpha : = u^{(1)}_{\frac{\alpha - r_p(\alpha)}{p}+1} \otimes 
u^{(2)}_{r_p(\alpha)}$. Since $u_\alpha \in S_H$ and $v_\beta \in S_H$ 
(by construction), it follows that $S \otimes R = \Gamma_H(u,v)$, where
$(u,v) \in S_H^{mp} \times S_H^{nq}$. Thus, $S \otimes R \in 
{\mathcal{Q}}_{mp,nq}(\F)$ (again, a consequence of 
\sref{Proposition}{prop:char_of_Q_m_n_for_both_fields}). 
\end{comment}
\noindent An important example of a matrix $A \in \M_{2^{m}}(\{-1,1\})$, which 
satisfies $\Vert A \Vert_{\infty, 1}^\R \leq 1$ and delivers 
$\sqrt{2}$ as a lower bound of $K_G^\R$ (cf. 
\sref{Proposition}{prop:ad_K_GR_2d_vs_K_GC_d}), and also plays a key role in 
the foundations of quantum mechanics and quantum information is the so-called 
\textit{Walsh-Hadamard transform} (also known as \textit{quantum gate} - cf. 
\cite{B2021}). In the following enlightening example, we extend the 
Walsh-Hadamard transform $H_m \in \M_{2^m}(\{-1,1\})$ to a complex 
Walsh-Hadamard transform $H_m^\C \in \M_{2^m}(\T)$ and disclose some surprising 
properties of that matrix. In particular, we will show that the (value of the) 
sign of any of the $4^{m}$ entries of the real Walsh-Hadamard transform can be 
specified precisely, in exactly $m$ calculation steps - for any $m \in \N$\,! 
To this end, recall the construction of the function $\pi_m : \Z \longrightarrow 
\{0,1\}^m$ (cf. \eqref{eq:Boolean_vector_function}) and put
\[
N_m(\nu, \mu) : = \langle\pi_m(\nu), \pi_m(\mu)\rangle_{\F_2^m} =
\sum_{i=1}^m b_i(\nu)b_i(\mu) = b_1(\nu)b_1(\mu) + N_{m-1}(f_2(\nu), f_2(\mu)),
\]
where $(\nu, \mu) \in [2^m] \times [2^m]$ and $N_0 : = 0$. 
$N_m(\nu, \mu)$ counts the number of all $i \in [m]$, such that 
$b_i(\nu)b_i(\mu) = 1$. In particular, $N_m(1, \mu) = 0$ for all 
$\mu \in [2^m]$ (since $b_1(1) = 0$). 
\begin{example}[\textbf{Real and complex Walsh-Hadamard transform}]
\label{ex:Hadamard_gate}
Let
\[
H_1 : = 
\tfrac{1}{\sqrt{2}}\begin{pmatrix}
1 & 1\\
1 & -1
\end{pmatrix} \text{ and } H_1^\text{op} : = R_2(i)H_1 = 
\tfrac{1}{\sqrt{2}}\begin{pmatrix}
-1 & 1\\
1 & 1
\end{pmatrix}\!.
\] 
For $m \in \N$, put
\[
H_{m+1} : = H_{m} \otimes H_{1} = H_{1} \otimes H_{m} = \tfrac{1}{\sqrt{2}}
\begin{pmatrix}
H_{m} & H_{m}\\
H_{m} & -H_{m}
\end{pmatrix}
\]
and
\[
H_{m+1}^\text{op} : = H_m \otimes H_1^\text{op} = H_1 \otimes H_m^\text{op} = 
\tfrac{1}{\sqrt{2}}
\begin{pmatrix}
H_{m}^\text{op} & H_{m}^\text{op}\\
H_{m}^\text{op} & -H_{m}^\text{op}
\end{pmatrix}\!. 
\]
Then the following properties are satisfied:
\begin{enumerate}
\item  
$(H_1)_{\alpha \beta} = \frac{1}{\sqrt{2}}\,(-1)^{(\alpha-1)(\beta-1)}$ for all 
$(\alpha, \beta) \in [2] \times [2]$. $H_1^\top = H_1 \in O(2)$ and
$(H_1^\text{op})^\top = H_1^\text{op} \in O(2)$. If $m \in \N_2$, then 
$(H_{m})^\top = H_{m} \in SO(2^m)$ and $(H_m^\text{op})^\top = H_m^\text{op} 
\in SO(2^m)$.
\item Let $m \in \N_2$ and $(\nu, \mu) \in [2^m] \times [2^m]$. Then
\begin{align}\label{eq:WH_matrix_sign_construction}
(H_m)_{\nu \mu} = \frac{1}{\sqrt{2}}\,(H_{m-1})_{f_2(\nu)\,f_2(\mu)}\cdot
(-1)^{b_1(\nu)\,b_1(\mu)} = \frac{1}{\sqrt{2^m}}\,
(-1)^{N_m(\nu, \mu)}\,,
\end{align}
In particular,
\[
(H_m)_{1 \mu} = (H_m)_{\mu 1} = \frac{1}{\sqrt{2^m}}\,.
\]
\item 
Let $m \in \N$. Then $H_m = \Re(H_m^\C) \in {\mathcal{Q}}_{2^m,2^m}(\R)$ and 
$H_m^\text{op} = \Im(H_m^\C) \in {\mathcal{Q}}_{2^m,2^m}(\R)$, where
\[
H_m^\C : = H_m + i\,H_m^\text{op} \in {\mathcal{Q}}_{2^m,2^m}(\C)\,.   
\]
In particular the matrix,
\[
H_1^\C = 
\begin{pmatrix}
1 & -i\\
i & 1
\end{pmatrix}H_1 =
\tfrac{1}{\sqrt{2}}\begin{pmatrix}
1-i & 1+i\\
1+i & -1+i
\end{pmatrix} = \overline{p}\,q^\top = \Gamma_\C(p,q)
\]
is of rank 1 and satisfies $\Vert H_{1}^\C \Vert_{\infty, 1} = \Vert p \Vert_1\,
\Vert q \Vert_1 = 4$, where $p : = (i,1)^\top \in \T^2$ and $q : = (\tfrac{1+i}{\sqrt{2}}, 
\tfrac{-1+i}{\sqrt{2}})^\top \in \T^2$ and $H_{m+1}^\C = H_1 \otimes H_m^\C$ for 
all $m \in \N$.  
\item For any $m \in \N$, $\Vert H_{m} \Vert_{\infty, 1} = 
\Vert H_m^\text{op} \Vert_{\infty, 1}$. The sequence $(\Vert H_{m} \Vert_{\infty, 1})_
{m \in \N}$ is non-decreasing. 
Moreover,
\begin{align}\label{eq:equivalence_of_infty_one_norms}
\sqrt{2}\,\Vert H_{m} \Vert_{\infty, 1} \leq \Vert H_{m+1} \Vert_{\infty, 1} \leq 
2\sqrt{2}\,\Vert H_{m} \Vert_{\infty, 1}
\end{align}
and
\begin{align}\label{eq:infty_one_norm_of_Quantum_gate_Kronecker_product}
(\sqrt{2})^m \leq \Vert H_{m} \Vert_{\infty, 1} \leq (\sqrt{2})^{3m-2}
\end{align}
for all $m \in \N$. In particular, $\Vert A^\text{Had}_{m} \Vert_{\infty, 1} 
\leq 1$ for all $m \in \N$, where 
\[
A^\text{Had}_{m} : = \frac{1}{(\sqrt{2})^{3m-2}}\,H_m = 
\frac{1}{2^{2m-1}}((\sqrt{2})^m\,H_m)  
\]
and
\begin{align}\label{eq:infty_one_norm_of_quantum_gate_sharp}
\frac{1}{\sqrt{2}}\,\Vert H_{1} \Vert_{\infty, 1} = 
\Vert A^\text{Had}_{1} \Vert_{\infty, 1} = 1.
\end{align}
\end{enumerate}
\end{example}
\begin{remark}
After some ``skillful searching'', a then simple calculation shows that also
\[
\Vert A^\text{Had}_{2} \Vert_{\infty, 1} = 1\,.
\]
If we namely consider the vectors $\widetilde{p} : = (1,1,-1,1)^\top \in \{-1,1\}^4$ 
and $\widetilde{q} : = (1,-1,1,1)^\top \in \{-1,1\}^4$, it follows that
\[
\text{tr}(A^\text{Had}_2 \widetilde{p} \widetilde{q}^\top) = 
\langle A^\text{Had}_2 \widetilde{p}, \widetilde{q}\rangle_{\R_2^4} = \frac{1}{8} 
\Bigg\langle\colvec{4}{\,\,2}{\!\!-2}{\,\,2}{\,\,2}, \colvec{4}{\,\,1}{\!\!-1}
{\,\,1}{\,\,1}\Bigg\rangle = 1\,.
\]
In particular, 
\begin{align}\label{eq:the_case_m_equal_2}
\big\vert (A^\text{Had}_{m}\,\widetilde{p})_i \big\vert = \frac{1}{2^m} 
\text{ for all } i \in [2^m]
\end{align}
(since $m = 2$). This naturally leads to the (open) question, whether 
$\Vert A^\text{Had}_{m}\Vert_{\infty, 1} = 1$ for all $m \in \N$ and whether 
\eqref{eq:the_case_m_equal_2} holds for all $m \in \N$. It seems that we cannot 
make use of induction on $m \in \N$ here. In fact, if $\nu \not= 2$ and $\Vert p\Vert_
{l_\infty^{2^\nu}} \leq 1$, then $\big\vert (A^\text{Had}_{\nu}\,p)_{i_0}\big\vert 
\not= \frac{1}{2^\nu}$ for some $i_0 \in [2^\nu]$\,! The case $\nu = 1$ follows 
from the fact that for any $a, b \in [-1,1]$, $\vert a + b \vert = 1 = 
\vert a - b \vert$ if and only if $\colvec{2}{a}{b} \in \Big\{\colvec{2}{0}{1}, 
\colvec{2}{0}{-1}, \colvec{2}{1}{0}, \colvec{2}{-1}{0}\Big\}$. In order to 
verify the claim for $\nu > 2$, assume by contradiction that there exist $m \in 
\N_3$ and $\widetilde{p} \in B_{l_\infty^{2^m}}$, such that $\frac{1}{2^m} = 
\big\vert (A^\text{Had}_{m}\,\widetilde{p})_i \big\vert = \frac{1}
{(\sqrt{2})^{3m-2}}\,\vert(H_m\,\widetilde{p})_i \vert = \frac{1}{2^{2m-1}}
\,\big\vert\sum\limits_{j=1}^{2^m}(-1)^{N_m(i,j)}\widetilde{p}_j \big\vert$ for 
all $i \in [2 ^m]$. Put $\widetilde{q} : = 2^m A^\text{Had}_{m} \widetilde{p}$. 
Then $\widetilde{q} \in \{-1,1\}^{2^m}$ (due to the assumption) and
\begin{align*}
1 &= 2^m\,\Vert A^\text{Had}_{m} \widetilde{p}\Vert_{l_2^{2^m}}^2 = 
\langle A^\text{Had}_{m} \widetilde{p}, \widetilde{q}\rangle_
{\R_2^{2^m}} = 2^m\,\langle \widetilde{p}, (A^\text{Had}_{m})^2\,
\widetilde{p}\rangle_{\R_2^{2^m}}\\ 
&= \frac{2^m}{2^{3m-2}}\,\Vert \widetilde{p}\Vert_{l_2^{2^m}}^2 \leq 
\frac{1}{2^{3m-2}}\,2^{2m} = \frac{1}{2^{m-2}}\,.
\end{align*}
On the other hand, $\frac{1}{2^{m-2}} \leq \frac{1}{2} < 1$ (since $m \geq 3$ by 
assumption), which is absurd. Observe that in any case 
$\big\vert(A^\text{Had}_{m}\,x)_i\big\vert \leq \frac{2^m}{2^{2m-1}} = 
\frac{2}{2^m}$ for all $i \in [2^m]$ and $x \in B_{l_\infty^{2^m}}$. 
Consequently,
\[
\big\Vert A^\text{Had}_{m}\,\frac{p}{2}\big\Vert_{l_\infty^{2^m}} < \frac{1}{2^m}
\]
for all $m \in \N_3$ and $p \in B_{l_\infty^{2^m}}$. Our conjecture is that there 
exists $\widetilde{m} \in \N_3$ such that $\Vert A^\text{Had}_
{\widetilde{m}}\Vert_{\infty, 1} < 1$.
\end{remark}
 To be more explicit, note e.g. that
\[
H_{2} = \frac{1}{2}
\begin{pmatrix}
1 & 1 & 1 & 1\\
1 & \!\!\!\!-1 & 1 & \!\!\!\!-1\\
1 & 1 & \!\!\!\!-1 & \!\!\!\!-1\\
1 & \!\!\!\!-1 & \!\!\!\!-1 & 1
\end{pmatrix},\hspace{0.2cm}
H_{2}^\text{op} = \frac{1}{2}
\begin{pmatrix}
-1 & 1 & \!\!\!\!-1 & 1\\
\,\,\,\,1 & 1 & 1 & 1\\
-1 & 1 & 1 & \!\!\!\!-1\\
\,\,\,\,1 & 1 & \!\!\!\!-1 & \!\!\!\!-1
\end{pmatrix}
\]
and
\[
H_{3} = \frac{1}{(\sqrt{2})^3}
\begin{pmatrix}
1 & 1 & 1 & 1 & 1 & 1 & 1 & 1\\
1 & \!\!\!\!-1 & 1 & \!\!\!\!-1 & 1 & \!\!\!\!-1 & 1 & \!\!\!\!-1\\
1 & 1 & \!\!\!\!-1 & \!\!\!\!-1 & 1 & 1 & \!\!\!\!-1 & \!\!\!\!-1\\
1 & \!\!\!\!-1 & \!\!\!\!-1 & 1 & 1 & \!\!\!\!-1 & \!\!\!\!-1 & 1\\
1 & 1 & 1 & 1 & \!\!\!\!-1 & \!\!\!\!-1 & \!\!\!\!-1 & \!\!\!\!-1\\
1 & \!\!\!\!-1 & 1 & \!\!\!\!-1 & \!\!\!\!-1 & 1 & \!\!\!\!-1 & 1\\
1 & 1 & \!\!\!\!-1 & \!\!\!\!-1 & \!\!\!\!-1 & \!\!\!\!-1 & 1 & 1\\
1 & \!\!\!\!-1 & \!\!\!\!-1 & 1 & \!\!\!\!-1 & 1 & 1 & \!\!\!\!-1\\
\end{pmatrix}\,.
\]
To see how smoothly and quickly \eqref{eq:WH_matrix_sign_construction} can be 
applied to $H_3$, let us perform a calculation for the two matrix entries 
$(H_3)_{6,4}$ and $(H_3)_{7,3}$ as an example. There are two ways to perform 
the calculation process. Either we go through each step of the recursion 
relation, or we count, step by step, how often the products 
$b_i(\nu)b_i(\mu)$ are equal to 1. So, either we proceed with 
\[
(H_3)_{6,4} = -\frac{1}{\sqrt{2}}\,(H_2)_{3,2} = -\frac{1}{\sqrt{2^2}}\,
(H_1)_{2,1} = \frac{1}{\sqrt{2^3}}(-1)
\]
and
\[
(H_3)_{7,3} = \frac{1}{\sqrt{2}}\,(H_2)_{4,2} = -\frac{1}{\sqrt{2^2}}\,
(H_1)_{2,1} = \frac{1}{\sqrt{2^3}}(-1)\,,
\]
or we count: $N_3(6,4) = 1$ (since $b_1(6)b_1(4) = 1, b_2(6)b_2(4) = b_1(3)b_1(2) 
= 0$ and $b_3(6)b_3(4) = b_2(3)b_2(2) = b_1(2)b_1(1) = 0$). Similarly, we obtain that 
$N_3(7,3) = 1$. However, since $\frac{1}{\sqrt{2^3}}(-1) = (H_3)_{7,3} \not= 
\frac{1}{\sqrt{2^3}} = \frac{1}{\sqrt{2^3}}(-1)^{(7-1)(3-1)}$, 
\cite[Exercise 2.33, (2.55)]{NC2000} seems to be wrong.
\begin{remark}[\textbf{CHSH inequalities}]\label{rem:CHSH_inequalities}
As we have seen (just by making use of elementary calculus on the 
real line), the following inequality holds
\begin{align}\label{eq:CHSH}
\vert x_1 y_1 + x_1 y_2 + x_2 y_1 - x_2 y_2 \vert \leq 2 \text{ for all } 
(x_1,x_2,y_1,y_2) \in [-1,1]^4,
\end{align}
which is equivalent to $A^\text{Had}_{1} = \tfrac{1}{\sqrt{2}}H_1 \in 
B_{\mathfrak{L}(l_\infty^2, l_1^2)}$. Even $A^\text{Had}_{1} \in 
S_{\mathfrak{L}(l_\infty^2, l_1^2)}$ holds (cf. 
\eqref{eq:infty_one_norm_of_quantum_gate_sharp}). We also know that 
$\Vert A^\text{Had}_{1} \Vert_{\infty, 1} \leq 1$ is equivalent to
\[
\vert \text{tr}(A^\text{Had}_{1} B)\vert \leq 1 \,\text{ for all }\, B \in 
B_{{\mathfrak{N}}(l_1^2, l_\infty^2)}
\]
(cf. \eqref{eq:trace_duality_fin_dim_case}). On the other hand, 
\[
\sqrt{2} = \frac{1}{\sqrt{2}}\,\vert{\text{tr}}(I_2)\vert
= \frac{1}{\sqrt{2}}\,\vert{\text{tr}}(H_1^2)\vert =
\vert {\text{tr}}(A^\text{Had}_{1} \widetilde{S})\vert > 1 
\text{ for some } \widetilde{S} \in {\mathcal{Q}}_{2,2} = 
B_{{\mathfrak{L}}_2(l_1^2, l_\infty^2)}
\]
(namely, $\widetilde{S} : = H_1$), implying again that 
$B_{{\mathfrak{N}}(l_1^2, l_\infty^2)}$ is strictly contained in
$B_{{\mathfrak{L}}_2(l_1^2, l_\infty^2)}$. Particularly, physicists, who 
are working in the foundations and philosophy of quantum mechanics recognise 
that \eqref{eq:CHSH} -- which are just inequalities between certain real numbers 
-- \textit{instantly imply} the famous CHSH inequalities. CHSH stands for 
John Clauser, Michael Horne, Abner Shimony, and Richard Holt, who introduced the 
inequalities (between expectation values) in \cite{CHSH1969} (cf. 
\url{https://www.nobelprize.org/prizes/physics/2022/clauser/facts/}) and used 
them as a means of proving Bell's theorem. In the 2-dimensional case (i.e., if 
$m=n=2$) the CHSH inequalities coincide with the so-called ``Bell inequalities'', 
assigned to the matrix $A^\text{Had}_{1}$. Somewhat vaguely, it is said that the 
matrix $\widetilde{S}  = H_1 \in {\mathcal{Q}}_{2,2}(\R)$ ``violates 
the Bell inequalities''. That ``violation'' implies that certain consequences of 
spatial entanglement in quantum mechanics can not be reproduced by classical 
probability theory in the sense of A. Kolmogorov (i.e., it cannot be reduced to 
``local hidden-variable theories'').
\end{remark}
\begin{remark}[\textbf{An application of $H_m$ in evolutionary biology}]
\label{rem:appl_in_evolutionary_biology}
The Walsh-Hadamard transform can even be found in evolutionary biology, 
specifically in relation to the challenge of reconstructing evolutionary trees 
from events several million years in the past. (cf. \cite{HP1989})! In order to
recognise this, we consider the orthogonal matrix $I_1^{(-)} : = 
\begin{pmatrix}
1 & 0\\
0 & -1
\end{pmatrix} \in O(2)$. A straightforward proof by induction shows that the 
matrix family $\{H^{(m)} : m \in \N\} = \{(1)\} \dotcup \{H^{(m)} : m \in \N_2\}$,
consisting of invertible matrices $H^{(m)} \in \M_{2^{m-1}}(\R)$, introduced in 
\cite{HP1989}, in fact can be represented as
\[
H^{(1)} : = (1)
\]
and
\[ 
H^{(m)} : = \begin{pmatrix}
H^{(m-1)} & -H^{(m-1)}\\
H^{(m-1)} & H^{(m-1)}
\end{pmatrix} = H^{(2)} \otimes H^{(m-1)} = \bigotimes\limits_{i=1}^{m-1} H^{(2)} =
\sqrt{2^{m-1}}\,H_{m-1}\,I_{m-1}^{(-)}
\]
if $m \in \N_2$, where $I_l^{(-)} : = I_1^{(-)} \otimes I_{l-1}^{(-)} = 
\bigotimes\limits_{i=1}^l I_1^{(-)} \in O(2^l)$ for all $l \in \N_2$.  
\end{remark}
\noindent It is quite instructive to compare 
\eqref{eq:QM_is_GT_constant_times_acx_hull_of_rank_1} to 
\eqref{eq:rep_of_c_times_QC_real_case} (real case), 
\eqref{eq:rep_of_c_times_QC_complex_case} (complex case) 
and the following
\begin{proposition}\label{prop:QM_is_GT_constant_times_acx_hull_of_rank_1_generalised}
Suppose there exist $c \in (1, \infty)$, a sequence $(r_\nu)_{\nu \in \N} \in 
B_{l_1(\F)}$, a probability space $(\Omega, \mathscr{F}, \P)$ and sequences 
$(P_{1, \nu})_{\nu \in \N}, \ldots, (P_{m, \nu})_{\nu \in \N}$, 
$(Q_{1, \nu})_{\nu \in \N}, \ldots, (Q_{n, \nu})_{\nu \in \N}$ of random 
variables which map into $S_\F$ $\P$-a.s., 
such that for all $m,n \in \N$, for all $\F$-Hilbert spaces $H$, 
for all $(u, v) \in S_H^m \times S_H^n$, and for all $(i,j) \in [m] \times [n]$ 
the following equality holds:
\begin{align*}
\langle v_j, u_i\rangle_H = \Gamma_H(u,v)_{ij} = c\,\sum_{\nu=1}^\infty r_\nu\,
\E_{\P}[\overline{P_{i,\nu}} Q_{j,\nu}]. 
\end{align*}
Then
\[
K_G^\F \leq c\,.
\]
\end{proposition}
\begin{comment}
\begin{align*}
\vert\text{tr}(A^\ast \Gamma_H(u,v))\vert &= 
c\big\vert \sum_{\nu=1}^\infty r_\nu \big(\sum_{i=1}^m\sum_{j=1}^n 
\overline{a_{ij}}\,\E_{\P}[\overline{P_{i,\nu}} Q_{j,\nu}]\big)\big\vert 
\leq c \sum_{\nu=1}^\infty \vert r_\nu \vert\,\big\vert\sum_{i=1}^m
\sum_{j=1}^n \overline{a_{ij}}\,
\E_{\P}[\overline{P_{i,\nu}} Q_{j,\nu}]\big\vert\\
&\leq c \sum_{\nu=1}^\infty \vert r_\nu \vert\,\E_{\P}\big[\big\vert\sum_{i=1}^m
\sum_{j=1}^n \overline{a_{ij}}\,\overline{P_{i,\nu}} Q_{j,\nu}\big\vert\big] = 
c \sum_{\nu=1}^\infty \vert r_\nu \vert\,\E_{\P}\vert\text{tr}
(A^\ast \Gamma_\F(\textbf{P}_\nu,\textbf{Q}_\nu))\vert\,,
\end{align*}
where the components of the random vectors $\textbf{P}_\nu : \Omega 
\longrightarrow {\F \cap \overline{\D}}^m$ and $\textbf{Q}_\nu : \Omega \longrightarrow 
{\F \cap \overline{\D}}^n$ are defined as $(\textbf{P}_\nu)_i : = 
P_{i,\nu}$ and $(\textbf{Q}_\nu)_j : = Q_{j,\nu}$. 
Since $\vert\text{tr}(A^\ast \Gamma_\F(\textbf{P}_\nu(\omega), 
\textbf{Q}_\nu(\omega)))\vert \leq \Vert A \Vert_{\infty, 1}$ for 
almost all $\omega \in \Omega$, it therefore follows that
\[
\vert{\text{tr}}(A^\ast \Gamma_H(u,v))\vert \leq 
c \sum_{\nu=1}^\infty \vert r_\nu \vert \Vert A \Vert_{\infty, 1} \leq 
c\Vert A \Vert_{\infty, 1}
\]
(since $\sum_{\nu=1}^\infty \vert r_\nu \vert \leq 1$, by assumption).
\end{comment}
\noindent Of particular interest is the value of $K_G^\R(2)$. A 
\textit{lower} bound is rather easy to detect: $\sqrt{2} \leq K_G^\R(2)$, 
implying the important fact that $K_G^\R > 1$. We just have to work with the 
Walsh-Hadamard transform $H_1 = \Gamma_{\R_2^2}(u_1, u_2, v_1, v_2) \in 
{\mathcal{Q}}_{2,2}(\R) \cap O(2)$ (cf. \sref{Example}{ex:Hadamard_gate}). To this 
end, consider again the symmetric matrix $A^\text{Had}_{1} = \tfrac{1}{\sqrt{2}}H_1$. 
Recall that $\Vert A^\text{Had}_{1} \Vert_{\infty, 1} = 1$ (due to 
\eqref{eq:infty_one_norm_of_quantum_gate_sharp}), and note that 
$\text{tr}(A^\text{Had}_{1} H_1) = \frac{1}{\sqrt{2}}\,\text{tr}(H_1^2) = 
\frac{1}{\sqrt{2}}\,\text{tr}(I_2) = \sqrt{2}$. Consequently, it follows that
\[
\sqrt{2} = \vert{\text{tr}}((A^\text{Had}_{1})^\top H_1)\vert \leq K_G^\R(2,2;2) \leq 
\min\{K_G^\R(2,2), K_G^\R(2)\} \leq K_G^\R.
\]
Much less trivial is the proof of the reverse direction, performed by Krivine 
in \cite{Kr1977, Kr1979}. Within the scope of 
\sref{Proposition}{prop:K_G_R_resp_K_G_C_as_sup_of_nuclear_norms} he namely 
represented - in the real 2-dimensional case - any $\Gamma_{\R_2^2}(u,v) \in 
{\mathcal{Q}}_{2,2}(\R)$ as matrix $\sum_{\nu=1}^\infty b_\nu\,T_\nu : l_1^n 
\longrightarrow l_\infty^m$ such that $\Vert\Gamma_{\R_2^2}(u,v)\Vert_
{\mathfrak{N}} = \Vert\sum_{\nu=1}^\infty b_\nu\,T_\nu\Vert_{\mathfrak{N}} \leq 
\sqrt{2}$ (which he called ``norme de la fonction $\cos(x-y)$ dans le produit 
tensoriel projectif $C[-\pi, \pi] \,\widehat{\otimes}\, C[-\pi, \pi]$'' 
in \cite{Kr1977}). Actually, the main building block in his proof is an 
intricate sophisticated representation of the function $\R \times \R \ni (x,y) 
\mapsto \cos(x-y)$ by convolution (cf. \autoref{thm:Krivine_main_step_to_KGR_of_2}). 
However, that representation allows us to provide a short, straightforward 
proof of \sref{Corollary}{cor:Krivine_KGR_of_2} - without the use of any 
tensor product structure.
\begin{theorem}[\textbf{Krivine, 1977}]\label{thm:Krivine_main_step_to_KGR_of_2}
Consider the probability space $([-\pi,\pi], {\mathcal{B}}([-\pi,\pi]), \mu)$, where 
$\mu : = \frac{1}{2\pi}\lambda_1{\big\vert}_{{\mathcal{B}}([-\pi,\pi])}$. Then 
there exist two functions $p, q : \R \longrightarrow [-1,1]$ and a sequence 
$(r_n)_{n \in \N}$ of real numbers such that $\sum_{n=1}^\infty \vert r_n \vert = 1$
and 
\begin{align*}
\cos(x-y) &= \sqrt{2}\,\sum_{n=1}^\infty r_n\E_{\mu}[p(n x-S)
q(n y-S)]\\ 
&= \frac{\sqrt{2}}{2\pi}\sum_{n=1}^\infty r_n \int_{[-\pi, \pi]} 
p(n x - \omega)q(n y - \omega)\lambda_1(\textup{d}\omega)
\end{align*}
for all $x, y \in \R$, where $\R \ni t \mapsto q(t) : = \text{sign}(\cos(t))$ 
and $[-\pi, \pi] \ni \omega \mapsto S(\omega) : = \omega$.
\end{theorem}
\begin{corollary}[\textbf{Krivine, 1977}]\label{cor:Krivine_KGR_of_2}
\[
K_G^\R(2) = \sqrt{2}\,.
\]
\end{corollary}
\begin{comment}
Consider the 2-dimensional standard Euclidean vector space $H : = \R_2^2$. Let 
$u = \colvec{2}{u_1}{u_2} \in S_H = \S^1$ and $v = \colvec{2}{v_1}{v_2} \in S_H = 
\S^1$. Since
\[
\langle u, v \rangle_H = \big\langle \colvec{2}{\cos(x)}{\sin(x)}, 
\colvec{2}{\cos(y)}{\sin(y)}\big\rangle_H = \cos(x-y)\,,
\]
for some $x,y \in \R$, \autoref{thm:Krivine_main_step_to_KGR_of_2} unveils 
as a special case of 
\sref{Proposition}{prop:QM_is_GT_constant_times_acx_hull_of_rank_1_generalised}, 
and the claim follows.
\end{comment}
\noindent \sref{Proposition}{prop:K_G_R_resp_K_G_C_as_sup_of_nuclear_norms}, 
together with \sref{Corollary}{cor:Krivine_KGR_of_2} help us to find a ``fitting'' 
relation between $K_G^\R(2d)$ and $K_G^\C(d)$. To this end, we firstly supplement 
and prove once again (for the sake of completeness) \cite[Corollary 2.14., (38)]
{FL2020}. 
\begin{proposition}\label{prop:ad_K_GR_2d_vs_K_GC_d}
Let $d, m, n \in \N$, $A \in \M_{m,n}(\R)$ and $B \in \M_{m,n}(\C)$. Then
\begin{align}\label{eq:real_and_complex_GT_norms_in_the_real_case}
\Vert A \Vert_{\R_2^d}^{\text{G}} \leq \Vert A \Vert_{\C_2^d}^{\text{G}} = 
\Vert A \Vert_{\R_2^{2d}}^{\text{G}} \leq K_G^\R(m,n; 2d)\,\Vert A \Vert_{\infty, 1}^\R\,.
\end{align}
In particular, 
\begin{align}\label{eq:real_and_complex_infty_1_norms_in_the_real_case}
\Vert A \Vert_{\infty, 1}^{\R} \leq \Vert A \Vert_{\infty, 1}^{\C} = 
\Vert A \Vert_{\R_2^2}^{\text{G}} \leq \sqrt{2}\,\Vert A \Vert_{\infty, 1}^{\R}
\end{align}
and
\begin{align}\label{eq:GT_norms_of_Re_and_Im_parts}
\Vert\Re(B)\Vert_{\C_2^d}^{\text{G}} \leq \Vert B \Vert_
{\C_2^d}^{\text{G}} \text{ and } \Vert\Im(B)\Vert_{\C_2^d}^{\text{G}} \leq 
\Vert B \Vert_{\C_2^d}^{\text{G}}\,.
\end{align}
Moreover,
\[
\Vert A^\text{Had}_{1} \Vert_{\infty, 1}^{\R} = 1 < \sqrt{2} = 
\Vert A^\text{Had}_{1} \Vert_{\infty, 1}^{\C}\,.
\]
\end{proposition}
\noindent Since
\[
\Vert A \Vert_{{\R_2^{2d}}}^{\text{G}} 
\stackrel{\eqref{eq:real_and_complex_GT_norms_in_the_real_case}}{=}
\Vert A \Vert_{{\C_2^d}}^{\text{G}} 
\leq K_G^\C(d)\,\Vert A\Vert_{\infty, 1}^{\C} 
\stackrel{\eqref{eq:real_and_complex_infty_1_norms_in_the_real_case}}{\leq} 
\sqrt{2}\,K_G^\C(d)\,\Vert A\Vert_{\infty, 1}^{\R}
\]
for all $m, n \in \N$ and $A \in \M_{m,n}(\R)$ and 
\begin{align*}
\Vert B \Vert_{{\C_2^d}}^{\text{G}} &= \Vert \Re(B) + i \Im(B)\Vert_{{\C_2^d}}^
{\text{G}} \leq \Vert \Re(B)\Vert_{{\C_2^d}}^{\text{G}} + \Vert \Im(B)\Vert_
{{\C_2^d}}^{\text{G}}\\ 
&\!\!\stackrel{\eqref{eq:real_and_complex_GT_norms_in_the_real_case}}{\leq} 
K_G^\R(2d)(\Vert\Re(B)\Vert_{\infty, 1}^{\R} + \Vert\Im(B)\Vert_{\infty, 1}^{\R}) 
\stackrel{\eqref{eq:GT_norms_of_Re_and_Im_parts}}{\leq} 2\,K_G^\R(2d)\,
\Vert B\Vert_{\infty, 1}^\C 
\end{align*}
for all $m, n \in \N$ and $B \in \M_{m,n}(\C)$, we obtain
\begin{corollary}\label{cor:K_GR_2d_vs_K_GC}
\[
\tfrac{1}{\sqrt{2}}\,K_G^\R(2d) \leq K_G^\C(d) \leq 2\,K_G^\R(2d) 
\text{ for all } d \in \N.
\]
In particular,
\begin{align}\label{eq:complex_versus_real_GT_constant}
\tfrac{1}{\sqrt{2}} K_G^\R \leq K_G^\C \leq 2\,K_G^\R\,.
\end{align}
\end{corollary}
\noindent Note that the implication \eqref{eq:complex_versus_real_GT_constant} 
contains \cite[Theorem 10.6]{J1987}. We do not know whether the second estimation 
could be improved to $K_G^\C(d) \stackrel{?}{\leq} \sqrt{2}\,K_G^\R(2d)$ for 
all $d \in \N$. In particular, since $\Re(B)$ and $\Im(B)$ do not commute, we 
do not know whether $(\Vert \Re(B) + i \Im(B)\Vert_{{\C_2^d}}^{\text{G}})^2 
\stackrel{?}{\leq} (\Vert\Re(B)\Vert_{{\C_2^d}}^{\text{G}})^2 + 
(\Vert\Im(B)\Vert_{{\C_2^d}}^{\text{G}})^2$ holds for all $m, n \in \N$ and 
$B \in \M_{m,n}(\C)$ (in analogy to $\vert x + iy\vert^2 = x^2 + y^2$ for all 
$x,y \in \R$).
%
\section{The Gaussian inner product splitting property}
 Next, we are going to disclose a crucial \textit{joint} multivariate 
Gaussian ``splitting property'' of inner products of vectors on the unit sphere 
of an arbitrary \textit{separable} $\F$-Hilbert space. That result (which should 
be compared to the construction of Gaussian Hilbert spaces or Moore's theorem on 
the characterisation of kernel functions (cf. \cite[Theorem 2.14]{PR2016})) 
runs like a thread throughout the whole paper, including its implementation in Theorem 
\ref{thm:h_fg_real_case} and Theorem 
\ref{thm:odd_completely_real_analytic_functions_at_0_vs_complex_h_fg}. It holds 
for both fields, $\F=\R$ and $\F = \C$, and plays a significant role, when we are 
looking for a specific Gaussian random structure in quantum correlation matrices. 
So, let $k, m, n \in \N$. Fix $S = \Gamma_H(u,v) = (u_i^\ast v_j)_{ij} 
\in Q_{m,n}(\F)$, where $(u,v) \in S_H^m \times S_H^n$ and $(i,j) \in [m] \times 
[n]$. Based on our analysis so far, if $\zeta_{ij} : = u_i^\ast v_j = 
\langle v_j, u_i\rangle_H$ is given, then we only know about the existence of a 
joint Gaussian random vector $\vc{{\textbf{Z}}_{ij}}{{\textbf{W}}_{ij}} \sim 
{\F}N_{2k}(0, \Sigma_{2k}(\zeta_{ij}))$. A priori, we cannot say whether 
it is even possible to allocate to $\zeta_{ij} = u_i^\ast v_j$ a joint 
Gaussian random vector of type $\vc{{\textbf{Z}}_{i}}{{\textbf{W}}_{j}} \sim 
{\F}N_{2k}(0, \Sigma_{2k}(\zeta_{ij}))$. In fact, our next cornerstone 
result reveals that such a ``joint Gaussian splitting of an inner product'' 
is guaranteed if we assume that $H$ is separable, $u_i \in S_H$ and $v_j \in 
S_H$. For this, we fix an arbitrary complete probability space 
$(\Omega, \mathscr{F}, \P)$ and construct a suitable random field.  
\begin{proposition}[\textbf{Inner product splitting}]\label{prop:correl_splitting}
Let $k \in \N$, $\F \in \{\R,\C\}$ and $H$ be a separable $\F$-Hilbert 
space. There exists a family $\{{\textbf{Z}}_x : x \in H\}$ of random vectors
${\textbf{Z}}_x \equiv (Z_x^{(1)},Z_x^{(2)}, \ldots, Z_x^{(k)})^\top$ in $\F^k$,
such that
\begin{align}\label{eq:jointly_Gaussian}
\vc{{\textbf{Z}}_x}{{\textbf{Z}}_y} \sim \F N_{2k}(0, C_{2k}(x,y)) 
\text{ for all } x, y \in H, 
\end{align}
where 
\[
C_{2k}(x,y) : = \begin{pmatrix}
\Vert x \Vert_H^2\,I_k & \langle x, y\rangle_H\,I_k\\
\langle y, x\rangle_H\,I_k & \Vert y \Vert_H^2\,I_k
\end{pmatrix}\!. 
\]
In particular,
\[
\E[\vert Z_x^{(\nu)} \vert^2] = \Vert x \Vert_H^2\,\text{ and }\,
\langle x, y \rangle_H = \E[Z_x^{(\nu)}\,\overline{Z_y^{(\nu)}}]\,\text{ for all }\,
\nu \in [k] \text{ and } x, y \in H\,.
\] 
If $w \in S_H$, then ${\textbf{Z}}_w \sim \F N_{k}(0, I_k)$ and
\begin{align}\label{eq:jointly_Gaussian_on_the_unit_sphere}
\vc{{\textbf{Z}}_u}{{\textbf{Z}}_v} \sim 
\F N_{2k}(0, \Sigma_{2k}(\langle u, v \rangle_H)) 
\text{ for all } u, v \in S_H. 
\end{align}
If $e_1, e_2 \in S_H$ are orthogonal, then 
\begin{align}\label{eq:jointly_Gaussian_on_the_clsd_unit_disk}
\vc{{\textbf{Z}}_\zeta}{{\textbf{Z}}_1} \sim 
\F N_{2k}(0, \Sigma_{2k}(\zeta)) \text{ for all }
\zeta \in \overline{\D},
\end{align}
where ${\textbf{Z}}_\zeta : = {\textbf{Z}}_{\zeta e_1 + 
\sqrt{1-\vert\zeta\vert^2} e_2}$ and ${\textbf{Z}}_1 : = {\textbf{Z}}_{e_1}$. 
In particular, $\vc{{\textbf{Z}}_{\langle u, v \rangle_H}}{{\textbf{Z}}_1} 
\stackrel{d}{=} \vc{{\textbf{Z}}_u}{{\textbf{Z}}_v}$ for all $u, v \in S_H$.
${\textbf{Z}}_x \in L^2(\Omega)^k$ for all $x \in H$, and $Tx : = 
{\textbf{Z}}_x$ defines a bounded linear operator $T \in \mathfrak{L}(H, 
L^2(\Omega)^k)$, such that $\frac{1}{\sqrt{k}}\,T$ is an isometry. 
For any $\nu \in [k]$, the family $\{Z_x^{(\nu)} : x \in H\}$ is an $H$-isonormal 
process. 
\end{proposition}
\noindent Let $\Sigma = (\sigma_{ij})_{(i,j) \in [n] \times [n]} \in C(n; \F)$ 
be an arbitrary correlation matrix. Then $\Sigma = \Gamma_{H_n}(w,w)$ for some 
$w \equiv (w_1, w_2, \ldots, w_n) \in S_{H_n}^n$, where $H_n : = \F_2^n$ 
(due to \sref{Lemma}{lem:charact_of_corr_matrices}), implying that $\sigma_{ij} = 
w_i^\ast w_j = \langle w_j, w_i\rangle_{H_n}$ for all $i, j \in [n]$. Consequently, 
\eqref{eq:jointly_Gaussian_on_the_unit_sphere}, applied to all pairs $w_i, w_j 
\in S_H$, immediately results in
\begin{corollary}[\textbf{Correlation matrix splitting}]\label{cor:corr_matrix_splitting}
Let $n \in \N$ and $\Sigma = (\sigma_{ij})_{(i,j) \in [n] \times [n]} \in 
C(n; \F)$. Let $k \in \N$. Then there exist $n$ $\F^k$-valued random vectors 
${\textbf{Z}}_1, {\textbf{Z}}_2, \ldots, {\textbf{Z}}_n$ such that 
\[
\vc{{\textbf{Z}}_i}{{\textbf{Z}}_j} \sim 
\F N_{2k}(0, \Sigma_{2k}(\sigma_{i j})) \text{ and } 
\frac{1}{\sqrt{k}}{\textbf{Z}}_i \in S_{H_k}
\] 
for all $(i,j) \in [n] \times [n]$, where $H_k : = L^2(\Omega)^k$. In particular, 
\[
\sigma_{i j} = \E[Z_i^{(\nu)}\,\overline{Z_j^{(\nu)}}] = 
\big\langle \frac{1}{\sqrt{k}}{\textbf{Z}}_i, \frac{1}{\sqrt{k}}{\textbf{Z}}_j
\big\rangle_{H_k} \text{ for all } (i,j) \in [n] \times [n] \text{ and } \nu \in [k].
\]
Moreover, for any $\F$-Hilbert space $H$, for any $u, v \in S_H$, there exist 
two joint Gaussian random vectors $\textbf{W}_u$ and $\textbf{W}_v$ in $\F^k$, 
such that $\frac{1}{\sqrt{k}}{\textbf{W}}_u \in S_{H_k}$, $\frac{1}{\sqrt{k}}
{\textbf{W}}_v \in S_{H_k}$, $W_u^{(\nu)} \in S_{L^2(\Omega)}$, $W_v^{(\nu)} \in 
S_{L^2(\Omega)}$ and
\begin{align*}
\langle v, u\rangle_H = \E[W_v^{(\nu)}\,\overline{W_{\scriptscriptstyle u}^{(\nu)}}\,] = 
\big\langle\frac{1}{\sqrt{k}}{\textbf{W}}_v, \frac{1}{\sqrt{k}}{\textbf{W}}_u
\big\rangle_{H_k}\,\text{ for all } \nu \in [k].
\end{align*}
\end{corollary}
\chapter{Powers of inner products of random vectors, uniformly distributed on the
sphere} 
\section{Gaussian \text{sign}-correlation}
 It seems to be the case that any rigorous proof of the Grothendieck 
inequality is built on two \textit{equalities}, namely the Grothendieck equality 
(if $\F = \R$ - cf. e.g. \cite{FLZ2018, K1992}, or the proof of 
\cite[Prop. 4.4.2]{DFS2008}) and the Haagerup equality (if $\F = \C$ - 
see \cite{FLZ2018, H1987, K1992}). In fact, if we reveal the 
inherent \textit{bivariate} Gaussian random structure, these equalities emerge as 
two special cases of the representation of a single Pearson correlation coefficient 
which applies likewise for the real case and the complex case 
(\sref{Corollary}{cor:unification_of_Grothendieck_and_Haagerup}). Rewritten in 
terms of real Gaussian random vectors (if $\F = \R$) and complex Gaussian 
random vectors (if $\F = \C$) namely, we firstly obtain a representation of the 
two equalities, indicating an already lurking common underlying probabilistic 
structure for both fields, $\R$ and $\C$ (cf. also
\eqref{eq:expected_l_2_norm_of_real_Gaussian} and 
\eqref{eq:expected_l_2_norm_of_complex_Gaussian}).
To this end, recall (cf., e.g., \cite{AAR1999} and \cite[Chapter II]{BEMOT1953}) 
that for any $a,b\in \C$, any $c \in \C\setminus\{-n : n \in \N_0\}$ and any 
$z \in \D$ the well-defined power series
\[
{}_2 F_1(a, b, c; z) : = \frac{\Gamma(c)}{\Gamma(a)\Gamma(b)}\sum_{n=0}^{\infty} 
\frac{\Gamma(a+n)\Gamma(b+n)}{\Gamma(c+n)}\,\frac{z^n}{n!} = 1 + 
\frac{\Gamma(c)}{\Gamma(a)\Gamma(b)}\sum_{n=1}^{\infty} 
\frac{\Gamma(a+n)\Gamma(b+n)}{\Gamma(c+n)}\,\frac{z^n}{n!}
\]
denotes the \textit{Gaussian hypergeometric function}. If in addition 
$\Re(c) > \Re(a+b)$ then the series converges absolutely on $\T$ and satisfies 
${}_2 F_1(a, b, c; 1) = \frac{\Gamma(c)\Gamma(c-a-b)}{\Gamma(c-a)\Gamma(c-b)}$ 
(Gauss Summation Theorem). 
Recall that $S_{\C^n} : = \{w \in \C^n : 
\Vert w \Vert_{{\C}_2^n} = 1\} = J_2^{-1}(\S^{2n-1})$ denotes the unit sphere 
in $\C^n$ (where $n \in \N$, of course).
\begin{thmGI*}
Let $n \in \N$ and $u, v \in \S^{n-1}$. Let $\textbf{X} \equiv 
(X_1, \ldots, X_n)^\top \sim N_n(0, I_n)$ be a standard-normally distributed 
real Gaussian random vector. Then
\begin{align}\label{eq:Grothendieck}
\begin{split}
\E[{\text{sign}}(u^\top \textbf{X}){\text{sign}}(v^\top \textbf{X})] &= 
\frac{2}{\pi}\,\arcsin(u^\top v)\\
&= \frac{2}{\pi}\,u^\top v\,\,
{}_2 F_1(\tfrac{1}{2}, \tfrac{1}{2}, \tfrac{3}{2}; (u^\top v)^2)\\
&= \E[\vert X_1 \vert]^2\,u^\top v\,\,
{}_2 F_1(\tfrac{1}{2}, \tfrac{1}{2}, \tfrac{3}{2}; (u^\top v)^2).
\end{split}
\end{align}
\end{thmGI*}
\begin{thmHI*}
Let $n \in \N$ and
$\textbf{Z} \equiv (Z_1, \ldots, Z_n)^\top \sim {\C}N_n(0, I_n)$ be a 
standard-normally distributed complex Gaussian random vector. Then

\scalebox{0.9}{
\vbox{
\begin{align}\label{eq:Haagerup}
\begin{split}
\E[{\text{sign}}(u^\ast \textbf{Z}){\text{sign}}(\overline{v^\ast \textbf{Z}})] 
&= \frac{\pi}{4}\,{\text{sign}}(u^\ast v)\big(\frac{1}{\pi}\int_{0}^{2\pi}
\arcsin(\vert u^\ast v \vert\cos(t))\cos(t)\,\textup{d}t\big)\\
&= \frac{\pi}{4}\,\text{sign}(u^\ast v)\,\vert u^\ast v \vert\,\,
{}_2 F_1(\tfrac{1}{2}, \tfrac{1}{2}, 2; \vert u^\ast v \vert^2)\\
&= \E[\vert Z_1 \vert]^2\,u^\ast v\,\,{}_2 F_1(\tfrac{1}{2},
\tfrac{1}{2}, 2; \vert u^\ast v \vert^2).
\end{split}
\end{align}
}}

\end{thmHI*}
\begin{remark}\label{rem:integral_rep_of_Haagerup}
The second equality in \eqref{eq:Haagerup} is a particular case of the equality
\begin{align}\label{eq:integral_repr_of_the_Haagerup_function}
\frac{1}{\pi}\int_{0}^{2\pi}\arcsin(x\cos(t))\cos(t)\,\textup{d}t = x \,
{}_2 F_1(\tfrac{1}{2}, \tfrac{1}{2}, 2; x^2) \text{ for all } x \in [-1,1],
\end{align}
implied by the Maclaurin series representation of the function $\arcsin$ and the well-known 
fact that $\int_{0}^{2\pi}\cos^{2(n+1)}(t)\,\textup{d}t = \frac{2\sqrt{\pi}}{(n+1)!}
\Gamma(n+\frac{3}{2}) = \frac{2\sqrt{\pi}}{\Gamma(n+2)}(n+\frac{1}{2})
\Gamma(n+\frac{1}{2})$ for all $n \in \N_0$.
\end{remark}
\noindent We can see clearly that both, \eqref{eq:Grothendieck} and 
\eqref{eq:Haagerup} do not depend on the choice of the dimension $n$. The 
reason for this is \sref{Lemma}{lem:correl_Gaussians_II}, but not the use 
of the sign function. If we namely fix an arbitrary random vector 
$\textbf{W} \sim {\F}N_n(0, I_n)$ and consider the matrix $A_{u,v} : = 
\begin{pmatrix}
u_1 & u_2 & \ldots & u_n\\
v_1 & v_2 &\ldots & v_n
\end{pmatrix} \in \M_{2, n}(\F)$, then $(u^\ast \textbf{W}, v^\ast \textbf{W})^\top 
= A_{u,v} \textbf{W} \sim N_2(0, \Sigma_2(u^\ast v))$ (due to 
\sref{Lemma}{lem:correl_Gaussians_II}). Hence, if $(S_1, S_2)^\top \sim 
N_2(0, \Sigma_2(u^\ast v))$ is given, then $\P_{(S_1, S_2)^\top} = 
\P_{(u^\ast \textbf{W}, v^\ast \textbf{W})^\top} = (A_{u,v})_\ast\,\P_{\textbf{W}}$. 
The change of variables formula therefore implies that (in particular) 
\textit{for any} choice of a.\,e. bounded functions $f, g \in L^\infty(\F)$, the 
following equality holds:
\begin{align}\label{eq:GT_and_Haag_unified}
\E[f(u^\ast \textbf{W})g(\overline{v^\ast \textbf{W}})] =
\E[(f \otimes \overline{g})\circ A_{u,v}(\textbf{W})] =
\int_{\F^2}f \otimes \overline{g}\,\textup{d}((A_{u,v})_\ast\,\P_{\textbf{W}}) =
\E[f(S_1)g(\overline{S_2})].
\end{align}
 Consequently, if we also include \sref{Lemma}{lem:charact_of_corr_matrices} 
(or the obvious fact that the mapping $S_{\F^n} \times S_{\F^n}
\ni (u,v) \mapsto u^\ast v \in \overline{\D} \cap \F$ is onto for any $n \in 
\N_2$) and \sref{Corollary}{cor:little_GT_and_Gaussian_structure}, then 
\eqref{eq:GT_and_Haag_unified}, applied to $f : = g : = \text{sign}$ implies 
that \eqref{eq:Grothendieck} and \eqref{eq:Haagerup} are special cases of an 
equality which ``just'' involves the function $\text{sign} : \R \longrightarrow 
\{-1,1\}$ and a $2$-dimensional Gaussian random vector, where the latter consists 
of two arbitrarily correlated random variables, though. Remembering 
the fact \eqref{eq:the_2_dim_corr_matrix}, we obtain: 
\begin{corollary}[\textbf{Gaussian $\text{sign}$-correlation coefficient}]
\label{cor:unification_of_Grothendieck_and_Haagerup}
Fix $\F \in \{\R, \C\}$. Let $\Sigma \in C(2; \F)$ and 
$(S_1, S_2) \sim {\F}N_2(0, \Sigma)$. Then $\Sigma = \Sigma_2(\zeta)$ for some 
$\zeta \in \F \cap \overline{\D}$, and the Pearson correlation coefficient 
between ${\text{sign}}(S_1)$ and ${\text{sign}}(S_2)$ is given by
\begin{align}\label{eq:GT_and_HAG_unification}
\E[{\text{sign}}(S_1){\text{sign}}(\overline{S_2})] 
= \E[\vert S_1 \vert]^2\,\zeta\,{}_2 F_1(\tfrac{1}{2},\tfrac{1}{2}, 
\tfrac{d_{\F}+2}{2}; \vert \zeta \vert^2) = 
\frac{1}{k_G^\F}\,\zeta\,{}_2 F_1(\tfrac{1}{2},\tfrac{1}{2}, 
\tfrac{d_{\F}+2}{2}; \vert \zeta \vert^2),
\end{align}
where $d_{\R} : = 1$ and $d_{\C} : = 2$.
\end{corollary}
\section{Integration over $\S^{n-1}$ and the Gamma function}
In fact, the Grothendieck equality as well as the Haagerup equality 
instantly unfold  as a \textit{special case} of a (much more general) 
result, where we explicitly describe all non-negative integer powers of an 
expectation of inner products of suitably correlated - real - random 
vectors, uniformly distributed on the unit sphere (cf. \autoref{thm:3F2_rep} 
and \sref{Proposition}{prop:2F1_rep}). In this regard, we possibly should point to 
the so-called ``kernel trick'', used also for the computation of inner products 
in high-dimensional feature spaces using simple functions defined on pairs of 
input patterns which is a crucial ingredient of support vector machines in 
statistical learning theory; i.e., learning machines that construct decision 
functions of sign type. This trick allows the formulation of nonlinear variants 
of any algorithm that can be cast in terms of inner products 
(cf. \cite[Chapter 5.6]{V2000}). 

 Firstly, it is quite helpful to understand the actual source of the 
values $\frac{2}{\pi}$ and $\frac{\pi}{4}$ (cf. \eqref{eq:the_one_dim_case}, 
\sref{Corollary}{cor:little_GT_and_Gaussian_structure}, \sref{Proposition}{prop:2F1_rep} 
and \cite[Chapter 8.7]{DF1993}).
\begin{lemma}\label{lem:double_factorial_fact}
Let $(b_n)_{n \in \N}$ be the sequence of real numbers, defined as
\begin{align}\label{eq:the_constants_b_n}
b_n : = \begin{cases} 1&\text{if } n \text{ is even}\\ 
\sqrt{\pi/2}&\text{if } n \text{ is odd\,\,\,.}\end{cases}
\end{align}
Then
\begin{align}\label{eq:Gamma_of_n_half}
\Gamma\big(\frac{n}{2}\big) = \frac{(n-2)!!}{\sqrt{2^{n-2}}}\,b_n \,\,\, 
\text{ for all } n \in \N .
\end{align}
In particular, $b_n = b_{2m+n}$ for all $m, n \in \N$, and
\begin{align}\label{eq:one_half_version-of_GHFs}
{}_2 F_1\big(\frac{k}{2}, \frac{l}{2},\frac{m}{2}; z\big) 
= \frac{(m-2)!!}{(k-2)!!\,(l-2)!!}
\sum_{n=0}^\infty \frac{(2(n-1)+k)!!\,(2(n-1)+l)!!}{(2n)!!\,(2(n-1)+m)!!}\,z^n 
\end{align}
for all $k, l, m \in \N$ and $z \in \D$. If in addition $m > k+l$, then 
\eqref{eq:one_half_version-of_GHFs} holds for any $z \in \T$.
\end{lemma}
\begin{comment}
Regarding the proof of \eqref{eq:Gamma_of_n_half}, we have to distinguish two 
cases; namely the even case (i.e., $n = 2l$ for some $l \in \N$) and the 
odd case (i.e., $n = 2k+1$ for some $k \in \N_0$). However, 
since $(2l - 2)!! = (2(l-1))!! = 2^{l-1}(l-1)!$ and $\Gamma(k + \frac{1}{2}) 
= \frac{(2k-1)!!}{2^{k}}\,\sqrt{\pi}$, \eqref{eq:Gamma_of_n_half} follows at once.
Equipped with \eqref{eq:Gamma_of_n_half}, the proof of the representation 
\eqref{eq:one_half_version-of_GHFs} just involves a few remaining basic algebraic 
transformations (including a multiple shortening of fractions), implied by the 
definition of Gaussian hypergeometric functions. 
\end{comment} 
\noindent We also need results about the real and complex Gaussian randomness 
structure, embedded in the Gamma function, which are of their own interest; 
built on an important link between the Gamma function and powers of absolute 
moments of standard normally distributed real random variables. To this end, 
we consider both, the real and the complex unit sphere as a probability space. 
Put
\begin{align*}
\begin{split}
\sigma_{n-1}(A) &: = 
\frac{\omega_n(A)}{\omega_n} = 
\frac{\Gamma(n/2)}{2 \pi^{n/2}}\,\omega_n(A) = 
\frac{\Gamma(n/2)}{2 \pi^{n/2}}\,n\,\lambda_n(\{r\xi: 0 < r \leq 1 \text{ and } 
\xi \in A\})\\
&\,\,= \frac{\Gamma(\tfrac{n}{2}+1)}{\pi^{n/2}}\,
\lambda_n(\{r\xi: 0 < r \leq 1 \text{ and } \xi \in A\}),
\end{split}
\end{align*}
where $A \in {\mathcal{B}}(\S^{n-1})$ and $\omega_n \equiv \omega_n(\S^{n-1}) = 
\frac{2 \pi^{n/2}}{\Gamma(n/2)}$ denotes the surface area of the unit sphere 
$\S^{n-1} \subseteq \R^n$. $\sigma_{n-1}$ denotes the rotation-invariant 
probability measure (Haar measure) on $\S^{n-1}$. Moreover, $\sigma_n^{\C} := 
(J_2^{-1})_\ast\,\sigma_{2n-1}$ denotes the surface area probability measure on 
the complex unit sphere $S_{\C^n}$. 
\begin{proposition}\label{prop:Gamma_fct_and_integration_over_spheres}
Let $n \in \N$, $X \sim N_1(0,1)$, $\textbf{X} \sim N_n(0, I_n)$, 
$Z \sim {\C}N_1(0,1)$, $\textbf{Z} \sim {\C}N_n(0, I_n)$ and 
$\textbf{Y} =\vc{\textbf{Y}_1}{\textbf{Y}_2} \sim N_{2n}(0, I_{2n})$.
Let $p,q \in \R$ such that $p > -1$ and $q > 0$. Then
\begin{enumerate}
\item
\begin{align}\label{eq:Gamma_function_abs_N_0_1_integral}
2\int_{0}^{\infty} s^p\,\gamma_1(\textup{d}s) = \E[{\vert X \vert}^p] = 
\frac{2^{p/2}}{\sqrt{\pi}}\,\Gamma(\tfrac{p+1}{2})
\end{align}
and
\[
\Gamma(q) = \frac{\sqrt{2\pi}}{2^{q}}\,\E[{\vert X \vert}^{2q-1}].
\]
In particular,
\begin{align}\label{eq:absolute_Gaussian_moments}
\E[{\vert X \vert}^k] = \frac{(k-1)!!}{b_k} = \begin{cases}(k-1)!! & 
\text{if } k \text{ is even}\\ 
\sqrt{\frac{2}{\pi}}\,(k-1)!! & \text{if } k \text{ is odd}\end{cases} 
\hspace{1em}\text{ for all } k \in \N_0\,,
\end{align}
where $b_k$ satisfies \eqref{eq:the_constants_b_n}.
\item
Let $n \geq 2$ and $f : \R^n \longrightarrow \R$, such that 
\[
f(rx) = r^p\,f(x) \text{ for all } (r,x) \in (0, \infty) \times \R^n\,.
\]
Then $f \in L^1(\R^n, \gamma_n)$ if and only if $f\big\vert_{\S^{n-1}} \in 
L^1(\S^{n-1},\sigma_{n-1})$, and
\begin{align}\label{eq:an_important_integral}
\E[f(\textbf{X})] &= \int_{\R^n}f(x)\gamma_n(\textup{d}x)\, = \, 2^{p/2}\,
\frac{\Gamma\big(\frac{n+p}{2}\big)}{\Gamma\big(\frac{n}{2}\big)}
\int\limits_{{\S}^{n-1}} f(\xi)\,\textup{d}\sigma_{n-1}(\xi)\,. 
\end{align}
\item Let $b : \C^n \longrightarrow \C$, such that 
\[
b(rz) = r^p\,b(z) \text{ for all } (r,z) \in (0, \infty) \times \C^n\,.
\]
Then $b \in L^1(\C^n, \gamma_n^{\C})$ if and only if $\Re(b)\circ 
J_2^{-1}\big\vert_{\S^{2n-1}} \in L^1(\S^{2n-1},\sigma_{2n-1})$ and 
$\Im(b)\circ J_2^{-1}\big\vert_{\S^{2n-1}} \in L^1(\S^{2n-1},\sigma_{2n-1})$, 
and
\begin{align}\label{eq:an_important_complex_integral}
\begin{split}
\E[b(\textbf{Z})] &= \int_{\C^n}b(z)\gamma_n^{\C}(\textup{d}z) = 
\frac{\Gamma(n + \tfrac{p}{2})}{(n-1)!} \int\limits_{S_{\C^n}} b(\zeta)\,
\textup{d}{\sigma}_n^{\C}(\zeta)\\
&= \frac{\Gamma(n + \tfrac{p}{2})}{(n-1)!}\Big(\int\limits_{{\S}^{2n-1}} 
\Re(b(y_1 + i\,y_2))\,\textup{d}\sigma_{2n-1}((y_1,y_2)) + i \int\limits_{{\S}^{2n-1}} 
\Im(b(y_1 + i\,y_2))\,\textup{d}\sigma_{2n-1}((y_1,y_2))\Big)\\
&= 2^{-p/2}\big(\E[\Re(b(\textbf{Y}_1 + i\,\textbf{Y}_2))] + 
i\,\E[\Im(b(\textbf{Y}_1 + i\,\textbf{Y}_2))]\big).
\end{split}
\end{align}
\end{enumerate}
\end{proposition}
\begin{corollary}\label{cor:the_random_norm_case}
Let $p \in (-1, \infty)$, $m, n \in \N$, $X \sim N_1(0,1)$, $\textbf{X} 
\sim N_n(0, I_n)$, $\textbf{Y} \sim N_{2n}(0, I_{2n})$, $Z \sim {\C}N_1(0,1)$ 
and $\textbf{Z} \sim {\C}N_n(0, I_n)$. Then
\begin{align}\label{eq:expected_l_2_norm_of_real_Gaussian}
\E[\Vert \textbf{X} \Vert_{\R_2^n}^p ] = 2^{p/2}\,
\frac{\Gamma(\frac{n+p}{2})}{\Gamma(\frac{n}{2})} =
\frac{\E[\vert X \vert^{n-1+p}]}{\E[\vert X \vert^{n-1}]} 
\end{align}
and
\begin{align}\label{eq:expected_l_2_norm_of_complex_Gaussian}
\E[{\Vert\textbf{Z}\Vert}_{{\C}_2^n}^p] = \frac{\Gamma(n + \frac{p}{2})}
{(n-1)\,!}  = 2^{-p/2}\,\E[{\Vert\textbf{X}\Vert}_{\R_2^{2n}}^p] = 
2^{-p/2}\,\frac{\E[\vert X \vert^{2n-1+p}]}{\E[\vert X \vert^{2n-1}]}\,. 
\end{align}
In particular,
\begin{align*}
\E[\Vert \textbf{X} \Vert_{\R_2^n}^m ] = a_n(m)\,\frac{(n-2+m)!!}{(n-2)!!}\,\, 
\text{ and } \,\, \E[{\Vert\textbf{Z}\Vert}_{{\C}_2^n}^m] = a_{2n}(m)\,2^{-m/2}\,
\frac{(2n-2+m)!!}{(2n-2)!!}
\end{align*}
and
\begin{align}\label{eq:the_one_dim_case}
\E[\vert X\vert] = \sqrt{\frac{2}{\pi}}\,\text{ and }\,\E[\vert Z\vert] = 
\sqrt{\frac{\pi}{4}} = \frac{\sqrt{\pi}}{2}\,,
\end{align}
where 
\begin{align}\label{eq:toggle_switch_function} 
a_n(m) : = \frac{b_{n+m}}{b_n} = \begin{cases}\sqrt{\pi/2} & 
\text{if } n \text{ is even and } m \text{ is odd}\\ 
\sqrt{2/\pi} & \text{if } n \text{ is odd and } m \text{ is odd}\\
1 & \text{if } m \text{ is even}
\end{cases}
\end{align}
and $b_n$ is defined as in \sref{Lemma}{lem:double_factorial_fact}.
\end{corollary}
\begin{comment}
We just have to apply \sref{Proposition}{prop:Gamma_fct_and_integration_over_spheres} 
to the function $\F^n \ni z \mapsto f_p(z) : = \Vert z \Vert^p$ (and to recall 
that by definition $\Vert \cdot \Vert_{\R_2^1} := \vert \cdot \vert$). 
\end{comment}
\noindent A further, very important special case (which allows an easy 
proof of \autoref{thm:3F2_rep}) arises if we consider the function 
$\R^n \ni x \mapsto \langle u, x\rangle_{l^2_n}^m = (u^\top x)^m$, where 
$m \in \N_0$ and $u \in \S^{n-1}$ are given.
\begin{corollary}\label{cor:spherical_integration_of_inner_products}
Let $p \in (-1, \infty), n \in \N_2, x \in \S^{n-1}$ and 
$Y \sim N_1(0,1)$. Let $f : \R \longrightarrow \R$ be a function such that 
$f \in L^1(\R, \gamma_1)$ and $f(ry) = r^p\,f(y)$ for all $(r,y) \in 
(0, \infty) \times \R$. Then $f(x^\top\cdot) \in L^1(\S^{n-1},\sigma_{n-1})$, 
and
\[
\int\limits_{{\S}^{n-1}} f(x^\top\,u)\,\sigma_{n-1}(\textup{d}u) = 
2^{-p/2}\frac{\Gamma(\frac{n}{2})}{\Gamma(\frac{n+p}{2})}\,\E[f(Y)]\,. 
\]
In particular,
\begin{align}\label{eq:the_power_case}
\int\limits_{{\S}^{n-1}} (x^\top\,u)^{m}\,\sigma_{n-1}(\textup{d}u) = 
\frac{1+(-1)^m}{2}\,\frac{\Gamma(\tfrac{m+1}{2})
\Gamma(\tfrac{n}{2})}{\sqrt{\pi}\,\Gamma(\tfrac{m+n}{2})} = 
\frac{1+(-1)^m}{2}\big(\frac{\Gamma(\frac{n}{2})}{\sqrt{\pi}\,\Gamma(\frac{n-1}{2})}\,
B(\tfrac{m+1}{2}, \tfrac{n-1}{2})\big)
\end{align}
for all $m \in \N_0$. Here, $(0, \infty) \times (0, \infty) \ni (x_1,x_2)
\mapsto B(x_1,x_2) : = \int_{0}^{1} t^{x_1-1}\,(1-t)^{x_2-1}\,\textup{d}t = 
\frac{\Gamma(x_1)\,\Gamma(x_2)}{\Gamma(x_1 + x_2)}$ denotes the real beta 
function. 
\end{corollary}
\begin{comment}
Since $\Vert x \Vert_{l^2_n} = 1$, it follows that $x^\top\,\textbf{Y} 
\stackrel{d}{=} Y \sim N_1(0,1)$ for any $\textbf{Y} \sim N_n(0, I_n)$. Thus, 
$\E[f(Y)] = \E[f(x^\top\,\textbf{Y})]$, so that we may apply 
\eqref{eq:an_important_integral} to the function 
$\R^n \ni x \mapsto f(x^\top\,\cdot)$. Regarding the 
particular case $\R \ni y \mapsto y^m$ (respectively, $\R^n \ni x \mapsto 
(x^\top\,\cdot)^m$), we only have to include \sref{Lemma}{lem:double_factorial_fact} 
and the well-known fact that the $m$-th moment of $Y \sim N_1(0,1)$ satisfies 
$\E[Y^m] = \frac{1+(-1)^m}{2}\,(m-1)!!$\,.
\end{comment}
\begin{lemma}\label{lem:Gamma_fct_and_constants_c_k}
Consider the sequence 
$(c_k)_{k \in \N}$, defined as
\[
c_k : = \frac{1}{\sqrt{{}_2 F_1(\tfrac{1}{2}, \tfrac{1}{2}, 
\frac{k+2}{2}; 1)}}\,.
\]
Let $\textbf{X} \sim N_k(0, I_k)$ and $\textbf{Z} \sim {\C}N_k(0, I_k)$. Then
\begin{align}\label{eq:the_sequence_lb_c_k_rb}
c_k = \sqrt{\frac{2}{k}}\,\frac{\Gamma(\frac{k+1}{2})}
{\Gamma(\frac{k}{2})} = a_k\,\frac{1}{\sqrt{k}}\,
\frac{(k-1)!!}{(k-2)\,!!} = \frac{1}{\sqrt{k}}\,
\E[\Vert \textbf{X} \Vert_{\R_2^k}] = \sqrt{\frac{2\pi}{k}}\,
\frac{\omega_{k}}{\omega_{k+1}}\,,
\end{align}
where $a_k : = a_k(1)$ satisfies \eqref{eq:toggle_switch_function}, $\omega_1 := 2$ 
and $\omega_m$ denotes the surface area of the unit sphere $\S^{m-1}$ $(m \in \N_2)$. 
In particular, $c_1^2 =\frac{2}{\pi}$, $c_2^2 =\frac{\pi}{4}$ and
\begin{align}\label{eq:the_sequence_lb_c_2k_rb}
c_{2k} = \frac{1}{\sqrt{k}}\,\frac{\Gamma(k+\frac{1}{2})}
{(k-1)!} = \frac{1}{\sqrt{k}}\,\sqrt{\frac{\pi}{4}}\,
\frac{(2k-1)!!}{(2k-2)\,!!} = \frac{1}{\sqrt{k}}\,
\E[\Vert \textbf{Z} \Vert_{\C_2^k}]\,.
\end{align}
Moreover, $0 < c_k < 1$ for all $k \in \N$, $\lim\limits_{k \to \infty} c_k = 1$, 
and
\[
\frac{1}{\sqrt{k}}\,\E[{\Vert\textbf{W}\Vert}_{{\F}_2^k}] = c_{\nu_k^\F} =
\frac{1}{\sqrt{{}_2 F_1(\tfrac{1}{2}, \tfrac{1}{2}, \frac{\nu_k^\F+2}{2}; 1)}} 
\stackrel{k \to \infty}{\longrightarrow} 1 \,\, \text{ for all } \textbf{W} 
\sim {\F}N_k(0, I_k)\,,
\]
where $\nu_k^\R : = k$ and $\nu_k^\C : = 2k$.
\end{lemma}
\noindent In the context of \autoref{thm:GT_Niemi_PSD_case}, the 
one-dimensional special cases of \eqref{eq:expected_l_2_norm_of_real_Gaussian} 
and \eqref{eq:expected_l_2_norm_of_complex_Gaussian}
disclose a unification of the real and complex Gaussian structure, encoded 
at least in the \textit{little} Grothendieck constant: 
\begin{corollary}\label{cor:little_GT_and_Gaussian_structure}
Let $\F \in \{\R, \C\}$. Then the little Grothendieck constant $k_G^\F$ can be written as
\[
k_G^\F = 
{}_2 F_1(\tfrac{1}{2}, \tfrac{1}{2}, \tfrac{d_\F+2}{2}; 1)
= \begin{cases}\pi/2&\text{if } \F = \R \text{ and } d_\R = 1\\ 
4/\pi &\text{if } \F = \C \text{ and } d_\C = 2\end{cases}.
\]
\end{corollary}
\begin{corollary}[\textbf{Krivine, 1979}]\label{cor:Krivine_s_integral_rep}
Let $f \in C([-1,1])$, $k \in \N_2$ and $u \in \S^{k-1}$. Then
\[
\int_{\S^{k-1}}f(\langle u, v \rangle_{\R_2^k})\textup{d}\sigma_{k-1}(v) = 
\int_{-1}^{1} f(t)\,\textup{d}\Q_k(t) = \sqrt{\frac{k-1}{2\pi}}\,c_{k-1}
\int_{-1}^{1} f(t)(1-t^2)^{\frac{k-3}{2}}\,\textup{d}t\,,
\]
where $\Q_k(\textup{d}s) : = \frac{\Gamma(k/2)}{\sqrt{\pi}\,\Gamma((k-1)/2)}\,
(1-t^2)^{\frac{k-3}{2}}\,\textup{d}s$ is a probability measure on $[-1,1]$.
\end{corollary}
\begin{comment}
We just have to link \eqref{eq:the_sequence_lb_c_k_rb} with the (unindexed) 
equality on page 27 of \cite{Kr1979}.
\end{comment}
\begin{remark}[\textbf{Absolutely $p$-summing operators and GT in matrix form}]
\label{rem:GT_vs_abs_summing_operators}
Recall that $T \in \mathfrak{L}(E,F)$ between Banach spaces $E$ and $F$ is called 
absolutely $p$-summing ($1 \leq p < \infty$) if there exists a constant $c \geq 0$ 
such that for all $n \in \N$ and $x_1, \ldots, x_n \in E$
\begin{align}\label{eq:absolutely_p_summing_operator}
\big(\sum_{i=1}^n \Vert Tx_i\Vert^p\big)^{1/p} \leq c\,w_p(x_1, \ldots, x_n)\,,
\end{align}
where $w_p(x_1, \ldots, x_n) : = \sup\limits_{\psi \in B_{E^\prime}}
\big(\sum_{i=1}^n \vert \langle x_i, \psi\rangle \vert^p\big)^{1/p}$. The 
$p$-summing norm $\Vert T \Vert_{{\mathfrak{P}}_p}$ is defined as the infimum of 
all constants $c \geq 0$ which satisfy \eqref{eq:absolutely_p_summing_operator} 
(cf., e.g., \cite[Chapter 11]{DF1993} or \cite[Chapter 2]{DJT1995}). Expressed 
in the terminology of absolutely $1$-summing operators, Grothendieck proved that 
his inequality in particular is \textit{equivalent} to
\begin{align}\label{eq:GT_and_absolutely_1_summing_operators}
\Vert T \Vert_{{\mathfrak{P}}_1} \leq K_G^\F\,\Vert T \Vert  
\end{align}
for all $\F$-Hilbert spaces $H$, $n \in \N$ and finite rank operators 
$T \in {\mathfrak{L}}(l_1^n, H)$ (cf. \cite{M1973, P1985} and 
\cite[Theorem 10.7]{J1987}). Actually, the proof of 
\eqref{eq:GT_and_absolutely_1_summing_operators} in the finite rank case 
is quite simple. It is based on the following two facts. Firstly, since 
$(l_1^n)^\prime \cong l_\infty^n$, it follows that for all $a_1, \ldots, a_m \in 
l_1^n$ ($a_i \equiv (a_{i1},\ldots,a_{in})^\top$), $\Vert A \Vert_{\infty, 1} = 
w_1(a_1, \ldots, a_m)$ and $\Vert A^\ast A \Vert_{\infty, 1} = 
(w_2(a_1, \ldots, a_m))^2$, where 
\[
A : = (a_1 \,\brokenvert\, a_2 \,\brokenvert\,\cdots\,\brokenvert\,a_m)^\top 
\equiv (a_{ij}) \in \M_{m,n}(\F).
\]
Secondly, since $H^\prime \cong H$ (Riesz), we obtain that 
\[
\sum_{i=1}^m \Vert T\,a_i\Vert_H = \sum_{i=1}^m \big\Vert \sum_{j=1}^n a_{ij}\,
T e_j\big\Vert_H = \text{tr}(A^\ast\Gamma_H(u,z_T)) \leq K_G^\F\,\Vert T \Vert \,
\Vert A \Vert_{\infty, 1} = K_G^\F\,\Vert T \Vert \,w_1(a_1, \ldots, a_m),
\]
for some $u \in B_H^m$ and $z_T \in H^n$, where the latter is defined as 
$(z_T)_j : = T e_j \in \Vert T \Vert\,B_H$ ($j \in [n]$).
Similarly, we obtain the well-known ${\mathfrak{P}}_2$-representation of 
the little Grothendieck inequality (``little GT''), which even is 
\textit{equivalent} to little GT in matrix form, since the use 
of $({\mathfrak{P}}_2, \Vert \cdot \Vert_{{\mathfrak{P}}_2})$ in fact naturally 
implies the emergence of the \textit{positive semidefinite} matrix 
$A^\ast A \in \M_{n}(\F)^+$:
\begin{align*}
\sum_{i=1}^m \Vert T\,a_i\Vert_H^2 &= 
\sum_{l=1}^n\sum_{j=1}^n \big(\sum_{i=1}^m \overline{a_{il}}\,a_{ij}\big)
\big\langle (z_T)_j, (z_T)_l\big\rangle_H = \text{tr}((A^\ast A)\Gamma_H(z_T,z_T))\\
&\leq k_G^\F\,\Vert T \Vert^2 \,\Vert A^\ast A \Vert_{\infty, 1} = 
k_G^\F\,\Vert T \Vert^2 (w_2(a_1, \ldots, a_m))^2,
\end{align*}
In other words,
\begin{align}\label{eq:little_GT_and_absolutely_2_summing_operators}
\Vert T \Vert_{{\mathfrak{P}}_2} \leq \sqrt{k_G^\F}\,\Vert T \Vert 
\end{align}
for all $\F$-Hilbert spaces $H$, $n \in \N$ and (finite rank operators) 
$T \in {\mathfrak{L}}(l_1^n, H)$. It is surprising that no attention seems to 
have been paid to the equivalence of the psd matrix version of little GT and 
the absolutely $2$-summing version of little GT so far. So, its worth to state 
it here. Moreover, if we combine \cite[Theorem 11.10]{DF1993}, 
\eqref{eq:expected_l_2_norm_of_real_Gaussian} and 
\eqref{eq:expected_l_2_norm_of_complex_Gaussian}, we can somewhat simplify the
representation of the $1$-summing norm of $Id_{\F_2^k}$ for 
any $k \in \N$. We namely have
\[
\Vert Id_{\R_2^k} \Vert_{{\mathfrak{P}}_1} = \sqrt{\frac{\pi}{2}}\,
\E[{\Vert\textbf{X}\Vert}_{\R_2^k}] = \sqrt{\frac{\pi}{2}}\,
a_k\,\frac{(k-1)!!}{(k-2)!!} \stackrel{\eqref{eq:the_sequence_lb_c_k_rb}}{=} 
\sqrt{k}\big(\sqrt{\frac{\pi}{2}}\,c_k\big) \,\, \text{ for all } k \in \N
\]
and
\[
\Vert Id_{\C_2^k} \Vert_{{\mathfrak{P}}_1} = \sqrt{\frac{4}{\pi}}\,
\E[{\Vert\textbf{Z}\Vert}_{\C_2^k}] = \frac{2}{\pi}\,\Vert Id_{\R_2^{2k}}\Vert_
{{\mathfrak{P}}_1} = \frac{(2k-1)!!}{(2k-2)!!} 
\stackrel{\eqref{eq:the_sequence_lb_c_2k_rb}}{=} 
\sqrt{k}\big(\sqrt{\frac{4}{\pi}}\,c_{2k}\big) \,\, \text{ for all } k \in \N\,,
\]
where $a_k : = a_k(1)$ satisfies \eqref{eq:toggle_switch_function}. Consequently, 
\sref{Lemma}{lem:Gamma_fct_and_constants_c_k} recovers \cite[Proposition 8.8]{J1987}, 
respectively \cite[Corollary 11.10]{DF1993} and reveals a link between norms of 
certain absolutely $1$-summing operators and values of Gaussian hypergeometric 
functions (since $\lim\limits_{\nu \to \infty}c_\nu = 1$):
\[
\lim\limits_{k \to \infty}\frac{1}{k}\Vert Id_{\F_2^k} \Vert_{{\mathfrak{P}}_1}^2 =
k_G^\F = {}_2 F_1\big(\frac{1}{2}, \frac{1}{2}, \frac{d_{\F} + 2}{2}; 1\big).
\]
That is,
\[
\lim\limits_{k \to \infty}\frac{1}{k}\Vert Id_{\R_2^k} \Vert_{{\mathfrak{P}}_1}^2 =
\frac{\pi}{2} = {}_2 F_1\big(\frac{1}{2}, \frac{1}{2}, 
\frac{3}{2}; 1\big)\,\text{ and }\,\lim\limits_{k \to \infty}\frac{1}{k}
\Vert Id_{\C_2^k} \Vert_{{\mathfrak{P}}_1}^2 = \frac{4}{\pi} = 
{}_2 F_1\big(\frac{1}{2}, \frac{1}{2}, 2; 1\big),
\]
implying a further facet of a link between (Euclidean norms of) Gaussian random 
vectors and the $1$-Banach ideal of absolutely $1$-summing operators.
\end{remark}
\noindent \sref{Proposition}{prop:Gamma_fct_and_integration_over_spheres} also
allows us to give a straightforward, simple proof of the following important 
well-known surface integral characterisation of the trace of a matrix: 
\begin{corollary}\label{cor:trace_as_surface_integral}
Let $\F \in \{\R, \C\}$, $m, n \in \N$, $A \in \M_{m,m}(\F), B \in \M_{m,n}(\F), 
C \in \M_{n,m}(\F)$ and $D \in \M_{n,n}(\F)$. Let $\textbf{Z} \sim 
{\F}N_m(0, I_m)$. 
Then the following properties hold:
\begin{enumerate}
\item
\[
{\text{tr}}(\Re(A)) + i\,{\text{tr}}(\Im(A)) = {\text{tr}}(A) = 
\E[f_A(\textbf{Z})], 
\] 
where $\F^m \ni z \mapsto f_A(z) : = z^\ast Az = {\text{tr}}(A zz^\ast)$.
\item 
\[
{\text{tr}}(A) = m\big(\int\limits_{{\S}^{m-1}} u^\top \Re(A)u\,
\textup{d}\sigma_{m-1}(u) + i\,\int\limits_{{\S}^{m-1}} u^\top \Im(A)u\,
\textup{d}\sigma_{m-1}(u)\big).
\]
In particular,
\[
\int\limits_{{\S}^{m+n-1}} \text{tr}(B \Gamma_\R(u,v))\,\textup{d}
\sigma_{m+n-1}(\vc{u}{v}) = \int\limits_{{\S}^{m+n-1}} u^\top Bv\,\textup{d}
\sigma_{m+n-1}(\vc{u}{v}) = 0.
\]
\end{enumerate}
\end{corollary}
\begin{comment}
(i) Since $\textbf{Z} \sim {\F}N_m(0, I_m)$, it follows that 
$\E[\textbf{Z}\textbf{Z}^\ast] = I_m$. Consequently, the linearity of
$\E$ implies that  
\begin{align*}
{\text{tr}}(A) = {\text{tr}}(A \E[\textbf{Z}\textbf{Z}^\ast]) =
\E[{\text{tr}}(A\textbf{Z}\textbf{Z}^\ast)] = 
\E[{\text{tr}}(\textbf{Z}^\ast A\textbf{Z})] = \E[f_A(\textbf{Z})].
\end{align*}
\\[0.2em]
\noindent (ii)  Obviously, we may assume that $A \in \M_m(\R)$. Since then
$f_A(rx) = r^2 f_A(x)$ for all $(r,x) \in (0, \infty) \times \R^m$, we may apply
\sref{Proposition}{prop:Gamma_fct_and_integration_over_spheres} to $p =2$, 
implying that 
\[
{\text{tr}}(A) = 2\,\frac{\Gamma\big(\tfrac{m}{2}+1\big)}
{\Gamma\big(\frac{m}{2}\big)}\int\limits_{{\S}^{m-1}} f_A(u)
\,\textup{d}\sigma_{m-1}(u) = m\int\limits_{{\S}^{m-1}} 
u^\top Au\,\textup{d}\sigma_{m-1}(u)
\]
(since $\Gamma\big(\tfrac{m}{2}+1\big) = \tfrac{m}{2}\, 
\Gamma\big(\tfrac{m}{2}\big)$).
Finally, if we put $w : = \vc{u}{v}$ and 
$\Delta(B) : = \frac{1}{2}\begin{pmatrix}
0 & B\\
B^\top & 0
\end{pmatrix} \in \M_{m+n, m+n}(\R)$, it obviously follows that 
$w^\top \Delta(B)w = u^\top Bv$ and $\text{tr}(\Delta(B)) = 0$. Consequently, we
obtain
\[
0 = \text{tr}(\Delta(B)) = (m+n)\int\limits_{{\S}^{m+n-1}} 
u^\top Bv\,\textup{d}\sigma_{m+n-1}(w).
\]  
\end{comment}
\section{Integrating powers of inner products of random vectors, uniformly 
distributed on $\S^{n-1}$}
Let us recall the (real, respectively complex) correlation matrices
\[
\Sigma_{2n}(z) : = 
\begin{pmatrix}
    I_n & z\,I_n \\
    \overline{z}\,I_n & I_n
  \end{pmatrix}
=
	\begin{pmatrix}
	1 & z\\
	\overline{z} & 1
	\end{pmatrix}
	\otimes
	I_n\,,
\]
where $\vert z \vert \leq 1$ and $n \in \N$ (cf. \eqref{eq:GT_corr_matrix} 
and \sref{Proposition}{prop:structure_of_Sigma_2n_zeta}). 
Moreover, if $\textbf{X} \sim N_n(0, I_n)$ and $A \in \mathcal{B}(\S^{n-1})$ is
an arbitrary Borel subset of the unit sphere $\S^{n-1} \subseteq \R^n$ $(n \geq 2)$, 
\eqref{eq:an_important_integral}, applied to the function $\R^n\setminus\{0\} \ni 
x \mapsto f_A(x) : = \ind_A(\frac{\!\!\!\!\!\!x}{\Vert x \Vert_{\R_2^n}})$ (and 
$p = 0$) implies that in particular
\[
\P\big(\frac{\!\!\!\!\!\textbf{X}}{\Vert\textbf{X}\Vert_{\R^n}} \in A\big) = 
\E[f_A(\textbf{X})] = \int\limits_{{\S}^{n-1}} \ind_A(\xi)\,\textup{d}\sigma_{n-1}(\xi) =
\sigma_{n-1}(A)
\]
(since $\textbf{X} \not= 0$ $\P$-a.s.).
Thus, we get again the well-known fact that $\frac{\!\!\!\!\!\!\textbf{X}}
{\Vert\textbf{X}\Vert_{\R^n}}$ is uniformly distributed on the unit sphere 
$\S^{n-1}$ if $\textbf{X} \sim N_n(0, I_n)$. Note that, in addition, for 
any $X \sim N_1(0,1)$, it is true that
\[
\P\big(\frac{X}{\vert X \vert} = \varepsilon \big) = \P(X > 0 \text{ and } 
\varepsilon = 1) + \P(X \leq 0 \text{ and } 
\varepsilon = -1) = \frac{1}{2}\,\text{ for all }\,\varepsilon \in \{-1,1\}\,,
\]
so that we could also say that the random variable $\frac{X}{\vert X \vert}$ is 
uniformly distributed on the ``unit sphere'' $S^0 : = \{-1,1\} \subseteq \R^1$ 
if $X \sim N_1(0,1)$.
\begin{theorem}\label{thm:3F2_rep}
Let $m, n \in \N$, $\rho \in (-1,1)$ and $\vc{\textbf{X}}{\textbf{Y}} = 
{\text{vec}}(X_1, \ldots, X_n, Y_1, \ldots, Y_n) \sim 
N_{2n}\big(0, \Sigma_{2n}(\rho)\big)$. 
\begin{enumerate}
\item If $m$ is odd then 

\scalebox{0.9}{
\vbox{
\[
\E\big[\big\langle\frac{\!\!\!\!\!\textbf{X}}{\Vert\textbf{X}\Vert_{\R_2^n}},
\frac{\!\!\!\!\!\textbf{Y}}{\Vert\textbf{Y}\Vert_{\R_2^n}}
\big\rangle^m_{\R_2^n} \big] = c_{\text{odd}}(m,n)\,\rho\,(1-\rho^2)^{\frac{n}{2}}\,
{}_3 F_2\big(\frac{n+1}{2}, \frac{n+1}{2}, \frac{m+2}{2}; \frac{3}{2}, 
\frac{m+n+1}{2}; \rho^2\big)\,,
\]
}}

where
\[
c_{\text{odd}}(m,n) : = 
\frac{m}{\sqrt{\pi}}\,\frac{\Gamma^2(\frac{n+1}{2})\,
\Gamma(\frac{m}{2})}{\Gamma(\frac{n}{2})
\Gamma(\frac{m+n+1}{2})}\,. 
\]
\item If $m$ is even then

\scalebox{0.9}{
\vbox{
\[
\E\big[\big\langle\frac{\!\!\!\!\!\textbf{X}}{\Vert\textbf{X}\Vert_{\R_2^n}},
\frac{\!\!\!\!\!\textbf{Y}}{\Vert\textbf{Y}\Vert_{\R_2^n}}\big\rangle^m_{\R_2^n}\big] = 
c_{\text{even}}(m,n)\,(1-\rho^2)^{\frac{n}{2}}\,{}_3 F_2\big(\frac{n}{2}, \frac{n}{2}, 
\frac{m+1}{2}; \frac{1}{2}, \frac{m+n}{2}; \rho^2\big)\,,
\]
}}

where
\[
c_{\text{even}}(m,n) : = 
\frac{1}{\sqrt{\pi}}\,\frac{\Gamma(\frac{n}{2})\,\Gamma(\frac{m+1}{2})}
{\Gamma(\frac{m+n}{2})}\,.
\] 
\end{enumerate}
In particular,
\begin{align}\label{eq:special_non_power_case}
\begin{split}
\E\big[\big\langle\frac{\!\!\!\!\!\textbf{X}}{\Vert\textbf{X}\Vert_{\R_2^n}},
\frac{\!\!\!\!\!\textbf{Y}}{\Vert\textbf{Y}\Vert_{\R_2^n}}\big\rangle_{\R_2^n}\big] 
&= c_n^2\,(1-\rho^2)^{\frac{n}{2}}\,\rho\,{}_2 F_1\big(\frac{n+1}{2}, 
\frac{n+1}{2}; \frac{n+2}{2}; \rho^2\big)\\
&= c_n^2\,\rho\,{}_2 F_1\big(\frac{1}{2}, \frac{1}{2}; \frac{n+2}{2}; \rho^2\big)\,,
\end{split}
\end{align}
where $c_n : = c_{\text{odd}}(1,n) = \sqrt{\frac{2}{n}}\,
\frac{\Gamma(\frac{n+1}{2})}{\Gamma(\frac{n}{2})}$.
\end{theorem}
\begin{remark}
The special case \eqref{eq:special_non_power_case} is contained in (the proof 
of) \cite[Lemma 2.1]{BOFV2014}. 
\end{remark}
\begin{remark}
Let $p \in (-1, \infty)$ and $f \in L^2(\gamma_1)$, satisfying $f(ry) = r^p\,f(y)$ 
for all $(r,y) \in (0, \infty) \times \R$. A natural question is, whether 
$I_f : = \E\big[\big\langle f\big(\frac{\!\!\!\!\!\textbf{X}}{
\Vert\textbf{X}\Vert_{\R_2^n}},\frac{\!\!\!\!\!\textbf{Y}}{
\Vert\textbf{Y}\Vert_{\R_2^n}}\big\rangle\big)_{\R_2^n} \big]$ can be similarly 
represented as $I = I_{f_m}$? In general, it seems that a closed-form representation 
of $I_f$ is not possible. However, our proof of \autoref{thm:3F2_rep} 
clearly reveals that
\begin{align*}
I_f &= \frac{\sqrt{\pi}}{2^{p/2}\,\Gamma(n/2)}\,(1-\rho^2)^{n/2}\,
\sum_{\nu=0}^\infty \frac{\Gamma^2\big(\frac{\nu+n}{2}\big)}{2^{\nu/2}\,
\Gamma(\frac{\nu+1}{2})\Gamma(\frac{\nu+2}{2})\Gamma(\frac{\nu+p+n}{2})}\,
\E[f(X) X^\nu]\,\rho^\nu\\
&= \E\big[f(X)\,\frac{\sqrt{\pi}}{2^{p/2}\,\Gamma(n/2)}\,(1-\rho^2)^{n/2}\,
\sum_{\nu=0}^\infty \frac{\Gamma^2\big(\frac{\nu+n}{2}\big)}{2^{\nu/2}\,
\Gamma(\frac{\nu+1}{2})\Gamma(\frac{\nu+2}{2})\Gamma(\frac{\nu+p+n}{2})}\,
(\rho X)^\nu\big],
\end{align*}
where $X \sim N_1(0,1)$. Observe that the factor $2^{\nu/2}$ in the 
denominator cancels out in the $f_m$-case (i.e., if $f = f_m$)! 
It originates from the integral 
\[
\int_{\S^{n-1}}f(\langle u, v \rangle_{l^2_n})\,
\langle u, v \rangle_{l^2_n}^{\nu}\,\textup{d}\sigma_{n-1}(v) = 2^{-p/2}\,
2^{-\nu/2}\,\frac{\Gamma(\frac{n}{2})}{\Gamma(\frac{\nu+p+n}{2})}\,\E[f(X) X^\nu]
\]
(cf. \sref{Corollary}{cor:spherical_integration_of_inner_products}). In 
general, that important reduction of the fraction cannot be maintained, though. 
Hence, if we represent the latter series as the sum of the even series 
part and the odd series part, we obtain
\[
I_f = I_f^\text{even} + I_f^\text{odd}\,,
\]
where
\[
I_f^\text{even} : = \E\big[f(X)\,\frac{\sqrt{\pi}}{2^{p/2}\,\Gamma(n/2)}\,
(1-\rho^2)^{n/2}\,\sum_{l=0}^\infty \frac{\Gamma^2\big(l+\frac{n}{2}\big)}
{\Gamma(l+\frac{1}{2})\Gamma(l+\frac{p+n}{2})}\,\frac{(\tfrac{1}{2} 
\rho^2 X^2)^l}{l!}\big]
\]
and
\[
I_f^\text{odd} : = \E\big[X f(X)\,\frac{\sqrt{\pi}}{2^{(p+1)/2}\,\Gamma(n/2)}\,
\rho (1-\rho^2)^{n/2}\,\sum_{l=0}^\infty \frac{\Gamma^2\big(l+\frac{1+n}{2}\big)}
{\Gamma(l+\frac{3}{2})\Gamma(l+\frac{p+n+1}{2})}\,\frac{(\tfrac{1}{2} 
\rho^2 X^2)^l}{l!}\big].
\]
A straightforward calculation shows that
\[
I_f^\text{even} = c_+(p,n) \sqrt{\frac{\pi}{2}}\frac{1}{(p-1)!!\,b_{p+1}}\,
(1-\rho^2)^{n/2}\E\big[f(X)\,{}_2 F_2\big(\frac{n}{2}, 
\frac{n}{2}; \frac{1}{2}, \frac{p+n}{2}; (\tfrac{1}{\sqrt{2}} 
\rho X)^2\big)\big]
\]
and
\[
I_f^\text{odd} = c_-(p,n) \sqrt{\frac{\pi}{2}}\,\frac{1}{p!!\,b_{p+2}}\,
\rho (1-\rho^2)^{n/2}\E\big[X f(X)\,{}_2 F_2\big(\frac{n+1}{2}, 
\frac{n+1}{2}; \frac{3}{2}, \frac{p+n+1}{2}; (\tfrac{1}{\sqrt{2}} 
\rho X)^2\big)\big]\,,
\]
where $b_{n}$ satisfies \eqref{eq:the_constants_b_n} ($n \in \{p+1, p+2\}$).
\end{remark}
\chapter{Completely correlation preserving functions} 
\section{Completely real analytic functions and the entrywise matrix functional 
calculus}
 Already while looking for the smallest upper bound of both, $K_G^\R$ 
and $K_G^\C$, we are lead to a deep interplay of different subfields of 
mathematics (both, pure and applied) including Gaussian harmonic analysis and 
{building blocks} of Malliavin calculus (Mehler kernel, Ornstein-Uhlenbeck 
semigroup, Hermite polynomials), integration over spheres in $\R^n$), 
complex analysis (analytic continuation and biholomorphic mappings, special 
functions), combinatorial analysis (inversion of Taylor series and ordinary 
partial Bell polynomials), matrix analysis (positive semidefinite matrices, 
block matrices) and multivariate statistics and high-dimensional Gaussian 
dependence modelling (correlation matrices, real and complex Gaussian random 
vectors, Gaussian measure). 

In particular, we have to look for those functions which map correlation matrices of any 
size and any rank entrywise into a correlation matrix of the same size again, by means of the 
so-called \textit{Hadamard product} of matrices:
\begin{definition}[\textbf{Hadamard product}]
Let $m,n \in \N$. Let $A = (a_{ij}) \in \M_{m,n}(\F)$ and $B = (b_{ij}) \in 
\M_{m,n}(\F)$. The Hadamard product $A \ast B \in \M_{m,n}(\F)$ is 
defined as
\[
(A \ast B)_{ij} : = a_{ij}\,b_{ij} \hspace{0.3cm} ((i,j) \in [m] \times [n])\,.
\] 
\end{definition}
\noindent The Hadamard product is sometimes called the \textit{entrywise product}, for 
obvious reasons, or the \textit{Schur product}, because of some early and basic results about 
the product obtained by Issai Schur (cf. \cite{HJ2013}). Like the usual matrix product, the 
distributive law also holds for the Hadamard product: $A\ast(B + C) = A \ast B + A \ast C$. 
Unlike the usual matrix product, the Hadamard product is commutative: $A \ast B = B \ast A$.
\begin{remark}
Very often, the Hadamard product is denoted by the symbol $\circ$. However, given our view, 
the symbolic notation $\circ$ perhaps could lead to a minor ambiguity, since quite regularly, 
$\circ$ denotes composition of mappings. This is why we adopt the symbol $\ast$ instead, used 
for the definition of the Schur product in \cite[p 29 ff]{P2002}.
\end{remark}
\begin{remark}[\textbf{Hadamard product as subordinated Kronecker product}]
\label{rem:Hadamard_subord_to_Kronecker}
Fix $m, n \in \N$ and $A, B \in \M_{m,n}(\F)$. There exists an interesting link 
between the Hadamard product $A \ast B$ and the Kronecker product $A \otimes B$, 
induced by \eqref{eq:multipl_of_Kronecker_products} and 
\eqref{eq:tensor_prod_rep_II}. In order to recognise this, let $(i,j) \in [m] 
\times [n]$ be given arbitrarily. Then
\begin{align*}
(A \ast B)_{ij} &\,\,\,= (e_i^\top A e_j)(e_i^\top B e_j) \stackrel{(!)}{=}
e_i^\top \otimes e_i^\top (A \otimes B) e_j \otimes e_j\\ 
&\stackrel{\eqref{eq:tensor_prod_rep_II}}{=} e_{(i-1)m+i}^{(m^2)\top}(A \otimes B)
e_{(j-1)n+i}^{(n^2)} = (A \otimes B)_{\Psi_m(i,i), \Psi_n(j,j)}\,.
\end{align*}
Consequently,
\[
A \ast B = (A \otimes B)_{\psi}\,,
\]
where $\psi(i,j) : =  (\Psi_m(i,i), \Psi_n(j,j)) = ((i-1)m+i, (j-1)n+j)$ for 
all $(i,j) \in [m]\times[n]$.
\end{remark}
\noindent If we combine the latter fact and 
\sref{Proposition}{prop:tensor_product_of_quant_corr_is_quant_corr}, we 
immediately obtain (cf. also \eqref{eq:Schur_product_of_corr_matrices}):
\begin{proposition}\label{prop:Schur_product_of_quant_corr_is_quant_corr}
Let $m, n \in \N$, $S \in \M_{m,n}(\F)$ and $R \in \M_{m,n}(\F)$.
If $S \in {\mathcal{Q}}_{m,n}(\F)$ and $R \in {\mathcal{Q}}_{m,n}(\F)$, then
$S \ast R \in {\mathcal{Q}}_{m,n}(\F)$. Note again that the Hadamard product 
is commutative.
\end{proposition}
\noindent Although, the proof of the following two facts are just a consequent 
application of the definition of the Schur product, they are of interest on their 
own, and they help us strongly to support, inter alia, a quick proof of Theorem 
\ref{thm:Schur_multiplication}. In this context, diagonal matrices play an important 
role: if $a = (a_1, \ldots, a_p)^\top \in \F^p$, then $D_a \in \M_p(\F)$ denotes 
the matrix, whose $(i,j)$'th entry is given by $\delta_{ij}\,a_j$. Moreover, we 
need the matrix $J^{(m,n)} \in \M_{m,n}(\F)$, whose $(i,j)$'th entry is given 
by $1$: $J^{(m,n)} : = \sum_{i=1}^m\sum_{j=1}^n e_i e_j^\top$. 
\begin{lemma}\label{lem:Gamma_F_of_Schur_product_of_vectors}
Let $m, n \in \N$, $B \in \M_{m,n}(\F)$ and $(x,y) \in \F^m \times \F^n$. Then
\[
\Gamma_\F(x, y) \ast B = \overline{x}y^\top \ast B 
= D_{x}^\ast B D_y\,.
\]
In particular,
\[
\Gamma_\F(y, y) \ast A \in \M_n(\F)^+ \text{ for all } A \in \M_n(\F)^+\,.
\] 
\end{lemma}
\begin{lemma}[\textbf{Hadamard product factor shifting}]
\label{lem:Frobenius_inner_product_and_Schur_mult}
Let $m, n \in \N$ and $A, B, C \in \M_{m,n}(\F)$. Then
\begin{align}\label{eq:Frobenius_inner_product_and_Schur_mult}
\langle A \ast \overline{B}, J^{(m, n)}\rangle_{F} = \langle A, B\rangle_{F} = 
\langle J^{(m, n)}, \overline{A} \ast B\rangle_{F}\,.
\end{align}
In particular, 
\[
\langle A \ast B, C\rangle_{F} = \langle B, \overline{A} \ast C \rangle_{F}\,. 
\] 
\end{lemma}
\begin{remark}
\sref{Lemma}{lem:Frobenius_inner_product_and_Schur_mult} implies that for any $A \in 
\M_n(\F)$ the adjoint of the ``Schur multiplier'' $S_A : 
{(\M_{m,n}(\F), \Vert \cdot \Vert_F) \longrightarrow (\M_{m,n}(\F), 
\Vert \cdot \Vert_F)}$, $B \mapsto A \ast B$ coincides with the ``Schur 
multiplier'' $S_{\overline{A}}$ (cf. \cite[p. 29 ff]{P2002}):
\[
S_A^\ast = S_{\overline{A}}\,.
\] 
\end{remark}
{
\noindent Fix $m, n \in \N$. The proof of \cite[Theorem 2.2]{FRR1987} reveals 
that in fact \textit{any} matrix $A \in \M_{m,n}(\F)$ is a Schur multiplier:
\[
\Vert S_A : \mathfrak{L}(l_2^n, l_2^m) \longrightarrow 
\mathfrak{L}(l_2^n, l_2^m)\Vert = \sup\{\Vert Auv^\top\Vert_\ast : 
(u, v) \in S_{\F^n} \times S_{\F^m}\} < \infty\,. 
\]
Here, $\Vert Auv^\top \Vert_\ast$ is the trace norm of the matrix 
$Auv^\top$ (cf. \sref{Remark}{rem:GT_vs_trace_class}). Consequently, if we link 
\sref{Proposition}{prop:char_of_Q_m_n_for_both_fields} and 
\eqref{eq:char_of_quantum_corr_matrices} with a further fundamental and deep 
result of Grothendieck, we obtain the following noteworthy representation of 
\textit{arbitrary} $m \times n$ matrices (cf. \cite[Theorem 4.2]{C2018}, 
\cite[Theorem 8.7]{P2002} and \cite[Theorem 5.1 and Proposition 5.4]{P2001}):
\begin{proposition}\label{prop:Schur_mult_and_quantum_corr_rep}
Let $m, n \in \N$ and $A \in \M_{m,n}(\F)$. 
Put $\lambda_A(m,n) : = \Vert A : l_1^n \longrightarrow l_\infty^m \Vert_{\mathfrak{L}_2}$. 
Then $\lambda_A(m,n) = \Vert S_A : \mathfrak{L}(l_2^n, l_2^m) \longrightarrow 
\mathfrak{L}(l_2^n, l_2^m)\Vert$, and
\[
A = \lambda_A(m,n)\,Q_A\,,
\]
where $Q_A \in {\mathcal{Q}}_{m,n}(\F)$.
\end{proposition}
Combining the latter result with the (not numbered) remark before 
\cite[Theorem 3.7]{P2002}, we obtain a further non-negligible
\begin{corollary}
Let $n \in \N$ and $B = (b_{ij}) \in \M_n(\F)^+$ be positive semidefinite. Then
\[
B = \max\limits_{i \in [n]}(b_{ii})\,Q_B\,,
\]
where $Q_B \in {\mathcal{Q}}_n(\F)$. In particular, if $b_{ii} = 0$ for all 
$i \in [n]$, then $B = 0$.
\end{corollary}
}
\noindent The usefulness of the structure of the Hadamard product is reflected 
in the Schur product theorem which states that the (closed, convex and self-dual) 
cone of all positive semidefinite matrices is stable under Schur multiplication. 
We give a short and completely self-contained proof:
\begin{theorem}[\textbf{Schur, 1911}]\label{thm:Schur_multiplication}
Let $n \in \N$. Let $A \in \M_n(\F)^+$ and $B \in \M_n(\F)^+$ be positive semidefinite. Then 
$A \ast B \in \M_n(\F)^+$. In particular,
\begin{align}\label{eq:Schur_product_of_corr_matrices}
C(n; \F) \ast C(n; \F) \subseteq C(n; \F).
\end{align} 
\end{theorem}
\begin{comment}
Fix $y \in \F^n$, and put $M : = \Gamma_\F(\overline{y},\overline{y}) \ast 
\overline{B} = y y^\ast \ast \overline{B}$. Because of 
\sref{Lemma}{lem:Gamma_F_of_Schur_product_of_vectors}, $M^\ast = M \in \M_n(\F)^+$ (since 
also $\overline{B} \in \M_n(\F)^+$). Consequently, 
\sref{Lemma}{lem:Frobenius_inner_product_and_Schur_mult} implies that
\begin{align*}
y^\ast(A \ast B)y &= \text{tr}(yy^\ast(A \ast B)) = \langle A \ast B, 
yy^\ast\rangle_{F} = \langle A, M\rangle_{F} = 
\text{tr}\big(M^{1/2} M^{1/2} A^{1/2} A^{1/2}\big)\\ 
&= \text{tr}\big((M^{1/2} A^{1/2})(A^{1/2} M^{1/2})\big) = 
\Vert A^{1/2} M^{1/2} \Vert_{F}^2 \geq 0\,,
\end{align*}
and the claim follows.
\end{comment}
\begin{definition}[\textbf{Entrywise functional calculus}]
Let $m,n \in \N$. Given $\emptyset \not= U \subseteq \F$, a function 
$f : U \longrightarrow \F$ and a matrix $A = (a_{ij}) \in \M_{m,n}(U)$ put
\[
f[A] : = (f(a_{ij})) \hspace{0.3cm} ((i,j) \in [m] \times [n])\,.
\]
\end{definition}
\noindent In particular, if $f(x) = \sum_{n=0}^\infty c_n\,x^n, x \in U$, where 
$c_n \in \F$ for all $n \in \N$, we have
\[
f[A]_{ij} = \sum_{n=1}^\infty c_n\,a_{ij}^n \text{ for all } (i,j) \in [m] \times [n]\,.
\]
Functions of the latter type, where $c_n \geq 0$ for all $n \in \N$ play a 
significant role, also with respect to an analysis of $K_G^\R$ and $K_G^\C$. 
This is particularly reflected in the next two results with respect to the real 
field, which however also play a key role in the complex case (cf. Theorem 
\ref{thm:h_b_c_complex_correl_case} and \sref{Corollary}{cor:Haagerup_incl}). 
Recall that a function $\psi : I \longrightarrow \R$, defined on an open interval 
$I \subseteq \R$, is \textit{absolutely monotonic}, if $\psi \in C^\infty(I)$ and 
$\psi^{(n)} \geq 0$ on $I$ for all $n \in \N_0$ (cf. e.g. 
\cite[Chapter 19]{K2022} and \cite[Chapter IV]{W1941} 
regarding a rigorous reassessment of that crucial and very rich concept, 
coined by S. N. Bernstein in 1929). Regarding a refresher of complex analysis 
of functions of one complex variable, we recommend to place the rich 
source \cite{SS2003} next to our paper.   
\begin{lemma}\label{lem:real_analytic_and_l1}
Let $r \in (0, \infty)$ and $\psi : (-r, r) \longrightarrow \R$ be real analytic. 
Suppose that $b : = \big(\frac{\psi^{(n)}(0)}{n!}\,r^{n}\big)_{n \in \N_0} \in l_1$. 
Then $\psi = \widetilde{\psi}\big\vert_{(-r,r)}$, where the complex, bounded 
function $\widetilde{\psi} : r \overline{\D} \longrightarrow \Vert b \Vert_1 
\overline{\D}$ is defined as 
\[
\widetilde{\psi}(z) : = \sum_{n = 0}^\infty \frac{\psi^{(n)}(0)}{n!}\,z^n 
\text{ for all } z \in r \overline{\D}.
\]
The real function $\psi_{\text{abs}} : [-r,r] \longrightarrow [-\Vert b \Vert_1, 
\Vert b \Vert_1], x \mapsto \sum_{n=0}^\infty \frac{\vert \psi^{(n)}(0)\vert}{n!}\,
x^{n}$ is continuous and bounded, such as the complex-valued function
$\widetilde{\psi}$. $\psi_{\text{abs}}\big\vert_{(-r, r)}$ 
is real analytic, and $\psi_{\text{abs}}\big\vert_{(0, r)}$ is absolutely 
monotonic on $(0, r)$. $\Vert b \Vert_1 = \psi_{\text{abs}}(r)$ and
\begin{align}\label{eq:psi_abs_estimation}
\vert \widetilde{\psi}(z)\vert \leq \psi_{\text{abs}}(\vert z\vert) 
\leq \psi_{\text{abs}}(r)\,\text{ for all }\, z \in r \overline{\D}\,.
\end{align} 
In particular, $\psi_{\text{abs}}$ can be extended to the continuous complex 
function $\widetilde{\psi}_{\text{abs}} : r \overline{\D} \longrightarrow 
\Vert b \Vert_1 \overline{\D}$, defined as
\[
\psi_{\text{abs}}(z) \equiv \widetilde{\psi}_{\text{abs}}(z) : = 
\sum_{n=0}^\infty \frac{\vert \psi^{(n)}(0)\vert}{n!}\,z^{n} \text{ for all } 
z \in r \overline{\D}.
\]
$\widetilde{\psi}\big\vert_{r \D}$ (respectively 
$\psi_{\text{abs}}\big\vert_{r \D}$) is the unique holomorphic extension of 
$\psi$ (respectively $\psi_{\text{abs}}\big\vert_{(-r, r)}$) on the domain $r \D$. 
\end{lemma}
\noindent Since that class of real analytic functions plays a recurring and 
decisive role in our paper (particularly for $r = 1$), and since real 
(and complex) analyticity of a function actually is a local property, it is 
justifiable to introduce the following definition: 
\begin{definition}\label{def:CRA}
Let $r \in (0, \infty)$. Put
\begin{align*}
W^\omega_+((-r, r)) : = \{\psi: \psi \in C^\omega((-r, r)) \text{ and } 
\big(\frac{\psi^{(n)}(0)}{n!} r^{n}\big)_{n \in \N_0} \in l_1\}. 
\end{align*}
Any element in $\psi \in W^\omega_+((-r, r))$ is said to be \textit{completely 
real analytic on $(-r, r)$ at $0$}.
\end{definition} 
\noindent Observe that by definition any function in $W^\omega_+((-r, r))$ 
coincides with its own Taylor series at $0$ ``completely'' 
(i.e., \textit{everywhere}) on its domain of definition $(-r, r)$ (and not 
``locally, around'' $0$ only). Moreover, due to \sref{Lemma}{lem:real_analytic_and_l1}, it 
follows that for any $\psi \in W^\omega_+((-r, r))$ also 
$\psi_{\text{abs}}{\big\vert}_{(-r, r)} \in W^\omega_+((-r, r))$, and
\begin{align}\label{eq:n_th_derivative_of_psi_abs_at_zero}
\psi_{\text{abs}}^{(n)}(0) = \vert \psi^{(n)}(0)\vert \text{ for all } n \in \N_0\,.
\end{align}
Let us explicitly highlight three facts, implied by 
\sref{Lemma}{lem:real_analytic_and_l1}. To this end, if $\alpha \in \T$ is given, we 
consider the biholomorphic function $M_\alpha : \C \longrightarrow \C$, defined 
as $M_\alpha(z) : = \alpha z$. Obviously, $M_\alpha^{-1} = M_{\overline{\alpha}}$ 
and $M_\alpha(\D) = \D$.  
\begin{remark}[\textbf{Sign condition}]\label{rem:odd_psi_abs_and_sign_condition}
Let $\psi$ and $\widetilde{\psi}$ be given as in 
\sref{Lemma}{lem:real_analytic_and_l1}. If $\psi$ is odd and $\text{sign}
(\psi^{(2n+1)}(0)) = (-1)^n$ for all $n \in \N_0$, then
\[
\psi_{\text{abs}}(z) = \frac{1}{i} \widetilde{\psi}(i z) = (M_{-i}\circ 
\widetilde{\psi} \circ M_i)(z) \text{ for all } z \in 
r \overline{\D}
\]
and
\[
\psi_{\text{abs}}^{(2n+1)}(\zeta) = (-1)^n\,{\widetilde{\psi}}^{(2n+1)}(i \zeta) 
\text{ for all } \zeta \in r \D \text{ and } n \in \N_0
\]
(since $i^{2n} = (-1)^n$ for all $n \in \N_0$). 
\end{remark}
\begin{remark}[\textbf{Wiener algebra}]\label{rem:disc_algebra_and_positive_Wiener_algebra}
Let $\psi$ be given as in \sref{Lemma}{lem:real_analytic_and_l1}. Assume that 
$r = 1$. Then $\widetilde{\psi}\big\vert_{\D}$ and $\psi_{\text{abs}}\big\vert_
{\D}$ are elements of
\[
W^+(\D) : = \{f : f(z) =\sum_{n=0}^\infty b_n z^n \text{ is holomorphic on } \D 
\text{ and satisfies } \sum_{n=0}^\infty \vert b_n\vert < \infty\}\,.
\]
$W^+(\D)$ is known as the \textit{Wiener algebra}. It is a unital commutative 
(complex) Banach algebra with respect to the norm $\Vert \sum_{n=0}^\infty b_n 
z^n \Vert_{W^+(\D)} : = \sum_{n=0}^\infty \vert b_n\vert$, where multiplication 
is defined as that one of analytic functions, via the Cauchy product formula 
(cf. \cite{M1989}). In particular, $\Vert \widetilde{\psi}\big\vert_{\D}\Vert_
{W^+(\D)} = \psi_{\text{abs}}(1) = \Vert\psi_{\text{abs}}\big\vert_{\D}\Vert_
{W^+(\D)}$. Since by construction of $W^+(\D)$ the series $f = \sum_{n=0}^\infty 
f_n$ is normally convergent in $\overline{\D}$ (with respect to the supremum 
norm), where $\overline{\D} \ni z \mapsto f_n(z) : = b_n z^n$, the series $f$ is 
uniformly convergent in $\overline{\D}$, and it follows that every element of 
$W^+(\D)$ can be continuously extended to an element of the unital commutative 
Banach algebra $A(\D)$, where 
\[
A(\D) : = \{g : g \in C(\overline{\D}) \text{ such that } g\big\vert_{\D} 
\text{ is holomorphic on } \D \}
\] 
is equipped with the supremum norm and the usual pointwise algebraic 
operations. $A(\D)$ denotes the well-known \textit{disc algebra} (cf. 
\cite[V.1., Example 4]{G2006} and \cite[Chapter III.E.]{W1991}). Thus, 
$W^+(\D)$ could be viewed as a subalgebra of the disc algebra $A(\D)$. Since 
$A(\D) \subseteq H^\infty \subseteq H^2 \subseteq H^1$, it is likely that 
a further link to the very rich theory of Hardy spaces might open up here 
(cf., e.g., \cite{G2006, K2002, R1980, R1987}) .
\end{remark}
\begin{remark}[\textbf{Inversion and complete real analyticity}]
\label{rem:inversion_and_CRA}
Like a common thread, the following highly non-trivial problem - which is 
decisive for computing the upper bounds of $K_G^\F$ - runs throughout our whole 
paper. Very generally formulated, let $\psi \in W^\omega_+((-1, 1))$ be odd and 
completely real analytic on $(-1,1)$ at $0$. Assume that $\psi = 
\psi_{\text{abs}}\big\vert_{(-1,1)}$ and that $\psi : (-1,1) \stackrel{\cong}
{\longrightarrow} (-1,1)$ is bijective. Assume further that $\psi^{-1}$ is real 
analytic on $(-1,1)$. Does $\psi^{-1} \in W^\omega_+((-1, 1))$ already apply; 
i.e., is then also $(\psi^{-1})_{\text{abs}}$ well-defined (cf., e.g., 
\sref{Lemma}{lem:hyperbolic_CCP_transform}, 
\sref{Corollary}{cor:h_ff_real_case_abs_value_crit}, \autoref{thm:odd_real_bdd_CCP}, 
\sref{Example}{ex:Krivine_recovered} and 
\sref{Example}{ex:Haagerup})? Both, \cite{BMMN2013} and \cite{H1987} present 
multi-page proofs, each one built on rather advanced complex analysis (including 
{intricate} holomorphic extensions), to answer this question to the positive 
for key functions $\psi \in W^\omega_+((-1, 1))$ related to computing the value 
of $K_G^\C$ and $K_G^\R$, respectively\,!
\end{remark}   
\section{Completely correlation preserving functions and Schoenberg's theorem}
\noindent Let $0 < c \leq r$ and $\psi$ be given as in 
\sref{Lemma}{lem:real_analytic_and_l1}. Assume that $\psi \not= 0$. Since 
$\psi \in W^\omega_+((-r,r))$, it follows that $\psi^{(n_0)}(0) \not= 0$ for 
some $n_0 \in \N_0$, implying that $\psi_{\text{abs}}(c) \geq 
\frac{\vert\psi^{(n_0)}(0)\vert}{n_0\,!}\,c^{n_0} > 0$. Thus, the function 
\begin{align}\label{eq:canonical_pre_CCP_transform}
(-1,1) \ni \rho \mapsto \psi_c(\rho) : = \frac{\psi(c\rho)}{\psi_{\text{abs}}(c)}
\end{align}
is well-defined, continuous, bounded and satisfies $0 \not= \psi_c \in 
W^\omega_+((-1,1))$. In Theorem 
\ref{thm:char_of_completely_real_analytic_functions_at_0_as_h_fg} we will shed 
light on the hidden structure of these functions $\psi_c$.
Let $n \in \N$. Since $(\psi_c)_{\text{abs}}(1) = 1$, it even follows that 
$(\psi_c)_{\text{abs}} = \frac{\psi_{\text{abs}}(c\,\cdot)}{\psi_{\text{abs}}(c)} 
: [-1,1] \longrightarrow [-1,1]$ transforms any real $n \times n$-correlation 
matrix entrywise into a real $n \times n$-correlation matrix; i.e.,
\begin{align}\label{eq:CCP_property}
(\psi_c)_{\text{abs}}[A] \in C(n; \R) \text{ for all } A \in C(n; \R), 
\text{ for all } n \in \N
\end{align} 
(due to \autoref{thm:Schur_multiplication}). 
Within the scope of our research, functions $f \in C([-1,1])$, satisfying
$f{\big\vert}_{(-1,1)} \in W^\omega_+((-1, 1))$ and 
$f{\big\vert}_{(-1,1)}(\rho) = \big(f{\big\vert}_{(-1,1)}\big)_{\text{abs}}(\rho)$ 
for all $\rho \in (-1,1)$, are of particular importance, due to the following 
version of a fundamental result of I. J. Schoenberg (cf. \cite[Theorem 16.2]{K2022}
and \cite{S1942}):
\begin{theorem}[\textbf{Schoenberg, 1942}]\label{thm:Schoenberg}
Let $f : [-1,1] \longrightarrow \R$ be a continuous function. Then the following statements 
are equivalent:
\begin{enumerate}
\item $f[A]$ is positive semidefinite for all $A \in \bigcup_{n=1}^\infty
\M_n([-1,1])^+$.
\item $f[\Sigma]$ is positive semidefinite for all 
$\Sigma \in \bigcup_{n=1}^\infty C(n; \R)$.
\item $f(x)$ equals a convergent series $\sum_{n=0}^\infty a_n\,x^n$ 
for all $x \in [-1,1]$, where $a_n \geq 0$ for all $n \in \N_0$.
\item $f{\big\vert}_{(-1,1)} \in W^\omega_+((-1, 1))$ and
$f = \big(f{\big\vert}_{(-1,1)}\big)_{\text{abs}}$.
\item $f{\big\vert}_{(-1,1)} \in W^\omega_+((-1, 1))$ and 
$f{\big\vert}_{(0,1)}$ is absolutely monotonic.
\item $f{\big\vert}_{(-1,1)}$ can be extended to a holomorphic function
$\D \ni z \mapsto \widetilde{f}(z) : = 
\sum_{n=0}^\infty a_n\,z^n$, where $a_n \geq 0$ for all $n \in \N_0$ and $f(1) =
\sum_{n=0}^\infty a_n < \infty$. 
\end{enumerate}
In particular, if (iii) or (vi) holds, then $a_n = \frac{f^{(n)}(0)}{n!} \geq 0$ 
for all $n \in \N_0$, and the series $\sum_{n=0}^\infty a_n\,x^n$ converges 
uniformly on $[0,1]$.
\end{theorem}
\noindent Fix $\F \in \{\R, \C\}$. Looking for common sources of the real and 
complex cases, we put 
$D_{\F} : = \{a \in \F : \vert a \vert < 1\}$, so that $D_{\R} = (-1,1)$ and
$D_{\C} = \D$. Property \eqref{eq:CCP_property} in 
\sref{Lemma}{lem:real_analytic_and_l1} deserves an autonomous and far-reaching
\begin{definition}[\textbf{Completely correlation preserving function}]
\label{def:nCP_and_CCP}
Fix $\F \in \{\R, \C\}$. Let $h : \overline{D_{\F}} \longrightarrow \F$ be a 
function.
\begin{enumerate}
\item Given $n \in \N$, $h$ is $n$-correlation-preserving (short: $n$-CP) if for any 
$n \times n$ correlation matrix $\Sigma \in C(n; \F)$ also $h[\Sigma] \in C(n; \F)$ is an 
$n \times n$ correlation matrix. 
\item $h$ is called completely correlation-preserving (short: CCP) if $h$ is 
$n$-correlation-preserving for all $n \in \N$.  
\end{enumerate}
\end{definition}
\noindent Since every CCP function $h$ is 2-CP, \sref{Definition}{def:nCP_and_CCP} 
directly implies that $h(\overline{D_{\F}}) \subseteq \overline{D_{\F}}$ and 
$h(1) = 1$. If $h : \overline{D_{\F}} \longrightarrow \F$ satisfies $h(1) > 0$, 
then $h[\Sigma] \in \M_n(\overline{D_{\F}})^+$ for all $\Sigma \in C(n; \F)$ if 
and only if $\frac{1}{h(1)}\,h$ is $n$-CP. If $h_1$ and $h_2$ are 
two CCP functions (respectively two $n$-CP functions), and if $\lambda \geq 0$, 
\sref{Definition}{def:nCP_and_CCP} immediately implies that also $h_2 \circ h_1$, 
$h_1 \circ h_2$ and $\lambda h_1 + (1-\lambda)h_2$ are CCP functions (respectively 
$n$-CP functions). Furthermore, since $h_1 h_2[A] = h_1[A] \ast h_2[A]$ for all 
matrices $A \in \M_{m,n}(\F), m, n \in \N$, it follows from 
\eqref{eq:Schur_product_of_corr_matrices} that also the product $h_1 h_2$ of two 
CCP functions (respectively two $n$-CP functions) again is CCP (respectively $n$-CP). 
Much less trivial is the fact (which involves the Grothendieck 
constant\,!) that in general, the inverse function $h^{-1}$ of an 
invertible CCP function $h$ is \textit{not} a CCP function (cf. 
\sref{Corollary}{cor:inverse_of_non_trivial_odd_CCP_is_not_CCP}, 
\sref{Remark}{rem:inverse_of_CCP_in_general_is_not_CCP} and 
\autoref{thm:odd_real_bdd_CCP}-\ref{(i)}), such as the 
inverse of $h : = \frac{2}{\pi}\arcsin$, given by $[-1,1] \ni y \mapsto 
h^{-1}(y) = \sin(\frac{\pi}{2} y) = \sum_{n=0}^\infty b_n \frac{y^n}{n !}$, 
where $b_n : = (-1)^{\floor{n/2}}\cdot\frac{1 - (-1)^n}{2}\cdot(\frac{\pi}{2})^n$. 

In the real case, 
\autoref{thm:Schoenberg} immediately leads to a full characterisation of 
continuous CCP functions, since:
\begin{theorem}\label{thm:CCP_via_Schoenberg}
Let $h : [-1,1] \longrightarrow \R$ be a continuous function. Then the following 
statements are equivalent:
\begin{enumerate}
\item $h$ is CCP.\label{(i)}
\item $h[A]$ is positive semidefinite for all $A \in \bigcup_{n=1}^\infty
\M_n([-1,1])^+$ and $h(1)=1$.\label{(ii)}
\item $h(x)$ has the unique series representation 
$h(x) = \sum_{n=0}^\infty a_n\,x^n$ for all $x \in [-1,1]$, where $a_n \geq 0$ for all 
$n \in \N_0$ and $\sum_{n=0}^\infty a_n = 1$.\label{(iii)}
\item $[-1,1] \ni \rho \mapsto h(\rho) = \E_\P[\rho^X] = 
\sum_{n=0}^\infty \P(X = n)\rho^n$ is the probability generating function 
of some discrete random variable $X$.\label{(iv)}
\item $h{\big\vert}_{(-1,1)} \in W^\omega_+((-1, 1))$, 
$h = \big(h{\big\vert}_{(-1,1)}\big)_{\text{abs}}$ and $h(1)=1$.\label{(v)}
\item $h{\big\vert}_{(-1,1)} \in W^\omega_+((-1, 1))$,  
$h{\big\vert}_{(0,1)}$ is absolutely monotonic and $h(1) = 1$.\label{(vi)}
\item $h{\big\vert}_{(-1,1)}$ can be extended to a holomorphic function
$\D \ni z \mapsto \widetilde{h}(z) : = 
\sum_{n=0}^\infty a_n\,z^n$, where $a_n \geq 0$ for all $n \in \N_0$ and 
$\sum_{n=0}^\infty a_n = 1$.\label{(vii)}
\item $h{\big\vert}_{(-1,1)}$ can be extended to a complex 
function $\widetilde{H} : \overline{\D} \longrightarrow \overline{\D}, z 
\mapsto \sum_{n=0}^\infty a_n\,z^n$, where $a_n \geq 0$ for all $n \in \N_0$ 
and $\sum_{n=0}^\infty a_n = 1$.\label{(viii)}
\end{enumerate}
If \textup{(iii)} or \textup{(vii)} or \textup{(viii)} is given, the series 
$\sum_{n=0}^\infty a_n\,z^n$ converges then absolutely on $\overline{\D}$ and 
$\sum_{n=0}^\infty a_n\,x^n$ converges uniformly on $[0,1]$. 
\end{theorem}
\begin{remark}[\textbf{CCP versus RKHS}]\label{rem:CCP_and_RKHS}
\autoref{thm:CCP_via_Schoenberg}-\ref{(iii)} also reveals a direct link to 
representing kernel Hilbert spaces (RKHS). \cite[Theorem 4.12.]{PR2016}, 
together with the remark on complexification in \cite[Chapter 5.1]{PR2016} 
namely implies that 
\[
(-1,1) \times (-1,1) \ni (s, t) \mapsto K_h(s,t) : = \sum_{n=0}^\infty a_n (st)^n =
\sum_{n=0}^\infty f_n(s) f_n(t)
\]
is a well-defined kernel function for a real RKHS $H(K_h)$ of real-analytic 
functions on $(-1,1)$ (cf. \cite[Definition 2.12]{PR2016}), where $f_n(\rho) : =
\sqrt{a_n}\,\rho^n$. The set of functions $f_n$ for which $a_n \not= 0$ is an 
orthonormal basis for $H(K_h)$.
\end{remark}
\noindent It is quite instructive to analyse how $n$-CP functions in particular 
relate to properties of Schoenberg's kernel functions $\mathcal{P}(\S^d)$, $d 
\in \N$. Schoenberg's approach was revisited by the geostatistician T. Gneiting 
in his impressive paper \cite{Gn2013}, where $\mathcal{P}(\S^d)$ is denoted as 
$\Psi_{d}$. So, let us recall the concept of a positive definite function on a 
metric space (which in particular is a kernel function (cf. 
\cite[Definition 2.12.]{PR2016}), originally coined by Schoenberg in \cite{S1942}. 
Supported by \sref{Lemma}{lem:charact_of_corr_matrices}, these kernel functions 
result in a further characterisation of real-valued $n$-CP functions. In order 
to recognise this observation, let $(X, d)$ be a metric space. Fix $d \in \N_2$. 
A function $f : [0, \infty) \longrightarrow \R$ is positive definite on $X$ if 
$f$ is continuous, and if \textit{for any $n \in \N$} and any $(x_1, x_2, \ldots, 
x_n) \in X^\nu$, the matrix $f[d(x_i, x_j)_{i,j=1}^n]$ is positive 
semidefinite (cf. \cite[Definition 16.15]{K2022} and \cite{Gn2013, S1942}). 
Consequently, if we apply Schoenberg's definition to the unit sphere $\S^{d-1}$, 
where the metric on $\S^{d-1}$ is induced by the geodesic distance $\S^{d-1} 
\times \S^{d-1} \ni (x,y) \mapsto \arccos(\langle x, y\rangle)$, 
\sref{Lemma}{lem:charact_of_corr_matrices} 
implies that any continuous function $f : [-1,1] \longrightarrow \R$ which 
satisfies $f \circ \cos : [0, \pi] \longrightarrow \R \in 
\mathcal{P}(\S^{d-1}) \equiv \Psi_{d-1}$ in particular is $d$-CP. 
Thus, \cite[Table 1]{Gn2013} shows us a wealth of non-trivial $d$-CCP functions, 
where $d \in [3]$. Similarly, we recognise that a continuous function $h : [-1,1] 
\longrightarrow \R$ is CCP if and only if $h \circ \cos \in \Psi_{\infty} = 
\bigcap_{d=1}^{\infty}\Psi_{d}$.    
 
If we combine these facts with \cite[Theorem 7]{Gn2013}, we are 
rewarded with a large class of (even) CCP functions. To this end, recall that a 
continuous function $\psi : [0, \infty) \longrightarrow \R$ is called 
\textit{completely monotonic} if $\psi\big\vert_{(0, \infty)} \in 
C^\infty((0, \infty))$ and $(-1)^n \psi^{(n)}(x) \geq 0$ for all $n \in \N_0$ 
and all $x > 0$ (cf., e.g., \cite[Definition 27.18]{K2022} and 
\cite[Definition 2c]{W1941}). Many explicit examples of completely 
monotonic functions are listed in \cite{MS2001}. They play a significant role in 
various subfields of probability theory including theory and applications of 
L\'{e}vy processes and infinite divisibility.
\begin{theorem}
Let $\psi : [0, \infty) \longrightarrow \R$ be completely monotonic and non-constant. 
Suppose that $\psi(0) = 1$, then 
\[
\psi\circ\arccos : [-1,1] \longrightarrow \R 
\text{ is a CCP function}.
\]
In particular, $\psi([0, \pi]) \subseteq [-1,1]$.
\end{theorem}
\noindent If we apply the complex analogue of Schoenberg's Theorem, coined by J.P.R. 
Christensen and P. Ressel in 1982 (cf. \cite[Theorem 16.7]{K2022}) to complex 
CCP functions, we obtain
\begin{theorem}\label{thm:complex_CCP_via_Christensen_and_Ressel}
Let $h : \overline{\D} \longrightarrow \C$ be a continuous complex function. 
Then the following statements are equivalent:
\begin{enumerate}
\item $h$ is CCP.
\item $h(z)$ has the unique series representation 
$h(z) = \sum_{k=0}^\infty\sum_{l=0}^\infty a_{k l}\,z^k\,\overline{z}^l$ for all 
$z \in \overline{\D}$, where $a_{k l} \geq 0$ for all $k, l \in \N_0$ and 
$\sum_{k=0}^\infty(\sum_{l=0}^\infty a_{k l}) = 1$.
\end{enumerate}
\end{theorem}
\noindent Observe that a significant implication of 
\autoref{thm:complex_CCP_via_Christensen_and_Ressel} is that in general 
\textit{complex} CCP functions are not holomorphic, respectively analytic. 
\autoref{thm:CCP_via_Schoenberg}, together with 
\sref{Lemma}{lem:real_analytic_and_l1} immediately implies how one can easily construct 
complex CCP functions out of real ones: 
\begin{remark}\label{rem:complex_CCP_construction}
Let $h = h_{f,f} : [-1,1] \longrightarrow [-1,1]$ be a real CCP function 
($\Vert f \Vert_{\gamma_k} = 1$) and $\widetilde{h} : 
\overline{\D} \longrightarrow \overline{\D}$ be as in 
\sref{Lemma}{lem:real_analytic_and_l1}. Then the following functions are all complex CCP
functions: $\widetilde{h}$, $\overline{\widetilde{h}}$ and 
$\widetilde{h}\,\cdot\,\overline{\widetilde{h}} = 
\vert \widetilde{h}\vert^2$.  
\end{remark}
\chapter{The real case: towards extending Krivine's approach}
\section{Some facts about real multivariate Hermite polynomials}
Another important estimation (even with upper bound $1$) which might help to 
support our search for a ``suitable'' CCP function which is different from the 
CCP function $\frac{2}{\pi}\arcsin$ in the real case, is included in the next 
two results. We start with the real case first. To this end, recall e.g. from 
\cite[Chapter 1.3]{B1998} that for fixed $k \in \N$ $\{H_\alpha : \alpha \in 
\N_0^k\}$ is an ortho\textit{normal} basis in $L^2(\gamma_k)$, where the 
$k$-variate Hermite polynomial $H_\alpha : \R^k \longrightarrow \R$ is defined 
by 
\begin{align}\label{eq:k_dim_Hermite_polynomials}
H_\alpha(x_1, x_2, \ldots, x_k) : = \prod\limits_{i=1}^k H_{\alpha_i}(x_i) =
\frac{(-1)^{\vert \alpha \vert}}{\sqrt{\alpha !}}
\frac{D^\alpha\varphi_k(x)}{\varphi_k(x)}\,. 
\end{align}
Here, $\varphi_k(x_1, x_2, \ldots, x_k) : = \prod\limits_{i=1}^k \varphi(x_i)$, 
$D^\alpha : = (\frac{\partial}{\partial x_1})^{\alpha_1}\cdots
(\frac{\partial}{\partial x_k})^{\alpha_k}$,
$\alpha ! : = \prod\limits_{i=1}^k \alpha_i !$ and $\vert \alpha \vert : = \sum_{i=1}^k \alpha_i$ 
($n \in \N_0$, $y \in \R, x = (x_1, \ldots, x_k)^\top \in \R^k$), and
\begin{align}\label{eq:one_dim_Hermite_pol}
\begin{split}
H_n(y) &:= \frac{1}{\sqrt{n!}}(-1)^n \exp\big(\frac{y^2}{2}\big)\,
\frac{\textup{d}^n}{\textup{d}y^n}\exp\big(-\frac{y^2}{2}\big) = 
\frac{(-1)^n}{\sqrt{n!}}\frac{\varphi^{(n)}(y)}{\varphi(y)}\\ 
&= \frac{1}{\sqrt{n!}}\,
\sum_{k=0}^{\floor{n/2}}H_{n,k}\,(-1)^k\,y^{n-2k} =
(-1)^n H_n(-y).
\end{split}
\end{align}
$H_n$ denotes the (probabilist's version of the) one-dimensional Hermite 
polynomial, where 
\[
H_{n,k} : = \binom{n}{2k}\,(2k-1)!! = \frac{n!}{2^k k!(n-2k)!}\,.
\] 
Thus,
\begin{align}\label{eq:Hermite_product}
H_\alpha\,\varphi_k = \frac{(-1)^{\vert \alpha \vert}}{\sqrt{\alpha !}}\,
\frac{D^\alpha\varphi_k}{\varphi_k}\,.
\end{align}
Since $(y + i\,z)^n = \sum_{l=0}^n \binom{n}{l}\,i^l\,z^l\,y^{n-l}$ and
$\E[X^m] = \frac{(-1)^m + 1}{2}\,(m-1)!!$ for all $m, n \in \N_0$, $y, z \in 
\R$ and $X \sim N_1(0,1)$, it follows that 
\[
\E[(y + i\,X)^n] = \sum_{l=0}^n \binom{n}{l}\,i^l\,\E[X^l]\,y^{n-l} =
\sum_{k=0}^{\floor{n/2}}\binom{n}{2k}\,(-1)^k\,(2k-1)!!\,y^{n-2k} =
{\sqrt{n!}}\,H_n(y)
\]
for all $n \in \N_0$, $y \in \R$ and $X \sim N_1(0,1)$. In particular, we
reobtain the numbers
\begin{align}\label{eq:Hermite_polynomials_at_zero}
H_{2l}(0) = (-1)^l\,\frac{1}{\sqrt{(2l)!}}\,(2l-1)!! = 
(-1)^l\sqrt{\frac{(2l-1)!!}{(2l)!!}} = (-1)^l\,\frac{\sqrt{(2l)!}}{2^l\,l!} 
\text{ and } H_{2l+1}(0) = 0
\end{align}
for all $l \in \N_0$ (cf. also \eqref{eq:Gauss_Copula_special_case_arcsin}). Recall 
that the generating function of the one-dimensional Hermite polynomials is 
given by (cf. \cite[(1.3.1)]{B1998})
\[
\R \times \R \ni (\lambda, x) \mapsto \exp\big(\lambda x - \frac{1}{2}\lambda^2\big) = 
\sum_{\nu=0}^\infty \frac{\lambda^\nu}{\sqrt{\nu!}}\,H_\nu(x) 
\]
implying that in the $k$-dimensional case the equality 
\begin{align}\label{eq:k_dim_generating_fct}
\exp\big(\lambda^\top\,x - \frac{1}{2}{\Vert\lambda\Vert}^2\big) =
\prod\limits_{i = 1}^k \exp\big(\lambda_i\,x_i - \frac{1}{2}\lambda_i^2\big) = 
\sum_{m \in \N_0^k}\frac{1}{\sqrt{m!}}\,H_m(x)\lambda^{m} =
\sum_{\nu=0}^\infty \big(\sum_{m \in C(\nu, k)}\frac{H_\alpha(x)}{\sqrt{m!}}\big)
\lambda^{m}
\end{align}
holds for all $\lambda, x \in \R^k$, where $C(\nu, k) : = 
\{m \in \N_0^k : \vert m \vert = \nu\}$. Namely, since $\vert H_m(x)\vert \leq (2\pi)^{k/4}
\exp\big(\tfrac{1}{4}\Vert x \Vert_2^2\big)$ for all $m \in \N_0^k$ and $x \in \R^k$ (see 
\cite[Lemma 1]{DILP2018}, respectively \cite[Proposition 2.3\,(i)]{IL2015}), the quotient 
test implies that $\sum_{\nu=0}^\infty\frac{1}{\sqrt{\nu!}}\,H_\nu(x)\lambda^\nu$ even 
converges absolutely. Hence, each one of the $k$ families 
$\big(\frac{1}{\sqrt{\nu!}}\,H_\nu(x_i)\lambda_i^\nu\big)_{\nu \in \N_0}$ is 
summable ($ i \in [k]$). Consequently, it follows that also the family 
$\big(\prod\limits_{i=1}^k\frac{1}{\sqrt{m_i!}}\,H_{m_i}(x_i)
\lambda_i^{m_i}\big)_{m \in \N_0^k} = \big(\frac{1}{\sqrt{m!}}\,
H_m(x)\lambda^{m}\big)_{m \in \N_0^k}$ is summable. 
\eqref{eq:k_dim_generating_fct} now follows from the reiteration of the double summation 
principle and the associativity formula for summable families. Consequently, 
we obtain (see also \cite[Lemma 1.3.2 (iii)]{B1998}):
\begin{align}\label{eq:k_dim_Taylor_coeff}
\sqrt{m!}\,H_m(x) = 
D^m \exp(\lambda^\top\,x - 
\tfrac{1}{2}{\Vert\lambda\Vert}^2)\big\vert_{\lambda = 0} =
\prod\limits_{\nu=1}^k (\tfrac{\partial}{\partial \lambda_\nu})^{m_\nu}
\exp(\lambda_\nu x - \tfrac{1}{2}\lambda_\nu^2)\big\vert_{\lambda_\nu = 0} 
\end{align}
for all $m \in \N_0^k$ and $x \in \R^k$. 

\noindent In fact, \eqref{eq:k_dim_Taylor_coeff} allows a further 
(purely analytic) proof of the Grothendieck \textit{equality} and its 
multi-dimensional generalisation. As a by-product, we provide a closed-form 
analytical representation of the multvariate distribution function of a 
Gaussian random vector $\textbf{X} \sim N_{2k}(0, \Sigma_{2k}(\rho))$ 
(cf. \sref{Proposition}{prop:closed_form_rep_of_multivariate_Gaussian_df}). In 
general, closed-form analytical representations of general multivariate Gaussian 
distribution functions are not available (cf., e.g., \cite{HS1980}). All that 
can be derived from the calculation of the following vital ``Fourier-Hermite 
coefficients'', which include $\text{sign}$ as a particular case (cf. 
\sref{Corollary}{cor:generalisation_of_Groth_equality} below):
\begin{theorem}\label{thm:GT_equality_and_beyond}
Let $k \in \N$, $m = (m_1, \ldots, m_k)^\top \in \N_0^k$, 
$a = (a_1, \ldots, a_k)^\top \in \R^k$ and 
$b = (b_1, \ldots, b_k)^\top \in \R^k$. Put $I_a : = 
\prod\limits_{i=1}^k [a_i, \infty)$ and $J_b : = \prod\limits_{i=1}^k (-\infty, b_i]$.
Then
\begin{align}\label{eq:Hermite_Fourier_coeff_2}
\int_{I_a}H_m\,\textup{d}\gamma_k = \langle \ind_{I_a}, H_m \rangle_{\gamma_k} = 
\prod\limits_{\stackrel{i=1}{m_i = 0}}^k (1-\Phi(a_i))\,
\prod\limits_{\stackrel{i=1}{m_i \not= 0}}^k \frac{1}{\sqrt{m_i}}
\varphi(a_i)H_{m_i-1}(a_i)
\end{align}
and
\begin{align}\label{eq:Hermite_Fourier_coeff_3}
\int_{J_b}H_m\,\textup{d}\gamma_k = \langle \ind_{J_b}, H_m \rangle_{\gamma_k} = 
(-1)^{c(k, m)}\prod\limits_{\stackrel{i=1}{m_i = 0}}^k \Phi(b_i)\,
\prod\limits_{\stackrel{i=1}{m_i \not= 0}}^k \frac{1}{\sqrt{m_i}}
\varphi(b_i)H_{m_i-1}(b_i),
\end{align}
where $c(k, m) : =  k - \sum_{i=1}^k\delta_{m_i\,0} \in \{0,1, \ldots, k\}$ 
counts the number of non-zero components of the vector $m$.  
\end{theorem}
\begin{corollary}\label{cor:generalisation_of_Groth_equality}
Let $k \in \N$, $m = (m_1, \ldots, m_k)^\top \in \N_0^k$, 
$a = (a_1, \ldots, a_k)^\top \in \R^k$ and $b = (b_1, \ldots, b_k)^\top \in 
\R^k$. Put $I_a : = \prod\limits_{i=1}^k [a_i, \infty)$ and $J_b : = 
\prod\limits_{i=1}^k (-\infty, b_i]$. Let $\chi_a : \R^k \longrightarrow \{-1,1\}$ and
$\psi_b : \R^k \longrightarrow \{-1,1\}$ be defined as
\[
\chi_a(x) : = 2\,\ind_{I_a}(x) - 1 
\text{ and } \psi_b(x) : = 1 - 2\,\ind_{J_b}(x) = 
-\chi_{-b}(-x).
\]
\begin{align}\label{eq:arcsin_CCP_via_Hermite_Fourier_coeff}
\langle \chi_a, H_0\rangle_{\gamma_k} = 2\prod\limits_{i=1}^k (1-\Phi(a_i)) - 1 
\text{ and } \langle \psi_b, H_0\rangle_{\gamma_k} = 1 - 2\prod\limits_{i=1}^k \Phi(b_i).
\end{align}
If $m \not=0$, then
\begin{align}\label{eq:arcsin_CCP_via_Hermite_Fourier_coeff_non_zero_case}
\begin{split}
{}&\langle \chi_a, H_m\rangle_{\gamma_k} = 
2\prod\limits_{\stackrel{i=1}{m_i = 0}}^k (1-\Phi(a_i))\,
\prod\limits_{\stackrel{i=1}{m_i \not= 0}}^k \frac{1}{\sqrt{m_i}}
\varphi(a_i)H_{m_i-1}(a_i)\\ 
&\text{and}\\ 
&\langle \psi_b, H_m\rangle_{\gamma_k} = 
2(-1)^{c(k,m)+1}\prod\limits_{\stackrel{i=1}{m_i = 0}}^k \Phi(b_i)\,
\prod\limits_{\stackrel{i=1}{m_i \not= 0}}^k \frac{1}{\sqrt{m_i}}
\varphi(b_i)H_{m_i-1}(b_i),
\end{split}
\end{align}
where $c(k, m) : =  k - \sum_{i=1}^k\delta_{m_i\,0} \in \{0,1, \ldots, k\}$. In 
particular,
\begin{align}\label{eq:Hermite_Fourier_coefficient_of_sign}
\langle\text{sign}, H_{2n}\rangle_{\gamma_1} = 0 \,\text{ and }\, 
\langle\text{sign}, H_{2n+1}\rangle_{\gamma_1} = (-1)^n \sqrt{\frac{2}{\pi}}\,
\frac{(2n-1)!!}{\sqrt{(2n+1)!}} = \frac{(-1)^n}{2^n}\,\sqrt{\frac{2}{\pi}}\,
\frac{(n+1)!}{\sqrt{(2n+1)!}}\,C_n
\end{align}
for any $n \in \N_0$, where $C_n : = \frac{1}{n+1}\binom{2n}{n}$ denotes the 
$n$'th Catalan number.
\end{corollary}
\noindent Moreover, because of the fact that $H_\nu(-x) = (-1)^\nu H_\nu(x)$ 
(respectively, $H_\nu(\vert x \vert) = ({\text{sign}}(x))^\nu H_\nu(x)$) for all 
$\nu \in \N_0$ and $x \in \R$, a $k$-fold multiplication of the values of 
one-dimensional Hermite polynomials implies that  
\[
H_n(-x) = (-1)^{\vert n \vert} H_n(x) \text{ for all } n \in \N_0^k \text{ and } 
x \in \R^k\,.
\]
Hence,
\begin{align}\label{eq:Fourier_coeff}
\langle g_\rho, H_n \rangle_{\gamma_k} =  \langle g, (H_n)_\rho \rangle_{\gamma_k}
= \rho^{\vert n \vert}\langle g, H_n \rangle_{\gamma_k} \text{ for all } g \in L^2(\gamma_k)
\text{ and } \rho \in \{-1,1\}\,,
\end{align}
where $\R^k \ni y \mapsto f_\rho(y) : = f({\text{sign}}(\rho)y) = f(\rho\,y)$ 
satisfies $\Vert f_\rho \Vert_{\gamma_k} = \Vert f \Vert_{\gamma_k}$ for $f \in 
\{g, H_n\}$ (which follows from a simple change of variables and the trivial 
fact that $\vert\rho^k\vert = 1$ for all $\rho \in \{-1,1\}$). Obviously, 
dependent on the smoothness structure of $g \in L^2(\gamma_k)$, it's quite a 
challenge to calculate the ``Fourier-Hermite coefficients'' 
$\langle g, (H_n)_\rho \rangle_{\gamma_k} = \int_{\R^k} g(x) H_n(x)
\gamma_k(\textup{d}x)$ explicitly. Regarding that particular problem, let us 
note a first important fact (cf. also 
\sref{Proposition}{prop:Sobolev_CCP_real_correl_case}):  
\begin{proposition}\label{prop:k_dim_Fourier_Hermite_coeff}
Let $k \in \N$, $\textbf{X} \sim N_k(0, I_k)$ and $g \in L^2(\gamma_k)$. Then the function 
$\R^k \ni \lambda \mapsto \E[g(\textbf{X}+\lambda)]$ is smooth around $0$. It 
satisfies
\begin{align}\label{eq:k_dim_Fourier_Hermite_coeff_transformation}
D^\alpha\,\E[g(\textbf{X}+\lambda)]\big\vert_{\lambda = 0} = \sqrt{\alpha!}\,
\langle g, H_\alpha \rangle_{\gamma_k} = \sqrt{\alpha!}\,
\E[g(\textbf{X})H_\alpha(\textbf{X})] \text{ for all } \alpha \in \N_0^k\,. 
\end{align}
In particular, if in addition $g$ is smooth, then
\[
\langle g, H_\alpha \rangle_{\gamma_k} = \E[g(\textbf{X})H_\alpha(\textbf{X})] = 
\frac{1}{\sqrt{\alpha!}}\,\E[D^\alpha\,g(\textbf{X})]
\]
for all $\alpha \in \N_0^k$.
\end{proposition}
\begin{comment}
Consider the function $\R^k \times \R^k \ni (x, \lambda) \mapsto 
F(x, \lambda) : = \exp(\langle x, \lambda \rangle_2 - \tfrac{1}{2}\,
\Vert \lambda \Vert_2)$. A straightforward application of the $k$-dimensional change of 
variables formula implies that $I_g(\lambda) : = \E[g(\textbf{X}+\lambda)] = 
\E[g(\textbf{X})F(\textbf{X}, \lambda)]$. Consequently, due to 
\eqref{eq:k_dim_generating_fct}, it follows that
\[
I_g(\lambda) = \sum_{\alpha \in \N_0^k}
\frac{\E[g(\textbf{X})H_\alpha^\ast(\textbf{X})]}
{\alpha!}\lambda^\alpha\,,
\]
where we put $H_\alpha^\ast(x) : = \sqrt{\alpha!}\,H_\alpha(x)$ for all 
$\alpha \in \N_0^k$ and $x \in \R^k$. 
Thus, the multi-indices version of the $k$-dimensional Taylor formula (cf. 
\cite[Theorem 2.8.4]{DK2004}) implies that
\[
D^\alpha\,\E[g(\textbf{X}+\lambda)]\big\vert_{\lambda = 0} =
D^\alpha\,I_g(\lambda)\big\vert_{\lambda = 0} = 
\E[g(\textbf{X})H_\alpha^\ast(\textbf{X})] =
\sqrt{\alpha!}\,\E[g(\textbf{X})H_\alpha(\textbf{X})].
\]
\end{comment}
\section{Real CCP functions and covariances: a Fourier-Hermite analysis 
approach}
Fix $f, g \in L^2(\R^k, \gamma_k), \nu \in \N_0$ and $k \in \N$. Put 
$C(\nu, k) : = \{n \in \N_0^k : \vert n \vert = \nu\}$ 
and 
\begin{align}\label{eq:coefficients_of_h_f_g}
p_\nu(f,g) : = \sum_{n \in C(\nu, k)}\langle f, H_n\rangle_{\gamma_k}\,
\langle g, H_n\rangle_{\gamma_k} = p_\nu(g,f)\,.
\end{align}
Since the inequality of arithmetic and geometric means in 
particular holds for any pair of elements, indexed by elements of the set of 
finitely many elements $C(\nu, k)$, it allows (the proof of) a direct transfer 
of H\"{o}lder's inequality by means of summation over $C(\nu, k)$, whence
\begin{align}\label{eq:implication_of_Hoelder_on_C_nu_k}
\vert p_\nu(f,g)\vert \leq \sqrt{p_\nu(f,f)}\sqrt{p_\nu(g,g)} 
\text{ for all } \nu \in \N_0\,. 
\end{align}
Since $\{C(\nu, k) : \nu \in \N_0\}$ obviously is a partition of the set $\N_0^k$ and 
$\{H_\alpha : \alpha \in \N_0^k\}$ is an orthonormal basis in $L^2(\gamma_k)$, 
H\"older's inequality again implies that
\begin{align}\label{eq:series_rep_of_real_h_fg}
[-1,1] \ni \rho \mapsto h_{f,g}(\rho) : = \sum_{n \in \N_0^k}
\langle f, H_n\rangle_{\gamma_k}\,\langle g, H_n\rangle_{\gamma_k}\,
\rho^{\vert n \vert} = \sum_{\nu = 0}^{\infty}p_\nu(f,g)\,\rho^\nu
\end{align}
converges absolutely, and
\begin{align}\label{eq:h_fg_is_bdd}
{
h_{f,g}(\rho) \leq \sum_{\nu = 0}^{\infty}\vert p_\nu(f,g)\vert \leq 
\Vert f \Vert_{\gamma_k} \Vert g \Vert_{\gamma_k}
}\,.
\end{align}
{Hence, $h_{f,g}$ is well-defined and bounded.} Note also that by construction 
$h_{f,g} = h_{g,f}$, $h_{f,g} = \frac{1}{4}(h_{f+g,f+g} - h_{f-g,f-g})$ (due to 
the polarisation equality) and $\vert h_{p,p}(\rho)\vert \leq 
\sum_{n \in \N_0^k}\langle p, H_n\rangle_{\gamma_k}^2\,
\vert\rho\vert^{\vert n \vert} \leq \Vert p \Vert_{\gamma_k}^2$ for all $p \in 
L^2(\R^k, \gamma_k)$ and $\rho \in [-1,1]$. 
Consequently, $h_{f,g}{\big\vert}_{(-1, 1)} \in W^\omega_+((-1, 1))$.  

Observe also that \textit{for any} $f = \sum_{n \in \N_0^k} a_n H_n \in 
L^2(\gamma_k)$ and \textit{any} $g = \sum_{n \in \N_0^k} b_n H_n \in 
L^2(\gamma_k)$, it follows that $a_n = x_n(f) : = \langle f, H_n \rangle_
{\gamma_k}$ and $b_n = x_n(g) = \langle g, H_n \rangle_{\gamma_k}$ for all 
$n \in \N_0^k$, implying that $p_\nu(f,g) = \sum_{n \in C(\nu, k)} a_n b_n$, 
whence
\[
h_{f,g}(\rho) = \sum_{\nu = 0}^\infty \big(\sum_{n \in C(\nu, k)} a_n b_n\big)\,
\rho^\nu \,\text{ for all }\, \rho \in [-1,1]
\]
and $\sum_{\nu = 0}^\infty\big(\sum_{n \in C(\nu, k)} a_n^2\big) = \sum_{n \in \N_0^k} 
a_n^2 = \Vert f \Vert_{\gamma_k}^2$. In particular,
\[
h_{f,f}(\rho) = \sum_{\nu = 0}^\infty \big(\sum_{n \in C(\nu, k)} a_n^2\big)\,
\rho^\nu \,\text{ for all }\, \rho \in [-1,1]\,.
\]
This immediately results in another important statement which should be 
compared with \cite[Proposition 5 and Th\'{e}or\`{e}me 3]{Kr1979} and 
\eqref{eq:entrywise_mapping_of_inner_products_to_inner_products}: 
\begin{proposition}\label{prop:inner_product_rep_of_h_f_g}
Let $k \in \N$ and $f, g \in S_{L^2(\gamma_k)}$. If $H$ is a real Hilbert space, 
then there exists a real Hilbert space $\H$ such that for any $u, v \in S_H$,
\[
h_{f,g}(\langle u, v \rangle_H) = \langle \psi_u(f), \psi_v(g) \rangle_{\H}\,,
\]
where $\psi_w : L^2(\gamma_k) \longrightarrow \H$ is a mapping which satisfies 
$\psi_w(S_{L^2(\gamma_k)}) \subseteq S_{\H}$ for any $w \in S_H$.
In particular, for any $m,n \in \N$, the following statements hold:
\[
h_{f,g}[S] \in \mathcal{Q}_{m,n} \text{ for all } S \in \mathcal{Q}_{m,n}
\]
and
\[
h_{f,g}[A] \in \M_n(\R)^+ \text{ for all } A \in \M_n(\R)^+\,.
\]
\end{proposition}
\noindent We will recognise soon that an \textit{additional boundedness assumption} 
on $f = \sum_{n \in \N_0^k} a_n H_n \in L^2(\gamma_k)$ itself is of utmost 
importance in relation to an approximation of the smallest upper bound of 
$K_G^\F$ (cf. \autoref{thm:GT_generalising_Krivine}, \autoref{thm:odd_real_bdd_CCP} 
and \autoref{thm:odd_complex_bdd_CCP_general_version}). To perform this highly 
non-trivial task, we have to look strongly for ``suitable'' $f = 
\sum_{n \in \N_0^k} a_n H_n \in L^2(\gamma_k)$ which in addition are bounded 
(a.s.); i.e., we have to look for some $M > 0$ such that (pointwise!) for 
($\gamma_k$-almost) all $x \in \R^k$,
\begin{align}\label{eq:boundedness_assumption}
\vert f(x)\vert \leq M
\end{align}
(as is the case with \eqref{eq:Hermite_rep_of_sign_in_L2_gamma_1}).  

We are now fully prepared to extend these important facts to one of 
our key results in this paper. In particular, we provide a multi-dimensional 
generalisation of the one-dimensional case $k=1$ (cf. 
\cite[Section 3.1]{BGGLT2019}) and specify a non-obvious 
tightening of the upper bound of $h_{f,g}$. Here, it should be noted that the 
inclusion of \sref{Proposition}{prop:char_of_the_Sigma_2n_prob_law} would allow to
view \autoref{thm:h_fg_real_case} as a straightforward simple implication of 
the key results in \cite[Chapter 11]{O'D2014}, including the consideration of 
\cite[Definition 11.10]{O'D2014} (cf. also \cite[Chapter 5.6.1]{B2019}, 
\cite[Lemma 2.2]{RS2009} and \sref{Remark}{rem:noise_stability}). However, we 
give a self-contained proof, built on the well-established Ornstein-Uhlenbeck 
semigroup (whose construction is recalled in the proof) and which sheds light 
also on the impact of the negative correlation case $-1 < \rho < 0$ in shape 
of an alternating sign change in the related power series.
\begin{theorem}\label{thm:h_fg_real_case}
Let $k \in \N, f, g \in L^2(\R^k, \gamma_k)$, $\textbf{S} \sim 
N_k(0, I_k)$, $\rho \in [-1,1]$ and $\vc{\textbf{X}}{\textbf{Y}} \sim 
N_{2k}(0, \Sigma_{2k}(\rho))$. Then the following properties hold:
\begin{enumerate}
\item
$h_{f,g}{\big\vert}_{(-1, 1)} \in W^\omega_+((-1, 1))$,
\begin{align}\label{eq:rep_of_h_fg}
h_{f,g}(0) = 
\E[f(\textbf{S})]\,\E[g(\textbf{S})] \text{ and }
h_{f,g}(\rho) = \E[f(\textbf{X})g(\textbf{Y})] = 
\text{Cov}(f(\textbf{X}), g(\textbf{Y})) + h_{f,g}(0).
\end{align}
In particular, 
$h_{f,g}(1) = \E[f(\textbf{S})g(\textbf{S})] = \langle f, g\rangle_{\gamma_k}$, 
$h_{f,f}(1) = \Vert f \Vert_{\gamma_k}^2$ and $h_{f,g}(-1) = 
\E[f(\textbf{S})g(\textbf{-S})] = \langle f, g_{-1}\rangle_{\gamma_k}$.  
\item $h_{f,g} : [-1,1] \longrightarrow \R$ is bounded and satisfies
\begin{align}\label{eq:boundedness_of_h}
\begin{split}
\vert h_{f,g}(\rho) \vert &{\leq \big(h_{f,g}{\big\vert}_{(-1, 1)}\big)_
{\text{abs}}(\vert \rho\vert) \leq \big(h_{f,g}{\big\vert}_{(-1, 1)}\big)_
{\text{abs}}(1)} {\leq \Vert f \Vert_{\gamma_k}\Vert g \Vert_{\gamma_k}} 
\end{split}
\end{align}
for all $\rho \in [-1,1]$.
\item 
If $\rho \in (-1,1)$, then 
\begin{align}\label{eq:h_fg_correl_representation}
\begin{split}
h_{f,g}(\rho) = \frac{1}{(2\pi)^k(1 - \rho^2)^{k/2}}\,\int_{\R^{k}}
\int_{\R^{k}}f(x)g(y)\exp\big(-\frac{\Vert x \Vert^2 + \Vert y \Vert^2 - 
2\rho\langle x,y \rangle}{2(1-\rho^2)}\big)\textup{d}^{k}x\,\textup{d}^{k}y.
\end{split}
\end{align} 
\item If $H$ is a real separable Hilbert space and $u, v \in S_H$, there 
exists a family $\{{\textbf{X}}_w : w \in S_H \}$ of $\R^k$-valued random 
vectors, such that $\vc{{\textbf{X}}_u}{{\textbf{X}}_v} \sim 
N_{2k}(0, \Sigma_{2k}(\langle u, v \rangle_H))$
and
\begin{align}\label{eq:quantum_corr_image}
h_{f,g}(\langle u, v \rangle_H) = 
\E[f({\textbf{X}}_u)g({\textbf{X}}_v)].  
\end{align}
\item If either $f$ or $g$ is odd, then $p_{2\nu}(f,g) = 0$ for all 
$\nu \in \N_0$ and
\[
{\text{cov}}(f(\textbf{X}), g(\textbf{Y})) = \E[f(\textbf{X})\,g(\textbf{Y})]
= h_{f,g}(\rho) = \rho\sum_{\nu = 0}^{\infty}p_{2\nu+1}(f,g) \rho^{2\nu}\,.
\]
In particular, $h_{f,g}$ is odd.
\end{enumerate}
\end{theorem}
\noindent Expressed in matrix notation, \autoref{thm:h_fg_real_case}-\ref{(iv)} 
leads directly to a further crucial result:
\begin{corollary}
Let $k, m, n \in \N$ and $f, g \in L^2(\gamma_k)$. Then the following matrix 
representations hold:
\begin{enumerate}
\item For any $S \equiv (s_{ij}) \in \mathcal{Q}_{m,n}$ there exist a 
random vector $\textbf{P}_f = (f(\textbf{X}_1), \ldots, f(\textbf{X}_m))^\top$ 
in $\R^m$ and a random vector $\textbf{Q}_g = (g(\textbf{Y}_1), \ldots, 
g(\textbf{Y}_n))^\top$ in $\R^n$, such that $(\textbf{X}_i, \textbf{Y}_j) \sim 
N_{2k}(0, \Sigma_{2k}(s_{ij}))$ for all $(i,j) \in [m] \times [n]$ and
\begin{align}\label{eq:splitting_in_mean_quantum_correl_case}
h_{f,g}[S] = \E[{\textbf{P}}_f {\textbf{Q}}_g^\top] =
\big(\E[f(\textbf{X}_i)g(\textbf{Y}_j)]\big)_{ij}. 
\end{align}
\item For any correlation matrix $\Sigma \equiv (\sigma_{ij}) \in 
C(n; \R)$ there exists a random vector $\textbf{R}_f = (f(\textbf{Z}_1), \ldots, 
f(\textbf{Z}_n))^\top$ in $\R^n\setminus\{0\}$, such that $(\textbf{Z}_i, 
\textbf{Z}_j) \sim N_{2k}(0, \Sigma_{2k}(\sigma_{ij}))$ for all $(i,j) \in [m] 
\times [n]$ and
\begin{align}\label{eq:splitting_in_mean_correlation_case}
h_{f,f}[\Sigma] = \E[\Theta_{f,f}] = \E[{\textbf{R}}_f {\textbf{R}}_f^\top] =
\big(\E[f(\textbf{Z}_i)f(\textbf{Z}_j)]\big)_{ij}\,,
\end{align}
where $\Theta_{f,f} : = {\textbf{R}}_f {\textbf{R}}_f^\top \in \M_n(\R)^+$ is a 
positive semidefinite random matrix of rank 1.
\end{enumerate}
\end{corollary}
\begin{corollary}\label{cor:h_fg_as_correl_coeff}
Let $k \in \N$, $f, g \in L^2(\R^k, \gamma_k)$, $r \in [-1,1]$, $\textbf{S} \sim 
N_k(0, I_k)$ and $\vc{\textbf{X}}{\textbf{Y}} \sim N_{2k}(0, \Sigma_{2k}(r))$. 
If $\Vert f \Vert_{\gamma_k} = 1$, $\Vert g \Vert_{\gamma_k} = 1$ and 
$\E[f(\textbf{S})] = \E[g(\textbf{S})] = 0$, then
\[
h_{f,g}(r) = \rho(f(\textbf{X}), g(\textbf{Y}))
\]
coincides with the Pearson correlation coefficient between the real random 
variables $f(\textbf{X})$ and $g(\textbf{Y})$.
\end{corollary}
{
Very recently, Krivine generalised the Grothendieck equality and constructed a 
related function $\Phi_{F,G} : [-1,1] \longrightarrow \{-1,1\}$ in the real case 
(cf. \cite{Kr2023}). In fact, if $\varepsilon > 0$ is fixed, 
\sref{Proposition}{prop:char_of_the_Sigma_2n_prob_law} implies that 
$\Phi_{F,G} = h_{f, g}$, where 
\[
\R^3 \ni (x_0, x_1, x_2)^\top \mapsto f(x_0, x_1, x_2) : = \text{sign}
\Big(\big\langle \colvec{2}{\cos(\sqrt{2}\varepsilon\,H_2(x_0))}
{\sin(\sqrt{2}\varepsilon\,H_2(x_0))}, \colvec{2}{x_1}{x_2}\big\rangle_{\R_2^2}\Big) =
\text{sign}(F(x_0, x_1, x_2))
\]
and
\[
\R^3 \ni (y_0, y_1, y_2)^\top \mapsto g(y_0, y_1, y_2) : = \text{sign}
\Big(\big\langle \colvec{2}{\cos(\sqrt{2}\varepsilon\,H_2(y_0))}
{-\sin(\sqrt{2}\varepsilon\,H_2(y_0))}, \colvec{2}{y_1}{y_2}\big\rangle_{\R_2^2}\Big) =
\text{sign}(G(y_0, y_1, y_2)).
\]
The function $\Phi_{F,G}$ in \cite{Kr2023} namely proves to be a beautiful 
example that fulfils the following statement (if $k = 1$ and $n = 2$):
\begin{proposition}\label{prop:Krivine_2023_generalised}
Let $k, n \in \N$, $\rho \in [-1,1]$ and $\vcfour{\textbf{X}_0}{\textbf{X}}
{\textbf{Y}_0}{\textbf{Y}}$ be a random vector, which maps into $\R^k \times 
\R^n \times \R^k \times \R^n \equiv \R^{2(k+n)}$, such that
$\vcfour{\textbf{X}_0}{\textbf{X}}{\textbf{Y}_0}{\textbf{Y}} \sim 
N_{2(k+n)}(0, \Sigma_{2(k+n)}(\rho))$. Let $F : \R^k \longrightarrow \S^{n-1} 
\subseteq \R^n$ and $G : \R^k \longrightarrow \S^{n-1} \subseteq \R^n$ be two 
arbitrary measurable functions. Put $\R^k\times \R^n \ni (x_0, x) \mapsto 
f(x_0, x) : = \text{sign}(\langle F(x_0), x\rangle_{\R_2^n})$ and 
$\R^k\times \R^n \ni (y_0, y) \mapsto g(y_0, y) : = 
\text{sign}(\langle G(y_0), y\rangle_{\R_2^n})$. Then

\scalebox{0.88}{
\vbox{
\begin{align*}
h_{f,g}(\rho) &= \frac{2}{\pi}\,
\E\big[\arcsin\big(\rho\,\langle F(\textbf{X}_0), G(\textbf{Y}_0)\rangle_
{\R_2^n}\big)\big]\\
&= \frac{2}{\pi}\int_{\R^{k}} \arcsin\big(\rho\,\langle F(x_0), G(y_0)\rangle_
{\R_2^n}\big)\,\exp\big(-\frac{\Vert x_0 \Vert^2 + \Vert y_0\Vert^2 - 
2\rho\langle x_0,y_0 \rangle}{2(1-\rho^2)}\big)
\frac{\textup{d}^{k}x_0\,\textup{d}^{k}y_0}{(2\pi)^k(1 - \rho^2)^{k/2}}.
\end{align*}
}}
\end{proposition}
}
In case of $f = g$ some important extra analytical facts emerge; 
particularly if in addition $f$ is odd (which is the relevant case for the topic 
of this paper):
\begin{theorem}\label{thm:h_ff_real_odd_case}
Let $k \in \N$ and $f \in L^2(\R^k, \gamma_k)$ be odd, such that 
$r : = \Vert f \Vert_{\gamma_k}^2 > 0$. Then the following properties hold:
\begin{enumerate}
\item $h_{f,f}{\big\vert}_{(-1, 1)} \in W^\omega_+((-1, 1))$. 
\item $h_{f,f}$ is a CCP function if and only if 
$r = 1$.
\item $h_{f,f} : [-1,1] \longrightarrow [-r, r]$ is an odd strictly 
increasing homeomorphism, which satisfies $h_{f,f}((-1,0)) = (-r,0)$ and 
$h_{f,f}((0,1)) = (0,r)$. Moreover,
\begin{align}\label{eq:derivative_of_CCP_at_zero}
0 \leq h_{f,f}^\prime(0) = \sum_{i=1}^k\big(\int_{\R^k}f(x)x_i
\gamma_k(\textup{d}x)\big)^2 = \int_{\R^{2k}} f(x)f(y)\langle x,y \rangle_{\R_2^k}
\gamma_{2k}(\textup{d}(x,y)).
\end{align}
\item If $h_{f,f}^\prime(0) > 0$, then $h_{f,f}^\prime(\rho) > 0$ for 
all $\rho \in (-1,1)$. In particular, $h_{f,f}^{-1}{\big\vert}_{(-r, r)} = 
\big(h_{f,f}{\big\vert}_{(-r,r)}\big)^{-1}$ is real analytic on $(-r, r)$ 
if and only if $h_{f,f}^\prime(0) > 0$. 
\item 
If $h_{f,f}^{-1}{\big\vert}_{(-r, r)} \in W^\omega_+((-r, r))$, then 
\begin{align}\label{eq:abs_h_f_f_inverse_estimation}
\big(h_{f,f}^{-1}{\big\vert}_{(-r, r)}\big)_{\text{abs}}(y) \geq 
\frac{1}{h_{f,f}^\prime(0)}\,y
\text{ for all } y \in [0, r]\,.
\end{align}
In particular, $s(r, f) : = \big(h_{f,f}^{-1}{\big\vert}_{(-r, r)}\big)_
{\text{abs}}(r) > 0$. If $\alpha > 0$ and $h_{f,f}^\prime(0) < \frac{r}{\alpha}$, 
then $s(r, f) > \alpha$.
Moreover,
\[
[-1,1] \ni t \mapsto \psi(t) : = \frac{1}{s(r, f)}\,
\big(h_{f,f}^{-1}{\big\vert}_{(-r, r)}\big)_{\text{abs}}(r\,t)
\]
is an odd, strictly increasing and homeomorphic CCP function. In particular, the
function $\big(h_{f,f}^{-1}{\big\vert}_{(-r, r)}\big)_{\text{abs}} : [-r,r] 
\stackrel{\cong}{\longrightarrow} [-s(r, f), s(r, f)]$ is strictly increasing, odd and 
homeomorphic, as well as its inverse. $\big(h_{f,f}^{-1}{\big\vert}_
{(-r, r)}\big)_{\text{abs}}(y) = s(r, f)\,\psi(\frac{y}{r})$ for all $y \in [-r,r]$ 
and $\big(\big(h_{f,f}^{-1}{\big\vert}_{(-r, r)}\big)_{\text{abs}}\big)^{-1}(x) = 
r\,\psi^{-1}(\frac{x}{s(r, f)})\, \text{ for all }\, x \in [-s(r, f), s(r, f)]$.
\end{enumerate}
\end{theorem}
\noindent A first implication for \textit{odd} CCP functions is a strong 
improvement of the boundedness condition in \autoref{thm:h_fg_real_case}-\ref{(ii)}; 
induced by the Schwarz lemma from complex analysis: 
\begin{proposition}\label{prop:odd_CCP_and_Schwarz_lemma}
Let $k \in \N$ and $f \in S_{L^2(\gamma_k)}$ be odd. Assume that $h_{f,f}(r) 
\not= r$ for some $r \in (0,1)$. Then
\begin{align}\label{eq:odd_CCP_and_Schwarz_lemma}
h_{f,f}(\rho) < \rho\,\text{ for all }\, \rho \in (0,1)\, \text{ and }\,
h_{f,f}(\tau) > \tau\,\text{ for all }\, \tau \in (-1,0).
\end{align}
Moreover, $0 \leq h_{f,f}^\prime(0) < 1$. If in addition, $h_{f,f}^{-1}
{\big\vert}_{(-1,1)} \in W^\omega_+((-1, 1))$, then
\begin{align}\label{eq:inverse_of_odd_CCP_and_CRA_implication_for_abs}
\big(h_{f,f}^{-1}{\big\vert}_{(-1,1)}\big)_{\text{abs}}(1) > 1\,.
\end{align}
\end{proposition}
\noindent All of a sudden we end up with an important general and non-obvious 
structural statement about \textit{odd} CCP functions; namely:
\begin{corollary}\label{cor:inverse_of_non_trivial_odd_CCP_is_not_CCP}
Let $k \in \N$ and $f \in S_{L^2(\gamma_k)}$ be odd. Assume that $h_{f,f}(r) 
\not= r$ for some $r \in (0,1)$. Then $h_{f,f}^{-1}$ is not a CCP function. 
\end{corollary}
\begin{comment}
We just have to link \autoref{thm:CCP_via_Schoenberg}-\ref{(v)} and 
\eqref{eq:inverse_of_odd_CCP_and_CRA_implication_for_abs}.
\end{comment}
\noindent We will recognise soon that \autoref{thm:h_ff_real_odd_case} is 
strongly linked with the value of the Grothendieck constant $K_G^\F$ (cf. Theorem 
\ref{thm:odd_real_bdd_CCP}). We namely obtain in a natural way another 
significant new definition; the so-called ``hyperbolic CCP transform''. To make 
this concept understandable, recall \sref{Lemma}{lem:real_analytic_and_l1} 
and \sref{Remark}{rem:odd_psi_abs_and_sign_condition} and reconsider the odd CCP 
function $\psi : = h_{\text{sign}, \text{sign}} = \frac{2}{\pi}\arcsin$; i.e., 
the Grothendieck function. Since $\psi^{-1} = \sin(\frac{\pi}{2}\,\cdot)$ on 
$[-1,1]$, it follows that $\big(\psi^{-1}{\big\vert}_{(-1,1)}\big)_{\text{abs}}
(\tau) = \sinh(\frac{\pi}{2}\tau) = \frac{1}{i}\sin(\frac{\pi}{2}\,i\tau)$ for 
all $\tau \in [-1,1]$. Note that $\big(\psi^{-1}{\big\vert}_{(-1,1)}\big)_
{\text{abs}}(1) = \sinh(\pi/2) > 1$, implying that $[-1,1] \subseteq 
(-\sinh(\pi/2), \sinh(\pi/2))$. Thus, 
\begin{align}\label{eq:hyp_trsfm_of_Krivine_fct}
\begin{split}
\frac{2}{\pi}\,\ln(y+\sqrt{y^2+1}) &= \frac{2}{\pi}\sinh^{-1}(y) = 
\big(\big(\psi^{-1}{\big\vert}_{(-1,1)}\big)_{\text{abs}} \big)^{-1}(y)\\
&= \frac{2}{\pi}\big(\frac{1}{i}\,\sin^{-1}(i\,y)\big) = \frac{1}{i}\,\widetilde{\psi}(iy)
\end{split}
\end{align}
for all $y \in [-1, 1]$. The Maclaurin series representation of the real CCP 
function $\psi$ therefore implies the Maclaurin series of $y \mapsto 
\frac{1}{i}\,\widetilde{\psi}(iy)$ on (the whole of) $\R$ is given by 
\[
\frac{1}{i}\,\widetilde{\psi}(iy) 
\stackrel{\eqref{eq:Gauss_Copula_special_case_arcsin}}{=} 
\frac{2}{\pi}\sum_{\nu = 0}^\infty (-1)^\nu\,\frac{((2\nu-1)!!)^2}{(2\nu+1)!} 
y^{2\nu+1} = -\frac{1}{3\pi} + \sum_{\stackrel{\nu = 0}{\nu \not= 1}}^\infty 
(-1)^\nu\,\frac{((2\nu-1)!!)^2}{(2\nu+1)!} y^{2\nu+1}\,. 
\]
Consequently, it follows that
\[
\big(\big(\psi^{-1}{\big\vert}_{(-1,1)}\big)_{\text{abs}} \big)^{-1}(\rho) = 
\frac{1}{i}\,\widetilde{\psi}(i\rho) < \psi(\rho) < \rho\,\text{ for all }\,\rho 
\in (0, 1)\,.
\]
\autoref{thm:h_ff_real_odd_case} and 
\sref{Proposition}{prop:odd_CCP_and_Schwarz_lemma}, together with 
\sref{Remark}{rem:odd_psi_abs_and_sign_condition} imply that the Grothendieck 
function is a special case of
\begin{lemma}[\textbf{Hyperbolic CCP transform}]
\label{lem:hyperbolic_CCP_transform}
Let $k \in \N$ and $\psi = h_{f,f}$, where $f \in S_{L^2(\gamma_k)}$ is odd. 
Assume that $\psi^{-1}{\vert}_{(-1,1)} \in W^\omega_+((-1,1))$ and $\psi(r) 
\not= r$ for some $r \in (0,1)$. Let the complex function $F : = 
$ \rlap{\raisebox{-0.9ex}{$\widesim{\phantom{\psi{\vert}_{(-1,1)}}}$}}
$\psi{\vert}_{(-1,1)}$ $: \overline{\D} \longrightarrow 
\overline{\D}$ be defined as in \sref{Lemma}{lem:real_analytic_and_l1}. Then
$F(\D) \subseteq \D$ and
\[
F(z) = \sum_{\nu=0}^\infty p_{2\nu+1}(f,f)\,z^{2\nu+1}\,\text{ for all }\, z \in 
\overline{\D}\,,
\]
where each $p_{2\nu+1}(f,f) \in [0, \infty)$ satisfies 
\eqref{eq:coefficients_of_h_f_g}. Put $s^\ast : = (\psi^{-1}{\vert}_
{(-1,1)})_{\text{abs}}(1)$. Then $s^\ast > 1$ and
\[
\psi^{\text{hyp}} : = ((\psi^{-1}{\vert}_{(-1,1)})_{\text{abs}})^{-1} : 
[-s^\ast, s^\ast] \stackrel{\cong}{\longrightarrow} [-1,1]
\]
is an odd, strictly increasing homeomorphism, which satisfies $[-1,1] \subseteq 
(-s^\ast, s^\ast)$, $\psi^{\text{hyp}}((0,1]) \subseteq (0,1)$ and 
$\psi^{\text{hyp}}([-1,0)) \subseteq (-1,0)$. Moreover, 
\begin{align}\label{eq:psi_hyp_at_1_vs_derivative_of_psi_at_zero}
0 < \psi^{\text{hyp}}(1) \leq \psi^\prime(0) < 1\,,
\end{align}
and the following two statements hold:
\begin{enumerate}
\item
\[
\vert\psi^{\text{hyp}}(\rho)\vert < \psi(\vert\rho\vert) < \vert\rho\vert \,
\text{ for all }\,\rho \in (-1, 1)\setminus\{0\}\,.
\]
\item If $\text{sign}((\psi^{-1}{\vert}_{(-1,1)})^
{(2n+1)}(0)) = (-1)^n$ for all $n \in \N_0$, then
\[
\psi^{\text{hyp}}(x) = \sum_{\nu=0}^\infty (-1)^\nu\,p_{2\nu+1}(f,f)\,x^{2\nu+1} =
\frac{1}{i}\,F(ix) = (M_{-i}\circ F\circ M_i)(x)\,\text{ for all }\, x \in [-1, 1]\,.
\]
{In particular, 
\begin{align}\label{eq:sign_condition_and_hyp_CCP_fct_value}
\psi^{\text{hyp}}(1) = \frac{1}{i}\,F(i). 
\end{align}
}
\end{enumerate}
\end{lemma}
\noindent Next, consider for example, the function $a : = \frac{1}{\sqrt{2}}(1 + H_2) : 
\R \longrightarrow \R$. Then $a \in S_{L^2(\gamma_1)}$ and $[-1,1] \ni \rho 
\mapsto h_{a,a}(\rho) = \frac{1}{2}(1 + \rho^2)$ defines an \textit{even} CCP 
function, such that $h_{a,a}(0) = \frac{1}{2} > 0$. With this example in mind, the 
assumption in the second part of the following result is not empty.
\begin{proposition}\label{prop:product_structure}
Let $k \in \N$, $\rho \in [-1,1]$ and $a,b, f, g \in L^2(\R^k, \gamma_k)$. Then
$a \otimes f \in L^2(\R^{2k}, \gamma_{2k})$ and $b \otimes g \in 
L^2(\R^{2k}, \gamma_{2k})$, and
\[
h_{a \otimes f, b \otimes g} = h_{a,b}\cdot h_{f,g}\,. 
\]
In particular, the product $h^{\text{ev}}\,h^{\text{odd}}$ of a continuous, even 
CCP function and a continuous, odd CCP function is an odd, strictly increasing 
homeomorphic CCP function. If $(h^{\text{odd}})^\prime(0) > 0$ and 
$h^{\text{ev}}(0) > 0$, then $(h^{\text{ev}}\,h^{\text{odd}})^{-1}{\big\vert}_
{(-1,1)}$ is real analytic. A product of two odd CCP functions is an even CCP 
function.
\end{proposition}
\noindent If we apply \autoref{thm:h_fg_real_case} to a pair of $k$-dimensional 
Hermite polynomials and recall that $\{H_\alpha : \alpha \in \N_0^k\}$ actually is an 
orthonormal basis in $L^2(\gamma_k)$ (cf. \eqref{eq:k_dim_Hermite_polynomials}), we 
immediately reobtain another remarkable property of Hermite polynomials (cf. 
\cite[Lemma 1.1.1]{N2009} and \cite[Proposition 11.33]{O'D2014}).
\begin{corollary}\label{cor:correlated_k_dim_Hermite_polynomials}
Let $k \in \N$, $\rho \in [-1, 1]$, $\alpha = (\alpha_1, \ldots, \alpha_k) \in \N_0^k$ and 
$\beta = (\beta_1, \ldots, \beta_k)^\top \in \N_0^k$. If 
$\vc{\textbf{X}}{\textbf{Y}} \sim N_{2k}(0, \Sigma_{2k}(\rho))$, then
\[
h_{H_\alpha, H_\beta}(\rho) = \E[H_\alpha(\textbf{X})\,H_\beta(\textbf{Y})] 
= \delta_{\alpha,\beta} \rho^{\vert\alpha\vert} = 
\prod\limits_{i=1}^k \delta_{\alpha_i, \beta_i}\,\rho^{\alpha_i}.
\]
\end{corollary}
\begin{remark}[\textbf{Noise stability}]\label{rem:noise_stability}
Fix $k \in \N$. Let $f \in L^2(\gamma_k)$ and $A, B$ Borel sets in $\R^k$. Let
$\rho \in [-1,1]\setminus\{0\}$ and $\vc{\textbf{X}}{\textbf{Y}} 
\sim N_{2k}(0, \Sigma_{2k}(\rho))$. Then $\text{sign}(\rho)\textbf{X} 
\sim N_k(0, I_k)$. Consequently, \eqref{eq:parametr_of_OU_sg} implies that
\begin{align}\label{eq:noise_stability_operator}
T_{\vartheta(\vert \rho \vert)} f_\rho(y) =
\E[f(\rho\,y + \sqrt{1-\rho^2}\,\textbf{X})] = 
\int_{\R^k}f(\rho\,y + \sqrt{1-\rho^2}\,x)\gamma_k(\textup{d}x) = U_\rho f(y) 
\end{align}
for all $y \in \R^k$, where $U_\rho$ is the Gaussian noise operator 
(cf. \cite[Definition 11.12]{O'D2014}). $U_\rho$ 
also is well-defined for $\rho = 0$, with constant value 
$U_0 f = \E[f(\textbf{X})]$. Observe that \eqref{eq:noise_stability_operator} implies 
that $U_\rho = T_{\vartheta(\vert \rho \vert)} M_\rho$, where the isometry 
$M_\rho = M_\rho^{-1} \in {\mathfrak{L}}(L^2(\gamma_k),L^2(\gamma_k))$ is given 
by $M_\rho f : = f_\rho$. A special case of \eqref{eq:OU_correl_rep} is given by
\[
\P(\textbf{X} \in A, \textbf{Y} \in B) = 
\E[{\ind_{A}}(\textbf{X}){\ind_{B}}(\textbf{Y})] = 
h_{{\indsm_{A}},\indsm_{B}}(\rho) = 
\langle \ind_{A}, U_\rho \ind_{B} \rangle_{\gamma_k}. 
\] 
Therefore, \autoref{thm:h_fg_real_case} (under inclusion of 
\sref{Proposition}{prop:char_of_the_Sigma_2n_prob_law}) encompasses the key concept of 
Gaussian noise stability, most commonly introduced as ${\textbf{Stab}}_\rho[f] : = 
\langle f, U_\rho f \rangle_{\gamma_k}$ ($\rho \in [-1,1]$, $f \in L^2(\gamma_k)$). 
Gaussian noise stability also comprises deep connections to geometry of minimal 
surfaces, hypercontractivity, isoperimetric inequalities, communication complexity 
and Gaussian copulas. A very comprehensive introductory processing of these 
topics can be found in \cite[Chapter 11]{O'D2014} and the references therein, 
including the seminal results of E. Mossel and J. Neeman.
\end{remark} 
\noindent In fact, it can be verified that $T_t$ (and hence $U_\rho = 
T_{\vartheta(\vert \rho \vert)} M_\rho$) is even a \textit{nuclear} operator, 
implying that each $T_t$ (and each $U_\rho$) in particular is a 
\textit{compact} Hilbert-Schmidt operator! More precisely, we have
\begin{proposition}\label{prop:OU_semigroup_elements_are_nuclear}
Let $k \in \N$, $t \geq 0$ and $\rho \in (-1,1)\setminus\{0\}$. 
Then both, $T_t \in \mathfrak{L}(L^2(\gamma_k), L^2(\gamma_k))$, and 
$U_\rho \in \mathfrak{L}(L^2(\gamma_k), L^2(\gamma_k))$ are Hilbert-Schmidt 
operators, satisfying
\[
\Vert T_t \Vert_{{\mathfrak{S}}_2} =  \frac{1}{(1-e^{-2t})^{k/2}} 
\text{ and } \Vert U_\rho \Vert_{{\mathfrak{S}}_2} = \frac{1}{(1 - \rho^2)^{k/2}}.
\]
$T_t$ as well as $U_\rho$ are even nuclear, and
\begin{enumerate} 
\item 
$\Vert T_t \Vert_{\mathfrak{N}} \leq \frac{1}{(1 - e^{-t})^k}$.
\item 
$\Vert U_\rho \Vert_{\mathfrak{N}} \leq \frac{1}{(1 - \vert\rho\vert)^k}$.
\end{enumerate}
\end{proposition}
\section{Examples of real CCP functions, Gaussian copulas and an extension of 
Stein's lemma}
As was to be expected, \autoref{thm:h_fg_real_case}
and \autoref{thm:h_ff_real_odd_case} give us first non-trivial examples, such 
as  
\[
h_{\text{sign},\text{sign}} : [-1,1] \longrightarrow [-1,1], x \mapsto
\frac{2}{\pi}\arcsin(x) = \frac{2}{\pi}\,x\,
{}_2 F_1\big(\frac{1}{2}, \frac{1}{2}, \frac{3}{2}; x^2 \big)
\] 
in the one-dimensional real case (implying the Grothendieck equality) and 
\[
h_{f_2, f_2} : [-1,1] \longrightarrow [-1,1], \tau \mapsto \frac{\pi}{4}  
\tau\,{}_2 F_1\big(\frac{1}{2}, \frac{1}{2}, 2; \tau^2 \big)\,
\]
in the one-dimensional complex case (implying the Haagerup equality), where 
\[
\R^2 \ni x \mapsto f_2(x) : = \begin{cases}\sqrt{2}\,
\frac{x_1}{\Vert x \Vert_2}&\text{if } x \not= 0\\ 
0 &\text{if } x = 0\end{cases}\,.
\]
(cf. \cite[Lemma 3.2. and Proof of Theorem 3.1]{H1987} and \sref{Example}{ex:Haagerup}). 
If we namely 
apply \autoref{thm:3F2_rep} to $m = 1$ and arbitrary $k \in \N$, we will 
recognise that \autoref{thm:h_fg_real_case} and Theorem 
\ref{thm:h_ff_real_odd_case} lead us to CCP functions 
$h_{f,f}$, where $f \in S_{L^2(\gamma_k)}$ is \textit{bounded a.e.}. These 
examples also include \cite[Lemma 2.1]{BOFV2014} as a special case. More 
precisely, we have:  
\begin{proposition}\label{prop:2F1_rep}
Let $k \in \N$, $\rho \in [-1,1]$ and $\vc{\textbf{X}}{\textbf{Y}} \equiv 
(X_1, \ldots, X_k, Y_1, \ldots, Y_k)^\top \sim N_{2k}(0, \Sigma_{2k}(\rho))$.
Consider the odd function
\[
\R^k \ni x \mapsto f_k(x) : = \begin{cases}\sqrt{k}\,
\frac{x_1}{\Vert x \Vert_{\R_2^k}}&\text{if } x \not= 0\\ 
0&\text{if } x = 0\end{cases}\,.
\]
Then $f_k \in S_{L^2(\gamma_k)} \cap L^\infty(\gamma_k)$ and 
$\Vert f_k \Vert_{\infty} = \sqrt{k}$. The function $h_{f_k,f_k} : [-1,1] 
\longrightarrow [-1,1]$ is an odd strictly increasing homeomorphism which is CCP 
and satisfies $h_{f_k,f_k}{\big\vert}_{(-1, 1)} \in W^\omega_+((-1, 1))$. Let 
$c_k$ be defined as in \sref{Lemma}{lem:Gamma_fct_and_constants_c_k}. Then
$0 < h_{f_k,f_k}^\prime(0) = c_k^2 < 1$. $h_{f_k,f_k}^{-1}{\big\vert}_{(-1, 1)}$ is 
real analytic, and
\begin{align}\label{eq:k_dim_elliptic_integral_and_h_ff}
\begin{split}
h_{f_k,f_k}(\rho) &= k\,\E\big[\frac{X_1\,Y_1}{\Vert\textbf{X}\Vert_{\R_2^k}\,
\Vert\textbf{Y} \Vert_{\R_2^k}}\big] = \E\big[\big\langle\frac{\!\!\!\!\!\textbf{X}}
{\Vert\textbf{X}\Vert_{\R_2^k}},\frac{\!\!\!\!\!\textbf{Y}}{\Vert\textbf{Y}\Vert_
{\R_2^k}}\big\rangle_{\R_2^k}\big]\\ 
&= c_k^2\,\rho\,{}_2 F_1\big(\frac{1}{2}, \frac{1}{2},\frac{k+2}{2}; \rho^2 \big) 
= c_k^2\,k\,!!\sum_{n=0}^\infty \frac{((2n-1)!!)^2}{(2n)!!\,(2n+k)!!}\,\rho^{2n+1} 
\end{split}
\end{align}
In particular, the function $[-1,1] \ni \rho \mapsto h_{f_k,f_k}(\rho) - 
c_k^2\,\rho$ is CCP as well.
If $\vert \rho \vert < 1$, then
\begin{align}\label{eq:k_dim_elliptic_integral_repr_of_h_ff}
h_{f_k,f_k}(\rho) = \sqrt{\frac{2k}{\pi}}\,c_k\,
\rho\,\int_0^1 \frac{(\sqrt{1-t^2})^{k-1}}{\sqrt{1-\rho^2\,t^2}}\,\textup{d}t\,.
\end{align}
Moreover, 
\[
\E\big[\frac{X_i\,Y_i}{\Vert \textbf{X}\Vert_{\R_2^k}\,\Vert\textbf{Y}\Vert_
{\R_2^k}}\big] = 
\frac{c_k^2}{k}\,\rho\,{}_2 F_1\big(\frac{1}{2}, \frac{1}{2},\frac{k+2}{2}; \rho^2 \big) 
\,\text{ and }\,\E\big[\frac{\!\!X_i^2}{\Vert\textbf{X} \Vert_{\R_2^k}^2}\big] = 
\frac{1}{k}
\]
for all $i \in [k]$. 
\end{proposition}
\begin{remark}
Let $k \in \N$ and $c_k$ be defined as in 
\sref{Lemma}{lem:Gamma_fct_and_constants_c_k}. Then \cite[Theorem 1]{BBT2011} 
in fact implies that
\[
\frac{\pi}{2}\,c_d^2 \leq K_G^\R(d) \,\text{ for all }\, d \in \N_3\,.
\]
\end{remark}
\begin{remark}[\textbf{Krivine rounding scheme reconsidered}] 
Fix $k \in \N$. Due to \sref{Proposition}{prop:2F1_rep}, the following set of 
real-valued functions is non-empty:
\begin{align}\label{eq:set_of_attainable_functions}
{\mathcal{G}}_k := \big\{f : f \in S_{L^2(\gamma_k)}, f \text{ is odd}, 
h_{f,f}^{\prime}(0) > 0\big\}. 
\end{align}
Let $f \in {\mathcal{G}}_k$. If in addition, $\vert f \vert = 1$ and 
$h_{f,f}^{-1}{\big\vert}_{(-1, 1)} \in W^\omega_+((-1, 1))$ 
\sref{Lemma}{lem:hyperbolic_CCP_transform} implies that $\{f\circ\sqrt{2}, f\circ\sqrt{2}\}$ 
is a Krivine rounding scheme, introduced in \cite[Definition 2.1]{BMMN2013} 
(since $h_{f,f}$ coincides with the function $H_{f\circ\sqrt{2}, f\circ\sqrt{2}}$, 
introduced in \cite[Definition 2.1]{BMMN2013} and $c(f\circ\sqrt{2}, 
f\circ\sqrt{2}) = h_{f,f}^{\text{hyp}}(1) \in (0, 1)$). The latter fact should 
be compared with \autoref{thm:odd_real_bdd_CCP}\,!
\end{remark}
\noindent Now it is no longer surprising that \sref{Proposition}{prop:2F1_rep} 
encompasses the Grothendieck equality and the Haagerup equality as particular cases. 
\begin{example}[\textbf{$k=1$ (Grothendieck)}]
Fix $\rho \in [-1,1]$. It is well-known that $\arcsin(\rho) = \rho\,
{}_2 F_1\big(\frac{1}{2}, \frac{1}{2},\frac{3}{2}; \rho^2 \big)$, whence
$h_{f_1, f_1}(\rho) = \frac{2}{\pi}\arcsin(\rho)$.
\end{example}
\begin{example}[\textbf{$k=2$ (Haagerup)}]
Let $n \in \N$, $u, v \in \C^n$ such that $\Vert u \Vert = 1$ and 
$\Vert v \Vert = 1$. Let $Z \sim {\C}N_n(0, I_n)$. The Haagerup equality 
\eqref{eq:Haagerup}, together with \eqref{eq:integral_repr_of_the_Haagerup_function} 
implies that
\[
\E[{\text{sign}}(u^\ast Z){\text{sign}}(\overline{v^\ast Z})] = 
{\text{sign}}(u^\ast v)\,h_{f_2, f_2}(\vert u^\ast v \vert).
\]
\end{example}
\noindent In the following, we shed some light on the underlying structure of 
the function $h_{\psi_a, \psi_b}$, where $\psi_p : = 1 - 2\,\ind_{J_p} \in 
S_{L^2(\gamma_k)}$ and $J_p : = \prod\limits_{i=1}^k (-\infty, p_i]$ for all $p \equiv 
(p_1, \ldots, p_k)^\top \in \R^k$ (introduced in 
\sref{Corollary}{cor:generalisation_of_Groth_equality}). Recall here that for any 
$n \in \N$, $\Sigma \in \M_n(\R)^+$ and $\textbf{X} \equiv 
(X_1, \ldots, X_{n})^\top \sim N_{n}(0, \Sigma)$, 
$\Phi_{0, \Sigma} : \R^n \longrightarrow [0,1]$ denotes the $n$-variate 
distribution function of $\textbf{X}$, i.e.,
\[
\Phi_{0, \Sigma}(x) : = F_{\textbf{X}}(x) = \P\big(\bigcap_{i=1}^n 
\{X_i \leq x_i\}\big) \text{ for all } x \equiv (x_1, \ldots, x_n)^\top 
\in \R^n\,.
\] 
\begin{proposition}\label{prop:closed_form_rep_of_multivariate_Gaussian_df}
Let $k \in \N$, $\rho \in [-1,1]$ and $\textbf{X} \equiv (X_1, \ldots, X_{2k})^\top 
\sim N_{2k}(0, \Sigma_{2k}(\rho))$. Let $a, b \in \R^k$ and $x \equiv 
(x_1, \ldots, x_{2k})^\top = \vc{a}{b}$. Then
\begin{align}\label{eq:rep_of_h_psi_a_psi_b}
h_{\psi_a, \psi_b}(\rho) = 
4 \Phi_{0, \Sigma_{2k}(\rho)}(x) + 1 - 2 \Phi_{0, I_k}(a) - 2 \Phi_{0, I_k}(b).
\end{align}
Furthermore,
\begin{align*}
\Phi_{0, \Sigma_{2k}(\rho)}(x) = \Phi_{0, I_k}(a)\Phi_{0, I_k}(b) + 
\sum_{\nu=1}^\infty d_\nu(x; k)\,\rho^\nu\,,
\end{align*}
where
\begin{align*}
d_\nu(x; k) : = \sum_{m \in C(\nu, k)}
\big(\prod\limits_{\stackrel{i=1}{m_i = 0}}^{2k}\,\Phi(x_i)
\prod\limits_{\stackrel{i=1}{m_i \not= 0}}^{2k} \frac{1}{\sqrt{m_i}}
\varphi(x_i)H_{m_i-1}(x_i)\big).
\end{align*}
\end{proposition}
\begin{comment}
Fix $\textbf{X} \equiv (X_1, \ldots, X_{2k})^\top \sim 
N_{2k}(0, \Sigma_{2k}(\rho))$ and $x = \vc{a}{b}$. Since $\psi_a(x_1)\,
\psi_b(x_2) = (1 - 2\,\ind_{J_a}(x_1))(1 - 2\,\ind_{J_b}(x_2)) = 
1 - 2\,\ind_{J_a}(x_1) - 2\,\ind_{J_b}(x_2) + 4\,\ind_{J_a}(x_1)\ind_{J_b}(x_2)$ 
for all $x_1, x_2 \in \R^k$ and $\textbf{X} = \vc{\textbf{X}_1}{\textbf{X}_2}$, 
where $\textbf{X}_1 : = (X_1, \ldots, X_{k})^\top \sim N_{k}(0, I_k)$ and 
$\textbf{X}_2 : = (X_{k+1}, \ldots, X_{2k})^\top \sim N_{k}(0, I_k)$, an 
application of \autoref{thm:h_fg_real_case} to the function 
$h_{\psi_a, \psi_b}$ implies that
\begin{align*}
F_{\textbf{X}}(x) &= \P\big(\bigcap_{i=1}^{2k} \{X_i \leq x_i\}\big) =
\E[\ind_{J_a}(\textbf{X}_1)\ind_{J_b}(\textbf{X}_2)] = 
\frac{1}{4}(h_{\psi_a, \psi_b}(\rho) - 1 + 2\Phi_{0, I_k}(a) + 2\Phi_{0, I_k}(b))\\
&= \frac{1}{4}\big((1 - 2\Phi_{0, I_k}(a))(1 - 2\Phi_{0, I_k}(b)) + 
\sum_{\nu=1}^\infty p_\nu(\psi_a, \psi_b) \rho^\nu - 1 + 2\Phi_{0, I_k}(a) + 
2\Phi_{0, I_k}(b)\big)\\
&= \Phi_{0, I_k}(a)\Phi_{0, I_k}(b) + \frac{1}{4}\sum_{\nu=1}^\infty 
\big(\sum_{m \in C(\nu, k)}\langle \psi_a, H_m\rangle_{\gamma_k}\,
\langle \psi_b, H_m\rangle_{\gamma_k}\big)\rho^\nu\,.
\end{align*}
The third equality is equivalent to \eqref{eq:rep_of_h_psi_a_psi_b}. 
\sref{Corollary}{cor:generalisation_of_Groth_equality} clearly finishes the 
proof.
\end{comment}
\noindent \sref{Proposition}{prop:closed_form_rep_of_multivariate_Gaussian_df} 
immediately gives us another significant example - which contains the Grothendieck 
equality as a special case again. Due to the famous Theorem of Sklar, it 
emerges from a lurking multivariate Gaussian copula; i.e., from a certain 
finite-dimensional distribution function with uniformly distributed marginals 
(cf., e.g., \cite{Oe2015} and the references therein). 
\begin{example}[\textbf{A lurking $\Sigma_{2k}(\rho)$-Gaussian copula}]
Let $k \in \N$, $m = (m_1, \ldots, m_k)^\top \in \N_0^k$, 
$\alpha = (\alpha_1, \ldots, \alpha_k)^\top \in (0,1)^k$ and 
$\beta = (\beta_1, \ldots, \beta_k)^\top \in 
(0,1)^k$. Let $a = (a_1, \ldots, a_k)^\top \in \R^k$ and 
$b = (b_1, \ldots, b_k)^\top \in \R^k$, where
$a_i : = \Phi^{-1}(\alpha_i)$ and $b_i : = \Phi^{-1}(\beta_i)$ for all $i \in 
[k]$. Then   
\[
h_{\psi_a,\psi_b}(\rho) = 1 - 2 \big(\prod\limits_{i=1}^{k}\alpha_i + 
\prod\limits_{i=1}^{k}\beta_i \big) + 4\,c_{\Sigma_{2k}(\rho)}(\alpha, \beta) 
\text{ for all } \rho \in [-1,1],
\]
where 
\[
(0,1)^k \ni (u_1, \ldots, u_k) \mapsto c_{\Sigma_{2k}(\rho)}(u_1, \ldots, u_k) : = 
\Phi_{0, \Sigma_{2k}(\rho)}(\Phi^{-1}(u_1),\ldots, \Phi^{-1}(u_k))
\] 
denotes the $2k$-dimensional Gaussian copula with respect to the correlation 
matrix $\Sigma_{2k}(\rho)$ (cf., e.g., \cite{Oe2015}). Consequently, if 
$(\alpha, \beta) = (\tfrac{1}{2},\tfrac{1}{2})$, it follows that $a = 0$, 
$b = 0$ and
\begin{align}\label{eq:Gauss_Copula_special_case}
h_{\psi_0,\psi_0}(\rho) = 1 - \tfrac{1}{2^{k-2}} + 
4 \Phi_{0, \Sigma_{2k}(\rho)}(0) = 1 - \tfrac{1}{2^{k-2}} + 
4\P\big(\bigcap_{i=1}^{2k} \{X_i \leq 0\}\big)
\end{align}
for all $\rho \in [-1,1]$ and $\textbf{X} = (X_1, \ldots, X_{2k})^\top 
\sim N_{2k}(0, \Sigma_{2k}(\rho))$ (since $\Phi(0) = \tfrac{1}{2}$). In fact, if 
$(\alpha, \beta) = (\tfrac{1}{2},\tfrac{1}{2})$ and $k=1$, then
$\psi_0 = \text{sign} - \ind_{\{0\}}$, implying that
\eqref{eq:Gauss_Copula_special_case} even reduces to the Grothendieck equality.
In order to recognise this, recall first that $(2l)! = 2^l\,l!\,(2l-1)!!$
for any $l \in \N_0$, whence
\begin{align}\label{eq:arcsin_part}
((2l-1)!!)^2 = \binom{2l}{l}\frac{(2l)!}{4^l}\,.
\end{align}
Consequently, the power series representation of the function 
$h_{\psi_0,\psi_0}$, together with 
\sref{Corollary}{cor:generalisation_of_Groth_equality} and 
\eqref{eq:Hermite_polynomials_at_zero} implies that 
\begin{align}\label{eq:Gauss_Copula_special_case_arcsin}
\begin{split}
4\P(X_1 \leq 0, X_{2}\leq 0) - 1 &= 
h_{\text{sign},\text{sign}}(\rho) = h_{\psi_0,\psi_0}(\rho)\\ 
&= \sum_{n=0}^\infty (\langle \psi_0, H_n \rangle_{\gamma_1})^2\,\rho^n = 
\frac{2}{\pi}\sum_{l=0}^\infty \frac{1}{2l+1} H_{2l}^2(0)\,\rho^{2l+1}\\
&\stackrel{\eqref{eq:Hermite_polynomials_at_zero}}{=}
\frac{2}{\pi}\sum_{l=0}^\infty ((2l-1)!!)^2 \,\frac{\rho^{2l+1}}{(2l+1)!} 
\stackrel{\eqref{eq:arcsin_part}}{=} \frac{2}{\pi}\arcsin(\rho) 
\end{split}
\end{align}
for all $\rho \in [-1,1]$ and $(X_1, X_{2})^\top \sim N_{2}(0, I_{2})$ 
(cf. also \cite{HS1980, M2013}). Again, we recognise that the Hermite expansion 
of the function $\text{sign}$ - in $L^2(\gamma_1)$ - is given by
\begin{align}\label{eq:Hermite_rep_of_sign_in_L2_gamma_1} 
\text{sign} = \psi_0 = \sqrt{\frac{2}{\pi}}\sum_{l=0}^\infty
(-1)^{l}\,\frac{(2l-1)!!}{\sqrt{(2l+1)!}}\,H_{2l+1}
\end{align}
(cf. \eqref{eq:Hermite_Fourier_coefficient_of_sign}).
\end{example}
\noindent A further interesting one-dimensional example ($k=1$) is given by the 
function $h_{\Phi, \Phi}$, where $\Phi = F_X$ is the (continuous) distribution 
function of a standard normally distributed random variable $X \sim N_1(0,1)$. 
To this end, recall that the (continuous) random variable $U : = \Phi(X) \sim 
U(0,1)$ is uniformly distributed on $[0,1]$, implying that $\Phi \in 
L^2(\gamma_1)$, with $\langle \Phi, 1 \rangle_{\gamma_1} = \E[U] = \frac{1}{2}$ 
and $\Vert \Phi \Vert_{\gamma_1}^2 = \E[U^2] = Var(U) + \E^2[U] = \frac{1}{12} 
+ \frac{1}{4} = \frac{1}{3}$ (cf. \cite[Remark 2.17.]{Oe2015}). In particular, 
$\E[\kappa(X)] = 0$, where $\kappa : = 2\sqrt{3}\Phi - \sqrt{3} = 
\sqrt{3}(2 \Phi - 1)$. Now, we are ready to prove
\begin{proposition}\label{prop:Phi_minus_one_half}
Let $X \sim N_1(0,1)$ and $\kappa : = \sqrt{3}\big(2\Phi - 1\big)$. Then the following 
properties hold:
\begin{enumerate}
\item 
\[
\frac{\textup{d}^n}{\textup{d}t^n}\E[\Phi(X+t)] = 2^{-n/2}\,\varphi^{(n-1)}
\big(\frac{t}{\sqrt{2}}\big) \text{ for all } n \in \N \text{ and } t \in \R\,.
\]
\item
\[
\E[\Phi(X+t)] = \Phi\big(\frac{t}{\sqrt{2}}\big) \text{ for all } n \in \N 
\text{ and } t \in \R\,.
\]
\item $\E[\Phi(X)\,H_{2n+1}(X)] = 
(-1)^n\,\sqrt{\frac{1}{2\pi}}\,\sqrt{\frac{1}{4^n}\frac{1}{2n+1}{\binom{2n}{n}}}\,
\sqrt{(\frac{1}{2})^{2n+1}}$ for all $n \in \N_0$. In particular, 
$\E[\kappa(X)\,H_{2n+1}(X)] = (-1)^n\,\sqrt{\frac{6}{\pi}}\,\sqrt{\frac{1}{4^n}
\frac{1}{2n+1}{\binom{2n}{n}}}\,\sqrt{(\frac{1}{2})^{2n+1}}$ for all $n 
\in \N_0$.
\item $h_{\Phi,\Phi}(\rho) = \frac{1}{4} + \frac{1}{2\pi}\,
\arcsin(\frac{\rho}{2}) \text{ and } h_{\kappa,\kappa}(\rho) = \frac{6}{\pi}
\,\arcsin(\frac{\rho}{2}) \text{ for all } \rho \in [-1,1]$. 
\item $h_{\sqrt{3}\,\Phi,\sqrt{3}\,\Phi} = 3\,h_{\Phi,\Phi} = 
\frac{3}{4}(1 + \frac{2}{\pi}\,\arcsin(\frac{1}{2}\cdot))$ is a strictly 
increasing homeomorphic CCP function which maps $[-1,1]$ onto 
$[\frac{1}{2}, 1]$ and is neither odd nor even. $h_{\kappa,\kappa}$ is an odd, 
strictly increasing homeomorphic CCP function. $h_{\sqrt{3}\,\Phi,\sqrt{3}\,
\Phi}^{-1}(t) = -2\cos(\frac{2\pi}{3} t)$ for all $t \in [\frac{1}{2}, 1]$, 
$h_{\kappa,\kappa}^{-1}(s) = 2\sin(\frac{\pi}{6} s)$ and 
$\big(h_{\kappa,\kappa}^{-1}\big)_{\text{abs}}(s) = 2\sinh(\frac{\pi}{6} s)$ 
for all $s \in [-1,1]$. In particular, $\big(h_{\kappa,\kappa}^{-1}\big)_
{\text{abs}}(1) > 1$. $h_{\kappa,\kappa}^{\text{hyp}}(\rho) = \frac{6}{\pi}
\sinh^{-1}(\frac{\rho}{2})$ for all $\rho \in [-1,1]$. 
\item $\kappa$ is odd, $\Vert \kappa \Vert_{\gamma_1} = 1$ and 
$h_{\kappa,\kappa}^{\prime}(0) = \frac{3}{\pi} > 0$. Moreover, $\kappa \in 
L^\infty(\gamma_1)$, and $\Vert\kappa\Vert_{\infty} = \sqrt{3}$.  
\end{enumerate}
\end{proposition}
\noindent We will soon realise that \autoref{thm:h_fg_real_case} actually 
reflects a characterisation of the class of all real continuous CCP functions (see Theorem 
\ref{thm:real_CCP_rep_thm}).

Within the scope of our analysis of the Grothendieck constants, we need 
\textit{invertible} CCP functions, implying that we may ignore even CCP functions, 
defined on $[-1,1]$ (since these are non-injective). However, thanks to 
\autoref{thm:CCP_via_Schoenberg}, $[-1,1] \ni \rho \mapsto \rho\,h(\rho)$ is an odd CCP 
function for any continuous, even CCP function $h$. In fact, we will recognise next that in particular 
any CCP function $h_{f,f}$, which satisfies $h_{f,f}(-1) < 1$ induces canonically 
an odd CCP function $h_{f^\times,f^\times}$, where $f^\times \in S_H$ is odd, $H : = 
L^2(\R^k, \gamma_k)$. In the complex case it seems that we even have to use odd 
functions, to avoid calculations with non-positive semidefinite entrywise absolute values of 
correlation matrices (cf. \eqref{eq:Haagerup_style} and 
\autoref{thm:odd_complex_bdd_CCP_general_version}-\ref{(ii)}). 
\begin{proposition}[\textbf{Odd CCP transform}]\label{prop:odd_CCP_transform}
Let $k \in \N$ and $f \in L^2(\R^k, \gamma_k)$. Consider the odd function 
$\R^k \ni x \mapsto f^{\text{odd}}(x) : = \frac{f(x)-f(-x)}{2}$ (the odd part of $f$). 
Then $f$ is odd if and only if $f^{\text{odd}} = f$. In particular, $(f^\text{odd})^
\text{odd} = f^{\text{odd}}$. $f$ is even if and only if $f^{\text{odd}} = 0$. Moreover, 
$f^{\text{odd}} \in L^2(\R^k, \gamma_k)$ and
\begin{enumerate}
\item
\begin{align*}
\Vert f^{\text{odd}} \Vert_{\gamma_k}^2 = \frac{h_{f,f}(1)- h_{f,f}(-1)}{2} \geq 0\,.
\end{align*}
\item $f$ is even $\lambda_k$-a.s. if and only if $h_{f,f}(-1) = 
h_{f,f}(1)$.
\item $f$ is odd $\lambda_k$-a.s. if and only if $h_{f,f}(-1) = 
-h_{f,f}(1)$. 
\item Assume that $h_{f,f}(1) > h_{f,f}(-1)$. Put
\begin{align}\label{eq:odd_CCP_transform}
f^\times : = \sqrt{\frac{2}{h_{f,f}(1)- h_{f,f}(-1)}}\,f^{\text{odd}}\,.
\end{align} 
Then $f$ is not even, $\Vert f^\times \Vert_{\gamma_k} = 1$, and $h_{f^\times, 
f^\times} = \frac{2}{h_{f,f}(1)- h_{f,f}(-1)}\,h_{f^{\text{odd}},
f^{\text{odd}}}$ is an odd CCP function.
\end{enumerate}
\end{proposition}
\begin{comment}
Due to \autoref{thm:h_fg_real_case}, respectively Theorem 
\ref{thm:h_ff_real_odd_case}, we just have to verify the non-trivial parts of the 
claims (i), (ii) and (iii).
\\[0.2em] 
\noindent (i) Fix $\textbf{X} \sim N_k(0, I_k)$. Then
\[
\Vert f^{\text{odd}} \Vert_{\gamma_k}^2 = \frac{1}{4}\,\E[(f^\circ)^2(\textbf{X})] =
\frac{1}{4}\big(\E[f^2(\textbf{X})] - 2\,\E[f(\textbf{X})f(\textbf{-X})] +  
\E[f^2(\textbf{-X})]\big).
\]
Since $\textbf{-X} \stackrel{d}{=} \textbf{X} \sim N_k(0, I_k)$ and 
$\vc{\textbf{X}}{\textbf{-X}} \sim N_{2k}(0, \Sigma_{2k}(-1))$, it follows from 
\autoref{thm:h_fg_real_case} that
\[
\Vert f^{\text{odd}} \Vert_{\gamma_k}^2 = \frac{1}{2}\big(\Vert f \Vert_{\gamma_k}^2 - 
\E[f(\textbf{X})f(\textbf{-X})]\big) = \frac{1}{2}(h_{f,f}(1)- h_{f,f}(-1)).
\]
(ii) follows directly from (i).
\\[0.2em] 
\noindent (iii) Consider the even function $f^{\text{ev}} : = 
f - f^{\text{odd}} \in L^2(\R^k, \gamma_k)$; i.e., the even part of $f$. By 
construction, $f^{\text{ev}}(x) = \frac{f(x) + f(-x)}{2}$ for all $x \in \R^k$. 
Similarly, as in the proof of (i), we therefore obtain
\[
\Vert f^{\text{ev}} \Vert_{\gamma_k}^2 = \frac{h_{f,f}(1) + h_{f,f}(-1)}{2}\,,
\]
wherefrom claim (iii) obviously follows.
\end{comment}
\begin{example}
Since $\arcsin(\frac{1}{2}) = \frac{\pi}{6}$, 
\sref{Proposition}{prop:Phi_minus_one_half}-\ref{(iv)} implies that 
$h_{\Phi, \Phi}(1) = \frac{1}{3}$ and $h_{\Phi, \Phi}(-1) = \frac{1}{6}$, whence
\[
\Phi^\times = \sqrt{3}(2\Phi-1) = \kappa\,.
\]
\end{example}
\begin{example}
Let $[-1,1] \ni \rho \mapsto C^{\text{Ga}}(\frac{1}{2}, \frac{1}{2}; \rho)$ 
denote the bivariate Gaussian copula with Pearson's correlation coefficient $\rho$ as 
parameter, evaluated at $(\frac{1}{2}, \frac{1}{2})$. Then
\[
C^{\text{Ga}}(\frac{1}{2}, \frac{1}{2}; \rho) = 
\Phi_{\Sigma_2(\rho)}(\Phi^{-1}(\frac{1}{2}), 
\Phi^{-1}(\frac{1}{2})) = \P(X \leq 0, Y \leq 0) = h_{g,g}(\rho)\,,
\]
where $\Phi_{\Sigma_2(\rho)}$ denotes the bivariate distribution function of the
random vector $(X,Y)^\top \sim N_{2}(0, \Sigma_{2}(\rho))$ and 
$g : = \ind_{(-\infty, 0]}$. Thus, $g^{\text{odd}} = -\frac{1}{2}\,\text{sign}$, 
$h_{g,g}(1) = \P(X \leq 0) = \Phi(0) = \frac{1}{2}$ and $h_{g,g}(-1) = 
\P(X \leq 0, -X \leq 0) = \P(X=0) = 0$. 
Consequently,
\[
g^\times = -\text{sign} = 2\,\ind_{(-\infty, 0)} - 1 + \ind_{\{0\}}\,,  
\]
and it follows that $h_{g^\times,g^\times} = \frac{2}{\pi}\arcsin$ (on $[-1,1]$).
\end{example}
\noindent Regarding an explicit calculation of the correlation coefficients
$\E[f(\textbf{X})g(\textbf{Y})]$, induced by a pair of non-linearly transformed 
Gaussian random vectors and given ``sufficiently smooth'' functions $f$ and $g$, 
we implement a few facts from the theory of distributions and test function 
spaces. A detailed in-depth introduction to these ``generalised functions'' and 
test function spaces and their analysis is provided by e.\,g. 
\cite{HT2008, J2001, R1991}. To this end, let $k \in \N$, $f \in 
L^1_{{\text{loc}}}(\R^k)$ (i.e., $f$ is locally integrable) and $\psi \in 
{\mathcal{S}}_k$ be an arbitrary test function, where ${\mathcal{S}}_k$ denotes 
the Schwartz space of rapidly decreasing test functions on $\R^k$. Put
\[
\langle \psi, \Lambda_f\rangle : = \int_{\R^k} f(x)\psi(x)\textup{d}^k{x}\,.
\]     
Observe that the latter symbol $\langle \cdot, \cdot \rangle$ denotes the 
duality bracket on ${\mathcal{S}}_k \times {\mathcal{S}}_k^\prime$ and not 
an inner product of Hilbert space elements. When there is no ambiguity, we 
adopt the common habit to identify the tempered distribution $\Lambda_f \in 
{\mathcal{S}}_k^\prime$ with $f \in L^1_{{\text{loc}}}(\R^k)$ itself. More 
generally, recall that tempered distributions and their derivatives are 
elements of the dual space ${\mathcal{S}}_k^\prime$ of ${\mathcal{S}}_k$, such 
as the Dirac delta distribution $\delta_0 = \delta$, defined via 
$\langle \psi, \delta \rangle : = \psi(0)$ for any $\psi \in {\mathcal{S}}_k$. 
Observe that there is no $g \in L^1_{{\text{loc}}}(\R^k)$ such that $\delta = 
\Lambda_g$. 

 Let us quickly recall how differentiation of tempered distributions is 
defined. It originates from a reiteration of the integration by parts formula 
and is also known as ``weak differentiation'' (cf. e.\,g. \cite[Chapter 6.12]{R1991}):
\begin{align}\label{eq:weak_derivative}
\langle \psi, D^n u \rangle : = (-1)^{\vert n \vert} 
\langle D^n \psi, u \rangle = (-1)^{\vert n \vert} \langle \psi, 
D^{n-1}(Du) \rangle
\end{align}
for all $\psi \in {\mathcal{S}}_k, u \in {\mathcal{S}}_k^\prime$ and $n \in 
\N_0^{k}$. Consequently, since $H_n\,\varphi_k \in {\mathcal{S}}_k$ for all 
$n \in \N_0^k$ and $D^n u \in {\mathcal{S}}_k^\prime$ for all $u \in 
{\mathcal{S}}_k$, \eqref{eq:Hermite_product} in particular implies:
\[
\langle H_n\,\varphi_k, u \rangle = \frac{1}{\sqrt{n!}}
\langle \varphi_k, D^n u \rangle \text{ for all } u \in {\mathcal{S}}_k^\prime 
\text{ and } n \in \N_0^k\,.
\]
Given an arbitrary compact subset 
$K \subseteq \R^k$ and $f \in L^2(\gamma_k)$, it follows that 
\[
\int_{K} \vert f \vert\,\textup{d}\lambda_k = 
\sqrt{2\pi}\int_{\R^k}\ind_{K}(x)\exp(\frac{1}{2}\Vert x \Vert^2)\,
\vert f(x) \vert \,\gamma_k(\textup{d}x) \leq c_K\,\Vert f \Vert_{\gamma_k}\,,
\]
where $c_K : = \sqrt{2\pi}\int_K e^{\Vert x \Vert^2}\,\gamma_k(\textup{d}x) = 
\sqrt{2\pi}\int_K \exp(\frac{1}{2}\Vert x \Vert^2)\lambda_k(\textup{d}x) 
< \infty$. Consequently, $L^2(\gamma_k) \subseteq L^1_{{\text{loc}}}(\R^k)$, 
implying that for any $\textbf{X} \sim N_k(0, I_k)$
\begin{align}\label{eq:Hermite_f_bracket}
\E[H_n(\textbf{X})f(\textbf{X})] = \langle H_n, f\rangle_{\gamma_k} = 
\langle H_n\,\varphi_k, \Lambda_f \rangle = 
\frac{1}{\sqrt{n!}}\langle \varphi_k, D^{n}\Lambda_f \rangle \text{ for all } f \in 
L^2(\gamma_k) \text{ and } n \in \N_0^k\,.
\end{align}
\begin{example}[\textbf{$n$'th weak derivative of $\Lambda_{\indsm_{[0, \infty)}}$}]
\label{ex:nth_deriv_of_Heaviside_step_fct}
Let $k =1$ and $n \in \N_0$. A very easy proof by induction on $n$, based on
\eqref{eq:weak_derivative}, together with \eqref{eq:one_dim_Hermite_pol} 
firstly reveals that  
\[
\langle\varphi, D^n\delta\rangle = (-1)^n\,\varphi^{(n)}(0) = \frac{\sqrt{n!}}
{\sqrt{2\pi}}\,H_n(0)\,.
\]
\eqref{eq:Hermite_polynomials_at_zero} therefore implies that
\[
\langle\varphi, D^{2l}\delta\rangle = (-1)^l\,\frac{1}{\sqrt{2\pi}}\,(2l-1)!!\,
\text{ for all } l \in \N_0\,.
\]
Similarly, \eqref{eq:weak_derivative} implies that $D\Lambda_{\indsm_{[0, \infty)}} 
= \delta$. Since $\text{sign} = 2\,\ind_{[0, \infty)}-1$ ($\gamma_1$-almost surely) and 
$\langle H_{2l+1}, 1\rangle_{\gamma_1} = \langle H_{2l+1}, H_0\rangle_{\gamma_1} =
\delta_{2l+1, 0} = 0$, it follows that
\begin{align*}
\langle H_{2l+1}, \text{sign}\rangle_{\gamma_1} &= 2\langle H_{2l+1}, 
\ind_{[0, \infty)}\rangle_{\gamma_1} \stackrel{\eqref{eq:Hermite_f_bracket}}{=} 
\frac{2}{\sqrt{(2l+1)!}}\langle\varphi, D^{2l+1}\Lambda_{\indsm_{[0, \infty)}}\rangle\\ 
&= \frac{2}{\sqrt{(2l+1)!}}\langle\varphi, D^{2l}\delta\rangle = (-1)^l\,
\sqrt{\frac{2}{\pi}}\,\frac{(2l-1)!!}{\sqrt{(2l+1)!}}
\end{align*}
for all $l \in \N_0$. This outcome -- which is solely based on a multiple 
weak differentiation of $\Lambda_{\indsm_{[0, \infty)}}$ -- should be compared now 
with (the derivation of) \eqref{eq:Gauss_Copula_special_case_arcsin}\,! 
\end{example}
\noindent \eqref{eq:Hermite_f_bracket} can be strongly simplified if $f$ is 
smooth (cf. also \sref{Proposition}{prop:k_dim_Fourier_Hermite_coeff}). More 
precisely, if $N \in \N_0, D^n f \in L^2(\gamma_k) \cap C(\R^k)$ for all 
$n \in \N_0^k$, satisfying $0 \leq \vert n \vert \leq N$, and if $\textbf{X} 
\sim N_k(0, I_k)$, then $D^{n}\Lambda_f = \Lambda_{D^{n} f}$ (cf. 
\cite[Chapter 6.13]{{R1991}}), and 
\begin{align}\label{eq:of_Stein_type}
\E[H_n(\textbf{X})f(\textbf{X})] = \langle H_n, f\rangle_{\gamma_k} = 
\frac{1}{\sqrt{n!}}\E[D^n f(\textbf{X})] \text{ for all } n \in \N_0^k\,, 
\text{ with } \vert n \vert \leq N\,. 
\end{align}
In fact, even more can be said. Recall that (by construction of Sobolev spaces) the latter 
equality also holds without the smoothness assumption, if we just assume that $f \in 
L^2(\gamma_k) \cap W_{{\text{loc}}}^{N, 1}(\R^k)$, so that in this case $D^n f$ denotes the 
distributional $n$'th derivative in $L_{{\text{loc}}}^1(\R^k)$ instead (cf. e.\,g.
\cite[Chapter 3.1]{HT2008}). Independent of an application of Theorem 
\ref{thm:h_fg_real_case} to the primary topic of our work, it implies further 
non-trivial consequences including a generalisation of Stein's Lemma (cf. 
\cite[Theorem 1 and Example 1]{LN2008}) and a certain ``decorrelation property'' 
of harmonic functions. We only have to combine \eqref{eq:of_Stein_type} and 
\autoref{thm:h_fg_real_case}, resulting at once in
\begin{proposition}\label{prop:Sobolev_CCP_real_correl_case}
Let $k\in \N$, $\rho \in [-1,1]$ and $\vc{\textbf{X}}{\textbf{Y}} \sim 
N_{2k}(0, \Sigma_{2k}(\rho))$. Let $f, g \in L^2(\gamma_k)$. Then
\[
{\text{cov}}(f(\textbf{X}), g(\textbf{Y})) = 
\sum_{\nu = 1}^\infty \big(\sum_{n \in C(\nu, k)}\frac{1}{n!}\,
\langle \varphi_k, D^{n}\Lambda_f \rangle\,\langle \varphi_k, D^{n}
\Lambda_g \rangle \big)\,\rho^{\nu}\,.
\]
If in addition $(f,g) \in W_{{\text{loc}}}^{N, 1}(\R^k) \times 
W_{{\text{loc}}}^{N, 1}(\R^k)$, or if $(D^n f, D^n g) \in (L^2(\gamma_k) \cap C(\R^k))
\times (L^2(\gamma_k) \cap C(\R^k))$ for all $n \in \N_0^k$, satisfying 
$0 \leq \vert n \vert \leq N$ for some $N \in \N$, then

\scalebox{0.88}{
\vbox{
\[
{\text{cov}}(f(\textbf{X}), g(\textbf{Y})) = 
\sum_{\nu=1}^N\big(\sum_{n \in C(\nu, k)}\frac{1}{n!}\,\E[D^n f(\textbf{X})]
\E[D^n g(\textbf{Y})]\big)\rho^{\nu} + 
\sum_{\nu = N+1}^\infty \big(\sum_{n \in C(\nu, k)}\frac{1}{n!}\,
\langle \varphi_k, D^{n}\Lambda_f \rangle\,\langle \varphi_k, D^{n}\Lambda_g \rangle \big)
\,\rho^{\nu}\,.
\]
}}

\noindent In particular, if $N=1$, then
\[
{\text{cov}}(f(\textbf{X}), Y_i) = 
\rho\,\E[\frac{\partial f}{\partial x_i}(\textbf{X})]
\]
for all $i \in [k]$, respectively
\[
\E[f(\textbf{X})\textbf{Y}] = \rho\,\E[\nabla f(\textbf{X})].
\]
If $N = 2$, then
\[
{\text{cov}}(f(\textbf{X}), \Vert\textbf{Y}\Vert^2) = 
\rho^2\,\E[\triangle f(\textbf{X})].
\] 
\end{proposition}
\noindent In the one-dimensional case (i.\,e., $k=1$), a direct application of 
\sref{Proposition}{prop:Sobolev_CCP_real_correl_case}, respectively 
\eqref{eq:of_Stein_type}, implies two further results which are of their own 
interest. In particular, our approach enables the provision of a quick and short 
proof of the following generalisation of Stein's Lemma to standard Gaussian 
random powers (cf. \cite[Theorem 1]{KM2022}):
\begin{corollary}
\label{cor:power_version_of_Stein_s_Lemma}
Let $m \in \N_0$, $f \in L^2(\gamma_1)$ and $X \sim N_1(0,1)$ be given. If 
$D^l f \equiv f^{(l)} \in L^2(\gamma_1) \cap C(\R)$ for all 
$0 \leq l \leq m$, then
\[
\E[f(X)X^m]  = \E[(-i)^m\,H_m(i D)(f(X))] = \sum_{\nu=0}^{\floor{m/2}}
\binom{m}{2\nu}(2\nu-1)!!\,\E[D^{m-2\nu} f(X)]\,.
\]
\end{corollary}
\begin{comment}
Fix $m \in \N_0$. Put $f_m(x) : = x^m$ ($x \in \R$). Firstly, recall the 
well-known result that
\[
\E[f_\nu(X)] = \E[X^\nu] = \frac{1+(-1)^{\nu}}{2}\,(\nu-1)!! \text{ for all } 
\nu \in \N_0\,,
\]
implying that 
\[
\langle f_m, H_l \rangle_{\gamma_1} \stackrel{\eqref{eq:of_Stein_type}}{=} 
\frac{1}{\sqrt{l!}}\,\E[D^l f_m(X)] =
\begin{cases}\frac{1+(-1)^{m-l}}{2}\,\sqrt{l!}\,\binom{m}{l}\,(m-l-1)!! 
&\text{if } l \leq m\\ 
0 &\text{if } l > m\end{cases}
\]
Consequently,
\[
\E[f(X)f_m(X)] = h_{f,f_m}(1) = \sum_{l=0}^m \langle f, H_l \rangle_{\gamma_1}\,
\langle f_m, H_l \rangle_{\gamma_1} = \sum_{l=0}^m \frac{1+(-1)^{m-l}}{2}
\E[D^l f(X)]\,\binom{m}{l}(m-l-1)!!\,.
\]
If $m$ is even, it therefore follows that
\[
\E[f(X)f_m(X)]  = \sum_{k=0}^{m/2} 
\E[D^{2k} f(X)]\,\binom{m}{2k}(m-2k-1)!! = \sum_{\nu=0}^{m/2} 
\E[D^{m-2\nu} f(X)]\,\binom{m}{2\nu}(2\nu-1)!!\,.
\]
Similarly, if $m$ is odd, we obtain
\[
\E[f(X)f_m(X)]  = \sum_{\nu=0}^{(m-1)/2}\E[D^{m-2\nu} f(X)]\,\binom{m}{2\nu}
(2\nu-1)!!\,.
\]
However, since $\binom{m}{2\nu}(2\nu-1)!! = H_{m,\nu}$ (due to 
\eqref{eq:one_dim_Hermite_pol}), the claim clearly follows.
\end{comment}
\begin{corollary}
Let $\nu \in \N_0$, $f, g \in L^2(\gamma_1)$ and $\rho \in (-1, 1)$.
Then
\[
h_{f,g}^{(\nu)}(\rho) = 
\sum_{n=0}^\infty \frac{1}{n!}\,\langle \varphi, D^{n}(D^{\nu}\Lambda_f) \rangle
\langle \varphi, D^{n}(D^{\nu}\Lambda_g) \rangle\rho^n\,.
\]
In particular, if in addition $(f, g) \in  C^\infty(\R^k) \times C^\infty(\R^k)$, then
\[
h_{f,g}^{(\nu)}(\rho) = h_{f^{(\nu)}, g^{(\nu)}}(\rho) \text{ for all } \rho \in (-1,1).
\]
\end{corollary}
\begin{comment}
Since $h_{f,g}$ is real analytic, $h_{f,g}^{(\nu)}$ is well-defined (where 
$h_{f,g}^{(0)} : = h_{f,g}$, of course). Due to \eqref{eq:Hermite_f_bracket} and Theorem 
\ref{thm:h_fg_real_case}, applied to $k=1$, the remaining part of the proof is a 
straightforward proof by induction on $\nu \in \N_0$. 
\end{comment}
\begin{remark}
Firstly, observe that the strong impact of the standard (multivariate) \textit{Gaussian} 
law, is reflected in \eqref{eq:Hermite_f_bracket}, primarily implied by the fact that the 
Gaussian density function $\varphi_k$ and hence each $H_n\,\varphi_k$, is rapidly decreasing: 
$H_n\,\varphi_k \in {\mathcal{S}}_k$ for all $n \in \N_0^k$ (since each $H_n$ is a 
polynomial). This fact and the structure of the (multivariate) Hermite polynomials (which 
itself is also induced by the structure of $\varphi_k$) namely enables a reiterated use of 
the integration by parts formula in the smooth case, respectively a use of the Sobolev space 
$W_{{\text{loc}}}^{N, 1}(\R^k)$ in the non-smooth case, implying the transition of each 
inner product $\langle H_n, f\rangle_{\gamma_k}$ on the Hilbert space $L^2(\gamma_k)$ into 
the duality bracket $\langle \varphi_k, D^n\,\Lambda_f \rangle$ on 
${\mathcal{S}}_k \times {\mathcal{S}}_k^\prime$. The latter, however, seems to be more 
convenient for performing specific computations.
\end{remark}

\section{Upper bounds of $K_G^\R$ and inversion of real CCP functions}
 Our next step is to embed Grothendieck's original approach as well as 
Krivine's improvement into a general framework. We are going to show that we 
may \textit{substitute} the Grothendieck function $h_{\text{sign},\text{sign}} = 
\frac{2}{\pi}\arcsin$ through \textit{invertible} CCP functions $h_{f,f}$, 
generated by \textit{bounded} functions $f : \R^k \longrightarrow \R$ (see 
\autoref{thm:GT_generalising_Krivine} and \autoref{thm:odd_real_bdd_CCP}). 
Various proofs of Krivine's main result ($K_G^\R \leq \frac{\pi}
{2 \ln(1 + \sqrt{2})} \approx 1.782$) are set out in detail in 
\cite[Section 5]{H1994}, \cite[proof of Lemma 10.5]{J1987} and \cite{Z2018}. We 
will recognise soon that inverses of odd CCP functions also lead to a further 
crucial construction of correlation matrices, lurking in the following important 
implication of \sref{Lemma}{lem:real_analytic_and_l1}:
%
\begin{theorem}\label{thm:char_of_completely_real_analytic_functions_at_0_as_h_fg}
Let $r > 0$, $0 < c \leq r$ and $0 \not= \psi \in W^\omega_+((-r, r))$. Then
$\psi_{\text{abs}}\big\vert_{[0, r]}$ is strictly increasing and 
$\psi_{\text{abs}}(c) > 0$. For any $k \in \N$ there exist 
$\alpha_k \equiv \alpha^{\psi, c}_k \in S_{L^2(\gamma_k)}$ and $\beta_k \equiv 
\beta^{\psi, c}_k \in S_{L^2(\gamma_k)}$, such that the following properties 
hold:
\begin{enumerate}
\item If $\textbf{X} \sim N_k(0, I_k)$, then 
$\sqrt{\psi_{\text{abs}}(c)}\,\E[\alpha_k(\textbf{X})] = 
\text{sign}(\psi(0))\sqrt{\vert \psi(0) \vert}$ and 
$\sqrt{\psi_{\text{abs}}(c)}\,\E[\beta_k(\textbf{X})] = 
\sqrt{\vert \psi(0)\vert}$. In particular, $\E[\alpha_k(\textbf{X})] = 
\text{sign}(\psi(0))\E[\beta_k(\textbf{X})]$.
\item If $c < r$, then $[-c,c] \subseteq (-r,r)$ and
\begin{align}\label{eq:rep_of_psi_vs_psi_abs_general_case}
\psi(c\rho) = \psi_{\text{abs}}(c)h_{\alpha_k,\beta_k}(\rho)
\text{ for all } \rho \in [-1,1].
\end{align}
In particular, $\psi(c) = \psi_{\text{abs}}(c)\langle \alpha_k, \beta_k\rangle_
{\gamma_k}$.
\item
If $c \leq r$, then
\[
\psi_{\text{abs}}(c\rho) = \psi_{\text{abs}}(c)\,h_{\alpha_k, \alpha_k}(\rho) =
\psi_{\text{abs}}(c)\,h_{\beta_k, \beta_k}(\rho) \text{ for all } \rho \in [-1,1].
\]
\item If $c < r$, $H$ is an arbitrary $\R$-Hilbert space, then there 
exists a $\R$-Hilbert space $\H$ such that for any $u, v \in S_H$
\begin{align}\label{eq:rep_of_psi_vs_psi_abs_inner_product_general_case}
\psi(c\langle u, v\rangle_H) = \psi_{\text{abs}}(c)\langle \psi_u(\alpha_k), 
\psi_v(\beta_k)\rangle_{\H} = \langle a_u, b_v\rangle_{\H}\,,
\end{align}
where for any $w \in S_H$, $\Vert a_w\Vert_{\H}^2 = \Vert b_w\Vert_{\H}^2 
= \psi_{\text{abs}}(c)$ and $\psi_w : L^2(\gamma_k) \longrightarrow \H$ is 
a mapping which satisfies $\psi_w(S_{L^2(\gamma_k)}) \subseteq S_{\H}$.
In particular, $\psi(c) = \langle a_w, b_w\rangle_{\H}$ for all $w \in S_H$.
\\[0.2em]
\noindent If $c \leq r$, then
\begin{align*}
\psi_{\text{abs}}(c\langle u, v\rangle_H) &= 
\psi_{\text{abs}}(c)\langle\psi_u(\alpha_k), \psi_v(\alpha_k)\rangle_
{\H} = \psi_{\text{abs}}(c)\langle \psi_u(\beta_k), \psi_v(\beta_k)\rangle_
{\H}\\
&= \langle a_u, a_v\rangle_{\H} = \langle b_u, b_v\rangle_{\H}\,,
\end{align*}
\end{enumerate}
\end{theorem}
\noindent If we combine \autoref{thm:CCP_via_Schoenberg}, Theorem 
\ref{thm:h_fg_real_case}, \sref{Proposition}{prop:odd_CCP_transform} and Theorem 
\ref{thm:char_of_completely_real_analytic_functions_at_0_as_h_fg}, (iii), a 
further significant characterisation of continuous odd CCP functions follows at 
once:
\begin{theorem}[\textbf{CCP representation theorem}]\label{thm:real_CCP_rep_thm}
Let $\psi : [-1,1] \longrightarrow \R$ be a continuous odd function and $k \in \N$. 
Then the following statements are equivalent:
\begin{enumerate}
\item $\psi$ is a CCP function.
\item $\psi = h_{f,f}$ for some odd $f \in S_{L^2(\gamma_k)}$. 
\end{enumerate}
\end{theorem}
\noindent 
\autoref{thm:char_of_completely_real_analytic_functions_at_0_as_h_fg}-\ref{(iv)}, 
together with \sref{Lemma}{lem:hyperbolic_CCP_transform} directly enriches us 
with another crucial result, related to a construction of real quantum 
correlation matrices which actually include lurking upper bounds of $K_G^\R$. 
{At this point we should especially recall the representation
\eqref{eq:sign_condition_and_hyp_CCP_fct_value} if the sign condition is 
present.}
\begin{corollary}\label{cor:h_ff_real_case_abs_value_crit}
Let $m, n, k \in \N$ and $f \in S_{L^2(\gamma_k)}$ be odd. Assume that 
$h_{f,f}^{-1}{\big\vert}_{(-1, 1)} \in W^\omega_+((-1, 1))$. Put 
\[
c(f) : = h_{f,f}^{\text{hyp}}(1)\,.
\] 
Then $c(f) \in (0,1)$ and
\[
h_{f,f}^{-1}[c(f)S] \in {\mathcal{Q}}_{m,n}\,
\text{ for all }\, S \in {\mathcal{Q}}_{m,n}\,.
\]
If 
$\Sigma =
\begin{pmatrix}
M & S \\
S^\top & N
\end{pmatrix} \in C(m+n; \R)$ is an arbitrary real $(m+n) \times (m+n)$ 
correlation matrix (with block elements $M \in C(m; \R), N \in C(n; \R)$ and 
$S \in {\mathcal{Q}}_{m,n}$), then
\begin{align}\label{eq:hyp_factor_in_quantum_correlation_matrix}
\begin{pmatrix}
(h_{f,f}^{-1})_{\text{abs}}[c(f)M] & h_{f,f}^{-1}[c(f)S]\\[0.7em]
h_{f,f}^{-1}[c(f)S^\top] & (h_{f,f}^{-1})_{\text{abs}}[c(f)N]
\end{pmatrix} \in C(m+n; \R)
\end{align}
again is an $(m+n) \times (m+n)$ correlation matrix with real entries.
\end{corollary}
\begin{comment} We just have to apply Theorem 
\ref{thm:char_of_completely_real_analytic_functions_at_0_as_h_fg}, (iv), 
including \eqref{eq:rep_of_psi_vs_psi_abs_inner_product_general_case}, to the 
function $\psi := h_{f,f}^{-1}{\big\vert}_{(-1, 1)}$ and the constant $c : = 
h_{f,f}^{\text{hyp}}(1)$ ($0 < c < 1$, due to 
\sref{Lemma}{lem:hyperbolic_CCP_transform}).
\end{comment}
\noindent If $\psi : D_\F \longrightarrow D_\F$ is an arbitrary CCP function 
and $\Sigma \in C(m+n; \F)$ is an arbitrary $m+n$-correlation matrix, then
\[
\psi[\Sigma] = 
\begin{pmatrix}
\psi[\Gamma_H(u,u)] & \psi[\Gamma_H(u,v)]\\[0.2em]
\psi[\Gamma_H(u,v)^\ast] & \psi[\Gamma_H(v,v)]
\end{pmatrix} \in C(m+n; \F)
\] 
again is a correlation matrix. Consequently, it follows that $\psi[\Gamma_H(u,u)] = 
\Gamma_K(x,x) \in C(m; \F)$, $\psi[\Gamma_H(v,v)] = 
\Gamma_K(y,y) \in C(n; \F)$ and $\psi[\Gamma_H(u,v)] = 
\Gamma_K(x,y) \in \mathcal{Q}_{m,n}(\F)$ for some Hilbert space $K$ over $\F$ 
and some $(x,y) \in S_K^m \times S_K^n$ (due to 
\sref{Corollary}{cor:block_matrix_structure_of_C_m_plus_n}). However, a 
considerably stronger result holds (cf. also 
\eqref{eq:hyp_factor_in_quantum_correlation_matrix}):
\begin{theorem}
\label{thm:real_hybrid_correlation_transform}
Let $m,n \in \N$, $A \in \M_m(\R)$, $S \in \M_{m,n}(\R)$ and $B \in \M_n(\R)$.
Consider the block matrix 
\[
M : = \begin{pmatrix}
A & S \\
S^\top & B
\end{pmatrix} \in \M_{m+n}(\R)\,.
\]
Let $0 < r < \infty$ and $f, g : (-r,r) \longrightarrow \R$ 
be two functions, such that $(-r,r) \ni x \mapsto f(x) = \sum_{\nu=0}^\infty 
a_\nu x^\nu \in W^\omega_+((-r,r))$ and $(-r,r) \ni x \mapsto g(x) = 
\sum_{\nu=0}^\infty b_\nu x^\nu \in W^\omega_+((-r,r))$. 
Let 
$q \in C^\omega((-r,r))$, such that
\begin{align}\label{eq:T0_C_omega_inequality}
\vert c_\nu \vert \leq \sqrt{\vert a_\nu\vert\,\vert b_\nu\vert} \text{ for all } 
\nu \in \N_0\,,
\end{align}
where $c_\nu : = \frac{q^{(\nu)}(0)}{\nu !}$.
Then 
$q \in W^\omega_+((-r,r))$, and the following properties 
hold:
\begin{enumerate}
\item If $M \in \M_{m+n}([-r,r])^+$ is positive semidefinite, and if any 
entry of the matrices $A, S$ and $B$ is an element of $[-r,r]$, then also  
\[
	\begin{pmatrix}
    f_{\text{abs}}[A] & \widetilde{q}[S]\\
    \widetilde{q}[S]^\top & g_{\text{abs}}[B]
	\end{pmatrix} \in \M_{m+n}([-r,r])^+
\]
is positive semidefinite, where $\widetilde{q}$ is defined as in 
\sref{Lemma}{lem:real_analytic_and_l1}. In particular, if $0 < c^\ast \leq r$ 
is a root of $f_{\text{abs}} - 1$, then
\[
\widetilde{q}[c^\ast\Gamma] \in \mathcal{Q}_{m,n} \text{ for all } 
\Gamma \in \mathcal{Q}_{m,n}\,.
\]
\item If $M \in C(m+n; \R)$ is a real correlation matrix, then also
\[
\begin{pmatrix}
h_{\alpha,\alpha}[A] & h_{\alpha,\beta}[S]\\
h_{\alpha,\beta}[S^\top] & h_{\beta,\beta}[B]
\end{pmatrix} \in C(m+n; \R)
\]
is a real correlation matrix for all $\alpha, \beta \in S_{L^2(\gamma_k)}$. 
In particular,
\begin{align}\label{eq:entrywise_mapping_of_inner_products_to_inner_products}
h_{\alpha,\beta}[\cdot] : {\mathcal{Q}}_{m,n} \longrightarrow {\mathcal{Q}}_{m,n} 
\text { for all } \alpha, \beta \in S_{L^2(\gamma_k)}.
\end{align}
\item For all $\alpha, \beta \in S_{L^2(\gamma_k)}$, for all Hilbert spaces $H$ and 
$u, v \in S_H$ there exist $d \in \N$ and $x, y \in \S^{d-1}$, such that
\[
h_{\alpha,\beta}(\langle u, v \rangle_H) = \langle x, y \rangle_{\R_2^d} = x^\top y\,.
\] 
\end{enumerate}
\end{theorem}
\noindent However, we will recognise that Theorem 
\ref{thm:real_hybrid_correlation_transform} cannot be fully transferred to the 
complex field, so that we have to distinguish carefully between the real case 
and the complex case here (cf. 
\sref{Proposition}{prop:complex_hybrid_correlation_transform}). If we link 
\eqref{eq:n_th_derivative_of_psi_abs_at_zero} and Theorem 
\ref{thm:real_hybrid_correlation_transform}, 
we obtain 
\begin{corollary}\label{cor:important_special_abs_case}
Let $m, n \in \N, A \in \M_{m,n}(\R)$ and
\[
\Sigma =
\begin{pmatrix}
M & S \\
S^\top & N
\end{pmatrix} \in C(m+n; \R)
\]
be an arbitrary real $(m+n) \times (m+n)$ correlation matrix (with block elements 
$M \in C(m; \R), N \in C(n; \R)$ and $S \in {\mathcal{Q}}_{m,n}$). Let 
$r > 0$ 
and $0 \not= \psi \in W^\omega_+((-r, r))$. If 
$0 < c \leq r$, then
\begin{align}\label{eq:abs_special_case}
\frac{1}{\psi_{\text{abs}}(c)}
\begin{pmatrix}
\psi_{\text{abs}}[c\,M] & \widetilde{\psi}[c\,S] \\[0.5em]
\widetilde{\psi}[c\,S]^\top & \psi_{\text{abs}}[c\,N]
\end{pmatrix} \in C(m+n; \R)
\end{align}
again is a correlation matrix with real entries. In particular,
\begin{align}\label{eq:entrywise_mapping_of_inner_products_to_inner_products_II}
\frac{1}{\psi_{\text{abs}}(c)}\widetilde{\psi}[c S] \in 
{\mathcal{Q}}_{m,n}\,\text{ for all }\,S \in {\mathcal{Q}}_{m,n}\,.
\end{align}
\end{corollary}
\noindent If we apply Theorem 
\ref{thm:char_of_completely_real_analytic_functions_at_0_as_h_fg}, (iv) and 
\eqref{eq:entrywise_mapping_of_inner_products_to_inner_products_II} to the 
\textit{inverses} of invertible functions in $W^\omega_+((-1,1))$, Bolzano's 
intermediate value theorem from calculus immediately implies a further crucial 
result. To this end, recall also \eqref{eq:psi_abs_estimation} and 
\eqref{eq:n_th_derivative_of_psi_abs_at_zero}.
\begin{theorem}[\textbf{Real inner product rounding}]
\label{thm:real_inner_product_rounding}
Let $\psi : [-1,1] \longrightarrow [-1,1]$ be a bijective real function. Assume 
that
\[
\psi{\vert}_{(-1,1)} \in W^\omega_+((-1,1)) \text{ and }
\psi^{-1}{\vert}_{(-1,1)} \in W^\omega_+((-1,1)).
\]
Then $\vert\psi^{-1}(0)\vert = (\psi^{-1}{\vert}_{(-1,1)})_
{\text{abs}}(0) \leq (\psi^{-1}{\vert}_{(-1,1)})_{\text{abs}}(1)$. 
Assume that
\[
\vert\psi^{-1}(0)\vert < 1 < (\psi^{-1}{\vert}_{(-1,1)})_
{\text{abs}}(1).
\]
Then there is a unique number $c^\ast \in (0,1)$, such that 
$(\psi^{-1})_{\text{abs}}(c^\ast) = 1$. Let $k \in \N$. There is $(\alpha, \beta) 
\in S_{L^2(\gamma_k)} \times S_{L^2(\gamma_k)}$ (dependent on $c^\ast$ and $k$) 
such that for any separable $\R$-Hilbert space $H$ and any $u, v \in S_H$, 
the following statements apply:
\begin{enumerate}
\item
\begin{align}\label{eq:real_inner_product_rep}
\langle u, v \rangle_H = \frac{1}{c^\ast} \psi(\rho_{u,v})
\end{align}
where $\rho_{u,v} : = h_{\alpha, \beta}(\langle u, v\rangle_H)$.
\item 
\[
c^\ast = \psi(\langle \alpha, \beta\rangle_{\gamma_1}).
\]
\item
Suppose that $\psi = h_{f,g}$ for some $\nu \in \N$ and $f, g \in L^2(\gamma_\nu)$.
If $\rho_{u,v} \in (-1,1)$, then
\[
c^\ast\,\langle u, v \rangle_H = 
\frac{1}{(2\pi)^{\nu}(1 - \rho_{u,v}^2)^{\nu/2}}\,
\int_{\R^{\nu}}\int_{\R^{\nu}}f(x)g(y)\exp\big(-\frac{\Vert x \Vert^2 + 
\Vert y \Vert^2 - 2\rho_{u,v} \langle x,y \rangle_2}{2(1-\rho_{u,v}^2)}\big) 
\textup{d}^{\nu}x\,\textup{d}^{\nu}y\,.
\]
If $m, n \in \N$ and $(x, y) \in S_H^m\times S_H^n$, then there exist $m+n$
$\R^k$-valued random vectors $\textbf{X}_1, \ldots, \textbf{X}_m, \textbf{Y}_1, 
\ldots, \textbf{Y}_n$, such that $\vc{\textbf{X}_i}{\textbf{Y}_j} \sim 
N_{2k}(0, \Sigma_{2k}(\rho_{ij}))$ for all $(i,j) \in [m] \times [n]$, 
and 
\begin{align}\label{eq:rep_of_c_times_QC_real_case}
\Gamma_H(x,y) = \frac{1}{c^\ast} h_{f,g}[R] =
\frac{1}{c^\ast}\,\E[{\textbf{P}}_f{\textbf{Q}}_g^\top] =
\frac{1}{c^\ast}\,\E[\Gamma_{\R}({\textbf{P}}_f, {\textbf{Q}}_g)],
\end{align}
where $({\textbf{P}}_f)_i : = f(\textbf{X}_i)$, 
$({\textbf{Q}}_g)_j : = g(\textbf{Y}_j)$, $\rho_{ij} : = 
h_{\alpha, \beta}(\langle x_i, y_j\rangle_H)$ and $R : = 
h_{\alpha, \beta}[\Gamma_H(x,y)] = (\rho_{ij})_{(i,j) \in [m] \times [n]}
\in \mathcal{Q}_{m,n}$. 
\item 
Moreover, $c^\ast \langle u, v\rangle_H \in (-1,1)$ for all $u, v \in S_H$ and
\begin{align}\label{eq:continuity_of_the_inverse_and_trace_equality}
\text{tr}(A^\top S) = \frac{1}{c^\ast}\,\text{tr}(A^\top\,\psi[Q_{c^\ast, \psi}])
\end{align}
for all $m,n \in \N$, for all $A \in \M_{m,n}(\R)$, for all $S \in 
\mathcal{Q}_{m,n}$, where $Q_{c^\ast, \psi} : = \psi^{-1}[c^\ast\,S] \in 
\mathcal{Q}_{m,n}$. 
\end{enumerate}
\end{theorem}
\noindent It is far from being trivial that it is possible to transfer Theorem 
\ref{thm:real_inner_product_rounding}
from the real field $\R$ to the complex field $\C$; at least if $f = g$ is odd 
(cf. \autoref{thm:complex_inner_product_rounding}).
In order to achieve this, we have to develop and implement certain non-trivial 
structural properties of the class of complex Hermite polynomials 
(cf. \autoref{thm:correlated_complex_Hermite_polynomials} and Theorem 
\ref{thm:odd_completely_real_analytic_functions_at_0_vs_complex_h_fg} below).

If $\nu \in \N_0$, $m, n \in \N$ and $A \in \M_{m,n}(\F)$, then $A^{\ast\,\nu}$ 
denotes the $\nu$-th entrywise power of the matrix $A$ (in terms of Schur 
multiplication), where $A^{\ast\,0} : = {\bf{1}}_m {\bf{1}}_n^\ast$ is the (rank 1) 
$m \times n$ matrix of all ones, where ${\bf{1}}_l : = (1, 1, \ldots, 1)^\top 
\in \F^l, l \in \N$ (by adapting the convention that $0^0 : = 1$ - cf. 
\cite[Remark 9.2]{K2022}). 

Now, we are fully prepared to embed both, Grothendieck's original estimation 
$K_G^\R \leq \sinh(\frac{\pi}{2}) \approx 2.301$, and Krivine's original estimation 
$K_G^\R \leq \frac{\pi}{2 \ln(1 + \sqrt{2})} \approx 1.782$ into a general framework. 
In particular, we will provide a further proof of the Grothendieck inequality itself 
(see \autoref{thm:GT_generalising_Krivine} below). 
Moreover, we will make a cute use of the little Grothendieck inequality, yielding 
a quite surprising outcome. To this end, we have to work with functions $h_{f,g}$, 
which are ``generated'' by \textit{bounded functions $f,g \in L^\infty(\R^k, \gamma_k) 
= L^\infty(\R^k, \lambda_k) \equiv L^\infty(\R^k)$, $k \in \N$}. 
\eqref{eq:boundedness_of_h} obviously implies that 
\[
\vert h_{f,g}(\rho)\vert \leq \Vert f \Vert_{\infty}\,\Vert g \Vert_{\infty} 
\,\text{ for all }\,f, g \in L^\infty(\R^k)\, \text{ and }\, \rho \in [-1,1]\,.
\]
\begin{remark}\label{rem:concepts_in_ML}
In statistical machine learning a function $f : \R^k \longrightarrow \{-1,1\}, 
k \in \N$ is a particular example of a mapping from some general domain of 
definition to the set $\{-1,1\}$, where the latter is known as ``concept'' 
(cf. e.g. \cite{S2009}). The concept $f$ obviously satisfies the condition 
$\Vert f \Vert_{\gamma_k} = 1 = \Vert f \Vert_\infty$.
\end{remark} 
\noindent We also need the following important result which holds for both fields 
interchangeably.
\begin{lemma}\label{lem:norm_estimation}
Let $\F \in \{\R, \C\}$ and $m, n \in \N$. Let ${\textbf{R}}$ be a 
$m$-dimensional random vector in $(\F \cap \overline{\D})^m$ and ${\textbf{S}}$ be a 
$n$-dimensional random vector in $(\F \cap \overline{\D})^n$, such that 
$\E[{\textbf{R}}{\textbf{S}}^\ast]$ exists. Then
\[
\Vert B \ast \E[{\textbf{R}}{\textbf{S}}^\ast]^{\ast\nu} \Vert^\F_{\infty, 1} 
\leq \Vert B \Vert^\F_{\infty, 1} \text{ for all } B \in \M_{m,n}(\F) \text{ and } \nu \in 
\N\,.
\]
\end{lemma}
\begin{theorem}\label{thm:GT_generalising_Krivine}
Let $k, m, n \in \N$ and $f, g \in L^\infty(\R^k)$. Put 
$r_\infty \equiv r_\infty(f,g) : = \Vert f \Vert_\infty\,\Vert g \Vert_\infty$. 
Then:
\begin{enumerate}
\item
\begin{align}\label{eq:property_of_quantum_corr_matrices}
\vert{\text{tr}}(A^\top h_{f,g}[S])\vert \leq r_\infty\,
\Vert A \Vert_{\infty, 1} \text{ for all } A \in \M_{m,n}(\R) \text{ and } S \in 
{\mathcal{Q}}_{m,n}
\end{align}
and
\begin{align}\label{eq:property_of_corr_matrices}
\vert{\text{tr}}(A^\top h_{f,f}[\Sigma])\vert \leq \Vert f \Vert_\infty^2\,
\max\limits_{x \in [-1,1]^n}\vert x^\top Ax \vert \text{ for all } A \in \M_n(\R) 
\text{ and } \Sigma \in C(n; \R).
\end{align}
\item Assume that $r_\infty > 0$ and $f$ or $g$ is odd. Let 
$h_{f,g} : [-1,1] \longrightarrow \R$ be continuous and injective. Then 
$h_{f,g}$ is a homeomorphism which is either strictly increasing or strictly 
decreasing and satisfies $h_{f,g}([-1,1]) = [-r,r]$, where $r \equiv r_k(f,g) 
: = \max\{-h_{f,g}(1), h_{f,g}(1)\} = \max\{-\langle f, g\rangle_{\gamma_k}, 
\langle f, g\rangle_{\gamma_k}\}$. Moreover, $0 < r \leq r_\infty$. Assume that 
$h_{f,g}^{-1}{\big\vert}_{(-r, r)} \in W^\omega_+((-r, r))$. Then the following 
two statements hold:
\begin{enumerate}
\item[\textup{(ii-1)}] If $\big(h_{f,g}^{-1}{\big\vert}_{(-r, r)}\big)_
{\text{abs}}(r) > 1$, then there is exactly one number $0 < c^\ast_k \equiv 
c^\ast_k(f,g) < r$, such that 
$\big(h_{f,g}^{-1}{\big\vert}_{(-r, r)}\big)_{\text{abs}}(c^\ast_k) = 1$ and 
\begin{align}\label{eq:upper_bound_real_case}
K_G^\R \leq \frac{r_\infty}{c^\ast_k}\,.
\end{align}
\item[\textup{(ii-2)}] If $ r = r_\infty$, 
then there exists a unique number $0 < \gamma^\ast_k \equiv \gamma^\ast_k(f,g) \in 
(0,r]$, such that  
\begin{align}\label{eq:real_GT_estimate_and_abs_fct}
K_G^\R = \big(h_{f,g}^{-1}{\big\vert}_{(-r,r)}\big)_{\text{abs}}
(\gamma^\ast_k) \leq \min\Big\{\frac{r}{c^\ast_k}\,,
\big(h_{f,g}^{-1}{\big\vert}_{(-r, r)}\big)_{\text{abs}}(r)\Big\}\,.
\end{align}
\end{enumerate}
\end{enumerate}
\end{theorem}
\noindent Inequality \eqref{eq:property_of_corr_matrices}, together with 
the equality \eqref{eq:d_1_psd_case} allows us to recover as special case 
(applied to $f = \text{sign}$) an interesting result of Y. Nesterov, yet 
without having to make use of their random hyperplane rounding technique 
(cf. \cite{N1998}):   
\begin{corollary}\label{cor:Nesterov_as_special_case}
Let $k, n \in \N$, $A \in \M_n(\R)^+$ be positive semidefinite and $f \in 
L^\infty(\R^k)$ be odd. Then
\begin{align*}
h_{f,f}^\prime(0)\sup\limits_{\Sigma \in C(n;\R)}{\text{tr}}(A \Sigma) 
&\leq \sup\limits_{\Sigma \in C(n;\R)}{\text{tr}}(A\,h_{f,f}[\Sigma]) \leq
\Vert f \Vert_\infty^2\,\sup\limits_{x \in \{-1,1\}^n}x^\top Ax \leq 
\Vert f \Vert_\infty^2\,\Vert A \Vert_{\infty, 1}\,.
\end{align*}
In particular,
\begin{align}\label{eq:upper_bound_of_derivative_of_CCP_at_zero} 
0 \leq h_{f,f}^\prime(0) \leq \frac{2}{\pi}\,\Vert f \Vert_{\infty}^2\,. 
\end{align}
\end{corollary}
\noindent In a similar vain, we directly obtain a further (and very short) proof 
for the value of the real little Grothendieck constant $k_G^\R$. To this end, we 
only have to combine \eqref{eq:property_of_quantum_corr_matrices} and 
\sref{Proposition}{prop:2F1_rep}:
\begin{corollary}[\textbf{Grothendieck, 1953}]
\label{rem:further_proof_re_real_little_GT_constant}
\[
k_G^\R = \frac{\pi}{2}.
\]
\begin{comment}
Let $A \in \M_n(\R)^+$ and $\Sigma \in C(n;\R)$ be arbitrary. We just have to 
work with the CCP function $[-1,1] \ni \rho \mapsto h_{f_1,f_1}(\rho) - 
c_1^2\,\rho = \frac{2}{\pi}\arcsin(\rho) - \frac{2}{\pi}\,\rho$. Because then
\[
0 \leq \tfrac{2}{\pi}\,\text{tr}(A\Sigma) = \text{tr}(A(\tfrac{2}{\pi}\,\Sigma)) \leq 
\text{tr}(A h_{f_1,f_1}[\Sigma]) 
\stackrel{\eqref{eq:property_of_quantum_corr_matrices}}{\leq} 
\Vert A\Vert_{\infty,1}\,.
\]
\end{comment}
\end{corollary}
\begin{remark}
A direct estimation of \eqref{eq:derivative_of_CCP_at_zero} leads to the upper bound 
$k\,\Vert f \Vert_{\infty}^2\, (\sqrt{\frac{2}{\pi}})^2 = 
\frac{2k}{\pi}\,\Vert f \Vert_{\infty}^2$, which, however, strongly depends 
on the dimension $k \in \N$. That upper bound, viewed as a function of $k$ is even strictly 
increasing. Our application of the little Grothendieck inequality implies the non-trivial 
result that \textit{for all $k \in \N$}, $h_{f,f}^\prime(0)$ actually is bounded above by 
$\frac{2}{\pi} \Vert f \Vert_\infty^2$ ``uniformly''. Moreover, if $f$ were an even function, 
then \eqref{eq:upper_bound_of_derivative_of_CCP_at_zero} would be trivial, since 
$h_{f,f}^\prime(0) = 0$. 
\end{remark}
\noindent \autoref{thm:GT_generalising_Krivine}-\ref{(i)} also implies a remarkable 
property of CCP functions. To this end, let $S \in \mathcal{Q}_{m,n}$ and 
$h = h_{f,f}$ be an arbitrary CCP function (not necessarily odd\,!), generated 
by some $f \in L^\infty(\R^k)$. Firstly, note that in any case,
\[
\frac{\pi}{2}\frac{1}{K_G^\R} 
\stackrel{\eqref{eq:known_upper_bounds_of_the_complex_GT_constant}}{<} 1 = 
\Vert f \Vert_{\gamma_k} \leq \Vert f \Vert_\infty \leq \Vert f \Vert^2_\infty\,.
\]
If also $h_{f,f}^{-1}$ were a CCP function, then $\widetilde{S} : = 
h_{f,f}^{-1}[S] \in \mathcal{Q}_{m,n}$ (due to 
\eqref{eq:entrywise_mapping_of_inner_products_to_inner_products}). Hence,
\[
\vert{\text{tr}}(A^\top S)\vert = 
\vert{\text{tr}}(A^\top h_{f,f}[\widetilde{S}])\vert \leq 
\Vert f \Vert^2_\infty\,\Vert A \Vert_{\infty, 1} \text{ for all } A \in 
\M_{m,n}(\R),
\] 	 
implying that $K_G^\R \leq \Vert f \Vert^2_\infty$. If we join the latter 
observation and \autoref{thm:CCP_via_Schoenberg}-\ref{(v)}, we obtain another 
interesting fact (which should be carefully compared with 
\sref{Corollary}{cor:inverse_of_non_trivial_odd_CCP_is_not_CCP}):
\begin{remark}\label{rem:inverse_of_CCP_in_general_is_not_CCP}
Let $k, m, n \in \R$ and $h = h_{f,f}$ be CCP for some $f \in S_{L^2(\gamma_k)} 
\cap L^\infty(\R^k)$. If $\Vert f \Vert_\infty < \sqrt{K_G^\R}$, then the 
inverse function $h_{f,f}^{-1}$ is not CCP. 
\end{remark} 
\noindent Next, we are going to summarise the remarkable properties of 
\textit{odd CCP functions} which we found so far and shed some light on an 
additional, quite surprising estimate for odd CCP functions $h_{f,f}$ which 
emerges if we assume in addition that $f \in S_{L^2(\gamma_k)}$ is 
\textit{bounded} (a.s.); i.e., if $f \in S_{L^2(\gamma_k)} \cap L^\infty(\R^k)$. 
In particular, if we combine \autoref{thm:h_ff_real_odd_case}-\ref{(ii)}, 
\ref{(v)} and \autoref{thm:GT_generalising_Krivine}-\ref{(ii)} in this 
case, we recover Grothendieck's upper bound as well as Krivine's upper bound 
at once. This follows from \sref{Example}{ex:Krivine_recovered}, which is a special 
case of our following key result for the real odd CCP case:
\begin{theorem}\label{thm:odd_real_bdd_CCP}
Let $k \in \N$ and $f \in S_{L^2(\gamma_k)} \cap L^\infty(\R^k)$ such that 
$\Vert f \Vert_{\gamma_k} = 1$. Then $\Vert f \Vert_\infty \geq 1$. Assume 
that $f$ is odd and $h_{f,f}^{-1}{\big\vert}_{(-1,1)} \in W^\omega_+((-1, 1))$. 
Then the following statements hold:
\begin{enumerate}
\item 
\[
\big(h_{f,f}^{-1}{\big\vert}_{(-1,1)}\big)_{\text{abs}}(y) \geq \frac{\pi}{2}\,
\frac{1}{\Vert f \Vert_\infty^2}\,y \text{ for all } y \in [0,1]\,.
\]
\item
\begin{align*}
K_G^\R \leq \frac{\Vert f \Vert_{\infty}^2}{h_{f,f}^\text{hyp}(1)}\,.
\end{align*}
{In particular, if $\text{sign}\big((h_{f,f}^{-1}{\vert}_{(-1,1)})^
{(2n+1)}(0)\big) = (-1)^n$ for all $n \in \N_0$, then
\[
K_G^\R \leq i\,\frac{\Vert f \Vert_{\infty}^2}{\widetilde{\psi_f}(i)}\,,
\]
where $\psi_f : = h_{f,f}{\big\vert}_{(-1,1)}$ and $\widetilde{\psi_f} : 
\overline{\D} \longrightarrow \overline{\D}$
is defined as in \sref{Lemma}{lem:real_analytic_and_l1}.
}
\item
Let $1 \leq c_\ast < K_G^\R$. If $\Vert f \Vert_\infty = 1$, then $0 < 
h_{f,f}^\text{hyp}(c_\ast) < 1$ and there is exactly one number $\gamma^\ast(f) 
\in (h_{f,f}^\text{hyp}(c_\ast) ,1]$, such that  
\begin{align}\label{eq:upper_bound_real_odd_unit_sphere_case}
K_G^\R = \big(h_{f,f}^{-1}{\big\vert}_{(-1,1)}\big)_{\text{abs}}(\gamma^\ast(f)) 
\leq \min\Big\{\frac{1}{h_{f,f}^\text{hyp}(1)},\, 
\big(h_{f,f}^{-1}{\big\vert}_{(-1,1)}\big)_{\text{abs}}(1)\Big\}\,.
\end{align}
\end{enumerate}
\end{theorem}
\begin{comment}
 Since $f \in L^\infty$, it follows that $f \in L^2(\gamma_k)$ and 
$1 = \Vert f \Vert_{\gamma_k} \leq \Vert f \Vert_\infty$. 
\\[0.2em]
(i) and (ii) The claimed inequalities directly follow from 
\eqref{eq:psi_abs_estimation}, \eqref{eq:abs_h_f_f_inverse_estimation}, 
\sref{Lemma}{lem:hyperbolic_CCP_transform}, 
{\eqref{eq:sign_condition_and_hyp_CCP_fct_value}},
\autoref{thm:GT_generalising_Krivine}, and 
\eqref{eq:upper_bound_of_derivative_of_CCP_at_zero}.
\\[0.2em]
\noindent (iii) Since $1 \leq c_\ast < K_G^\R \leq \big(h_{f,f}^{-1}{\big\vert}_
{(-1,1)}\big)_{\text{abs}}(1)$ (due to \autoref{thm:GT_generalising_Krivine}), 
it follows from \sref{Lemma}{lem:hyperbolic_CCP_transform} that $0 < 
h_{f,f}^\text{hyp}(1) \leq h_{f,f}^\text{hyp}(c_\ast) < \gamma^\ast(f) : = 
h_{f,f}^\text{hyp}(K_G^\R) \leq 1$. (iii) now follows from (ii).
\end{comment}
\noindent If we only assume that $f \in L^\infty(\R^k)\setminus\{0\}$, then an 
application of \autoref{thm:odd_real_bdd_CCP} to $\frac{1}{\Vert f \Vert_
{\gamma_k}}f \in S_{L^2(\gamma_k)}$ directly leads to
\begin{corollary}\label{cor:odd_real_bdd_CCP_general_version}
Let $k \in \N$ and $f \in L^\infty(\R^k)\setminus\{0\}$. Then $0 < r \equiv 
r_k(f) : = \Vert f\Vert_{\gamma_k}^2 < \infty$ and $\Vert f \Vert_\infty \geq 
\sqrt{r}$. Assume that $f$ is odd and $h_{f,f}^{-1}{\big\vert}_{(-r,r)} \in 
W^\omega_+((-r, r))$. Then $\big(h_{f,f}^{-1}{\big\vert}_{(-r,r)}\big)_
{\text{abs}}(r) > 1$, and the following statements hold:
\begin{enumerate}
\item
\[
\big(h_{f,f}^{-1}{\big\vert}_{(-r,r)}\big)_{\text{abs}}(y) \geq \frac{\pi}{2}\,
\frac{1}{\Vert f \Vert_\infty^2}\,y \,\text{ for all }\, y \in [0,r]\,.
\]
\item
\begin{align}\label{eq:upper_bound_real_general_bdd_case}
r\,K_G^\R \leq \frac{\Vert f \Vert_{\infty}^2}{h_{\frac{f}{\sqrt{r}},
\frac{f}{\sqrt{r}}}^\text{hyp}(1)}\,.
\end{align}
{In particular, if $\text{sign}\big((h_{f,f}^{-1}
{\big\vert}_{(-r,r)})^{(2n+1)}(0)\big) = (-1)^n$ for all $n \in \N_0$, then
\[
K_G^\R \leq i\,\frac{\Vert f \Vert_{\infty}^2}{\widetilde{\psi_f}(i)}\,,
\]
where $\psi_f : = h_{f,f}{\big\vert}_{(-1,1)}$ 
and $\widetilde{\psi_f} : \overline{\D} \longrightarrow r\,\overline{\D}$ is 
defined as in \sref{Lemma}{lem:real_analytic_and_l1}.
}
\item
Let $1 \leq c_\ast < K_G^\R$. If $\Vert f \Vert_\infty = \sqrt{r}$, then 
$0 < h_{\frac{f}{\sqrt{r}},\frac{f}{\sqrt{r}}}^\text{hyp}(c_\ast) < 1$ and there 
is exactly one number number $\gamma^\ast(f) \in (h_{\frac{f}{\sqrt{r}},
\frac{f}{\sqrt{r}}}^\text{hyp}(c_\ast) ,1]$, such that  
\begin{align*}
K_G^\R = \big(h_{f,f}^{-1}{\big\vert}_{(-r,r)}\big)_{\text{abs}}(r\,
\gamma^\ast(f)) \leq \min\Big\{\frac{1}{h_{\frac{f}{\sqrt{r}},
\frac{f}{\sqrt{r}}}^\text{hyp}(1)},\, \big(h_{f,f}^{-1}{\big\vert}_{(-r,r)}\big)_
{\text{abs}}(r)\Big\}\,.
\end{align*}
\end{enumerate}
\end{corollary}
\noindent The proof of \autoref{thm:odd_real_bdd_CCP}-\ref{(i)} shows us that 
here we may circumvent the rather strong assumption of $h_{f,f}^{-1}{\big\vert}_
{(-r,r)}$ being completely real analytic on $(-r,r)$ at $0$; at least in the 
following sense (cf. \sref{Example}{ex:Naor_et_al_s_result} as application for 
this):
\begin{proposition}\label{prop:thm_odd_real_bdd_CCP_without_CRA_condition}
Let $k \in \N$ and $f \in L^\infty(\R^k)\setminus\{0\}$. Then $0 < r \equiv 
r_k : = \Vert f\Vert_{\gamma_k}^2 < \infty$ and $\Vert f \Vert_\infty \geq 
\sqrt{r}$. Assume that $f$ is odd, $h_{f,f}^\prime(0) > 0$ and
\begin{align}\label{eq:root_condition}
\sum_{n=0}^\infty \frac{\vert(h_{f, f}^{-1})^{(2n+1)}(0)\vert}{(2n+1)!}\,
(c^\ast)^{2n+1} = 1\,,
\end{align}
for some $c^\ast \in (0, r]$. Then $h_{f,f}^{-1}{\big\vert}_{(-c^\ast,c^\ast)} \in 
W^\omega_+((-c^\ast, c^\ast))$, and the following 
statements hold:
\begin{enumerate}
\item
\[
\big(h_{f,f}^{-1}{\big\vert}_{(-c^\ast, c^\ast)}\big)_{\text{abs}}(y) \geq 
\frac{\pi}{2}\,\frac{1}{\Vert f \Vert_\infty^2}\,y \,\text{ for all } \, y \in 
[0, c^\ast] \,\text{ and }\,\big(h_{f,f}^{-1}{\big\vert}_{(-c^\ast,c^\ast)}\big)_
{\text{abs}}(c^\ast) = 1\,.
\]
\item
\begin{align}\label{eq:entrywise_mapping_of_inner_products_to_inner_products_III}
h_{f,f}^{-1}[c^\ast\,S] \in {\mathcal{Q}}_{m,n}\,\text{ for all }\, S \in
{\mathcal{Q}}_{m,n}\,. 
\end{align}
\item 
\begin{align*}
K_G^\R \leq \frac{\Vert f \Vert_{\infty}^2}{c^\ast}\,.
\end{align*}
\end{enumerate}
\end{proposition}
\begin{comment}
(i) The assumption clearly implies that $h_{f,f}^{-1}{\big\vert}_{(-c^\ast,c^\ast)} 
\in W^\omega_+((-c^\ast, c^\ast))$ and $\big(h_{f,f}^{-1}{\big\vert}_
{(-c^\ast, c^\ast)}\big)_{\text{abs}}(c^\ast) = 1$.
\\[0.2em]
\noindent (ii) and (iii) Observe that in particular the restriction $h_{f,f}^{-1}
{\big\vert}_{[-c^\ast, c^\ast]} : [-c^\ast, c^\ast] \longrightarrow [-1,1]$ is 
continuous. Consequently, if $S \in \mathcal{Q}_{m,n}$, then 
$S_0 : = h_{f,f}^{-1}[c^\ast\,S] \in \mathcal{Q}_{m,n}$ (due to 
\eqref{eq:entrywise_mapping_of_inner_products_to_inner_products_II}), and (ii) 
follows. Let $m, n \in \N$, $A \in \M_{m,n}(\R)$ and $S \in \mathcal{Q}_{m,n}$ 
be arbitrarily given. Since $S_0 \in \mathcal{Q}_{m,n}$, we therefore obtain
\begin{align*}
\vert{\text{tr}}(A^\top S) \vert &= \frac{1}{c^\ast}\,
\vert{\text{tr}}(A^\top h_{f,f}[S_0])\vert 
\stackrel{\eqref{eq:property_of_quantum_corr_matrices}}{\leq} \frac{1}{c^\ast}\,
\Vert A \Vert_{\infty, 1}\,,
\end{align*}
and (iii) follows as well.
\end{comment}  
\begin{example}[\textbf{Grothendieck and Krivine}]\label{ex:Krivine_recovered}
Once again, we consider the CCP function $\psi : = \frac{2}{\pi}\arcsin$. Recall 
that $\psi = h_{f,f}$, where $f : = \text{sign} \in S_{L^\infty(\R)} \cap 
S_{L^2(\gamma_k)}$. Due to \eqref{eq:hyp_trsfm_of_Krivine_fct}, it 
follows that
\[
\psi^{\text{hyp}}(1) = \frac{2}{\pi}\,\sinh^{-1}(1) = \frac{2}{\pi}\,
\ln(1 + \sqrt{2}).
\]
Hence, 
\[
K_G^\R \stackrel{\eqref{eq:upper_bound_real_odd_unit_sphere_case}}{\leq}
\min\Big\{\frac{1}{\psi^{\text{hyp}}(1)},\,\big(\psi^{-1}{\big\vert}_
{(-1,1)}\big)_{\text{abs}}(1)\Big\} \leq \frac{\pi}{2\ln(1 + \sqrt{2})} 
(\approx 1.78221) \leq \sinh(\frac{\pi}{2}) (\approx 2.30129)
\] 
precisely reflects Krivine's upper bound of 
$K_G^\R$ as well as Grothendieck's (larger) upper bound of 
$K_G^\R$\,! 
\end{example}
\begin{example}
Consider the CCP function $\kappa : = \sqrt{3}(2\phi-1)$ (cf. 
\sref{Proposition}{prop:Phi_minus_one_half}). Due to Theorem 
\ref{thm:odd_real_bdd_CCP}, (i), respectively 
\eqref{eq:upper_bound_real_general_bdd_case}, we obtain the following (weaker) 
estimation:
\[
K_G^\R \leq \frac{\pi}{6\,\ln(\tfrac{1}{2}(1 + \sqrt{5}))}\cdot 
\Vert \kappa\Vert_{\infty}^2 \leq \frac{\pi}{2\,\ln(\tfrac{1}{2}(1 + \sqrt{5}))} 
(\approx 3.26425)\,.
\]
We highly recommend the readers to check whether this estimation can be improved, 
if more generally the function $\kappa_\alpha : = \alpha (2\phi-1)$ is considered, 
where $0 < \alpha < \sqrt{3}$ is given (instead of the CCP function $\kappa = 
\kappa_{\sqrt{3}}$)! Observe also that $\Vert \kappa\Vert_{\infty} = \sqrt{3} 
\not=1 = \Vert \kappa\Vert_{\gamma_k}$.
\end{example}
\begin{example}
Fix $k \in \N_3$ and consider the function $\psi : = h_{f_k,f_k}$, introduced in
\sref{Proposition}{prop:2F1_rep}. \textit{Assume that} $\psi^{-1} \in 
W^\omega_+((-1, 1))$. $\psi$ then satisfies all assumptions, listed in 
\autoref{thm:odd_real_bdd_CCP}, and it follows that $0 < c^\ast_k : = 
h_{f_k,f_k}^{\text{hyp}}(1) < 1$ satisfies
\begin{align}\label{eq:upper_bound_real_CCP_case}
K_G^\R \leq \frac{k}{c^\ast_k}
\end{align}
(since $\Vert f_k \Vert_\infty^2 = k$). However, observe that the sequence 
$\big(\frac{k}{c^\ast_k}\big)_{k \in \N}$ is not bounded and hence cannot 
converge (since $0 < c_k^\ast \leq 1$). Moreover, in contrast to the previous 
two examples, we do not know whether also $\big(h_{f_k,f_k}^{-1}{\big\vert}_
{(-1,1)}\big)_{\text{abs}}$ can be represented in a closed analytical form; one 
of the major open problems in our search for the smallest upper bound of 
$K_G^\R$ (cf. \sref{Section}{sec:RP_1} and \sref{Example}{ex:Haagerup}, 
where the latter includes the approximation of the constant $c_k^\ast \approx 
0.71200$ if $k=2$). A straightforward, yet a bit laborious calculation with 
fractions, based on the table \eqref{eq:First_7_Taylor_coeff_of_the_inverse} 
(similarly to the special case $k=2$, studied in \sref{Example}{ex:Haagerup}), 
yields that the Maclaurin series of $h_{f_k,f_k}^{-1}(s) \equiv 
\sum_{\nu = 0}^\infty \beta_{2\nu+1}(k)\,s^{2\nu+1}$ can e.g. be approximated 
by the Taylor polynomial of degree 7 as:
\begin{align}
\begin{split}
h_{f_k,f_k}^{-1}(s) &= \beta_1(k)s + \beta_3(k)s^3 + \beta_5(k)s^5 + 
\beta_7(k)s^7 + o(\vert s\vert^7)\\
&= \frac{1}{c_k^2}\,s + \big(-\frac{1}{c_k^6}\frac{1}{2(k+2)}
\,s^2\big)\big(s + \frac{3}{4 c_k^4}\,\frac{k-2}{(k+2)(k+4)}\,s^3\\ 
&+ \frac{3}{8 c_k^8}\,\frac{9(k+2)^2-48(k+2)+64}{(k+2)^2\,(k+4)\,(k+6)}\,s^5\big) + 
o(\vert s\vert^7). 
\end{split} 
\end{align}
Observe also that $\beta_5(k) : = -\frac{3}{8 c_k^{10}}\,\frac{k-2}{(k+2)^2\,(k+4)} = 0$ 
if and only if $k = 2$ and that $\beta_5(k) \leq 0$ if and only if $k \geq 2$. 
Similarly, since the single (local) minimum of the function $\R \ni x \mapsto 
g(x) : = 9x^2-48x+64$ is attained at $x_\ast : = \frac{8}{3}$ and $g(x_\ast) = 0$, 
it follows that $\beta_7(k) \leq 0$.
\end{example}
\begin{example}
\label{ex:Naor_et_al_s_result}
{Our approach also can be applied to slightly modify the proof of the 
strongest result to date, namely that $K_G^\R < \frac{\pi}{2 \ln(1 + \sqrt{2})}$.} 
To this end, firstly observe that a simple change 
of variables reveals that \textit{for any} $f,g \in L^2(\gamma_k)$, the 
corresponding generalised function $H_{f\circ\sqrt{2}, g\circ\sqrt{2}}$, listed 
in \cite[Definition 2.1]{BMMN2013}, satisfies
\begin{align}\label{eq:Krivine_RS_vs_h_f_g} 
H_{f\circ\sqrt{2}, g\circ\sqrt{2}} = h_{f,g}\big\vert_{(-1,1)}
\end{align}
(due to \eqref{eq:h_fg_correl_representation}). So, $H_{f\circ\sqrt{2}, 
g\circ\sqrt{2}}$ is well-defined on $(-1,1)$. The Grothendieck function 
$\frac{2}{\pi}\arcsin$ is then generalised in \cite{{BMMN2013}} to the 
complex-valued function 
\begin{align}\label{eq:Naor_et_al_s_function}
F_{p, \eta} : = (1-p)H_0 + p H_\eta\,, 
\end{align}
where $0 \leq p \leq 1$, $0 \leq \eta < 1$ and $H_\eta : \S \longrightarrow \C$ 
is defined as in \cite[(41)]{BMMN2013}, where $\S : = \{z \in \C : 
\vert\Re(z)\vert < 1\}$. Independent of any complex analysis a priori, the 
construction of $H_\eta$, together with \eqref{eq:Krivine_RS_vs_h_f_g} implies 
that 
\[
H_\eta\big\vert_{\S \cap \R} = H_\eta\big\vert_{(-1,1)} = 
H_{g_\eta\circ\sqrt{2}, g_\eta\circ\sqrt{2}} 
\stackrel{\eqref{eq:Krivine_RS_vs_h_f_g}}{=} h_{g_\eta, g_\eta}\big\vert_{(-1,1)}\,.
\]  
Here, the odd function $g_\eta : \R^2 \longrightarrow \R$ is defined as 
$g_\eta(x) : = \sqrt{\frac{\pi}{2}}\,f_\eta(\frac{1}{\sqrt{2}}x)$, where 
the 2-dimensional sign concept $f_\eta : \R^2 \longrightarrow \{-1,1\}$ 
satisfies \cite[(40)]{BMMN2013}. Consequently, the function 
$\frac{2}{\pi}\,H_\eta\big\vert_{(-1,1)}$ actually is a restriction of the 
well-defined odd - and hence invertible - CCP function $\frac{2}{\pi}\,
h_{g_\eta, g_\eta} = h_{\sqrt{2/\pi}\,g_\eta, \sqrt{2/\pi}\,g_\eta}$ to the 
open interval $(-1,1)$ (due to \autoref{thm:real_CCP_rep_thm}).

 Let $\rho \in (-1,1)$ and $\vc{\textbf{X}}{\textbf{Y}} \sim 
N_{4}(0, \Sigma_{4}(\rho))$. Since $g_0(x) = \sqrt{\frac{\pi}{2}}\,
{\text{sign}}(x_2)$ for all $x = (x_1, x_2)^\top \in \R^2$, we may apply 
\sref{Proposition}{prop:char_of_the_Sigma_2n_prob_law} to the partitioned 
Gaussian random vector $\vc{\textbf{X}}{\textbf{Y}}$, and it follows that
\[
H_0(\rho) = h_{g_0, g_0}(\rho) = \frac{\pi}{2}\,
\E_{\tiny{\P_{\vc{\textbf{X}}{\textbf{Y}}}}}[\text{sign}(X_2)\,\text{sign}(Y_2)] =
\arcsin(\rho),
\]
which {slightly shortens the} proof of \cite[Lemma 4.3]{BMMN2013}. In particular, 
for any $p \in [0,1]$,  
\begin{align}\label{eq:convex_comb_of_CCPs}
\frac{2}{\pi}\,F_{p, \eta}\big\vert_{(-1,1)} = \psi_{p, \eta}\big\vert_{(-1,1)}
\end{align}
emerges as a restriction of the odd and hence \textit{strictly increasing, 
homeomorphic} CCP function $\psi_{p, \eta} : = (1-p)\frac{2}{\pi}\arcsin + 
p\,\frac{2}{\pi}\,h_{g_\eta, g_\eta}$ on $(-1,1)$ (due to 
\autoref{thm:CCP_via_Schoenberg} and \autoref{thm:h_ff_real_odd_case}-\ref{(iii)}).
Observe that $\psi_{p, \eta}^\prime(0) \geq (1-p)\frac{2}{\pi}$, implying 
that $\psi_{p, \eta}^{-1}\big\vert_{(-1,1)} \in C^\omega((-1,1))$ if $p < 1$ 
(due to \autoref{thm:h_ff_real_odd_case}-\ref{(iv)}). Because of the non-trivial 
result \cite[Theorem 5.1]{BMMN2013} (including its technically demanding proof) 
it follows the existence of $(p_0, \eta_0) \in (0,1) \times (0,1)$ and 
$c^\ast \in (\frac{2}{\pi}\ln(1+\sqrt{2}), \frac{9}{5\pi})$, such that 
\[
\sum_{n=0}^\infty \frac{\vert(\psi_{p_0, \eta_0}^{-1})^{(2n+1)}(0)\vert}
{(2n+1)!}\,(c^\ast)^{2n+1} = 1\,.
\]
(Given the outcome of \cite[Theorem 5.1]{BMMN2013}, we only have to set 
$c^\ast : = \frac{2\gamma}{\pi}$, where $\ln(1+\sqrt{2}) < \gamma < \frac{9}{10}$
satisfies \cite[(63)]{BMMN2013}.) Hence, $\psi_{p_0, \eta_0}^{-1}\big\vert_
{(-c^\ast,c^\ast)} \in W^\omega_+((-c^\ast, c^\ast))$. 
Since in particular $\psi_{p_0, \eta_0}^{-1}\big\vert_{[-c^\ast,c^\ast]}$ is 
continuous, we may apply 
\sref{Proposition}{prop:thm_odd_real_bdd_CCP_without_CRA_condition}. Consequently, 
if $S \in \mathcal{Q}_{m,n}$, then $S_0 : = \psi_{p_0, \eta_0}^{-1}[c^\ast\,S] 
\in \mathcal{Q}_{m,n}$ (due to 
\eqref{eq:entrywise_mapping_of_inner_products_to_inner_products_III}). 
Observe that for any $0 \leq \eta < 1$ the CCP function $\frac{2}{\pi}\,
h_{g_\eta, g_\eta} = h_{\sqrt{2/\pi}\,g_\eta, \sqrt{2/\pi}\,g_\eta}$ actually 
originates from the \textit{bounded} function $\sqrt{\frac{2}{\pi}}\,
g_\eta \in S_{L^\infty} \cap S_{L^2(\gamma_2)}$, such as $\frac{2}{\pi}\arcsin = 
h_{\text{sign}, \text{sign}}$. Let $S \in \mathcal{Q}_{m,n}$ be arbitrarily 
given. Since $S_0 \in \mathcal{Q}_{m,n}$, we therefore obtain
\begin{align*}
\vert{\text{tr}}(A^\top S) \vert &= \frac{1}{c^\ast}\,
\vert{\text{tr}}(A^\top \psi_{p_0, \eta_0}[S_0])\vert = \frac{1}{c^\ast}\,
\big\vert{\text{tr}}(A^\top \big((1-p)\frac{2}{\pi}\arcsin[S_0] + 
p\,\frac{2}{\pi}\,h_{g_\eta, g_\eta}[S_0]\big)\big\vert\\
&\stackrel{\eqref{eq:property_of_quantum_corr_matrices}}{\leq} 
\frac{1}{c^\ast}\,((1-p)\Vert A \Vert_{\infty, 1} + p\,\Vert A \Vert_{\infty, 1}) 
= \frac{1}{c^\ast}\,\Vert A \Vert_{\infty, 1}\,,
\end{align*}
whence 
\[
K_G^\R \leq \frac{1}{c^\ast} < \frac{\pi}{2\,\ln(1+\sqrt{2})}\,.
\]
In summary, given the - crucial - result \cite[Theorem 5.1]{BMMN2013},
our general framework can be applied here as well, leading to a {slightly 
modified} proof of \cite[Theorem 1.1]{BMMN2013}.
At this point, we would like to highlight {another, very recent and 
partially altered proof of \cite[Theorem 1.1]{BMMN2013} -- provided by 
Krivine again (cf. \cite[Theorem 1]{Kr2023})}. 
\end{example}
Despite the quite remarkable outcome of \autoref{thm:odd_real_bdd_CCP}, 
we should observe that its practical implementation seems to be quite difficult 
(at least without sufficiently large computer power). Primarily, as we already 
have seen, this is due to the following facts:
\begin{enumerate}
\item Either we have to know a closed form representation (or at least a ``close'' 
approximation of the Maclaurin series) of $h_{f,f}$, $h_{f,f}^{-1}$ and 
$\big(h_{f,f}^{-1}\big)_{\text{abs}}$ if $f \in L^\infty(\R^k)$ is given (such 
as is the case for $k=1$ and $f : = \text{sign} \in L^\infty(\R^1) \cap 
S_{L^2(\gamma_1)}$), or we have to check whether $h = h_{f,f}$ ``originates'' 
from some $f \in L^\infty(\R^k) \cap S_{L^2(\gamma_k)}$, if the functions 
$h$, $h^{-1}$ and $\big(h^{-1}\big)_{\text{abs}}$ are known to us.
\item However, already in the one-dimensional case (i.e., for $k=1$) that search 
requires rather complex calculation techniques, respectively some very helpful 
knowledge about Hermite polynomials. Moreover, if $k$ increases, we are confronted 
with a ``curse of combinatorial dimensionality'', since for any $\nu \in \N_0$ 
and $n \in \N$ it can be easily shown by induction on $k$ that the set 
$C(\nu, k) : = \{n \in \N_0^k : \vert n \vert = \nu\}$ which determines the 
structure of $h_{f,f}$ (cf. \autoref{thm:h_fg_real_case}) actually consists 
of $\binom{\nu+k-1}{k-1} = \frac{(\nu+k-1)!}{\nu !\,(k-1)!}$ elements. In 
particular, already $C(\nu, 2)$ consists of $\nu+1$ elements ($\nu \in \N_0$). 
For example, to determine $C(2, 11)$ explicitly, we would have to know all of 
its $66$ elements!
\end{enumerate}
We will recognise how deep actually we are confronted with a ``curse of 
combinatorial dimensionality'' if only a Maclaurin series representation 
of the function $h_{f,f}$ is given to us (cf. \sref{Section}{sec:RP_1}).
\chapter{The complex case: towards extending Haagerup's approach}
\section{Multivariate complex CCP functions and their relation to the real 
case}
\label{sec:Intro_to_the_fct_h_b_b}
Next, we are going to transfer the main results in the previous 
{chapter} from the real field $\R$ to the complex field $\C$. In order to 
achieve this, we have to implement non-trivial structural properties of the 
class of \textit{complex} Hermite polynomials (cf. Theorem 
\ref{thm:correlated_complex_Hermite_polynomials} and Theorem 
\ref{thm:odd_completely_real_analytic_functions_at_0_vs_complex_h_fg}).
The complex versions of \autoref{thm:h_fg_real_case} and 
\autoref{thm:h_ff_real_odd_case} also allow a generalisation of the Haagerup 
equality by transition from the $\Sigma_2(\zeta)$-correlated couple of 
two complex one-dimensional signum functions ${\text{sign}} : \C 
\longrightarrow \T$ to a $\Sigma_{2k}(\zeta)$-correlated couple of two 
(possibly different) arbitrary square-integrable functions $b, c : \C^k 
\longrightarrow \C$, where $k \in \N$ can be arbitrarily large 
(cf. \sref{Corollary}{cor:Haagerup_incl} and \sref{Example}{ex:Haagerup}). By 
definition (cf. \cite[Lemma 3.2. and Proof of Theorem 3.1]{H1987}), the 
complex ${\text{sign}}$-function is given as 
\[
{\text{sign}}(z) : = 
\begin{cases}\frac{z}{\vert z \vert}&\text{if } z \in \C^\ast\\ 
0&\text{if } z = 0
\end{cases}\,.
\]
We introduce the following helpful symbolic constructions and shortcuts. Fix 
$k,l \in \N$. Let $b : \C^k \longrightarrow \C$ and $c : \C^l \longrightarrow 
\C$ be two functions. Put
\[
b {}_k\!\otimes_l c : = (b\circ P_k)(c\circ Q_l),
\]
where 
\[
P_k : =
\begin{pmatrix}
    1 & 0 & \ldots & 0 & 0 & 0 & \ldots & 0 \\
    0 & 1 & \ldots & 0 & 0 & 0 & \ldots & 0 \\
		\vdots & \vdots & \ddots & \vdots & \vdots & \vdots & \ddots & \vdots\\
		0 & 0 & \ldots & 1 & 0 & 0 & \ldots & 0
  \end{pmatrix} = 
	(I_k \,\brokenvert\, 0) \in \M_{k, k+l}(\C)
\] 
and
\[Q_l : = 
\begin{pmatrix}
    0 & 0 & \ldots & 0 & 1 & 0 & \ldots & 0\\
    0 & 0 & \ldots & 0 & 0 & 1 & \ldots & 0\\
		\vdots & \vdots &  \ddots & \vdots & \vdots & \vdots & \ddots & \vdots\\
		0 & 0 & \ldots & 0 & 0 & 0 & \ldots & 1
  \end{pmatrix} = 
	(0 \,\brokenvert\, I_l) \in \M_{l,k+l}(\C),
\]
implying that
\begin{align}\label{eq:vec_splitting}
b\,{}_k\!\otimes_l c({\text{vec}}(z\, \brokenvert \,w)) = 
b(z)c(w) \text{ for all } (z, w) \in \C^k 
\times \C^l\,.
\end{align}
Given the construction of $b\,{}_k\!\otimes_l c$ we may unambiguously shorten 
it simply to $b \otimes c$ (and suppress the listing of the dimensions of the 
domains of definition of $b$, respectively $c$). Let $d : \C^k \longrightarrow 
\C$ be a given function. Recall the induced functions 
$r(d) : \R^{2k} \longrightarrow \R$ and $s(d) : \R^{2k} \longrightarrow \R$, 
defined as $r(d) : = \Re(d) \circ \frac{1}{\sqrt{2}}\,J_2^{-1}$ and 
$s(d) : = \Im(d) \circ \frac{1}{\sqrt{2}}\,J_2^{-1}$ (cf. 
\eqref{eq:real_valued_fcts_r_b_and_s_b}). In order to facilitate reading, we 
put $d_\alpha(z) : = d({\text{sign}}(\alpha)z)$, where $\alpha \in \C$ and 
$z \in \C^k$ (implying that $d_0 = d(0)$). Consequently, the construction of 
$b \otimes \overline{c}$ implies that
\begin{align*}
r(b \otimes \overline{c})(\vc{\vc{x_1}{x_2}}{\vc{y_1}{y_2}}) &=
r(b)(\vc{x_1}{y_1})r(c)(\vc{x_2}{y_2})\\
&+ s(b)(\vc{x_1}{y_1})s(c)(\vc{x_2}{y_2})
\end{align*}
and
\begin{align*}
s(b \otimes \overline{c})(\vc{\vc{x_1}{x_2}}{\vc{y_1}{y_2}}) &=
s(b)(\vc{x_1}{y_1})r(c)(\vc{x_2}{y_2})\\ 
&- r(b)(\vc{x_1}{y_1})s(c)(\vc{x_2}{y_2})
\end{align*}
for all $x_1, x_2, y_1, y_2 \in \R^k$. In other words: 
\begin{align}\label{eq:product_case_I}
r(b \otimes \overline{c}) = (r(b) \otimes r(c))\circ G + 
(s(b) \otimes s(c))\circ G
\end{align}
on $\R^{4k}$, and
\begin{align}\label{eq:product_case_II}
s(b \otimes \overline{c}) = (s(b) \otimes r(c))\circ G - 
(r(b) \otimes s(c))\circ G\,,
\end{align}
on $\R^{4k}$, where again $G = G^\top = G^{-1} \in O(4n)$ is the matrix, 
introduced in \eqref{eq:orthog_matrix}. Finally, $d$ is odd (respectively, even) 
if and only if both, $r(d)$ and $s(d) = r(-i\,d)$ are odd 
(respectively, even).
\begin{lemma}\label{lem:facts}
Let $k \in \N$, $\rho \in [-1,1]$, $\alpha \in \C$, $\theta \in \R$, 
$c, d : \C^k \longrightarrow \C$ and $\vc{\textbf{X}}{\textbf{Y}} \sim 
N_{4k}(0, \Sigma_{4k}(\rho))$. Then
\begin{enumerate}
\item $r(d_{\alpha}) = r(d) \circ R_2({\text{sign}}(\alpha)I_k)$ and
$r(d_{\theta}) = r(d) \circ {\text{sign}}(\theta)I_{2k}$. 
\item If $\alpha \not= 0$, then $r(d) \in L^2(\R^{2k}, \gamma_{2k})$ 
if and only if $r(d_{\alpha}) \in L^2(\R^{2k}, \gamma_{2k})$. In this case, the
norms coincide: $\Vert r(d) \Vert_{\gamma_{2k}} = 
\Vert r(d_{\alpha}) \Vert_{\gamma_{2k}}$.
\item If $\alpha \not= 0$, $h_{r(d), r(d)}(\rho) = 
\E[r(d)(\textbf{X})r(d)(\textbf{Y})] = 
\E[r(d_{\alpha})(\textbf{X})r(d_{\alpha})(\textbf{Y})] = 
h_{r(d_{\alpha}), r(d_{\alpha})}(\rho)$.
\item $h_{r(c_{r}), r(d)}(\vert r\vert) = h_{r(c),r(d)}(r) 
\text{ for all } r \in [-1,1]$.
\end{enumerate}
\end{lemma}
\noindent Fix $\vc{\textbf{Z}}{\textbf{W}} \sim {\C}N_{2k}(0, \Sigma_{2k}(\zeta))$, 
where $\zeta \in \overline{\D}\setminus\{0\}$ and $k \in \N$. Let 
$b,c : \C^k \longrightarrow \C$, such that $b \otimes \overline{c} \in L^1(\C^{2k}, 
\P_{\vc{\textbf{Z}}{\textbf{W}}})$. We put
\[
h^\C_{b,c}(\zeta) : = h^\C(b,c; \zeta) : = \E\big[b(\textbf{Z})
\overline{c(\textbf{W})}\big] = \overline{h^\C_{c,b}(\overline{\zeta})}. 
\]
As in the real case (see \eqref{eq:splitting_in_mean_quantum_correl_case}), the 
joint multivariate Gaussian splitting property 
\eqref{eq:jointly_Gaussian_on_the_unit_sphere} of inner products of vectors on 
the unit sphere and \sref{Lemma}{lem:correl_Gaussians_I}-\ref{(ii)} imply the 
important observation that for any separable $\C$-Hilbert space $H$, for any 
$m, n \in \N$, and for any $(u, v) \in S_H^m \times S_H^n$, we have
\begin{align}\label{eq:h_bb_Schur_image_of_complex_Gamma_u_v}
h^\C_{b,c}[\Gamma_H(u,v)] = \E\big[\overline{\textbf{R}_c}\,
{\textbf{S}_b}^\top\big] = \E[\Gamma_\C(\textbf{R}_c,\textbf{S}_b)],
\end{align}
where the $\C^m$-valued random vector $\textbf{R}_c$ and the $\C^n$-valued random 
vector $\textbf{S}_b$ are defined as $(\textbf{R}_c)_i : = c(\textbf{Z}_{u_i})$ 
and $(\textbf{S}_b)_j : = b(\textbf{Z}_{v_j})$, respectively ($(i,j) \in 
[m] \times [n]$). Fix $f, g \in L^2(\R^{k}, \gamma_{k})$ and put
\[
H_{f,g} : = h_{f,f} + h_{g,g}.
\]
\eqref{eq:h_fg_correl_representation} implies a concrete integral representation of the 
function $H_{f,g}{\big\vert}_{(-1,1)}$, which particularly plays an important role in the 
complex case (cf. 
\sref{Corollary}{cor:Haagerup_incl}):

\scalebox{0.90}{
\vbox{
\begin{align}\label{eq:H_fg_integral_representation}
H_{f,g}(\rho) = \frac{1}{(2\pi)^k(1 - \rho^2)^{k/2}}\,\int_{\R^{k}}\int_{\R^{k}}
\big\langle\binom{f(x)}{g(x)}, \binom{f(y)}{g(y)}\big\rangle_{\R_2^2}
\exp\big(-\frac{\Vert x \Vert^2 + \Vert y \Vert^2 - 2\rho\langle x,y \rangle}
{2(1-\rho^2)}\big) \textup{d}^{k}x\,\textup{d}^{k}y
\end{align}
}}

\noindent for all $\rho \in (-1,1)$. 
\autoref{thm:h_fg_real_case} further implies 
that $H_{f,g}(\rho) = \sum_{\nu = 0}^{\infty}a_\nu\,\rho^\nu$ for all $\rho \in [-1,1]$, where 
$a_\nu : = p_\nu(f,f) + p_\nu(g,g) \geq 0$ for all $\nu \in \N_0$. In particular, 
$H_{f,g}(1) = \Vert f \Vert^2_{\gamma_{2k}} + \Vert g \Vert^2_{\gamma_{2k}}$, implying that 
$H_{f,g}$ is real analytic on $(-1,1)$, continuous on $[-1,1]$ and absolutely monotonic on 
$[0,1]$. Moreover, $H_{f,g}$ is bounded, and 
\begin{align}\label{eq:H_f_g_is_bdd}
\vert H_{f,g}(\rho) \vert \leq \Vert f \Vert^2_{\gamma_{2k}} + 
\Vert g \Vert^2_{\gamma_{2k}} = H_{f,g}(1)
\end{align}
for all $\rho \in [-1,1]$. Hence, if $f \not=0$ or $g \not= 0$, it follows that 
\begin{align}\label{eq:rep_of_H_f_g}
H_{f,g} = H_{f,g}(1)\,\psi_{f,g},
\end{align}
where the function $\psi_{f,g} : = \frac{H_{f,g}}{H_{f,g}(1)} : [-1,1] \longrightarrow 
[-1,1]$ is CCP (due to \autoref{thm:CCP_via_Schoenberg}). In particular,
$H_{f,g}\Big\vert_{(-1,1)} \in W^\omega_+((-1, 1))$. Thus, if $b \in 
L^2(\C^k, \gamma_k^\C)$, then \eqref{eq:complex_Gaussian_L2} implies 
that $r(b) \in L^2(\R^{2k}, \gamma_{2k})$ and $s(b) \in L^2(\R^{2k},\gamma_{2k})$, 
and 
\begin{align}\label{eq:H_Re_Im_at_1}
H_{r(b),s(b)}(1) = \Vert r(b) \Vert^2_{\gamma_{2k}} + \Vert s(b) \Vert^2_{\gamma_{2k}} = 
\Vert b \Vert^2_{\gamma_k^\C}.
\end{align}
If - in addition - $b$ is odd and 
$\Vert b \Vert^2_{\gamma_k^\C} > 0$, then also $r(b)$ and $s(b)$ are odd functions 
(such as $H_{r(b),s(b)}$), satisfying $r(b) \not= 0$ or $s(b) \not= 0$. 
Hence, we may apply \autoref{thm:h_ff_real_odd_case} to the well-defined odd CCP function 
$\psi_{r(b),s(b)}$, and it follows that $H_{r(b),s(b)} : [-1,1] \longrightarrow 
[-\Vert b \Vert^2_{\gamma_k^\C}, \Vert b \Vert^2_{\gamma_k^\C}]$ is a strictly increasing 
homeomorphism, such as the inverse function $(H_{r(b),s(b)})^{-1} : 
[-\Vert b \Vert^2_{\gamma_k^\C}, \Vert b \Vert^2_{\gamma_k^\C}] \longrightarrow [-1,1]$ (due 
to \eqref{eq:H_Re_Im_at_1} and \eqref{eq:rep_of_H_f_g}). Similarly, if we - further - 
assume that $H_{r(b),s(b)}^{\prime}(0) > 0$, \autoref{thm:h_ff_real_odd_case}-\ref{(iv)} shows that 
also $(H_{r(b),s(b)})^{-1}$ is real analytic on $(-\Vert b \Vert^2_{\gamma_k^\C}, 
\Vert b \Vert^2_{\gamma_k^\C})$. So, we could apply \sref{Lemma}{lem:real_analytic_and_l1} 
to the real analytic function $(H_{r(b),s(b)})^{-1}
{\big\vert}_{(-\Vert b \Vert^2_{\gamma_k^\C}, \Vert b \Vert^2_{\gamma_k^\C} )}$ (if 
the assumptions are given) to check the existence of 
$(H_{r(b),s(b)}^{-1})_{\text{abs}} : [-\Vert b \Vert_{\gamma_k^\C}^2, 
\Vert b \Vert_{\gamma_k^\C}^2] \longrightarrow \R$ then (a crucial assumption 
in \autoref{thm:complex_inner_product_rounding}). Equipped with these facts, 
we arrive at the complex version of \autoref{thm:h_fg_real_case}:
\begin{theorem}\label{thm:h_b_c_complex_correl_case}
Let $k \in \N$, $\zeta \in \overline{\D}$, $\textbf{L} \sim 
{\C}N_{k}(0, I_k)$ and $\vc{\textbf{Z}}{\textbf{W}} 
\sim {\C}N_{2k}(0, \Sigma_{2k}(\zeta))$. Let $b \in 
L^2(\C^k, \gamma_k^\C)$ and $c \in L^2(\C^k, \gamma_k^\C)$. 
Then 
$b \otimes \overline{c} \in L^1(\C^k \times \C^k, \P_{\vc{\textbf{Z}}{\textbf{W}}})$. 
If $\zeta = 0$, then
\begin{align}\label{eq:CCP_complex_independence_case}
\begin{split}
\overline{h^\C_{c,b}(0)} &= h^\C_{b,c}(0) = \E[b(\textbf{L})]\,
\E[\overline{c}(\textbf{L})]\\
&= h_{r(b), r(c)}(0) + h_{s(b), s(c)}(0) + 
i\,(h_{s(b), r(c)}(0) - h_{r(b), s(c)}(0)).
\end{split}
\end{align}
If $\zeta \not= 0$, then
\begin{align}\label{eq:CCP_complex_correl_case}
\overline{h^\C_{c,b}(\overline{\zeta})} = h^\C_{b,c}(\zeta) = 
h_{r(b_\zeta), r(c)}(\vert \zeta \vert) + 
h_{s(b_\zeta), s(c)}(\vert \zeta \vert) + 
i\,(h_{s(b_\zeta), r(c)}(\vert \zeta \vert) - 
h_{r(b_\zeta), s(c)}(\vert \zeta \vert)).
\end{align}
In particular,
\begin{align}\label{eq:CCP_complex_correl_case_II}
0 \leq h^\C_{b_{\overline{\zeta}},b}(\zeta) = H_{r(b), s(b)}(\vert \zeta \vert)
\text{ and } h^\C_{b,c}(1) = \langle b, c \rangle_{\gamma_k^\C}.
\end{align}
$h^\C_{b,c} : \overline{\D} \longrightarrow \C$ is bounded and satisfies
\begin{align}\label{eq:boundedness_of_h_complex_case}
\vert h^\C_{b,c}(\zeta) \vert \leq \Vert b \Vert_{\gamma_k^\C}\,
\Vert c \Vert_{\gamma_k^\C} \text{ for all } \zeta \in \overline{\D}.
\end{align}
\end{theorem}
\noindent By taking into account that 
${\text{sign}}(\zeta)\,\vert \zeta \vert^{2\nu+1} = 
\zeta\cdot\zeta^{\nu}\cdot\overline{\zeta}^{\;\nu}$ for all $\zeta \in \C$ and 
$\nu \in \N_0$, \autoref{thm:h_b_c_complex_correl_case} also leads to a 
straightforward generalisation of the Haagerup function 
(cf. \cite[Proof of Theorem 3.1]{H1987} and \sref{Example}{ex:Haagerup}). 
To this end, we introduce a class of complex-valued functions which could be 
viewed as a transfer of the class of all odd real-valued functions to the 
complex field and contains the complex signum function 
$\text{sign} : \C \longrightarrow \C$ as element.
\begin{definition}\label{def:circ_symm_fcts}
Let $\F \in \{\R, \C\}$ and $k \in \N$. A function $b : \F^k \longrightarrow \F$ 
is circularly symmetric if 
\[
b(\alpha z) = \alpha b(z) \text{ for all } (\alpha, z) \in S_\F \times 
\F^k\,.
\]
The set of all circularly symmetric functions is denoted by $CS_k(S_\F)$.
\end{definition}
\noindent \sref{Definition}{def:circ_symm_fcts} obviously implies that 
$CS_k(S_\R) = CS_k(\{-1,1\})$ coincides with the set of all odd real 
functions from $\R^k$ to $\R$ and that $CS_k(S_\C) = CS_k(\T)$. 
Moreover, $(CS_k(S_\F), \circ, \text{id})$ is a monoid (i.e., a semigroup, 
with unit element), where the binary operation $\circ$ is given by the 
composition of functions.     
\begin{corollary}\label{cor:Haagerup_incl}
Let $k \in \N$ and $\zeta \in \overline{\D}$. Let $b, c \in 
L^2(\C^k, \gamma_k^\C)$. Suppose that $b \in CS_k(\T)$. Then
\begin{enumerate}
\item
\begin{align}\label{eq:h_bc_b_odd_circular}
h^\C_{b,c}(\zeta) = {\text{sign}}(\zeta)(h_{r(b), r(c)}(\vert \zeta \vert) + 
h_{s(b), s(c)}(\vert \zeta \vert) + i\,(h_{s(b), r(c)}(\vert \zeta \vert) - 
h_{r(b), s(c)}(\vert \zeta \vert)))\,.
\end{align}

\begin{align}\label{eq:Haagerup_style}
h^\C_{b,b}(\zeta) = {\text{sign}}(\zeta)H_{r(b), s(b)}(\vert\zeta\vert) =
\zeta\sum_{\nu = 0}^\infty (p_{2\nu+1}(r(b), r(b)) + 
p_{2\nu+1}(s(b), s(b)))\zeta^{\nu}\,\overline{\zeta}^{\;\nu}\,.
\end{align}
In particular, we have:
\begin{enumerate}
\item[\textup{(i-1)}] $h_{b,b}^\C\big\vert_{[-1,1]} = H_{r(b), s(b)} = h_{r(b), r(b)} +
h_{s(b), s(b)}$ and $h_{b,b}^\C\big\vert_{(-1,1)} \in W^\omega_+((-1, 1))$.
\item[\textup{(i-2)}] If $\vert \zeta \vert = 1$, then $h^\C_{b,b}(\zeta) = \zeta\,
\Vert b \Vert_{\gamma_k^\C}^2$. 
\item[\textup{(i-3)}] If $\zeta \in \D$ and $\vc{\textbf{X}}{\textbf{Y}} \sim 
N_{2k}(0, \Sigma_{2k}(\vert\zeta\vert))$, then

\scalebox{0.88}{
\vbox{
\begin{align}\label{eq:H_b_b_integral_representation}
\begin{split} 
h^\C_{b,b}(\zeta) &= 
\frac{{\text{sign}}(\zeta)}{(2\pi)^{2k}(1 - \vert\zeta\vert^2)^k}\,
\int_{\R^{2k}}\int_{\R^{2k}}\big\langle\binom{r(b)(x)}{s(b)(x)}, 
\binom{r(b)(y)}{s(b)(y)}\big\rangle_{\R_2^2}
\exp\big(-\frac{\Vert x \Vert^2 + \Vert y \Vert^2 - 2\vert\zeta\vert\,
\langle x,y \rangle}{2(1-\vert\zeta\vert^2)}\big)\,\textup{d}^{2k}x\,
\textup{d}^{2k}y\\
&= {\text{sign}}(\zeta)(\E[r(b)(\textbf{X})r(b)(\textbf{Y})] + 
\E[s(b)(\textbf{X})s(b)(\textbf{Y})]).
\end{split}
\end{align}
}}
\end{enumerate}
\item $h^\C_{b,b}$ is bounded, and
\begin{align}\label{eq:boundedness_of_h_bb_complex_case}
\vert h^\C_{b,b}(\zeta) \vert \leq \Vert b \Vert_{\gamma_k^\C}^2 
\text{ for all } \zeta \in \overline{\D}.
\end{align}
\item  If $b \not= 0$, then $\frac{1}{\Vert b \Vert^2_{\gamma_k^\C}}\,
H_{r(b), s(b)}$ as well as $\frac{1}{\Vert b \Vert^2_{\gamma_k^\C}}\,h^\C_{b,b}$ 
are CCP functions. Both, $H_{r(b), s(b)} : [-1,1] \longrightarrow 
[-\Vert b \Vert_{\gamma_k^\C}^2, \Vert b \Vert_{\gamma_k^\C}^2]$ and 
$h^\C_{b,b} : \overline{\D} \longrightarrow 
\Vert b \Vert_{\gamma_k^\C}^2\overline{\D}$ are circularly symmetric 
homeomorphisms. $H_{r(b), s(b)}$ is strictly increasing and satisfies
$H_{r(b), s(b)}((-1,0)) = (-\Vert b \Vert_{\gamma_k^\C}^2,0)$ and 
$H_{r(b), s(b)}((0,1)) = (0,\Vert b \Vert_{\gamma_k^\C}^2)$. Moreover,
$h^\C_{b,b}(\D) = \D$ and  
\begin{align}\label{eq:inverse_fct_in_the_complex_case}
(h^\C_{b,b})^{-1}(w) = {\text{sign}}(w)H_{r(b), s(b)}^{-1}(\vert w \vert) 
\text{ for all } w \in \Vert b \Vert_{\gamma_k^\C}^2\overline{\D}\,.
\end{align}
\item $h_{r(b),r(b)}^\prime(0) > 0$ or $h_{s(b),s(b)}^\prime(0) > 0$ 
if and only if $H_{r(b), s(b)}^{-1}{\vert}_{(-\Vert b \Vert_{\gamma_k^\C}^2, 
\Vert b \Vert_{\gamma_k^\C}^2)} = (H_{r(b), s(b)}{\big\vert}_{(-1,1)})^{-1}$
is real analytic on $(-\Vert b \Vert_{\gamma_k^\C}^2, \Vert b \Vert_{\gamma_k^\C}^2)$.
\end{enumerate}
\end{corollary}
\begin{remark}
Already the trivial example $b : = 1$, $c:=1$ and $\zeta : = i$ shows us that the additional 
assumption $b \in CS_k(\T)$ in \sref{Corollary}{cor:Haagerup_incl} cannot be dropped.
\end{remark}
\noindent Recall \eqref{eq:real_L_2_as_complex_L2_product_case}, including the 
construction of the function $g^\C$ therein. Due to 
\sref{Proposition}{prop:product_structure} and \sref{Corollary}{cor:Haagerup_incl} 
we obtain 
\begin{remark}\label{rem:real_embedded_in_complex}
Let $k \in \N$ and $f \in L^2(\R^k, \gamma_k)$.  
Then $r(f^\C) = f \otimes 1$ and $s(f^\C) = 0$. Moreover,
\[
h_{f,f} = h_{f \otimes 1, f \otimes 1} = H_{r(f^\C), s(f^\C)}\,.
\]
If in addition $f^\C \in CS_k(\T)$, then $h_{f,f} = h^\C_{f^\C, f^\C}{\big\vert}_
{[-1,1]}$\,.
\end{remark}
\noindent At this point, it is very useful to recall Theorem 
\ref{thm:real_hybrid_correlation_transform} and its proof, where we also 
implemented Abel's theorem on \textit{real} power series. From complex analysis 
it is well-known that in general Abel's theorem on power series in this form 
does not hold for $\F = \C$. Moreover, already the structure of the somewhat 
``simpler'' complex-valued odd functions $h_{b,c}$ (cf. 
\eqref{eq:h_bc_b_odd_circular} and \eqref{eq:Haagerup_style}) seemingly does not 
allow a transfer of \autoref{thm:real_hybrid_correlation_transform}-\ref{(ii)} 
and \autoref{thm:real_hybrid_correlation_transform}-\ref{(iii)} to the complex 
case (due to Theorem \ref{thm:complex_CCP_via_Christensen_and_Ressel}). However, 
since \cite[Lemma 1.1.13, statement 2]{GMS2015} also holds also for the complex 
field, the proof of \autoref{thm:real_hybrid_correlation_transform} can be easily 
adapted, as well as the proof of \sref{Corollary}{cor:important_special_abs_case}, 
so that at least the following versions of a complex hybrid correlation transform 
can be stated at once:
\begin{proposition}\label{prop:complex_hybrid_correlation_transform}
Let $m,n \in \N$, $0 < r < \infty$ and
\[
M : = \begin{pmatrix}
    A & S \\
    S^\ast & B
  \end{pmatrix} \in M_{m+n}(r\,\overline{\D})^+
\]
be positive semidefinite, where any entry of the matrices $A, S$ and $B$ is an 
element of $r\,\overline{\D}$. Let $f, g : (-r,r) \longrightarrow \R$ be two 
functions, such that $(-r,r) \ni x \mapsto f(x) = \sum_{\nu=0}^\infty a_\nu x^\nu 
\in W^\omega_+((-r,r))$ and $(-r,r) \ni x \mapsto g(x) = \sum_{\nu=0}^\infty b_\nu 
x^\nu \in W^\omega_+((-r,r))$. Let $q: r\D \longrightarrow \C$ be a holomorphic 
function, such that 
\begin{align}\label{eq:hol_T0_C_omega_inequality}
\vert c_\nu \vert \leq \sqrt{\vert a_\nu\vert\,\vert b_\nu\vert} \text{ for all } 
\nu \in \N_0\,,
\end{align}
where $c_\nu : = \frac{q^{(\nu)}(0)}{\nu !}$.
Then $q_r \in W^+(\D)$, where $q_r(\zeta) : = q(r \zeta)$ for all $\zeta \in \D$. 
Put $r\,\overline{\D} \ni z \mapsto \tilde{q}(z) : = \sum_{\nu = 0}^\infty c_\nu 
z^\nu$, $r\,\overline{\D} \ni z \mapsto f_{\text{abs}}(z) : = \sum_{\nu = 0}^\infty 
\vert a_\nu\vert z^\nu$ and $r\,\overline{\D} \ni z \mapsto g_{\text{abs}}(z) : = 
\sum_{\nu = 0}^\infty \vert b_\nu\vert z^\nu$. Then $\widetilde{q_r}(\frac{z}{r}) 
= \tilde{q}(z)$ for all $z \in r\,\overline{\D}$ and 
$\tilde{q}\big\vert_{r\,\D} = q$, where $\widetilde{q_r} \in A(\D)$ is the 
continuous extension of $q_r$. Moreover, the following properties hold:
\begin{enumerate}
\item 
\[
	\begin{pmatrix}
    f_{\text{abs}}[A] & \tilde{q}[S]\\
    \tilde{q}[S]^\ast & g_{\text{abs}}[B]
	\end{pmatrix}
\in M_{m+n}(r\,\overline{\D})^+
\]
is positive semidefinite.
\item If $0 < c_\ast \leq r$ is a root of $f_{\text{abs}} - 1$, then
\[
\tilde{q}[c_\ast\,\Gamma] \in \mathcal{Q}_{m,n}(\C) \text{ for all } 
\Gamma \in \mathcal{Q}_{m,n}(\C)\,.
\]
\end{enumerate}
\end{proposition}
\noindent Thus, an application of \sref{Lemma}{lem:real_analytic_and_l1} 
immediately leads us to
\begin{corollary}\label{cor:important_special_complex_abs_case}
Let $m, n \in \N, A \in \M_{m,n}(\C)$ and
\[
\Sigma =
\begin{pmatrix}
M & S \\
S^\top & N
\end{pmatrix} \in C(m+n; \C)
\]
be an arbitrary complex $(m+n) \times (m+n)$ correlation matrix (with block elements 
$M \in C(m; \C), N \in C(n; \C)$ and $S \in {\mathcal{Q}}_{m,n}(\C)$).
Let $r > 0$ and 
$0 \not=\psi \in W^\omega_+((-r, r))$. Put $r\,\overline{\D} \ni z 
\mapsto \tilde{\psi}(z) : = \sum_{\nu = 0}^\infty a_\nu z^\nu$, where $a_\nu : = 
\frac{\psi^{(\nu)}(0)}{\nu !}$. If $0 < c \leq r$, then
\[
\frac{1}{\psi_{\text{abs}}(c)}
\begin{pmatrix}
\psi_{\text{abs}}[c\,M] & \tilde{\psi}[c\,S] \\[0.5em]
\tilde{\psi}[c\,S^\ast] & \psi_{\text{abs}}[c\,N]
\end{pmatrix} \in C(m+n; \C)
\]
again is a correlation matrix with complex entries. In particular,
\begin{align}\label{eq:complex_entrywise_mapping_of_inner_products_to_inner_products}
\frac{1}{\psi_{\text{abs}}(c)}\widetilde{\psi}[c S] \in 
{\mathcal{Q}}_{m,n}(\C)\,\text{ for all }\,S \in {\mathcal{Q}}_{m,n}(\C)\,.
\end{align}
\end{corollary}
\section{On complex bivariate Hermite polynomials}
 Our next aim is to reveal in detail that it is possible to transfer 
the content of \autoref{thm:GT_generalising_Krivine} from the real case to 
the complex one, while maintaining our constructive proof (including the 
intended avoidance of the tensor product language). However, we cannot simply 
copy the proof of \autoref{thm:GT_generalising_Krivine}. Nevertheless, 
we are going to unfurl that in fact it is possible to transfer 
\eqref{eq:upper_bound_real_case} and \eqref{eq:real_GT_estimate_and_abs_fct} from 
the real case to the complex one. To this end, we are going to work with a particular 
case of the rich class of complex bivariate Hermite polynomials, first considered 
by K. It\^{o} while working with complex multiple Wiener integrals (cf. 
\cite{I1953}). Similarly to the real case, we need to verify a convenient 
correlation property of a random version of these polynomials, which to the 
best of our knowledge have not been published before (see Theorem 
\ref{thm:correlated_complex_Hermite_polynomials} below). A detailed introduction 
to complex Hermite polynomials (which would exceed the topic of this monograph 
by far) can be studied in \cite{G2013, I2015}. 
Firstly, we have to recall the following general construction:
\begin{definition}[\textbf{Complex Hermite polynomial}]
Let $m, n \in \N_0$ and $z,w \in \C$. The complex Hermite polynomial $H_{m,n} : \C^2 
\longrightarrow \C$ is defined as
\[
H_{m,n}(z,w) : = \frac{1}{\sqrt{m! n!}}\sum_{j=0}^{m \wedge n}(-1)^j j! \binom{m}{j}
\binom{n}{j}z^{m-j}\,w^{n-j}\,.
\] 
\end{definition}
\noindent Within the scope of our work, we need the particular case of 
It\^{o}'s complex Hermite polynomials 

\scalebox{0.95}{
\vbox{
\begin{align*}
\C \ni z &\mapsto H_{m,n} \circ \kappa(z) := H_{m,n}(z, \overline{z})\\
&= i^{m+n} \sum_{j=0}^m \sum_{k=0}^n i^{j+k} (-1)^{j+n} \sqrt{\binom{m}{j}\binom{n}{k}}
s(j,k)H_{j+k}(\sqrt{2}\Re(z))\,s(m-j, n-k)H_{m-j+n-k}(\sqrt{2}\Im(z)),
\end{align*}
}}
where $\kappa(z) : = \vc{z}{\overline{z}}$ and $s(\nu, \mu): = 
\sqrt{\tfrac{(\nu+\mu)!}{\nu\,!\mu !}}$ for all $\nu, \mu \in \N_0$, including 
the following statements, which we give without proof (cf. \cite{I1953, I2015}).
\begin{theorem}\label{thm:exp_gen_fct_of_complex_Hermite_polynomials}
Let $m, n \in \N_0$ and $z,w \in \C$. Then
\begin{enumerate}
\item $\{H_{m,n} \circ \kappa : m, n \in \N_0\}$ is an orthonormal basis in the complex 
Hilbert space $L^2(\gamma_1^\C)$.
\item The exponential generating function of $\{H_{m,n} \circ \kappa : m, n \in \N_0\}$ 
is given as
\[
\sum_{m, n=0}^\infty H_{m,n}(z, \overline{z}) \frac{u^m}{\sqrt{m}!}\frac{v^n}{\sqrt{n}!} =
\exp(uz + v\overline{z}-uv) \text{ for all } u, v \in \C.
\] 
\end{enumerate}
\end{theorem}
\noindent \sref{Lemma}{lem:gen_of_CF_and_MGF}-\ref{(ii)} allows us to transfer 
(the special case $k=1$ of) \sref{Corollary}{cor:correlated_k_dim_Hermite_polynomials} 
to the complex case. More precisely, we have
\begin{theorem}\label{thm:correlated_complex_Hermite_polynomials}
Let $m,n,\nu, \mu \in \N_0$ and $\zeta \in \overline{\D}$. If 
$\vc{Z}{W} \sim \C N_{2}(0, \Sigma_{2}(\zeta))$, then
\[
\E\big[H_{m,n}(Z, \overline{Z})\overline{H_{\nu,\mu}(W, \overline{W})}\,\big] = 
\delta_{m,\nu}\,\delta_{n,\mu}\,\zeta^{m}\,\overline{\zeta}^{\,n} =
\delta_{(m,n), (\nu,\mu)}\,\zeta^{m}\overline{\zeta}^{\,n}.
\]
\end{theorem}
\section{Upper bounds of $K_G^\C$ and inversion of complex CCP functions}
 Equipped with the complex Hermite polynomials and Theorem 
\ref{thm:correlated_complex_Hermite_polynomials}, it is possible to transfer 
\autoref{thm:char_of_completely_real_analytic_functions_at_0_as_h_fg} from 
the real field $\R$ to the complex field $\C$; at least in the odd case. 
We ``just'' have to construct the mappings $\alpha^{\psi, c}_1 \in 
S_{L^2(\gamma_1^\C)}$ and $\beta^{\psi, c}_1 \in S_{L^2(\gamma_1^\C)}$ properly. 
\begin{theorem}\label{thm:odd_completely_real_analytic_functions_at_0_vs_complex_h_fg}
Let $k \in \N$ and $0 < c \leq 1$. Let $0 \not= \psi \in W^\omega_+((-1, 1))$ be odd.
Then there exist $\alpha \equiv \alpha_{\psi, c}, 
\beta \equiv \beta_{\psi, c} \in S_{L^2(\gamma_1^\C)}$, which satisfy 
the following properties:
\begin{enumerate}
\item $\E[\alpha(Z)] 
= \E[\beta(Z)] = 0$ for all $Z \sim {\C}N_1(0,1)$. 
\item If 
$c\,\zeta \in \D$, then 
\[
{\text{sign}}(\zeta)\psi(c\vert\zeta\vert) = 
\psi_{\text{abs}}(c)\,h_{\alpha, \beta}^\C(\zeta)
\]
and
\[
{\text{sign}}(\zeta)\psi_{\text{abs}}(c\vert\zeta\vert) = 
h_{\alpha, \alpha}^\C(\zeta)\psi_{\text{abs}}(c) = 
h_{\beta, \beta}^\C(\zeta)\psi_{\text{abs}}(c). 
\]
In particular,
\begin{align}\label{eq:psi_c_vs_psi_abs_c_complex_case}
\psi(c) = \psi_{\text{abs}}(c)\,\langle \alpha, \beta\rangle_{\gamma_1^\C}\,.
\end{align}
\item If $c \not= 1$ and $H$ is an arbitrary $\C$-Hilbert 
space, then
\[
{\text{sign}}(\langle u, v \rangle_H)\psi(c\vert\langle u, v\rangle_H\vert) = 
\psi_{\text{abs}}(c)\,h_{\alpha, \beta}^\C(\langle u, v\rangle_H) 
\]
and
\[
{\text{sign}}(\langle u, v \rangle_H)
\psi_{\text{abs}}(c\,\vert\langle u, v\rangle_H\vert) = \psi_{\text{abs}}(c)\,
h_{\alpha, \alpha}^\C(\langle u, v\rangle_H) = \psi_{\text{abs}}(c)\,
h_{\beta, \beta}^\C(\langle u, v\rangle_H)
\]
for all $u, v \in S_H$.
\end{enumerate}
\end{theorem}
\noindent Recall \sref{Corollary}{cor:Haagerup_incl} including the structure of 
$h^\C_{b,b}$, strongly built on the odd homeomorphic real CCP function 
$\frac{1}{\Vert b \Vert_{\gamma_k^\C}^2}\,H_{r(b), s(b)} : [-1,1] \longrightarrow 
[-1, 1]$. Thus, if we link \eqref{eq:Haagerup_style}, 
\eqref{eq:inverse_fct_in_the_complex_case} and Theorem 
\ref{thm:odd_completely_real_analytic_functions_at_0_vs_complex_h_fg} 
(where the latter is applied to the function $\psi : = 
\big(\frac{1}{\Vert b \Vert_{\gamma_k^\C}^2}\,H_{r(b),s(b)}\big)^{-1} = 
H_{r(b),s(b)}^{-1}(\Vert b \Vert_{\gamma_k^\C}^2\,\cdot)$), we  
are able to prove  
\begin{theorem}[\textbf{Complex inner product rounding}]
\label{thm:complex_inner_product_rounding}
Let $k \in \N$ and $b \in S_{L^2(\gamma_k^\C)}$. If $b \in CS_k(\T)$ and 
$(H_{r(b),s(b)})^{-1}{\big\vert}_{(-1,1)} \in W^\omega_+((-1,1))$, then 
$\big(H_{r(b),s(b)}^{-1}\big)_{\text{abs}}(1) > 1$, and there exist $\alpha_b 
\in L^2(\gamma_1^\C)$ and $\beta_b \in L^2(\gamma_1^\C)$, satisfying $0 < 
\langle\alpha_b, \beta_b\rangle_{\gamma_1^\C} < 1$, such that for all 
$\C$-Hilbert spaces $H$ and $u, v \in S_H$ the following properties are satisfied:
\begin{enumerate}
\item
\begin{align}\label{eq:complex_inner_product_rep}
\langle u, v \rangle_H = 
\frac{1}{c^\ast} h_{b,b}^\C(\zeta_{u,v}(b)),  
\end{align}
where $0 < c^\ast \equiv c^\ast(b) : = H_{r(b), s(b)}^{\text{hyp}}(1) < 1$ and 
$\zeta_{u,v}(b) : = h_{\alpha_b, \beta_b}^\C(\langle u, v\rangle_H) \in \D$.
\item 
\[
c^\ast = H_{r(b),s(b)}(\langle \alpha_b, \beta_b \rangle_{\gamma_1^\C}) = 
h_{b,b}^\C(\langle \alpha_b, \beta_b \rangle_{\gamma_1^\C}).
\]
\item
If $\vc{\textbf{X}}{\textbf{Y}} \sim 
N_{2k}(0, \Sigma_{2k}(\vert\zeta_{u,v}(b)\vert))$, then

\scalebox{0.87}{
\vbox{
\begin{align*}
c^\ast\,\langle u, v \rangle_H &= 
\frac{{\text{sign}}(\zeta_{u,v})}{(2\pi)^{2k}(1 - \vert\zeta_{u,v}\vert^2)^k}\,
\int_{\R^{2k}}\int_{\R^{2k}}\big\langle\binom{r(b)(x)}{s(b)(x)}, 
\binom{r(b)(y)}{s(b)(y)}\big\rangle_{\R_2^2}
\exp\big(-\frac{\Vert x \Vert^2 + \Vert y \Vert^2 - 2\vert\zeta_{u,v}(b)\vert\,
\langle x,y \rangle}{2(1-\vert\zeta_{u,v}(b)\vert^2)}\big) \textup{d}^{2k}x\,
\textup{d}^{2k}y\\[0.5cm]
&= {\text{sign}}(\zeta_{u,v})
\big(\E[r(b)(\textbf{X})r(b)(\textbf{Y})] + 
\E[s(b)(\textbf{X})s(b)(\textbf{Y})]\big),
\end{align*}
}}
\item
If $m, n \in \N$ and $(z, w) \in S_H^m\times S_H^n$, then there exist $m+n$ 
$\C^k$-valued random vectors $\textbf{Z}_1, \ldots, \textbf{Z}_m, \textbf{W}_1, 
\ldots, \textbf{W}_n$, such that $\vc{\textbf{Z}_i}{\textbf{W}_j} \sim 
{\C}N_{2k}(0, \Sigma_{2k}(\zeta_{ij}(b)))$ for all $(i,j) \in 
[m] \times [n]$, and 
\begin{align}\label{eq:rep_of_c_times_QC_complex_case}
\Gamma_H(z,w) = 
\frac{1}{c^\ast}\,\E[\overline{{\textbf{P}}_b}{\textbf{Q}}_b^\top]\,, 
\end{align}
where $({\textbf{P}}_b)_i : = b(\textbf{Z}_i)$, 
$({\textbf{Q}}_b)_j : = b(\textbf{W}_j)$ and 
$\zeta_{ij}(b) : = h_{\alpha_b, \beta_b}^\C(\langle z_i, w_j\rangle_H) \in \D$, 
$(i,j) \in [m] \times [n]$.
\end{enumerate}
\end{theorem}
\noindent Our next result shows that in fact also \autoref{thm:odd_real_bdd_CCP}, 
respectively \sref{Corollary}{cor:odd_real_bdd_CCP_general_version} can be 
transferred to the complex case:
\begin{theorem}\label{thm:odd_complex_bdd_CCP_general_version}
Let $k, m, n \in \N$ and $b, c \in L^\infty(\C^k)$. Then $0 \leq r \equiv r_k(b) : = 
\Vert b \Vert_{\gamma_k^\C}^2 < \infty$ and $\Vert b \Vert_{\infty} \geq 
\sqrt{r}$. Moreover, the following statements hold:
\begin{enumerate}
\item
\[
\vert{\text{tr}}(A^\ast h^\C_{b,c}[S])\vert \leq \Vert b \Vert_\infty\,
\Vert c \Vert_\infty\,\Vert A \Vert^\C_{\infty, 1} \text{ for all } 
A \in \M_{m,n}(\C) \text{ and } S \in {\mathcal{Q}}_{m,n}(\C). 
\]
\item Assume that $b \in CS_k(\T)\setminus\{0\}$ and 
$H_{r(b),s(b)}^{-1}{\big\vert}_{(-r, r)} \in W^\omega_+((-r, r))$. Put
$c^\ast \equiv c^\ast(b) : = H_{\frac{r(b)}{\sqrt{r}}, 
\frac{s(b)}{\sqrt{r}}}^\text{hyp}(1)$. Then $c^\ast \in (0,1)$, and the 
following properties are satisfied:
\begin{enumerate}
\item[\textup{(ii-1)}] 
\begin{align}\label{eq:upper_bound_complex_general_bdd_case}
r\,K_G^\C \leq \frac{\Vert b \Vert_{\infty}^2}{c^\ast}\,.
\end{align}
{In particular, if $\text{sign}\big((H_{\frac{r(b)}{\sqrt{r}}, 
\frac{s(b)}{\sqrt{r}}}^{-1}{\vert}_{(-1,1)})^
{(2n+1)}(0)\big) = (-1)^n$ for all $n \in \N_0$, then
\[
r\,K_G^\C \leq i\,\frac{\Vert b \Vert_{\infty}^2}{F_b(i)}\,,
\]
where $F_b : = $ \rlap{\raisebox{-0.3ex}{$\widesim{\phantom{\psi{\vert}_{(-1,1)}}}$}}
$H_{\frac{r(b)}{\sqrt{r}}, \frac{s(b)}{\sqrt{r}}}{\vert}_{(-1,1)}$ $: 
\overline{\D} \longrightarrow \overline{\D}$ is defined as in 
\sref{Lemma}{lem:real_analytic_and_l1}.
}
\item[\textup{(ii-2)}] 
Let $1 \leq \kappa_\ast < K_G^\C$. If $\Vert b \Vert_\infty = \sqrt{r}$, then 
$0 < H_{\frac{r(b)}{\sqrt{r}},\frac{s(b)}{\sqrt{r}}}^\text{hyp}(\kappa_\ast) < 
1$ and there is exactly one number number $\gamma^\ast \equiv \gamma^\ast(b) \in 
(H_{\frac{r(b)}{\sqrt{r}},\frac{s(b)}{\sqrt{r}}}^\text{hyp}(\kappa_\ast) ,1]$, 
such that  
\begin{align}\label{eq:upper_bound_complex_odd_bdd_coinciding_norms_case}
K_G^\C = \big(H_{r(b), s(b)}^{-1}{\big\vert}_{(-r,r)}\big)_{\text{abs}}(r\,
\gamma^\ast) \leq \min\big\{\frac{1}{c^\ast},\, \big(H_{r(b),s(b)}^{-1}
{\big\vert}_{(-r,r)}\big)_{\text{abs}}(r)\big\}\,.
\end{align}
\end{enumerate} 
\end{enumerate}
\end{theorem}
\noindent Recall again from \autoref{thm:complex_CCP_via_Christensen_and_Ressel} 
that in general complex CCP functions need not be analytic. Due to this fact and 
the structure of the complex-valued (circularly symmetric) functions $h_{b,b}$, 
which is built on absolute values and signs of complex numbers (cf. 
\sref{Corollary}{cor:Haagerup_incl}, including \eqref{eq:Haagerup_style}), it 
seems that \sref{Proposition}{prop:thm_odd_real_bdd_CCP_without_CRA_condition} 
cannot be easily transferred to the complex case (if at all).
\begin{example}[\textbf{Haagerup function}]\label{ex:Haagerup}
Both, \sref{Corollary}{cor:Haagerup_incl} and \sref{Proposition}{prop:2F1_rep} 
enable us to recover quickly a direct power series representation of Haagerup's 
CCP function $h_{b,b} : \D \longrightarrow \D$, where $b : = {\text{sign}}$. 
(Retranslated into the terminology of Haagerup, $h_{b,b} = \Phi$ and 
$H_{r(b), s(b)} = \varphi$ (cf. \cite[Lemma 3.5 and Theorem 3.1]{H1987})). 
To this end, let $0 \not= \zeta \in \overline{\D}$ and $\rho : = \vert \zeta \vert$. 
As shown on \cite[p. 200]{H1987}, the function $h_{b,b}$ can then be written in 
terms of ${\text{sign}}(\zeta)$ and the two complete elliptic integrals $E(\rho)$ 
and $K(\rho)$. However, in our approach (by which Haagerup's CCP function is 
obtained as a special case), we don't have to work with elliptic integration. 
Firstly note that 
\[
r(b)(x_1, x_2) = r(b)(x) = \frac{x_1}{\Vert x \Vert_2} \text{ and } s(b)(x) = r(b)(x_2, x_1)
= \frac{x_2}{\Vert x \Vert_2}
\]  
for all $x = \vc{x_1}{x_2} \in \R^2\setminus{\{0\}}$. Thus, if $\tau \in [-1,1]$ and
$\vc{\textbf{X}}{\textbf{Y}} \sim N_{2k}(0, \Sigma_{2k}(\tau))$, $k :=2$ are given, then
\begin{align*}
H_{r(b), s(b)}(\tau) &= h_{r(b), r(b)}(\tau) + h_{s(b), s(b)}(\tau)\\ 
&= \E[r(b)(\textbf{X})r(b)(\textbf{Y})] + 
\E[s(b)(\textbf{X})s(b)(\textbf{Y})]\\ 
&= \E[r(b)(X_1, X_2)r(b)(Y_1,Y_2)] + 
\E[r(b)(X_2, X_1)r(b)(Y_2,Y_1)]\\ 
&= \E\big[\frac{X_1\,Y_1}{\Vert\textbf{X}\Vert_2\,\Vert\textbf{Y}\Vert_2}\big] + 
\E\big[\frac{X_2\,Y_2}{\Vert\textbf{X}\Vert_2\,\Vert\textbf{Y}\Vert_2}\big]\\
&= \E\big[\big\langle\frac{\textbf{X}}{\Vert\textbf{Y}\Vert_2},
\frac{\textbf{Y}}{\Vert\textbf{X}\Vert_2}\big\rangle\big]\\
&= \frac{\pi}{4}\,\tau\,{}_2 F_1\big(\frac{1}{2}, \frac{1}{2}, 2; \tau^2 \big) = 
h_{f_2,f_2}(\tau),
\end{align*}
whereby the equalities in the last line follow from \sref{Proposition}{prop:2F1_rep} 
(applied to $k=2$). An application of \sref{Corollary}{cor:Haagerup_incl} 
(or the original proof of \cite[Theorem 3.1]{H1987}) therefore implies
\[
{\text{sign}}(\overline{\zeta})h_{b,b}(\zeta) = H_{r(b), s(b)}(\rho) = 
\frac{\pi}{4}\,\rho\,{}_2 F_1\big(\frac{1}{2}, \frac{1}{2}, 2; \rho^2 \big) = 
h_{f_2,f_2}(\rho).
\] 
Hence,
\[
h_{b,b}(\zeta) = {\text{sign}}(\zeta)\,\E\big[\big\langle
\frac{\textbf{X}}{\Vert\textbf{X}\Vert_2},\frac{\textbf{Y}}{\Vert\textbf{X}\Vert_2}
\big\rangle\big] =
\frac{\pi}{4}\,\zeta\,{}_2 F_1\big(\frac{1}{2}, \frac{1}{2}, 2; 
\vert\zeta\vert^2 \big) = {\text{sign}}(\zeta)\,h_{f_2,f_2}(\vert\zeta\vert)\,.
\]
A highly non-trivial part in \cite{H1987} consists of a multi-page proof of 
the fact that $h_{f_2,f_2}^{-1} = H_{r(b), s(b)}^{-1} \in W^\omega_+((-1, 1))$, 
implying that $\big(h_{f_2,f_2}^{-1}\big)_{\text{abs}}$ is well-defined 
(cf. \cite[Lemma 2.6]{H1987}). The Maclaurin series of $\big(h_{f_2,f_2}^{-1}\big)_
{\text{abs}}$ can e.g. be approximated by the Taylor polynomial of degree 7 as:
\begin{align}\label{eq:inversion_of_Haagerup_s_real_fct_part}
\big(h_{f_2,f_2}^{-1}\big)_{\text{abs}}(s) = \frac{4}{\pi}\,s + 
\frac{8}{\pi^3}\,s^3 + 0\cdot s^5 + \frac{16}{\pi^7}\,s^7 + o(\vert s\vert^7).
\end{align}
Hence (cf. \sref{Proposition}{prop:2F1_rep} and 
\eqref{eq:upper_bound_complex_odd_bdd_coinciding_norms_case}), 
\begin{align}\label{eq:abs_CCP_inverse_in_case_of_f_2}
\frac{K_G^\R}{\sqrt{2}} 
\stackrel{\eqref{eq:known_upper_bounds_of_the_complex_GT_constant}}{\leq} 
K_G^\C 
&\stackrel{\eqref{eq:upper_bound_complex_odd_bdd_coinciding_norms_case}}{\leq}
\min\big\{\frac{1}{h_{f_2,f_2}^{\text{hyp}}(1)}, \big(h_{f_2,f_2}^{-1}\big)_
{\text{abs}}(1)\big\}\\[0.5em]
&\hspace{3mm}\leq \big(h_{f_2,f_2}^{-1}\big)_{\text{abs}}(1) = \frac{4}{\pi} + 
\frac{8}{\pi^3} + \frac{16}{\pi^7} + o(1) \approx 1,53655 + o(1).
\end{align}
These facts follow from \cite{H1987}, respectively 
\eqref{eq:k_dim_elliptic_integral_and_h_ff} and 
\eqref{eq:First_7_Taylor_coeff_of_the_inverse}, where the latter is applied to
\[
\alpha_\nu : =
\begin{cases} 0 &\text{if } \nu \text{ is even}\\  
 \frac{\pi}{2}\frac{((\nu-2)!!)^2}{((\nu-1)!!)^2\,(\nu+1)} &\text{if } 
\nu \text{ is odd}
\end{cases}
\text{ and } \alpha_\nu^\times : = \frac{\alpha_\nu}{\alpha_1} = 
\frac{\alpha_\nu}{c_2^2} = \frac{4}{\pi}\,\alpha_\nu\,.
\]
Already a numerical calculation of the single root $0 < c^\ast < \frac{\pi}{4}$ 
of the polynomial $s \mapsto \frac{4}{\pi}\,s + \frac{8}{\pi^3}\,s^3 + 
\frac{16}{\pi^7}\,s^7 - 1$ leads to the number $\frac{1}{c^\ast} \approx 1.40449$. 
The latter outcome should now be compared with 
the result of Haagerup in \cite{H1987}\,.
\end{example}
\chapter{A summary scheme of the main result}
\label{chp:Summary}
 To highlight and summarise our approach, it completely suffices to list 
in detail the single steps and assumptions in the form of a ``flowchart'', possibly 
leading to a computer-aided approach regarding the implementation of an 
approximation to the lowest upper bound of the Grothendieck constant $K_G^\F$ 
as a next step. \textit{Very likely, supercomputers are required to perform 
these approximations. That (technical) implementation would go however far beyond 
the scope of our groundwork; especially since we have no access to equipment of 
this type.} 
\begin{center}
\setlength{\fboxsep}{25pt}
\fbox{\parbox{0.99\columnwidth}{
\begin{center}
Fix $\F \in \{\R, \C\}$ and $k \in \N$.
\end{center}
\begin{enumerate}
\item[(SIGN)] Choose a function 
$0 \not=b : \F^k \longrightarrow \F$, which satisfies the following conditions: 
\begin{enumerate}
\item $b$ is circularly symmetric 
(\textit{cf. \sref{Definition}{def:circ_symm_fcts}}\,); 
\item $b \in L^\infty(\F^k)$. 
\end{enumerate}
\item[(CCP)] Consider the function $b^\circ_\F : \C^k \longrightarrow \C$, 
defined as 
\[
b^\circ \equiv b^\circ_\F: = \begin{cases} \big(\frac{b}
{\Vert b \Vert_{\gamma_{k}}}\!\big)^\C 
&\text{if } \F = \R\\
\frac{b}{\Vert b \Vert_{\gamma_k^\C}} &\text{if } \F = \C\,.
\end{cases}
\] 
Then $\Vert b^\circ\Vert_{\gamma_k^\C} = 1$ and $1 \leq \Vert b^\circ\Vert_\infty = 
\frac{\Vert b\Vert_\infty^\F}{\Vert b\Vert_{\gamma_k^\F}}$ (\textit{cf. 
\eqref{eq:real_L_2_as_complex_L2_product_case}}). 
Construct its allocated homeomorphic real CCP function $H_{r(b^\circ), 
s(b^\circ)} = h_{r(b^\circ), r(b^\circ)} + h_{s(b^\circ), s(b^\circ)}$ 
(\textit{cf. \sref{Corollary}{cor:Haagerup_incl}-\ref{(iii)} and 
\sref{Remark}{rem:real_embedded_in_complex}}\,).
\item[(CRA)] Assume that ${H_{r(b^\circ), s(b^\circ)}^{-1}
{\big\vert}_{(-1, 1)}} \in W^\omega_+((-1, 1))$ (\textit{cf. 
\sref{Definition}{def:CRA}}\,).
\end{enumerate}
}}
\end{center}
An application of \sref{Lemma}{lem:hyperbolic_CCP_transform}, 
\sref{Corollary}{cor:odd_real_bdd_CCP_general_version} and Theorem 
\ref{thm:odd_complex_bdd_CCP_general_version} consequently leads to the 
following result which holds for both, $\R$ and $\C$ simultaneously:
\begin{center}
\setlength{\fboxsep}{25pt}
\fbox{\parbox{0.99\columnwidth}{
Assume that (SIGN), (CCP) and (CRA) are satisfied. Then 
{$\big(H_{r(b^\circ), s(b^\circ)}^{-1}{\big\vert}_{(-1, 1)}\big)_
{\text{abs}}(1) > 1$}. Put $c^\ast_\F : = H_{r(b^\circ), s(b^\circ)}^\text{hyp}(1)$ 
{(\textit{cf. \sref{Lemma}{lem:hyperbolic_CCP_transform}})}. Then 
$0 < c^\ast_\F < 1$, and the following statements hold:
\begin{enumerate}
\item  
\begin{align*}
K_G^\F \leq \frac{1}{c^\ast_\F}\,\Vert b^\circ\Vert_\infty^2\,.
\end{align*}
{In particular, if $\text{sign}\big((H_{r(b^\circ), s(b^\circ)}^{-1}
{\big\vert}_{(-1, 1)})^{(2n+1)}(0)\big) = (-1)^n$ for all $n \in \N_0$, then
\[
K_G^\F \leq i\,\frac{\Vert b^\circ\Vert_\infty^2}{\widetilde{\psi_b}(i)}\,,
\]
where $\psi_b : = H_{r(b^\circ), s(b^\circ)}{\big\vert}_{(-1, 1)}$ and 
$\widetilde{\psi_b} : \overline{\D} \longrightarrow 
\overline{\D}$
is defined as in \sref{Lemma}{lem:real_analytic_and_l1}.}
\item Let $1 \leq \kappa_\ast < K_G^\F$. If $\Vert b \Vert_\infty^\F = 
\Vert b \Vert_{\gamma_k^\F}$, then $0 < H_{r(b^\circ), s(b^\circ)}^\text{hyp}
(\kappa_\ast) < 1$ and there is exactly one number $\gamma^\ast_\F \equiv 
\gamma^\ast_\F(b) \in (H_{r(b^\circ), s(b^\circ)}^\text{hyp}(\kappa_\ast) ,1]$, 
such that  
\begin{align*}
K_G^\F = \big(H_{r(b^\circ), s(b^\circ)}^{-1}{\big\vert}_{(-1,1)}\big)_{\text{abs}}
(\gamma^\ast_\F) \leq \min\Big\{\frac{1}{c^\ast_\F},\, \big(H_{r(b^\circ), s(b^\circ)}^{-1}
{\big\vert}_{(-1,1)}\big)_{\text{abs}}(1)\Big\}\,.
\end{align*} 
\end{enumerate}
}}
\end{center}
 Again, we recognise that Maclaurin series representation (or at least 
its approximation by the Taylor polynomial of a given degree) and 
Maclaurin series inversion of CCP functions play the key role regarding the 
search for the lowest upper bound of the Grothendieck constant $K_G^\F$. 
Unfortunately, a closed form representation of the coefficients of the inverse of 
a Taylor series runs against a well-known combinatorial complexity issue (due to 
the presence of ordinary partial Bell polynomials as building blocks of these 
coefficients - cf. \sref{Section}{sec:RP_1} below for details), which in 
general does not allow a closed form representation of these coefficients. The 
inverse of the real function factor of the Haagerup function is one such 
example. It is given by
\[
[-1,1] \ni \tau \mapsto H_{r(\text{sign}), s(\text{sign})}(\tau) = 
h_{f_2, f_2}(\tau) = \frac{\pi}{4}\,\tau\,{}_2 F_1\big(\frac{1}{2}, 
\frac{1}{2}, 2; \tau^2\big)
\] 
in the complex case (cf. \cite{H1987}, Remark on page 216 and 
\sref{Example}{ex:Haagerup}), as opposed to the Grothendieck function 
\[
[-1,1] \ni \rho \mapsto H_{r(\text{sign}^\C), s(\text{sign}^\C)}(\rho) = 
h_{\text{sign}, \text{sign}}(\rho) = \frac{2}{\pi}\arcsin(\rho)
\] 
in the real case.
\chapter{Concluding remarks and open problems}\label{chp:future_research_topics}
Not very surprisingly, the long-standing, intensive and technically 
quite demanding attempts to compute the - still not available - value of the real 
and complex Grothendieck constants (an open problem since 1953) leads to further 
projects and open problems, such as the following ones; addressed in 
particular to researchers who also wish to get a better understanding of the 
reasons underlying these topics.   
\section{Open problem 1: Grothendieck constant versus Taylor series 
inversion}\label{sec:RP_1}
 Only between 2011 and 2013 it was shown that $K_G^\R$ is strictly smaller than 
Krivine's upper bound, stating that $K_G^\R < \frac{\pi}{2\ln(1+\sqrt{2})}$ 
(cf. \cite{BMMN2013} and \sref{Example}{ex:Naor_et_al_s_result}). Consequently, 
in the real case $\text{sign}$ is not the ``optimal'' function to choose 
(answering a question of H. K\"{o}nig to the negative - cf. \cite{K2001}). So, 
if we wish to reduce the value of the upper bound of the real Grothendieck 
constant we have to look for functions $b : \R^k \longrightarrow \R$ 
which are different from $\text{sign} : \R \longrightarrow \{-1,1\}$. 
However, these functions are required to satisfy any of the conditions in the 
workflow, listed in \sref{Chapter}{chp:Summary}. In particular, we have to 
look for both, the Fourier-Hermite coefficients of the Taylor series 
(respectively the approximating Taylor polynomial) of 
$H_{r(b),s(b)} = h_{r(b), r(b)} + h_{s(b), s(b)}$ and the coefficients of the 
Taylor series of both, the \textit{inverse} function 
$H_{r(b),s(b)}^{-1}$ and $\big(H_{r(b), s(b)}^{-1}\big)_{\text{abs}}$. 
It is well-known that the latter task increases rapidly in computational 
complexity if we want to calculate such Taylor coefficients of a higher degree, 
leading to the involvement of highly non-trivial combinatorial aspects, 
concretised by the use of partitions of positive integers and partial 
exponential Bell polynomials as part of the Taylor coefficients of the inverse 
Taylor series (a thorough introduction to this framework including the related 
Lagrange-B\"urmann inversion formula is given in \cite{C1974, K1951}).

 To reveal the origin of these difficulties let us focus 
on the one-dimensional real case. Let $b \in L^2(\R, \gamma_1)$ be given. Assume 
that $\alpha_0 : = h_{b,b}(0) = 0$ (which is the case if $b$ were odd).
First recall from \eqref{eq:series_rep_of_real_h_fg} that 
\[
h_{b,b}(\rho) = \sum_{n=1}^\infty \langle b, H_n \rangle_{\gamma_1}^2\,\rho^n
\]
for all $\rho \in [-1,1]$, where for $n \in \N$ and $x \in \R$
\[
H_n(x) : = \frac{1}{\sqrt{n!}} (-1)^n \exp\big(\frac{x^2}{2}\big)\,
\frac{\textup{d}^n}{\textup{d}x^n}\exp\big(-\frac{x^2}{2}\big) 
\]
denotes the (probabilistic version of the) $n$-th Hermite polynomial. Put 
$\alpha_n : = \langle b, H_n \rangle_{\gamma_1}^2$. If $h_{b,b}'(0) = \alpha_1 > 0$, 
we know that the real analytic function $h_{b,b}\big\vert_{(-1,1)}$ is invertible 
around $0 = h_{b,b}(0)$. Its inverse is also expressible as a power series there; 
i.e., around $0$, $\big(h_{b,b}\big\vert_{(-1,1)}\big)^{-1}$ is real analytic, 
too. Hence, given the assumption (CRA), listed in \sref{Chapter}{chp:Summary}, 
it follows that $g_b(y) : = h_{b,b}^{-1}(y) = \sum_{n=1}^\infty \beta_n \,y^n$ 
for all $y \in [-1,1]$, where $\beta_1 = \frac{1}{\alpha_1}$ and
\begin{align}\label{eq:coeff_of_the_inverse_function_Taylor_series__gen_case}
\begin{split}
\beta_n &= \frac{1}{n}\,\sum_{k=1}^{n-1}\frac{1}{\alpha_1^{n+k}}(-1)^{k}\,
\binom{n-1+k}{k}\,B_{n-1,k}^{\circ}(\alpha_2, \alpha_3, \ldots, \alpha_{n-k+1})\\
&= \frac{1}{n \alpha_1^n}\,\sum_{k=1}^{n-1}(-1)^{k}\,\binom{n-1+k}{k}\,
B_{n-1,k}^{\circ}(\frac{\alpha_2}{\alpha_1}, \frac{\alpha_3}{\alpha_1}, \ldots, 
\frac{\alpha_{n-k+1}}{\alpha_1})
\end{split}
\end{align}
for all $n \in \N_2$. In this context,
\[
B_{n,k}^{\circ}(x_1, x_2, \ldots, x_{n-k+1}) : =
\sum_{\nu \in P(n, k)}\binom{k}{\nu_1, \nu_2, \ldots, \nu_{n-k+1}}\,
\prod\limits_{i=1}^{n-k+1} x_i^{\nu_i} = \sum_{\nu \in P(n, k)}k!\frac{x^\nu}{\nu!}
\] 
denotes the ordinary partial Bell polynomial, where the multinomial coefficient 
$\binom{\sum_{i=1}^{n-k+1}\nu_i}{\nu_1, \nu_2, \ldots, \nu_{n-k+1}} : = 
\frac{\big(\sum_{i=1}^{n-k+1}\nu_i\big)!}{\prod\limits_{i=1}^{n+1-k}{\nu_i}!}$ 
represents the number of ways of allocating $\sum_{i=1}^{n-k+1}\nu_i$ distinct 
objects into $n-k+1$ distinct bins, with $\nu_i$ objects in the $i$'th bin and 
$P(n, k)$ indicates the set of all multi-indices $\nu \equiv 
(\nu_1, \nu_2, \ldots, \nu_{n+1-k}) \in \N_0^{n-k+1}$ ($k \leq n$) which 
satisfy the Diophantine equations 
\[
\sum_{i=1}^{n+1-k} \nu_i = k \text{ and } \sum_{i=1}^{n+1-k}i\,\nu_i = n\,;
\]
i.e., summation is extended over all partitions of the number $n$ into $k$ positive 
(non-zero) integers (cf. e.g. \cite{C1974, C2011, MMKS2019, T1992}). Observe that 
(dependent on the choice on $n$ and $k$, of course) these Diophantine equations 
may have an extremely large, if not even an unmanageable set of solutions! 
Already that definition implies the well-known and important fact that
\begin{align}\label{eq:calculation_rule_for_exp_partial_Bell_polynomials}
B_{n,k}^{\circ}(a b x_1, a b^2 x_2, \ldots, a b^{n-k+1} x_{n-k+1}) = a^k\,b^n 
B_{n,k}^{\circ}(x_1, x_2, \ldots, x_{n-k+1}),
\end{align}
for all $a, b, x_1, \ldots, x_{n-k+1} \in \C$. If the element $(x_1, \ldots, x_{n-k+1}) \in 
\N_0^{n-k+1}$ consists of at most 2 non-zero elements, $x_{i_1}$ and $x_{i_2}$, say 
($i_1 < i_2$), we only have to sum over the set of all $\nu \in P(n,k)$, such that $\nu_i = 0$
for all $i \notin \{i_1, i_2\}$ (since $0^0 = 1$). In this case, the Diophantine equations 
reduce to a simple 2-dimensional linear equation system. The latter has a solution 
$(\nu_{i_1}, \nu_{i_2}) \in \N_0^2$ if and only if
\begin{align}\label{eq:2-dim_Bell_case}
i_1 \leq \frac{n}{k} \leq i_2 \text{ and } i_2 - i_1 \text{ divides both, } k i_2 - n \in \N_0 
\text{ and } n - k i_1 \in \N_0\,.
\end{align}
The solution $(\nu_{i_1}, \nu_{i_2}) \in \N_0^2$ is then unique and given as
\[
\nu_{i_1} = \frac{k i_2 - n}{i_2 - i_1} \text{ and } \nu_{i_2} = \frac{n - k i_1}{i_2 - i_1}\,.
\] 
Hence, if \eqref{eq:2-dim_Bell_case} is satisfied, we immediately recognise that
\[
B_{n,k}^{\circ}(0, 0, \ldots, 0, x_{i_1},0, 0, \ldots, 0, x_{i_2},0, 0, \ldots, 0) =
\binom{k}{\nu_{i_1}}x_{i_1}^{\nu_{i_1}}x_{i_2}^{\nu_{i_2}}\,.
\]
In particular, if $k \leq n \leq 2k$, we reobtain the well-known special case
\[
B_{n,k}^{\circ}(x,y, 0, \ldots 0) = \binom{k}{2k-n}x^{2k-n}y^{n-k} \text{ for all } x,y 
\in \C\,.
\]
In the literature, one frequently finds the so-called \textit{exponential partial Bell 
polynomials $B_{n,k}$}, characterised as (cf. e.g. \cite[Remark on page 136]{C1974})
\begin{align}\label{eq:exponential_partial_Bell_vs_partial_Bell}
B_{n,k}(x_1, \ldots, x_{n-k+1}) : = \frac{n!}{k!}B_{n,k}^{\circ}\big(\frac{x_1}{1!},
\frac{x_2}{2!}, \ldots, \frac{x_{n-k+1}}{(n-k+1)!}\big),
\end{align}
implying that any result about $B_{n,k}$ can be directly converted into a result 
about $B_{n,k}^{\circ}$ and conversely. For example, if $n > k$, then 
\cite[formula (3l)]{C1974} transforms very pleasantly and reiteratively into
\begin{align*}
B_{n,k}^{\circ}(x_1, x_2, \ldots x_{n-k+1}) &= \sum_{l=k-\alpha(n,k)}^{k-1}
\binom{k}{l} x_1^l \,B_{n-k,k-l}^{\circ}(x_2, \ldots, x_{n-2k+l+1})\\
&= x_1^{k-\alpha(n,k)}\sum_{l=0}^{\alpha(n,k)-1}\binom{k}{\alpha(n,k)-l}\,
x_1^l \,B_{n-k,\alpha(n,k)-l}^{\circ}(x_2, \ldots, x_{(n-2k)^+ + l + 2})
\end{align*}
for all $x_1, x_2, \ldots, x_{n-k+1} \in \C$, where $\alpha(n,k) : = \min\{n-k, k\}$.
Another application of \eqref{eq:exponential_partial_Bell_vs_partial_Bell} implies 
the well-known fact that \eqref{eq:coeff_of_the_inverse_function_Taylor_series__gen_case} 
is equivalent to
\begin{align}\label{eq:coeff_of_the_inverse_function_Taylor_series__gen_case_II}
\big(h_{b,b}^{-1}\big)^{(n)}(0) = \sum_{k=1}^{n-1}\frac{1}{\delta_1^{n+k}}
(-1)^{k}\,\frac{(n-1+k)!}{(n-1)!}\,
B_{n-1,k}(\delta_2, \delta_3, \ldots, \delta_{n-k+1}),
\end{align} 
where $\delta_l : = \tfrac{h_{b,b}^{(l)}(0)}{l}, l \in [n-k+1]$ (cf. 
\cite[Theorem E on p. 150]{C1974}).

 Regarding an explicit recursive construction of these polynomials in 
full generality, yet without having to know the sets $P(n, k)$ beforehand, we 
recall the important fact that any ordinary partial Bell polynomial 
$B_{n,k}^\circ$ actually arises as a (kind of) discrete convolution of two 
ordinary partial Bell polynomial series. More precisely, we have:
\begin{lemma}\label{lem:Bell_polynomials_as_convolutions}
Let $m \in \N_0$, $k \in \N$, $n \in \N_k$ and $x_1, \ldots, x_{n-k+1} \in \C$. Then 
the following equalities are satisfied:
\begin{enumerate}
\item $B_{m,0}^\circ(x_1, \ldots, x_{m+1}) = \delta_{m0}$ and
\begin{align*}
B_{n,k}^\circ(x_1, \ldots, x_{n-k+1}) &= \sum_{i=k-1}^{n-1} x_{n-i}\,
B_{i,k-1}^\circ(x_1, \ldots, x_{i-k+2}) = \sum_{i=1}^{n-k+1} 
x_{i}\,B_{n-i,k-1}^\circ(x_1, \ldots, x_{n-k+2-i}).
\end{align*}
\item
\begin{align*}
n B_{n,k}^\circ(x_1, \ldots, x_{n-k+1}) &= 
k\sum_{i=k-1}^{n-1} (n-i) x_{n-i}\,
B_{i,k-1}^\circ(x_1, x_2, \ldots, x_{i-k+2})\\ 
&= k\sum_{i=1}^{n-k+1} i x_i\,
B_{n-i,k-1}^\circ(x_1, x_2, \ldots, x_{n-k+2-i}).
\end{align*}
\end{enumerate}
\end{lemma}
\begin{comment}
(i) Follows from  \cite[page 136, formula (3k)]{C1974}, respectively 
\cite[page 366, formula (13)]{OS1992}. 
\\[0.2em]
(ii) See \cite[formula (2.3) and its equivalent (unnumbered) representation on 
page 1546 (line 3)]{C2011}. 
\end{comment} 
\noindent \sref{Lemma}{lem:Bell_polynomials_as_convolutions} obviously implies 
the following multiple-sum representation of the ordinary partial Bell 
polynomials: 
\begin{align}\label{eq:Bell_closed_form}
B_{n,k+1}^{\circ}(x_1, x_2, \ldots, x_{n-k}) &=
\sum_{i_1=k}^{n-1}\sum_{i_2=k-1}^{i_1-1}\ldots\sum_{i_k=1}^{i_{k-1}-1}x_{n-i_1}\,
\prod\limits_{\nu=2}^k x_{i_{\nu-1}-i_\nu}\,x_{i_k}
\end{align}
for all $k \in \N$, $n \in \N_{k+1}$ and $x_1, x_2, \ldots, x_{n-k} \in \C$. 
For the convenience of the readers, we list a few examples of ordinary partial Bell 
polynomials that can be displayed in closed form. For review, we 
refer to the widely comprehensive table of these polynomials on page 309 of 
\cite{C1974} (enumerating all polynomials $B_{n,m}^{\circ}$ for which $10 \geq n 
\geq m \geq 1$). To this end, fix $k \in \N$ and $x_1, \ldots, x_{k+1} \in \C$. 
Then
\begin{enumerate}
\item $B_{0,0}(x_1) = 1$ and $B_{k,0}^{\circ}(x_1, \ldots, x_{k+1}) = 0$.
\item $B_{k,1}^{\circ}(x_1, \ldots, x_k) = x_k$ and $B_{k,k}(x_1) = x_1^k$.
\item $B_{k+1,k}^{\circ}(x_1, x_2) = k\,x_1^{k-1}\,x_2$.
\item $B_{k+2,k}^{\circ}(x_1, x_2, x_3) = \colvec{2}{k}{2}\,x_1^{k-2}\,x_2^2 + 
k\,x_1^{k-1}\,x_3$ \,if\, $k \geq 2$.
\item $B_{k+3,k}^{\circ}(x_1, x_2, x_3, x_4) = \colvec{2}{k}{3}\,x_1^{k-3}\,
x_2^3 + k(k-1)\,x_1^{k-2}\,x_2\,x_3 + k\,x_1^{k-1}\,x_4$ \,if\, $k \geq 3$.
\item $B_{k+4,k}^{\circ}(x_1, \ldots, x_5) = \colvec{2}{k}{4}\,x_1^{k-4}\,
x_2^4 + \colvec{2}{k}{3}\,x_1^{k-3}(3\,x_2^2\,x_3) + \colvec{2}{k}{2}\,
x_1^{k-2}(x_3^2 + 2\,x_2\,x_4) + k\, x_1^{k-1}\,x_5$ \,if\, $k \geq 4$.
\end{enumerate}
Moreover, we have
\[
B_{k,2}^\circ(x_1, \ldots, x_{k-1}) = 
\sum_{i=1}^{k-1} x_i\,x_{k-i} \,\text{ if }\, k \geq 2\,.
\]
Since the Taylor series of the inverse of the ``standardised'' 
Taylor series $\sum_{n=1}^\infty \alpha_n^\times\,\rho^n = \frac{1}{\alpha_1}\,
h_{b,b}(\rho)$ of the function $h_{b,b}$ obviously is given by $\big(\frac{1}{\alpha_1}\,
h_{b,b}\big)^{-1}(y) = h_{b,b}^{-1}(\alpha_1\,y) = \sum_{n=1}^\infty 
(\beta_n\,\alpha_1^n)y^n$ for all $y \in [-1,1]$, where $\alpha_n^\times : = 
\frac{\alpha_n}{\alpha_1}$ ($n \in \N$), it follows that the $n$-th Taylor 
series coefficient $\beta_n^\times$ of the Taylor series of $\big(\frac{1}{\alpha_1}\,
h_{b,b}\big)^{-1}$ is given by $\beta_n^\times = \beta_n\,\alpha_1^n$. 
Consequently,
\begin{align}\label{eq:beta_n_Bell_pol_rep}
\beta_n\,\alpha_1^n = \beta_n^\times = \frac{1}{n}\,\sum_{k=1}^{n-1}\,(-1)^k
\binom{n-1+k}{k}\,B_{n-1,k}^{\circ}(\alpha_2^\times, \alpha_3^\times, 
\ldots, \alpha_{n-k+1}^\times)
\end{align}
(due to \eqref{eq:coeff_of_the_inverse_function_Taylor_series__gen_case}).
In the odd case, i.e., if in addition $\alpha_{2n} = 0$ for all $n \in \N$, 
the intrinsic combinatorical complexity of 
\eqref{eq:coeff_of_the_inverse_function_Taylor_series__gen_case} can be even 
further reduced, possibly allowing a non-negligible saving of computing time 
(see \sref{Corollary}{cor:coefficients_in_the_odd_case}). 

 We explicitly list $\beta_1, \beta_2, \beta_3, \beta_4, \beta_5, 
\beta_6$ and $\beta_7$ in full generality. To this end, as discussed above, 
if $\alpha_1 \not= 0$, we have to consider the Taylor coefficents 
$\alpha_n^\times : = \frac{\alpha_n}{\alpha_1}, n \in \N$ of the 
``standardised'' Taylor series $\frac{1}{\alpha_1}\,h_{b,b}$. Note again that 
$\alpha_1^\times = 1$. 
\eqref{eq:coeff_of_the_inverse_function_Taylor_series__gen_case} therefore 
implies that
\begin{align}\label{eq:First_7_Taylor_coeff_of_the_inverse}
\begin{split}
\beta_1\,\alpha_1 &= 1 = \alpha_1^\times\\
\beta_2\,\alpha_1^2 &= -\alpha_2^\times\\ 
\beta_3\,\alpha_1^3 &= -\alpha_3^\times + 2(\alpha_2^\times)^2\\
\beta_4\,\alpha_1^4 &= -\alpha_4^\times + 5\alpha_2^\times
\alpha_3^\times - 5(\alpha_2^\times)^3\\
\beta_5\,\alpha_1^5 &= -\alpha_5^\times + 6\alpha_2^\times
\alpha_4^\times + 3(\alpha_3^\times)^2 - 21(\alpha_2^\times)^2
\alpha_3^\times + 14(\alpha_2^\times)^4\\
\beta_6\,\alpha_1^6 &= -\alpha_6^\times + 7\alpha_2^\times
\alpha_5^\times + 7\alpha_3^\times\alpha_4^\times - 
28\alpha_2^\times(\alpha_3^\times)^2 - 28(\alpha_2^\times)^2
\alpha_4^\times + 84(\alpha_2^\times)^3\alpha_3^\times - 
42(\alpha_2^\times)^5\\
\beta_7\,\alpha_1^7 &= -\alpha_7^\times + 8\alpha_2^\times
\alpha_6^\times + 8\alpha_3^\times \alpha_5^\times + 
4(\alpha_4^\times)^2 - 36(\alpha_2^\times)^2\alpha_5^\times - 
72\alpha_2^\times\alpha_3^\times\alpha_4^\times - 
12(\alpha_3^\times)^3 + 120(\alpha_2^\times)^3\alpha_4^\times\\ 
&\hspace{0.45cm}+ 180(\alpha_2^\times)^2(\alpha_3^\times)^2 - 
330(\alpha_2^\times)^4\alpha_3^\times + 132(\alpha_2^\times)^6\,.
\end{split}
\end{align}
In fact, if we make use of the key result, listed in \cite{W1951}, paired with 
the general Theorem \cite[p. 222]{CENN2002}, we are able to present a (purely linear 
algebraic and algorithmic) representation of each coefficient 
$\beta_n$, which avoids an explicit use of ordinary partial Bell polynomials (where no closed 
form seems to be available). To the best of our knowledge, in this context, 
that representation had not been published before. Instead of ordinary partial Bell 
polynomials, we have to calculate determinants of leading principal submatrices. Of course, 
the computational complexity induced by the increasing size of partial Bell 
polynomials transforms into the rapidly increasing computing time, induced by the increasing 
size of the determinants including the need to sum proper parts of determinants of 
\text{different} size. However, that summation is a recurrence relation (see 
\eqref{eq:recursive_constr_of_dets} and the instructive \sref{Example}{ex:det_calculus}). 

 Firstly, if $\alpha_1 \not= 0$, an enhancement of \cite{W1951} 
reveals the following explicit representation of each coefficient $\beta_n$ 
($n \in \N_2$):
\begin{align}\label{eq:matrix_rep_of_beta_ns}
\begin{split}
\beta_n &= \frac{(-1)^{n-1}}{n!\,\alpha_1^n}\,
\det(A_n \ast T_n(1, \alpha_2^\times, \alpha_3^\times, \ldots, 
\alpha_n^\times))\\
&= \frac{(-1)^{n-1}}{n!\,\alpha_1^{2n-1}}\,\det(A_n \ast 
T_n(\alpha_1, \alpha_2, \ldots, \alpha_n))\\ 
&= \frac{1}{n!\,\alpha_1^{2n-1}}\,\det(-(A_n \ast 
T_n(\alpha_1, \alpha_2, \ldots, \alpha_n))),
\end{split}
\end{align}
where $\ast$ again denotes the Hadamard product and the matrices $A_n \in 
\M_{n-1}(\R)$ and $T_n \equiv T_n(x_1, x_2, \ldots, x_n) \in 
\M_{n-1}(\R)$ ($x_1, \ldots, x_n \in \C$) are respectively defined as 
\[
A_n : = 
\begin{pmatrix}
n      & 1 			  & 0        & 0       & 0      & \cdots & 0   & 0\\
2n     & n+1  	  & 2        & 0       & 0      & \cdots & 0   & 0\\
3n     & 2n+1 	  & n+2      & 3       & 0      & \cdots & 0   & 0\\
4n     & 3n+1 	  & 2n+2     & n+3     & 4      & \cdots & 0   & 0\\
\vdots & \vdots   & \vdots   & \vdots  & \vdots & \vdots & \vdots   & \vdots\\
(n-3)n & (n-4)n+1	& (n-5)n+2 & \cdots  & \cdots & \cdots & n-3        & 0\\
(n-2)n & (n-3)n+1	& (n-4)n+2 & \cdots  & \cdots & \cdots & n+(n-3)  & n-2\\
(n-1)n & (n-2)n+1 & (n-3)n+2 & \cdots  & \cdots & \cdots & 2n+(n-3) & n+(n-2)
\end{pmatrix}
\]
and
\[
T_n \equiv T_n(x_1, x_2, \ldots, x_n) : = 
\begin{pmatrix}
x_2     & x_1 	  & 0        & 0        & 0        & \cdots   & 0        & 0\\
x_3     & x_2  	  & x_1 & 0        & 0        & \cdots   & 0        & 0\\
x_4     & x_3 	  & x_2 & x_1 & 0        & \cdots   & 0        & 0\\
x_5     & x_4 	  & x_3 & x_2 & x_1 & \ddots   & 0        & 0\\
\vdots       & \vdots       & \vdots   & \vdots   & \vdots   & \ddots   & \ddots   & \vdots\\
x_{n-2} & x_{n-3}	& \cdots   & \cdots   & x_3 & x_2 & x_1 & 0\\
x_{n-1} & x_{n-2}	& \cdots   & \cdots   & \cdots   & x_3 & x_2 & x_1\\
x_n     & x_{n-1} & x_{n-2} & \cdots & \cdots & x_4 & x_3 & x_2
\end{pmatrix}\,. 
\]
Obviously, the Toeplitz matrix $T_n(x_1, x_2, \ldots, x_n)$ is well-defined 
for any $x_1, \ldots, x_n \in \C$. Observe that the appearance of the rather 
uncommon factor $\frac{1}{\alpha_1^{2n-1}} = \frac{1}{\alpha_1^{n}}\cdot 
\frac{1}{\alpha_1^{n-1}}$ in \eqref{eq:matrix_rep_of_beta_ns} actually 
originates from the simple, yet important, transformation
\begin{align}\label{eq:WLOG_alpha_1_equals_1}
T_n \equiv T_n(\alpha_1, \alpha_2, \ldots, \alpha_n) = \alpha_1\,
T_n(1, \alpha_2^\times, \alpha_3^\times, \ldots, \alpha_n^\times).
\end{align}
More precisely, if $n \in \N$ and $(i,j) \in [n-1] \times [n-1]$, we have:
\[
(A_n)_{ij} : = \begin{cases}
0 &\text{if } j \geq i+2\\
i &\text{if } j = i+1\\ 
(i-j+1)n + j-1 &\text{if } j \leq i
\end{cases}
\text{ and }
(T_n)_{ij} : = \begin{cases}
0 &\text{if } j \geq i+2\\
\alpha_1 &\text{if } j = i+1\\  
x_{i-j+2} &\text{if } j \leq i
\end{cases}.
\]
Comparing \eqref{eq:coeff_of_the_inverse_function_Taylor_series__gen_case} and 
\eqref{eq:matrix_rep_of_beta_ns}, it follows that for all $\alpha_1 \in \C^\ast$ and 
$\alpha_2, \ldots, \alpha_n \in \C$
\[
\det((A_n \ast T_n(\alpha_1, \ldots, \alpha_n))) =
\alpha_1^{n-1}\sum_{k=1}^{n-1}(-1)^{n-1+k}\,
\frac{(n-1+k)!}{k!}\,B_{n-1,k}^{\circ}(\alpha_2^\times, \alpha_3^\times, 
\ldots, \alpha_{n-k+1}^\times).
\]
Hence,
\[
\det((A_n \ast T_n(1, \alpha_2^\times, \alpha_3^\times, \ldots, 
\alpha_n^\times))) = \sum_{k=1}^{n-1}(-1)^{n-1+k}\,\frac{(n-1+k)!}{k!}\,
B_{n-1,k}^{\circ}(\alpha_2^\times, \alpha_3^\times, \ldots, 
\alpha_{n-k+1}^\times).
\]
\begin{remark}[\textbf{Connection to Apostol's approach in \cite{A2000}}]
In fact, 
it can be shown that
\[
\det\big(-\big(A_n \ast T_n\big(x_1, \frac{x_2}{2!}, 
\frac{x_3}{3!}, \ldots, \frac{x_n}{n!}\big)\big)\big) = 
P_n(x_1, \ldots, x_n) \text{ for all } x_1, \ldots, x_n \in 
\C
\]
precisely coincides with the function $P_n$, introduced in \cite{A2000} (due to 
the convolution representation and the partial derivative structure of ordinary 
partial Bell polynomials)! 
In particular,
\begin{align*}
P_n(1, x_2, \ldots, x_n) &= 
\sum_{k=1}^{n-1}(-1)^{k}\,
\frac{(n-1+k)!}{k!}\,B_{n-1,k}^{\circ}\big(\frac{x_2}{2!}, 
\frac{x_3}{3!}, \ldots, \frac{x_{n-k+1}}{(n-k+1)!}\big)\\
&= \det\big(-\big(A_n \ast T_n\big(1, \frac{x_2}{2!}, 
\frac{x_3}{3!}, \ldots, \frac{x_n}{n!}\big)\big)\big).
\end{align*}
\end{remark}
\noindent Due to \eqref{eq:WLOG_alpha_1_equals_1}, it follows that  
\[
A_n \ast T_n(\alpha_1, \ldots, \alpha_n) = \alpha_1 (A_n \ast 
T_n(1, \alpha_2^\times, \ldots, \alpha_n^\times)) = \alpha_1 
B_n[n-1](\alpha_2^\times, \ldots, \alpha_n^\times),
\]
where
for any $p \in [n-1]$ and $x_1, \ldots, x_p \in \C$, the matrix $B_n[p] 
\equiv B_n[p](x_1, \ldots, x_p) \in \M_p(\C)$ is defined as
\[
\scalebox{0.80}{$
B_n[p] : =
\begin{pmatrix}
 n x_1      & 1			  & 0        & 0       & 0      & \cdots & 0   & 0\\
2n x_2     & (n+1)x_1  	  & 2        & 0       & 0      & \cdots & 0   & 0\\
3nx_3     & (2n+1)x_2 	  & (n+2)x_1      & 3       & 0      & \cdots & 0   & 0\\
4n x_4     & (3n+1)x_3 	  & (2n+2)x_2     & (n+3)x_1     & 4      & \cdots & 0   & 0\\
\vdots & \vdots   & \vdots   & \vdots  & \vdots & \vdots & \vdots   & \vdots\\
(p-2)n x_{p-2} & ((p-3)n+1)x_{p-3}	& ((p-4)n+2)x_{p-4} & \cdots  & \cdots & \cdots & p-2        & 0\\
(p-1)n x_{p-1} & ((p-2)n+1)x_{p-2}	& ((p-3)n+2)x_{p-3} & \cdots  & \cdots & \cdots & 
(n+(p-2))x_1  & p-1\\
p n x_p & ((p-1)n+1)x_{p-1} & ((p-2)n+2)x_{p-2} & \cdots  & \cdots & \cdots & 
(2n+(p-2))x_2 & (n+(p-1)x_1)
\end{pmatrix}
$}
\]
$B_n[p]$ therefore denotes the $p$-th leading principal submatrix of the 
matrix $B_n[n-1]$ if $p \in [n-1]$ (cf. e.g. \cite[0.7.1]{HJ2013}). 
More precisely, if $p \in [n-1]$ and $i, j \in [p]$, then:
\[
B_n[p](x_1, x_2,\ldots, x_p)_{ij} : =
\begin{cases}
0 &\text{if } j \geq i+2\\
i &\text{if } j = i+1\\  
((i-j+1)n + j-1)x_{i-j+1} &\text{if } j \leq i\,.
\end{cases}
\]
Equipped with all $p$ leading principal submatrices $B_n[1], 
B_n[2], \ldots, B_n[p-1]$ of $B_n[p]$ and the ``incipient matrix'' 
$B_n[0] : = (1)$, we may apply the main (unnumbered) theorem on page 222 of 
\cite{CENN2002} to the matrix $B_n[p](x_1, \ldots, x_p)$, and it follows that
\begin{align}\label{eq:recursive_constr_of_dets}
\det(B_n[p](x_1, \ldots, x_p)) &= 
(p-1)!\sum_{k=1}^{p}(-1)^{p-k}\,\frac{pn - (k-1)(n-1)}{(k-1)!}\,
x_{p-k+1}\,\det({B_n[k-1]}).
\end{align}
Consequently, if $p = n-1$, it follows that

\scalebox{0.90}{
\vbox{
\begin{align}\label{eq:det_of_B_n}
\begin{split}
\det(A_n \ast T_n(\alpha_1, \ldots, \alpha_n)) 
&= \alpha_1^{n-1}\,\det(B_n[n-1](\alpha_2^\times, \ldots, 
\alpha_n^\times))\\ 
&= \alpha_1^{n-1}\,(n-2)!\sum_{k=1}^{n-1}(-1)^{n-1-k}\,\frac{(n-1)(n-k+1)}{(k-1)!}\,
\alpha_{n-k+1}^\times\,\det(B_n[k-1](\alpha_2^\times, \ldots, 
\alpha_k^\times))\\
&= \alpha_1^{n-1}\,(n-1)!\,\sum_{k=1}^{n-1}(-1)^{n-k-1}\,\frac{n-k+1}{(k-1)!}\,
\alpha_{n-k+1}^\times\,\det(B_n[k-1](\alpha_2^\times, \ldots, 
\alpha_k^\times)).
\end{split}
\end{align}
}}

Thus, 
\begin{align*}
\beta_n\,\alpha_1^n &\stackrel{\eqref{eq:matrix_rep_of_beta_ns}}{=}
\frac{1}{n}\sum_{k=1}^{n-1}(-1)^{k}\,\frac{n-k+1}{(k-1)!}\,
\alpha_{n-k+1}^\times\,\det(B_n[k-1](\alpha_2^\times, \ldots, 
\alpha_k^\times))\,.
\end{align*}
 In particular, if $n = 2m+1 \in \N_3$ is odd ($m \in \N$) and 
$\alpha_{2l} : = 0$ for all $l\in \N$, then
\begin{align*}
\beta_n\,\alpha_1^n &= 
\frac{1}{n}\sum_{\stackrel{k=1}{k {\text{ odd}}}}^{n-1}(-1)^{k}\,
\frac{n-k+1}{(k-1)!}\,\alpha_{n-k+1}^\times\,
\det(B_n[k-1](0, \alpha_3^\times, 0 \ldots, 0, \alpha_{k}^\times))\\
&= -\alpha_{2m+1}^\times -\frac{1}{2m+1} 
\sum_{r=1}^{m-1}\frac{2(m-r)+1}{(2r)!}\,\alpha_{2(m-r)+1}^\times\,
\det(B_{2m+1}[2r](0, \alpha_3^\times, 0 \ldots, 0, 
\alpha_{2r+1}^\times)).
\end{align*}
Note also that \eqref{eq:recursive_constr_of_dets} implies that
\begin{align}\label{eq:recursive_constr_of_dets_in_the_odd_case}
\det(B_{2m+1}[2r](0, x_1, 0, x_2, 0, \ldots, 0, x_r)) = 
- 2(2r-1)!\,\big((2m+1)r\,x_r + \sum_{k=1}^{r-1} p_k(r,m)\,x_{r-k}\big)
\end{align}
for all $m \in \N$, $r \in [m]$ and $x_1, \ldots, x_r \in \C$, where
\[
p_k(r,m) : = \frac{2m(r-k) + r}{(2k)!}\,\det({B_{2m+1}[2k]}
(0, x_1, 0 \ldots, 0, x_k) \hspace{0.5cm} 
(k \in [r-1]).
\]
\begin{example}\label{ex:det_calculus}
Fix $m \in \N$. Assume for simplification that $\alpha_1 = 1$ $\alpha_{2l} = 0$ 
for all $l\in \N$. Then $\alpha_n^\times = \alpha_n$ for all $n \in \N$. If 
$r \in \{1, 2\}$, the calculation of $\det({B_{2m+1}[2r]})$ is very 
straightforward:
\[
\det(B_{2m+1}[2](0, \alpha_3)) = -2(2m+1)\alpha_3 = 
2!(-1)^1\binom{2m+1}{1}\alpha_3 = 2!(-1)^1\binom{2m+1}{1}B_{1,1}^\circ(\alpha_3),
\]
and
\begin{align*}
\det(B_{2m+1}[4](0, \alpha_3, 0, \alpha_5)) 
&\stackrel{\eqref{eq:recursive_constr_of_dets_in_the_odd_case}}{=} 
-12\,\big(2(2m+1)\alpha_5 + 
\frac{2m+2}{2}\,\det({B_{2m+1}[2]})\alpha_3\big)\\
&= -12\,\big(2(2m+1)\alpha_5 + 
\frac{2m+2}{2}\,\big(-2(2m+1)\alpha_3\big)\alpha_3\big)\\
&= 4!\big(-\binom{2m+1}{1}\alpha_5 + \binom{2m+2}{2}\alpha_3^2\big)\\
&= 4!\big((-1)^1\binom{2m+1}{1}B_{2,1}^\circ(\alpha_3, \alpha_5) + 
(-1)^2\binom{2m+2}{2} B_{2,2}^\circ(\alpha_3)\big).
\end{align*}
If $r =3$, a little more calculation effort is required, also triggered by a 
significant transformation step:

\scalebox{0.88}{
\vbox{
\begin{align*}
\det(B_{2m+1}[6](0, \alpha_3, 0, \alpha_5, 0 , \alpha_7)) 
&\stackrel{\eqref{eq:recursive_constr_of_dets_in_the_odd_case}}{=} 
-240\,\big(3(2m+1)\alpha_7 + 
\frac{4m+3}{2}\,\det({B_{2m+1}[2]})\alpha_5 + 
\frac{2m+3}{4!}\,\det({B_{2m+1}[4]})\alpha_3\big)\\
&= -240\,\big(3\binom{2m+1}{1}\alpha_7{-(4m+3)\binom{2m+1}{1}\alpha_3\alpha_5}\\  
&\hspace{0.45cm}{+(2m+3)}\big({-\binom{2m+1}{1}\alpha_5} + 
\binom{2m+2}{2}\alpha_3^2\big){\alpha_3}\Big)\\
&= -240\,\Big(3\binom{2m+1}{1}\alpha_7 {- 3\,\binom{2m+2}{2}(2\alpha_3\alpha_5)} + 
(2m+3)\binom{2m+2}{2}\alpha_3^3\Big)\\
&= 6!\,\Big(-\binom{2m+1}{1}\alpha_7 + \binom{2m+2}{2}\,(2\alpha_3\alpha_5) - \binom{2m+3}{3}
\alpha_3^3\Big)\\
&= 6!\,\Big((-1)^1\binom{2m+1}{1}B_{3,1}^\circ(\alpha_3, \alpha_5, \alpha_7) + (-1)^2\,
\binom{2m+2}{2}\,B_{3,2}^\circ(\alpha_3, \alpha_5)\\ 
&\hspace{1cm}+ (-1)^3\,\binom{2m+3}{3}\,B_{3,3}^\circ(\alpha_3)\Big).
\end{align*}
}}

 In fact, the emerging structure can be kept in the case of $r=4$, since

\scalebox{0.88}{
\vbox{
\begin{align*}
\det({B_{2m+1}[8]})(0, \alpha_3, 0, \alpha_5, 0, \alpha_7, 0, \alpha_9) 
&\stackrel{\eqref{eq:recursive_constr_of_dets_in_the_odd_case}}{=} 
-2\cdot 7!\,\Big(4(2m+1)\alpha_9 + 
\frac{6m+4}{2}\,\det({B_{2m+1}[2]})\alpha_7\\ 
&\hspace{1.1cm}+\frac{4m+4}{4!}\,\det({B_{2m+1}[4]})\alpha_5 + \frac{2m+4}{6!}\,
\det({B_{2m+1}[6]})\alpha_3\Big)\\
&\hspace{0.2cm}= 8!\,\Big(-\binom{2m+1}{1}\alpha_9 + \binom{2m+2}{2}
(2\alpha_3\alpha_7 + \alpha_5^2) - \binom{2m+3}{3}(3\alpha_3^2\alpha_5)\\ 
&\hspace{1.1cm}+ \binom{2m+4}{4}\alpha_3^4\Big)\\
&\hspace{0.2cm}= 8!\,\Big((-1)^1\,\binom{2m+1}{1}B_{4,1}^\circ(\alpha_3, 
\alpha_5, \alpha_7, \alpha_9) + (-1)^2\,\binom{2m+2}{2}B_{4,2}^\circ(\alpha_3, 
\alpha_5, \alpha_7)\\ 
&\hspace{1.1cm} + (-1)^3\,\binom{2m+3}{3}B_{4,3}^\circ(\alpha_3, \alpha_5) + 
(-1)^4\,\binom{2m+4}{4}B_{4,4}^\circ(\alpha_3)\Big).
\end{align*}
}}
\end{example}
\noindent A relentless focus on \sref{Example}{ex:det_calculus} therefore leads 
to a non-obvious simplification of \eqref{eq:beta_n_Bell_pol_rep} which 
reduces the analysis of complex partition sets $P(2m, k)$ and related non-trivial ordinary 
partial Bell polynomials $B_{2m,k}^\circ(0, \alpha_3, 0, \alpha_5, 0, \alpha_7, 0, \ldots)$ 
to that one of partition sets $P(m, l)$ and related ``fully occupied'' ordinary partial Bell 
polynomials $B_{m,l}^\circ(\alpha_3, \alpha_5, \alpha_7, \ldots)$. In fact, the 
following result holds:
\begin{proposition}\label{prop:determinant_and_Bell_polynomials}
Let $(x_n)_{n \in \N}$ be an arbitrary sequence of complex numbers. Let $m \in 
\N$ and $r \in [m]$.
Then
\[
\det(B_{2m+1}[2r](0,x_1, 0, x_2, \ldots, 0, x_r)) = 
(2r)!\sum_{l=1}^{r}(-1)^l\binom{2m+l}{l}\,
B_{r,l}^{\circ}(x_1, x_2, \ldots, x_{r-l+1}).
\]
\end{proposition}
\noindent Altogether, \sref{Lemma}{lem:Bell_polynomials_as_convolutions}, 
 \sref{Proposition}{prop:determinant_and_Bell_polynomials}
and \eqref{eq:double_sum_representation}, together with our previously mentioned 
analysis of the structure of $\beta_{2m+1}$, imply the following two fundamental 
results:
\begin{theorem}\label{thm:determinants_vs_Bell_polynomials}
Let $m \in \N_2$ and $x_1, \ldots, x_{m-1} \in \C$. Then
\begin{align*}
{}&\hspace{0.45cm}\sum_{r=1}^{m-1}\frac{2(m-r)+1}{(2r)!}\,x_{m-r}\,
\det\big({B_{2m+1}[2r](0, x_1, 0, x_2, \ldots, 
0, x_r)}\big)\\
&= \frac{1}{(2m)!}\sum_{r=2}^{m}(-1)^{r-1}\,\frac{(2m+r)!}{r!}\,
B_{m,r}^{\circ}(x_1, x_2, \ldots, x_{m-r+1}).
\end{align*}
\end{theorem}
\begin{comment}
By consecutive application of 
\sref{Proposition}{prop:determinant_and_Bell_polynomials} and 
\eqref{eq:double_sum_representation}, it follows that
\begin{align*}
{}&\hspace{0.45cm}\sum_{r=1}^{m-1}\frac{2(m-r)+1}{(2r)!}\,x_{m-r}\,
\det\big({B_{2m+1}[2r](0, x_1, 0, x_2, \ldots, 0, x_r)}\big)\\
&= \sum_{r=1}^{m-1}\sum_{l=1}^{r} (-1)^l\,(2(m-r) + 1)\,\,
\binom{2m+l}{l}\,x_{m-r}\, 
B_{r,l}^{\circ}(x_1, x_2, \ldots, x_{r-l+1})\\
&= \sum_{r=1}^{m-1} (-1)^r\,\binom{2m+r}{r} \sum_{l=r}^{m-1} (2(m-l) + 1)\,
x_{m-l}\,B_{l,r}^{\circ}(x_1, x_2, \ldots, x_{l-r+1}).
\end{align*}
However, the latter expression equals
\begin{align*}
{}&\hspace{0.45cm}\sum_{r=1}^{m-1} (-1)^r\,2\,\binom{2m+r}{r}\sum_{l=r}^{m-1}
(m-l)\,x_{m-l}\,B_{l,r}^{\circ}(x_1, x_2, \ldots, x_{l-r+1})\\
&+\sum_{r=1}^{m-1} (-1)^r\,\binom{2m+r}{r}\sum_{l=r}^{m-1}\,x_{m-l}\,
B_{l,r}^{\circ}(x_1, x_2, \ldots, x_{l-r+1}).
\end{align*}
Thus, we can now apply \sref{Lemma}{lem:Bell_polynomials_as_convolutions} 
to both summands, and a final shift of the index $r$ in the single remaining 
sum ($\sum_{r=1}^{m-1} a_{r+1} = \sum_{r=2}^{m} a_r$) clearly finishes the proof. 
\end{comment}
\begin{corollary}\label{cor:coefficients_in_the_odd_case}
Let $f(\rho) = \sum_{n=0}^\infty \alpha_{2n+1}\,\rho^{2n+1}$ be an odd real 
analytic function, convergent on $(-r,r) \subseteq \R$, where $r > 0$ denotes 
the radius of convergence of $f$. Assume that $f^{\prime}(0) = \alpha_1 \not= 
0$, implying that $f$ is invertible around $0$. Consider the real analytic odd 
inverse function $f^{-1} : V \longrightarrow \R$, where $V$ is an open 
neighbourhood of $f(0)=0$. If $f^{-1}(y) = \sum_{m=0}^\infty \beta_{2m+1}\,
y^{2m+1}$ for all $y \in V$, then $\beta_1 = \frac{1}{\alpha_1}$ and

\scalebox{0.95}{
\vbox{
\begin{align*}
\beta_{2m+1} &= -\frac{1}{\alpha_1^{2m+1}}\Big(\alpha_{2m+1}^\times +  
\frac{1}{2m+1}\sum_{r=1}^{m-1}\frac{2(m-r)+1}{(2r)!}\,\alpha_{2(m-r)+1}^\times\,
\det\big({B_{2m+1}[2r](0, \alpha_3^\times, 0, \alpha_5^\times, \ldots, 
0, \alpha_{2r+1}^\times)}\big)\Big)\\
&= \frac{1}{(2m+1)!\,\alpha_1^{2m+1}}\sum_{r=1}^{m}(-1)^r\,\frac{(2m+r)!}{r!}\,
B_{m,r}^{\circ}(\alpha_3^\times, \alpha_5^\times, \ldots, 
\alpha_{2(m-r+1)+1}^\times)
\end{align*}
}}
for all $m \in \N$, where $\alpha_{2\nu+1}^\times : = 
\frac{\alpha_{2\nu+1}}{\alpha_1}$ ($\nu \in \N$).
\end{corollary}
\noindent As a little, yet illuminating exercise, we recommend the readers to 
perform the rather quick calculation of the first 3 Taylor coefficents of the 
inverse of the odd function $f : = 3\sinh$, say, by applying 
\sref{Corollary}{cor:coefficients_in_the_odd_case}\,! The outcome could be 
double-checked by means of 
\url{https://en.wikipedia.org/wiki/Inverse_hyperbolic_functions#Series_expansions}. 
Moreover, because of \cite[Theorem 2]{T1992} the coefficients $\beta_{2m+1}$ 
satisfy the following, interesting recurrence relation: 
\begin{align*}
\beta_{2m+1}\alpha_1^{2m+1} &= -\sum_{r=0}^{m-1} (\beta_{2r+1}\alpha_1^{2r+1})\,
B_{2m+1,2r+1}^{\circ}(\alpha_1^\times, 0, \alpha_3^\times, 0, \alpha_5^\times, 
\ldots, 0, \alpha_{2(m-r)+1}^\times)\\
&= -\sum_{r=0}^{m-1} \beta_{2r+1}\,
B_{2m+1,2r+1}^{\circ}(\alpha_1, 0, \alpha_3, 0, \alpha_5, \ldots, 0, 
\alpha_{2(m-r)+1})
\end{align*}
for all $m \in \N$.

\noindent In a nutshell, we recognise that already in the one-dimensional case, 
at least two hard open problems appear. On the one hand we need to know the 
explicit value of the Fourier-Hermite coefficients (cf. 
\sref{Proposition}{prop:k_dim_Fourier_Hermite_coeff}) 
\[
\sqrt{n!}\,\langle b, H_n \rangle_{\gamma_1} = 
\sqrt{n!}\,\E[b(X)H_n(X)] = \frac{\textup{d}^n}{\textup{d}t^n}
\E[b(X+t)]\Big\vert_{t=0}\,,
\]
where $X \sim N_1(0,1)$. On the other hand, we have to look for a closed form 
expression of the coefficients $\beta_n$ (if it were available at all), where 
the latter involves the complex recursive structure of ordinary partial Bell 
polynomials or related determinants. For example (keeping the Haagerup function 
in mind - cf. \sref{Example}{ex:Haagerup} and 
\sref{Remark}{rem:integral_rep_of_Haagerup}), our question of the value of
\[
\frac{\pi^k}{4^k}\,B_{n,k}^{\circ}\Big(\frac{1}{8}, \frac{3}{64}, \frac{25}{1024}, 
\ldots, \frac{((2(n-k)+1)!!)^2}{((2(n-k+1))!!)^2\,(n-k+2)}\Big)
\]
very recently lead to an in depth-analysis, published in \cite{DLQ2022}. It 
appears to us that in general one cannot use proofs by standard induction on 
$n \in \N$ to verify statements about Bell polynomials. The Noetherian Induction 
Principle seems to be more appropriate here (as we have seen for example, in the 
proof of \sref{Proposition}{prop:determinant_and_Bell_polynomials}).
In this context, we would like to draw attention to another recently published 
paper, where the authors point to similar difficulties including the formulation 
of related - open - problems (cf. \cite{LNQY2020}). Moreover, the solved examples 
in \cite{LNQY2020} show the large combinatorial barriers which we have to resolve 
while working with (partial) Bell polynomials.
 
Keeping these problems and barriers in mind, the following research topics and 
problems - which actually do not require any knowledge of the Grothendieck 
inequality - arise naturally:
\begin{enumerate}
\item[(RP1)] Continue to investigate the structure of partial Bell polynomials; 
possibly under inclusion of the use of supercomputers and related 
computer algebra systems. 
\item[(RP2)] Develop a software package which puts 
\sref{Corollary}{cor:coefficients_in_the_odd_case} into practice. 
\item[(RP3)] Look for an explicit analytic expression for the \textit{inverse} 
function of the main building block of the Haagerup function; i.e., the inverse 
of the strictly increasing odd function
\[
[-1,1] \ni x \mapsto x\,{}_2 F_1(\frac{1}{2}, \frac{1}{2}, 2; x^2)
\]
(if available\,!), where as usual ${}_2 F_1(a,b,c; \cdot)$ denotes the classic 
Gaussian hypergeometric function (cf. \sref{Example}{ex:Haagerup}). Obviously, 
the inverse of $[-1,1] \ni x \mapsto x \cdot {}_2 F_1(\frac{1}{2}, \frac{1}{2}, 
\frac{3}{2}; x^2) = \arcsin(x) = \sin^{-1}(x)$ is the function $\sin$. However, 
what about the inverses of (invertible) functions ${}_2 F_1(a,b,c; \cdot)$ in 
general? Do we have to work with elliptic integrals here? What part does the 
Jacobi elliptic function play in this? A complex-analytic approach to a 
part of this problem using contour integration is given in \cite[Chapter 5, 
including Theorem 5.6.18]{B2019}. 
\end{enumerate}
\section{Open problem 2: Interrelation between the Grothendieck 
inequality and copulas}
 If we thouroughly overhaul the CCP function $[-1,1] \ni \rho \mapsto 
\frac{2}{\pi}\arcsin(\rho)$ we recognise that some knowledge of Gaussian copulas (i.e., 
finite-dimensional multivariate distribution functions of univariate marginals generated by 
the distribution function of Gaussian random vectors - cf., e.g., \cite{N2006, 
Oe2015, Sk1959, UADF2017}) and (the probabilistic version) of the Hermite 
polynomials might become very fruitful regarding our indicated search for different ``
suitable'' CCP functions. $[-1,1] \ni \rho \mapsto \psi(\frac{1}{2},\frac{1}{2}; t) = 
\frac{2}{\pi}\arcsin(\rho)$ namely reveals as a special case of the CCP function

\scalebox{0.85}{
\vbox{ 
\[
[-1,1] \ni \rho \mapsto \psi(p,p; \rho) = \frac{1}{c(p)}\sum_{n=1}^\infty \frac{1}{n}\,
H_{n-1}^2(\Phi^{-1}(p))\,\rho^{n} = \frac{1}{2\pi\,p(1-p)}\,\exp(-(\Phi^{-1}(p))^2)\,
\rho\sum_{n=0}^\infty \frac{1}{n+1}\,H_n^2(\Phi^{-1}(p))\,\rho^{n}\,, 
\]
}}

\noindent where $0 < p < 1$ and
\[
c(p) : = \frac{p(1-p)}{\varphi^2(\Phi^{-1}(p))} = 2\pi\,p(1-p)\exp((\Phi^{-1}(p))^2) 
= \sum_{n=0}^\infty \frac{1}{n+1}\,H_{n}^2(\Phi^{-1}(p))\,.
\]
If we put 
\[
b_p(x) : = {\text{sign}}(x - \Phi^{-1}(p)) = 
2\,\ind_{[\Phi^{-1}(p), \infty)}(x) - 1 = 1 - 2\,\ind_{(-\infty, \Phi^{-1}(p))}(x) \in 
\{-1,1\}\,,
\] 
where $x \in \R$, then the tetrachoric series expansion of the bivariate Gaussian copula (cf. 
\cite{AF2014, G1963, M2013}) implies the following generalisation of the Grothendieck 
equality:
\begin{eqnarray*}
h_p(\rho) \, : = \, h_{b_p}(\rho) : & = & \E[b_p(X)\,b_p(Y)]
\, = \, (2p-1)^2 + \frac{2}{\pi}\,\exp(-(\Phi^{-1}(p))^2)\,
\sum_{n=1}^\infty \frac{1}{n}\,H_{n-1}^2(\Phi^{-1}(p))\rho^{n}\\
& = & (2p-1)^2 + 4p(1-p)\psi(p,p; \rho)\,.
\end{eqnarray*}
Due to our construction of $\psi(p,p; \cdot)$ the latter is clearly equivalent to
\[
\rho(b_p(X),b_p(Y)) = \psi(p,p; \rho)
\]
for all $p \in (0,1)$, $\rho \in [-1,1]$ and $(X,Y) \sim N_2(0, \Sigma_2(\rho))$, where
$\rho(b_p(X),b_p(Y))$ denotes Pearson's correlation coefficient between the random variables
$b_p(X)$ and $b_p(Y)$. Unfortunately,
\[
h_p(\rho) = \psi(p,p; \rho) \text{ for all } \rho \in [-1,1] \text{ if and only 
if } p = \frac{1}{2}\,.
\]
These facts clearly lead to further research problems; namely:
\begin{enumerate}
\item[(RP4)] Prove whether there are $p \in (-1,1)\setminus\{\frac{1}{2}\}$ and functions 
$\chi_p : \R \longrightarrow \{-1,1\}$ such that $\psi(p,p; \rho) = 
h_{\chi_p}(\rho) = \E[\chi_p(X)\,\chi_p(Y)]$ for all $\rho \in [-1,1]$ 
and $(X,Y) \sim N_2(0, \Sigma_2(\rho))$, so that the condition (SIGN) of our workflow is 
satisfied for $h_{\chi_p}$.
\item[(RP5)] Generalise the above approach (which is built on the tetrachoric 
series of the bivariate Gaussian copula) to the $n$-variate case, where $n \in 
\N_3$.
\item[(RP6)] Verify whether the above approach can be transferred to the complex 
case. Could we then similarly generalise the Haagerup equality?
\item[(RP7)] If (RP4), respectively (RP5) holds, prove whether the  
condition (CRA) of the scheme holds. If this were the case, calculate 
(respectively approximate numerically) the related upper bound of $K_G^\R$. 
Include supercomputers and computer algebra systems if necessary. 
\end{enumerate}
\section{Open problem 3: Non-commutative dependence structures in quantum 
mechanics and the Grothendieck inequality}
Even a mathematical modelling of \textit{non-commutative dependence} 
in quantum theory and its applications to quantum information and quantum 
computing is strongly linked with the existence of the real Grothendieck 
constant $K_G^\R$. 

The latter can be very roughly adumbrated as follows: the experimentally proven non-
Kolmogorovian (non-commutative) nature of the underlying probability theory of quantum 
physics leads to the well-known fact that in general a normal state of a composite quantum 
system cannot be represented as a convex combination of a product of normal states of the 
subsystems. This phenomenon is known as \textit{entanglement} or \textit{quantum 
correlation}. The Einstein-Podolsky-Rosen paradox, the violation of Bell's inequalities 
(limiting \textit{spatial} correlation) and the Leggett-Garg inequalities (limiting 
\textit{temporal} correlation) in quantum mechanics and related theoretical and experimental 
research implied a particular focus on a deeper understanding of this type of correlation - and hence to the 
\textit{modelling of a specific type of dependence} of two (ore more) quantum observables in 
a composite quantum system, measured by two (or more) space-like separated instruments, each 
one having a classical parameter (such as the orientation of an instrument which measures the 
spin of a particle). In this context, a Leggett-Garg inequality (LGI) could be viewed 
as a ``Bell inequality in time''. The transition probability function, i.\,e., 
the \textit{joint} probability distribution of observables in some fixed state 
of the system (considered as a function of the aforementioned parameters) may 
violate Bell's inequalities and is therefore not realisable in ``classical'' 
(commutative) physics. The surprising fact, firstly recognised by B. S. 
Tsirel'son (cf. \cite[Ch. 11.2]{AS2017} and \cite{P2018, T1980, T1993}), 
is that also this - experimentally verified - gap is an implication of the 
existence of the real Grothendieck constant $K_G^\R > 1$ (also known as 
\textit{Tsirel'son bound})! In other words, $K_G^\R$ indicates ``how non-local 
quantum mechanics can be at most''.    

Already in the classical Kolmogorovian model, i.\,e., in the framework of 
probability space triples $(\Omega, \mathscr{F}, \P)$, a rigorous description 
of tail dependence - which \textit{exceeds} the standard dependence measure, 
given by Pearson's correlation coefficient, is a challenging task. To disclose 
(and simulate) the geometry of dependence one has to determine finite-dimensional 
multivariate distribution functions of univariate marginals, hence 
\textit{copulas}. In the description of research problem 2 we have seen that 
Gaussian copulas are lurking in the Grothendieck equality. More precisely, we 
have (cf. \cite{S1889}):  
\begin{ex*}[\textbf{Stieltjes, 1889}]
Let $\rho \in [-1,1]$. Let $X, Y \sim N_1(0,1)$ such that $\E[XY] = \rho$. Then
\[
\E[{\text{sign}}(X){\text{sign}}(Y)] = {4\,C^{\text{Ga}}(\tfrac{1}{2}, 
\tfrac{1}{2}; \rho ) - 1} = {\frac{2}{\pi}\arcsin(\rho)} = 
\frac{2}{\pi}\arcsin(\E[XY]),
\]
where $[-1,1] \ni \rho \mapsto C^{\text{Ga}}(\tfrac{1}{2}, \tfrac{1}{2}; \rho)$ 
denotes the bivariate Gaussian copula with Pearson's correlation coefficient 
$\rho$ as parameter, evaluated at $(\frac{1}{2}, \frac{1}{2})$. 
\end{ex*}
\noindent Always keeping in mind a \textit{non-commutative} version of the 
Grothendieck inequality (cf. \cite{P2012, T1980, T1993}), our conjecture is 
that copulas in function spaces play a non-negligible role here.
Unfortunately, compared to the finite-dimensional setting, the advent of the 
latter confronts us with hard problems. For example, by no means it 
is clear how marginals can be defined in an infinite-dimensional measurable 
vector space. If $X$ is a random variable in a separable Hilbert space $H$, 
projections onto an orthonormal basis $(\langle X, e_n \rangle)_{n \in \N}$ are 
reasonable candidates. This case was studied in \cite{HR2017}. If in addition 
the space considered is a reproducing kernel Hilbert space of functions, over 
$[0, 1]$ say, an equally natural option for marginals would be function 
evaluations $\{X(t) : t \in [0,1]\}$. Here, a new framework is required, 
including the preparation of a general concept of marginals for measurable 
vector spaces (cf. \cite{BDS2022}). Consequently, we get the following problems:
\begin{enumerate}
\item[(RP8)] Look for objects like ``non-commutative copulas'', leading to a 
search for ``non-commutative distribution functions in measurable vector spaces'', 
including a non-commutative version of the famous result of Sklar (cf. \cite{Oe2015} 
and the references therein).
\item[(RP9)] Define a ``multivariate'' spectral theory of 
\textit{non-commuting} normal operator tuples and introduce non-commutative 
tail dependency measures in non-commutative $C^\ast$-algebras and operator 
spaces.
\end{enumerate}
\noindent Let us close \sref{Chapter}{chp:future_research_topics} briefly with 
the following ``blue-sky'' research questions, which appear quite naturally 
and are completely unanswered. Can we improve the approximation results in the 
commutative case if we remove the underlying Gaussian structure in the 
Grothendieck inequality (for both fields, $\R$ and $\C$) and implement 
tail dependent distribution functions instead (such as the generalised extreme 
value (GEV) distribution)? What about infinitely divisible probability 
distributions in general? It is very likely that the use of correlation matrices 
and CCP functions, including linked \textit{Gaussian} copula approaches, would 
no longer suffice (just as it is the case with Brownian motion which is in fact 
a particular case of a L\'{e}vy process, yet without jumps). So, could even general 
semimartingale techniques help to improve the approximations (cf., e.g., 
\cite{EN2019, KR2020})?

\begin{ack*} I would like to thank Professor Feng \underline{Qi} (Henan 
Polytechnic University, China and independent researcher, Dallas, TX 75252, 
USA) for a very helpful correspondence online; particularly with respect 
to subtleties underlying the topics of \sref{Section}{sec:RP_1}. I also 
would like to thank the anonymous helpful non-virtual assistant for a very 
helpful electronic correspondence, which in particular has expanded the list of 
references to include further significant sources.
\end{ack*}

\noindent\dci{I declare that I do have no known competing financial interests 
or personal relationships that could have appeared to influence the work 
presented in this paper.}


\end{document}